\documentclass[a4paper,11pt]{article}

\usepackage{amsmath,amsfonts,amssymb,mathtools,amsthm}
\usepackage{epsfig,mathrsfs,enumerate,graphics}
\usepackage[usenames,dvipsnames]{color}
\usepackage{dsfont}
\usepackage{graphicx}
\usepackage{wrapfig}
\usepackage{caption}
\usepackage{subcaption}

\usepackage[all]{xy} \entrymodifiers={!!<0pt,0.7ex>+}

\usepackage[hyperindex=true, 
					colorlinks=false]{hyperref}

\newcommand{\ens}{\enspace}

\def\al{\alpha}
\def\be{\beta}
\def\ga{\gamma}
\def\Ga{{\Gamma}}
\def\de{\delta}

\def\De{\Delta}
\def\eps{{\varepsilon}}
\def\ka{\kappa}
\def\la{\lambda}

\def\om{\omega}

\def\Om{\Omega}

\def\sig{{\sigma}}
\def\Sig{{\Sigma}}

\def\th{{\theta}}
\def\Th{\Theta}
\newcommand{\ph}{\varphi}

\def\ze{{\zeta}}

\newcommand{\demi}{\frac{1}{2}}
\newcommand{\dem}{\tfrac{1}{2}}
\newcommand{\tiers}{\frac{1}{3}}

\newcommand{\ov}{\overline}

\newcommand{\lan}{\langle}
\newcommand{\ran}{\rangle}

\newcommand{\dist}{\operatorname{dist}}

\newcommand{\id}{\operatorname{id}}
\newcommand{\ID}{\mathop{\hbox{{\rm Id}}}\nolimits}

\newcommand{\I}{{\mathrm i}}
\newcommand{\dd}{{\mathrm d}}
\newcommand{\ex}{{\mathrm e}}
\newcommand{\ee}{\mathrm e}
\newcommand{\eel}{{\underline{\mathrm e}}}
\newcommand{\Log}{\mathrm{Log\hspace{1.5pt}}}
\newcommand{\Llog}{\mathrm{\gL\hspace{-.9pt} og\hspace{1.5pt}}}

\newcommand{\pa}{\partial}

\newcommand{\ii}{^{-1}}
\newcommand{\ic}{^{\circ(-1)}}
\newcommand{\ti}{\tilde}

\newcommand{\IM}{\mathop{\Im m}\nolimits}
\newcommand{\RE}{\mathop{\Re e}\nolimits}

\newcommand{\ie}{{i.e.}\ }
\newcommand{\cf}{{cf.}\ }
\newcommand{\eg}{{e.g.}\ }
\newcommand{\resp}{{resp.}\ }
\newcommand{\wrt}{{with respect to}}
\newcommand{\lhs}{{left-hand side}}
\newcommand{\rhs}{{right-hand side}}

\newcommand{\dst}{\displaystyle}

\newcommand{\C}{\mathbb{C}}
\newcommand{\Clog}{\ti\C}
\newcommand{\D}{\mathbb{D}}

\newcommand{\N}{\mathbb{N}}
\newcommand{\Q}{\mathbb{Q}}
\newcommand{\R}{\mathbb{R}}

\newcommand{\Z}{\mathbb{Z}}

\def\cA{\mathcal{A}}
\def\cB{\mathcal{B}}

\def\cL{\mathcal{L}}

\def\cN{\mathcal{N}}

\def\cV{\mathcal{V}}

\newtheorem{thm}{Theorem}[section]
\newtheorem{lemma}[thm]{Lemma}
\newtheorem{cor}[thm]{Corollary}

\theoremstyle{definition}

\newtheorem{rem}[thm]{Remark}
\newtheorem{Def}[thm]{Definition}
\newtheorem{nota}[thm]{Notation}
\newtheorem{exo}{Exercise}[section]
\newtheorem{exa}[exo]{Example}

\newcounter{parag}[section]
\renewcommand{\theparag}{{\thesection.\arabic{parag}}}
\newcommand{\parag}{\medskip \refstepcounter{parag} \noindent{{\bf\theparag}\ } }

\newcounter{parage}

\newcounter{paraga}

\newcommand{\pp}[1]{^{(#1)}}

\DeclarePairedDelimiter\abs{\lvert}{\rvert}%
\DeclarePairedDelimiter\norm{\lVert}{\rVert}%

\def\om{\omega}

\newcommand{\val}{\operatorname{val}}
\newcommand{\ord}{\operatorname{ord}}
\newcommand{\coeff}{\operatorname{coeff}}

\newcommand{\dDe}[1]{
  \raisebox{.23ex}{
  {$\stackrel{\raisebox{-.23ex}{$\scriptscriptstyle\bullet$}}\Delta_{
  \raisebox{-.23ex}[1ex][0ex]{$\scriptstyle#1$}}$}}
}

\newcommand{\dDep}[1]{
  \raisebox{.23ex}{
  {$\stackrel{\raisebox{-.23ex}{$\scriptscriptstyle\bullet$}}\Delta_{
  \raisebox{-.23ex}{$\scriptstyle#1$}}^{
  \raisebox{-.85ex}{$\scriptstyle +$}}$}}
}

\newcommand{\dDem}[1]{
  \raisebox{.23ex}{
  {$\stackrel{\raisebox{-.23ex}{$\scriptscriptstyle\bullet$}}\Delta_{
  \raisebox{-.23ex}{$\scriptstyle#1$}}^{
  \raisebox{-.85ex}{$\scriptstyle -$}}$}}
}

\newcommand{\dDepm}[1]{
  \raisebox{.23ex}{
  {$\stackrel{\raisebox{-.23ex}{$\scriptscriptstyle\bullet$}}\Delta_{
  \raisebox{-.23ex}{$\scriptstyle#1$}}^{
  \raisebox{-.85ex}{$\scriptstyle\pm$}}$}}
}

\newcommand{\hta}[1]{{\stackrel{\raisebox{-.23ex}{$\scriptscriptstyle\wedge$}}{#1}}}
\newcommand{\ch}[1]{{\stackrel{\raisebox{-.23ex}{$\scriptscriptstyle\vee$}}{#1}}}
\newcommand{\tr}[1]{{\stackrel{\raisebox{-.23ex}{$\scriptscriptstyle\triangledown$}}{#1}}}

\newcommand{\chfsub}[1]{%
\stackrel{\raisebox{-.23ex}{$\scriptscriptstyle\vee$}}{f}\hspace{-.35em}\vphantom{f}_{#1}
}

\newcommand{\htb}[1]{\raisebox{-.23ex}{${\stackrel{
            \raisebox{-.18ex}{$\scriptscriptstyle\wedge$}
          }{#1}
     }$}}

\newcommand{\chb}[1]{\raisebox{-.23ex}{${\stackrel{
            \raisebox{-.23ex}{$\scriptscriptstyle\vee$}
          }{#1}
     }$}}
\newcommand{\chbph}[1]{ \raisebox{-.23ex}{${\stackrel{
            \raisebox{-.23ex}{$\scriptscriptstyle\vee$}
          }{\ph}
     }$}^{\raisebox{-.63ex}{$\scriptstyle#1$}} }
\newcommand{\trb}[1]{\raisebox{-.23ex}{${\stackrel{
            \raisebox{-.23ex}{$\scriptscriptstyle\triangledown$}
          }{#1}
     }$}}
\newcommand{\htn}[1]{\raisebox{.20ex}{${\stackrel{
            \raisebox{-.20ex}{$\scriptscriptstyle\wedge$}
          }{#1}
     }$}}
\newcommand{\chn}[1]{\raisebox{.20ex}{${\stackrel{
            \raisebox{-.20ex}{$\scriptscriptstyle\vee$}
          }{#1}
     }$}}
\newcommand{\trn}[1]{\raisebox{.20ex}{${\stackrel{
            \raisebox{-.20ex}{$\scriptscriptstyle\triangledown$}
          }{#1}
     }$}}
\newcommand{\trigaE}[1]{\raisebox{.20ex}[1.85ex][0ex]{${\stackrel{
            \raisebox{-.20ex}{$\scriptscriptstyle\triangledown$}
          }{E}
     }$}^{I,\ga}_{#1}}
\newcommand{\strn}[1]{\raisebox{0ex}[1ex][0ex]{${\stackrel{
            \raisebox{-.20ex}{$\scriptscriptstyle\triangledown$}
          }{\scriptstyle#1}
     }$}}

\newcommand{\chg}[1]{\raisebox{-.31ex}{${\stackrel{
            \raisebox{-.23ex}{$\scriptscriptstyle\vee$}
          }{#1}
     }$}}

\newcommand{\SING}{\operatorname{SING}}
\newcommand{\ANA}{\operatorname{ANA}}

\newcommand{\simp}{^{\mathrm{simp}}}

\newcommand{\gGR}{\ti\gG^{\mathrm{RES}}}

\newcommand{\gGS}{\ti\gG\simp}

\newcommand{\bem}{\mbox{}^\flat\hspace{-.8pt}}

\newcommand{\sing}{\operatorname{sing}}

\newcommand{\var}{\operatorname{var}}
\newcommand{\cont}{\operatorname{cont}}
\newcommand{\Al}{\operatorname{\cA}}

\def\Bbibitem#1#2{\bibitem[#1]{#2}}

\newcommand{\defeq}{\coloneqq} 

\newcommand{\col}{\colon\thinspace}          

\newcommand{\zcz}{z\ii\C[[z\ii]]}

\newcommand{\gC}{\mathscr C}       
\newcommand{\gD}{\mathscr D}       
\newcommand{\gG}{\mathscr G}       
\newcommand{\gO}{\mathscr O}       
\newcommand{\gP}{\mathscr P}       
\newcommand{\gR}{\mathscr R}       
\newcommand{\gS}{\mathscr S}       
\newcommand{\gL}{\mathscr L}       
\newcommand{\gV}{\mathscr V}       

\newcommand{\eith}{\ee^{\I\th}}

\newcommand{\eul}{^{\raisebox{-.23ex}{$\scriptstyle\mathrm E$}}}
\newcommand{\Eul}{^{\raisebox{-.83ex}{$\scriptstyle\mathrm E$}}}
\newcommand{\poin}{^{\raisebox{-.23ex}{$\scriptstyle\mathrm P$}}}
\newcommand{\Poin}{^{\raisebox{-.83ex}{$\scriptstyle\mathrm P$}}}
\newcommand{\hurw}{^{\raisebox{-.23ex}{$\scriptstyle\mathrm H$}}}
\newcommand{\Hurw}{^{\raisebox{-.83ex}{$\scriptstyle\mathrm H$}}}

\newcommand{\Ddem}{\D_{\rho/2}}

\newcommand{\est}{\emptyset}
\newcommand{\bB}{\text{\boldmath{$B$}}}
\newcommand{\etab}{\text{\boldmath{$\eta$}}}
\newcommand{\bDe}{\text{\boldmath{$\De$}}}

\newcommand{\SSig}{ \{0,\om\} \cup \Sig_\om }

\newcommand{\ds}{
  \raisebox{-0.08ex}{
  {$\stackrel{\raisebox{-0.1ex}{$\scriptscriptstyle\bullet$}}\sig$}}
}

\newcommand{\dl}[1]{
  \raisebox{0.21ex}{
  {$\stackrel{\raisebox{-0.23ex}{$\scriptscriptstyle\bullet$}}\ell_{
  \raisebox{-.23ex}{$\scriptstyle#1$}}$}}
}

\newcommand{\dll}[2]{
  \raisebox{0.21ex}{$
       \stackrel{\raisebox{-0.23ex}{$\scriptscriptstyle\bullet$}}
       \ell_{\raisebox{-.23ex}{$\scriptstyle#1$}}^{
               \raisebox{-.6ex}{$\,\scriptstyle#2$}
        }
   $}
}

\newcommand{\DDe}[2]{\dDe{#1 \raisebox{-.6ex}{$\scriptstyle|\strn E_{#2}$}}}
\newcommand{\DDep}[2]{\dDep{#1 \raisebox{-.6ex}[0ex][0ex]{$\scriptstyle|\strn E_{#2}$}}}

\newcommand{\Ep}{ \raisebox{.20ex}[0ex][0ex]{${\stackrel{
            \raisebox{-.20ex}{$\scriptscriptstyle\triangledown$}
          }{E}
     }$}' }
\newcommand{\Epp}{ \raisebox{.20ex}[0ex][0ex]{${\stackrel{
            \raisebox{-.20ex}{$\scriptscriptstyle\triangledown$}
          }{E}
     }$}'' }

\newcommand{\omp}{\om'\in\Om'\cap d}
\newcommand{\ompp}{\om''\in\Om''\cap d}
\newcommand{\sigp}{\sig'\in\Om'\cap d}
\newcommand{\sigpp}{\sig''\in\Om''\cap d}

\newcommand{\bigop}{
              \stackrel{
                   \raisebox{-.20ex}{$\scriptscriptstyle\wedge$}
                    }{\dst\bigoplus_{\om\in\Om\cap d}} 
           }
\newcommand{\bigopp}{
              \stackrel{
                   \raisebox{-.20ex}{$\scriptscriptstyle\wedge$}
                    }{\dst\bigoplus_{\om'\in\Om'\cap d}} 
           }
\newcommand{\bigoppp}{
              \stackrel{
                   \raisebox{-.20ex}{$\scriptscriptstyle\wedge$}
                    }{\dst\bigoplus_{\om''\in\Om''\cap d}} 
           }
\newcommand{\bigpmN}{
              \stackrel{
                   \raisebox{-.20ex}{$\scriptscriptstyle\wedge$}
                    }{\dst\bigoplus_{\om\in \pm 2\pi\I\N}} 
           }

\newcommand{\beglabel}[1]{\begin{equation}	\label{#1}}
\newcommand{\elabel}{\end{equation}}

\newcommand{\eqrefa}[1]{(\ref{#1}a)}   
\newcommand{\eqrefb}[1]{(\ref{#1}b)}   

\newcommand{\DD}[2]{{\De\hspace{-.65em}\raisebox{.4ex}{$\scriptscriptstyle /$}\hspace{.3em}^{#1}_{#2}}}

\newcommand{\conc}{\raisebox{.25ex}{$\centerdot$}}

\newcommand{\imp}{\ens\Longrightarrow\ens}
\newcommand{\Imp}{\quad\Longrightarrow\quad}

\newcommand{\upp}{^\textnormal{up}}
\newcommand{\low}{^\textnormal{low}}
\newcommand{\uplow}{^\textnormal{up/low}}

\begin{document}

\thispagestyle{empty}
\begin{center}
{\bf \huge
Introduction to $1$-summability\\[1ex] and resurgence}\\[3ex]

David Sauzin
\end{center}
\vspace{1.2cm}

\begin{abstract}

This text is about the mathematical use of certain divergent power series.

The first part is an introduction to $1$-summability.
The definitions rely on the formal Borel transform and the Laplace transform 
along an arbitrary direction of the complex plane.
Given an arc of directions, if a power series is $1$-summable in that
arc, then one can attach to it a Borel-Laplace sum, i.e. a holomorphic
function defined in a large enough sector and asymptotic to that power
series in Gevrey sense.

The second part is an introduction to \'Ecalle's resurgence theory.
A power series is said to be resurgent when its Borel transform is convergent 
and has good analytic continuation properties: there may be singularities but 
they must be isolated.
The analysis of these singularities, through the so-called alien calculus, allows
one to compare the various Borel-Laplace sums attached to the same resurgent
$1$-summable series.
In the context of analytic difference-or-differential equations, this sheds light on 
the Stokes phenomenon.

A few elementary or classical examples are given a thorough treatment (the Euler
series, the Stirling series, a less known example by Poincar\'e).
Special attention is devoted to non-linear operations: $1$-summable
series as well as resurgent series are shown to form algebras which
are stable by composition. 
As an application, the resurgent approach to the classification of
tangent-to-identity germs of holomorphic diffeomorphisms in the
simplest case is included.
An example of a class of non-linear differential equations giving rise
to resurgent solutions is also presented.
The exposition is as self-contained as can be, requiring only some
familiarity with holomorphic functions of one complex variable.

\end{abstract}

\pagebreak

\setcounter{tocdepth}{1}
\tableofcontents

\pagebreak





This text grew out of a course given in a CIMPA school in Lima in 2008
and further courses taught at the Scuola Normale Superiore di Pisa between
2008 and 2010.
We tried to make it as self-contained as possible, so that it can be
used with undergraduate students, assuming only on their part some
familiarity with holomorphic functions of one complex variable.
The first part of the text
(Sections~\ref{sec:Prol}--\ref{secgrponesummdiffeos}) is an
introduction to $1$-summability. 
The second part (Sections~\ref{sec_resurfunct}--\ref{sec:Bridge}) is
an introduction to \'Ecalle's resurgence theory.

Throughout the text, we use the notations
\[
\N = \{0,1,2, \ldots\}, \qquad \N^* = \{1,2,3,\ldots\}
\]
for the set of non-negative integers and the set of positive integers.


\vspace{1.2cm}

\centerline{\Large\sc Introduction}
\addcontentsline{toc}{part}{\sc Introduction}

\vspace{.3cm}


\section{Prologue}	\label{sec:Prol}


\parag
At the beginning of the second volume of his \emph{New methods of celestial
mechanics} \cite{Poincare}, H.~Poincar\'e dedicates two pages to elucidating
\emph{``a kind of misunderstanding between geometers and astronomers about the meaning
of the word \emph{convergence}''}.
He proposes a simple example, namely the two series
\begin{equation}	\label{eqsimplePoinc}
\sum\frac{1000^n}{n!} \quad\textrm{and}\quad 
\sum\frac{n!}{1000^n}.
\end{equation}
He says that, for geometers (\ie mathematicians), the first one is convergent because the
term for $n=1.000.000$ is much smaller than the term for $n=999.999$, whereas
the second one is divergent because the general term is unbounded
(indeed, the $(n+1)$-th term is obtained from the $n$th one by multiplying
either by $1000/n$ or by $n/1000$).
On the contrary, according to Poincar\'e, astronomers will consider the
first series as divergent because the general term is an increasing function
of~$n$ for $n\le 1000$, and they will consider the second one as convergent
because the first $1000$ terms decrease rapidly.

He then proposes to reconcile both points of view by clarifying the role that
divergent series (in the sense of geometers) can play in the approximation of
certain functions.
He mentions the example of the classical Stirling series, for which the absolute
value of the general term is first a decreasing function of~$n$ and then an
increasing function; this is a divergent series and still, Poincar\'e says,
\emph{``when stopping at the least term one gets a representation of Euler's
gamma function, with greater accuracy if the argument is larger''}.
This is the origin of the modern theory of asymptotic expansions.\footnote{%
In fact, Poincar\'e's observations go even beyond this, in direction of least
term summation for Gevrey series, but we shall not discuss the details of all
this in the present article; the interested reader may consult \cite{Ramis},
\cite{gazRamisA}, \cite{gazRamisB}.}


\parag
In this text we shall focus on formal series given as power series expansions,
like the Stirling series for instance, rather than on numerical series. Thus, we
would rather present Poincar\'e's simple example~\eqref{eqsimplePoinc} in the
form of two formal series
\begin{equation}	\label{eqPoincEuler}
\sum_{n\ge0}\frac{1000^n}{n!} t^n \quad\textrm{and}\quad 
\sum_{n\ge0} \frac{n!}{1000^n} t^n,
\end{equation}
the first of which has infinite radius of convergence, while the second has zero
radius of convergence.
For us, \emph{divergent series} will usually mean a formal power series with
zero radius of convergence.

Our first aim in this text is to discuss the Borel-Laplace summation process as a way of
obtaining a function from a (possibly divergent) formal series,
the relation between the original formal series and this function being a
particular case of asymptotic expansion of Gevrey type.
For instance, this will be illustrated on Euler's gamma function and the Stirling series
(see Section~\ref{sec:Stirling}).
But we shall also describe in this example and others the phenomenon for which J.~\'Ecalle coined
the name \emph{resurgence} at the beginning of the 1980s and give a brief
introduction to this beautiful theory.


\section{An example by Poincar\'e}	\label{secExPoin}


Before stating the basic definitions and introducing the tools with which we
shall work, we want to give an example of a divergent formal series~$\ti\phi(t)$
arising in connection with a holomorphic function~$\phi(t)$ (later on, we shall
come back to this example and see how the general theory helps to understand
it).
Up to changes in the notation this example is taken from Poincar\'e's discussion
of divergent series, still at the beginning of \cite{Poincare}.

Fix $w\in\C$ with $0<|w|<1$ and consider the series of functions of the complex variable~$t$
\begin{equation}	\label{eqdefphiP}
\phi(t) = \sum_{k\ge0} \phi_k(t), \qquad
\phi_k(t) = \frac{w^k}{1 + k t}.
\end{equation}
This series is uniformly convergent in any compact subset of $U = \C^* \setminus
\big\{-1, -\demi, -\tiers, \ldots\big\}$, as is easily checked,
thus its sum~$\phi$ is holomorphic in~$U$.

We can even check that~$\phi$ is meromorphic in~$\C^*$ with a simple pole at
every point of the form $-\frac{1}{k}$ with $k\in\N^*$:
Indeed, $\C^*$ can be written as the union of the open sets
\[ \Om_N = \{ t\in \C \mid |t|>1/N \} \]
for all $N\ge1$;
for each~$N$, the finite sum $\phi_0+\phi_1+\cdots+\phi_N$
is meromorphic in~$\Om_N$ with simple poles at
$-1,-\demi,\ldots,-\frac{1}{N-1}$, on the other hand the functions~$\phi_k$ are
holomorphic in~$\Om_N$ for all $k\ge N+1$, with
$ \dst |\phi_k(t)| \le \frac{|w|^k}{k |t+1/k|} \le
\left(\frac{1}{N}-\frac{1}{N+1}\right)\ii \frac{|w|^k}{k} $,
whence the uniform convergence and the holomorphy in~$\Om_N$ of 
$\phi_{N+1}+\phi_{N+2}+\cdots$
%
follow, and consequently the meromorphy of~$\phi$.

We now show how this function~$\phi$ gives rise to a divergent formal series when~$t$ approaches~$0$.
For each $k\in\N$, we have a convergent Taylor expansion at the origin
\[
\phi_k(t) = \sum_{n\ge0} (-1)^n w^k k^n t^n
\quad \textrm{for $|t|$ small enough.}
\]
Since for each $n\in\N$ the numerical series
\begin{equation}	\label{eqdefbn}
b_n = \sum_{k\ge0} k^n w^k
\end{equation}
is convergent, one could be tempted to recombine the (convergent) Taylor expansions of the~$\phi_k$'s as
$\dst \phi(t) \textrm{``$=$''} \sum_k \left( \sum_n (-1)^n w^k k^n t^n \right)
\textrm{``$=$''} \sum_n (-1)^n \left(\sum_k k^n w^k \right) t^n
$,
which amounts to considering the well-defined formal series
\begin{equation}	\label{eqdefformalphiP}
\ti\phi(t) = \sum_{n\ge0} (-1)^n b_n t^n
\end{equation}
as a Taylor expansion at~$0$
for~$\phi(t)$.
But it turns out that \emph{this formal series is divergent!}

Indeed, the coefficients~$b_n$ can be considered as functions of the complex
variable $w = \ex^s$, for $w$ in the unit disc or, equivalently, for $\RE s<0$; we
have
$b_0 = (1-w)\ii = (1-\ex^s)\ii$ and
$b_n = \left(w\frac{d}{d w}\right)^n b_0 = \left(\frac{d}{d s}\right)^n b_0$.
Now, if~$\ti \phi(t)$ had nonzero radius of convergence, there would exist $A,B>0$ such
that $|b_n|\le A B^n$ and
the formal series 
\begin{equation}	\label{eqdefBorphiP}
F(\ze) = \sum (-1)^n b_n \frac{\ze^n}{n!}
\end{equation}
 would have infinite radius of convergence,
whereas, recognizing the Taylor formula of~$b_0$ \wrt\ the variable~$s$, we see
that
$F(\ze) = \sum (-1)^n \frac{\ze^n}{n!} \left(\frac{d}{d s}\right)^n b_0 = (1-\ex^{s-\ze})\ii$
has a finite radius of convergence ($F(\ze)$ is in fact the Taylor expansion
at~$0$ of a meromorphic function with poles on $s+2\pi\I\Z$, thus this radius of
convergence is $\dist(s,2\pi\I\Z)$).

Now the question is to understand the relation between the divergent formal series
$\ti\phi(t)$ and the function~$\phi(t)$ we started from.
We shall see in this course that the Borel-Laplace summation is a way of going from~$\ti\phi(t)$
to~$\phi(t)$, that $\ti\phi(t)$ is the asymptotic expansion of~$\phi(t)$ as
$|t|\to0$ in a very precise sense and we shall explain what resurgence means in this example.

\begin{rem}	\label{remdefGevrey}
We can already observe that the moduli of the
coefficients~$b_n$ satisfy
\begin{equation}	\label{ineqdefGevrey}
|b_n| \le A B^n n!, \qquad n\in\N,
\end{equation}
for appropriate constants $A$ and~$B$ (independent of~$n$).
Such inequalities are called \emph{$1$-Gevrey estimates} for the formal series
$\ti\phi(t) = \sum b_n t^n$.
For the specific example of the coefficients~\eqref{eqdefbn},
inequalities~\eqref{ineqdefGevrey} 
can be obtained by reverting the last piece of reasoning: since the meromorphic
function $F(\ze)$ is holomorphic for $|\ze| < d = \dist(s,2\pi\I\Z)$ and 
$b_n = (-1)^n F^{(n)}(0)$,
the Cauchy inequalities yield~\eqref{ineqdefGevrey} with any $B > 1/d$.
\end{rem}

\begin{rem}
The function~$\phi$ we started with is not holomorphic
(nor meromorphic) in any neighbourhood of~$0$, because of the accumulation at
the origin of the sequence of simple poles $-1/k$;
it would thus have been quite surprising to find a positive radius of convergence
for~$\ti\phi$.
\end{rem}


\vspace{1.2cm}

\centerline{\Large\sc The differential algebra $\C[[z\ii]]_1$}

\bigskip

\centerline{\Large\sc and the formal Borel transform}
\addcontentsline{toc}{part}{\sc The differential algebra $\C[[z\ii]]_1$ and the
formal Borel transform} 

\vspace{.3cm}


\section{The differential algebra $\big(\C[[z\ii]],\pa\big)$}


\parag
It will be convenient for us to set $z=1/t$ in order to ``work at~$\infty$''
rather than at the origin.  
At the level of formal expansions, this simply means that we shall deal with
expansions involving non-positive integer powers of the indeterminate.
We denote by 
\[
\C[[z\ii]] = \bigg\{
\ph(z) = \sum_{n\ge0} a_n z^{-n}, \;
\textrm{with any $a_0,a_1,\ldots \in \C$} 
\bigg\}
\]
the set of all these formal series.
This is a complex vector space, and also an algebra when we take into account
the Cauchy product:
\[
\Big( \sum_{n\ge0} a_n z^{-n} \Big) \Big( \sum_{n\ge0} b_n z^{-n} \Big)
= \sum_{n\ge0} c_n z^{-n}, \qquad
c_n = \sum_{p+q=n} a_p b_q.
\]

The natural derivation %
\begin{equation}
\pa = \dfrac{\dd\,\,}{\dd z}
\end{equation} 
makes it a \emph{differential algebra};
this simply means that we have singled out a $\C$-linear map which satisfies the Leibniz rule
\begin{equation}	\label{eqLeibnizrule}
\pa(\ph\psi) = (\pa\ph)\psi + \ph(\pa\psi),
\qquad \ph,\psi \in \C[[z\ii]].
\end{equation}
If we return to the variable~$t$ and define $D = -t^2\dfrac{\dd\,\,}{\dd t}$, we
obviously get an isomorphism of differential algebras between
$\big(\C[[z\ii]],\pa\big)$ and
$\big(\C[[t]],D\big)$
by mapping $\sum a_n z^{-n}$ to $\sum a_n t^n$.


\parag
The standard valuation (or `order') on $\C[[z\ii]]$ is the map
\beglabel{eqstdval}
\val \col \C[[z\ii]] \to \N\cup\{\infty\}
\elabel
defined by $\val(0)=\infty$ and
$\val(\ph) = \min\{n\in\N \mid a_n\neq0\}$ 
for $\ph = \sum a_n z^{-n} \neq 0$.

For $\nu\in\N$, we shall use the notation 
\begin{equation}	\label{eqdefznuCzii}
z^{-\nu}\C[[z\ii]] =
\bigg\{ \sum_{n\ge\nu} a_n z^{-n}, \;
\textrm{with any $a_\nu,a_{\nu+1},\ldots \in \C$} \bigg\}.
\end{equation}
This is precisely the set of all $\ph\in\C[[z\ii]]$ such that $\val(\ph)\ge\nu$.
In particular, from the viewpoint of the ring structure, $\mathfrak{I}=z\ii\C[[z\ii]]$ is the
maximal ideal of $\C[[z\ii]]$; its elements will often be referred to as
``formal series without constant term''.

Observe that
\begin{equation}	\label{ineqvald}
\val(\pa\ph) \ge \val(\ph) + 1,
\qquad \ph\in\C[[z\ii]],
\end{equation}
with equality if and only if $\ph\in\zcz$.


\parag	\label{parKrull}
It is an exercise to check that the formula
\beglabel{eqdefdistK}
d(\ph,\psi) = 2^{-\val(\psi-\ph)}, \qquad \ph,\psi\in\C[[z\ii]],
\elabel
defines a distance and that $\C[[z\ii]]$ then becomes a complete metric space.
The topology induced by this distance is called the Krull topology or the
topology of the formal convergence (or the $\mathfrak{I}$-adic topology).
It provides a simple way of using the language of topology to describe certain
algebraic properties.

We leave it to the reader to check that a sequence $(\ph_p)_{p\in\N}$ of formal
series is a Cauchy sequence if and only if, for each $n\in\N$, the
sequence of the $n$th coefficients is stationary:
there exists an integer $\mu_n$ such that $\coeff_n(\ph_p)$ is the
same complex number $\al_n$ for all $p\ge\mu_n$.
The limit $\dst\lim_{p\to\infty} \ph_p$ is then simply $\dst\sum_{n\ge0} \al_n z^{-n}$.
(This property of formal convergence of $(\ph_p)$ has no
relation with any topology on the field of coefficients, except with
the discrete one).

In practice, we shall use the fact that a series of formal series $\sum
\ph_p$ is formally convergent if and only if there is a sequence of integers 
$\nu_p \xrightarrow[p\to\infty]{} \infty$
such that $\ph_p \in z^{-\nu_p} \C[[z\ii]]$ for all~$p$.
Each coefficient of the sum $\ph=\sum\ph_p$ is then given by a finite sum:
the coefficient of~$z^{-n}$ in~$\ph$ is 
$\dst \coeff_n(\ph) = \sum_{p\in M_n} \coeff_n(\ph_p)$,
where $M_n = \{ p \mid \nu_p \le n \}$.


\begin{exo}
  Check that, as claimed above, \eqref{eqdefdistK} defines a distance which makes
  $\C[[z\ii]]$ a complete metric space; check that the subspace $\C[z\ii]$ of
  polynomial formal series is dense.
Show that, for the Krull topology, $\C[[z\ii]]$ is a topological ring (\ie addition and multiplication
are continuous) but not a topological $\C$-algebra (the scalar multiplication
is not).
Show that $\pa$ is a contraction for the distance~\eqref{eqdefdistK}.
\end{exo}


\parag
As an illustration of the use of the Krull topology, let us define the
\emph{composition operators} by means of formally convergent series.

Given $\ph,\chi\in\C[[z\ii]]$, we observe that $\val(\pa^p \ph) \ge \val(\ph) + p$
(by repeated use of~\eqref{ineqvald}), hence $\val(\chi^p \, \pa^p \ph) \ge \val(\ph)
+ p$ and the series
\beglabel{eqdefcircIDchi}
\ph \circ (\id+\chi) \defeq \sum_{p\ge0} \frac{1}{p!} \chi^p \, \pa^p \ph
\elabel
is formally convergent. Moreover
\beglabel{ineqvalcirc}
\val\big( \ph \circ (\id+\chi) \big) = \val(\ph).
\elabel
We leave it as an exercise for the reader to check that, for fixed~$\chi$,
the operator $\Th \col \ph \mapsto \ph \circ (\id+\chi)$ is a continuous
automorphism of algebra (\ie a $\C$-linear invertible map, continuous for the
Krull topology, such that $\Th(\ph\psi)= (\Th\ph)(\Th\psi)$).

A particular case is the shift operator
\begin{equation}	\label{eqdefshift}
T_c \col
\C[[z\ii]] \to \C[[z\ii]], \quad
\ph(z) \mapsto \ph(z+c)
\end{equation}
with any $c\in\C$
(the operator~$T_c$ is even a differential algebra automorphism, \ie
an automorphism of algebra which commutes with the differential~$\pa$).

The counterpart of these operators in $\C[[t]]$ via the change of indeterminate
$t=z\ii$ is
$\phi(t) \mapsto \phi(\frac{t}{1+ct})$ for the shift operator and, more
generally for the composition with $\id+\chi$,
$\phi \mapsto \phi\circ F$ with 
$F(t) = \frac{t}{1+t G(t)}$, $G(t)=\chi(t\ii)$.
See Sections~\ref{secGermsholdiffeos}--\ref{secinversgrpgG} for more on the composition of formal series
at~$\infty$ (in particular for associativity).

\begin{exo}
[\textbf{Substitution into a power series}]
Check that, for any $\ph(z) \in \zcz$, the formula
\[
H(t) = \sum_{p\ge0} h_p t^p 
\mapsto
H\circ\ph(z) \defeq \sum_{p\ge0} h_p \big(\ph(z)\big)^p
\]
defines a homomorphism of algebras from $\C[[t]]$ to $\C[[z\ii]]$, \ie a
linear map~$\Th$ such that $\Th 1=1$ and $\Th(H_1 H_2)=
(\Th H_1)(\Th H_2)$ for all $H_1,H_2$.
\end{exo}

\begin{exo}    \label{exoAllAlgEndo}
Put the Krull topology on~$\C[[t]]$ and use it to define the
composition operator $C_F \col \phi\mapsto\phi\circ F$ for any $F\in t\C[[t]]$;
check that $C_F$ is an algebra endomorphism of~$\C[[t]]$.
Prove that any algebra endomorphim~$\Th$ of~$\C[[t]]$ is of this form.
(Hint: justify that
$\phi \in t\C[[t]] \iff \forall\al\in\C^*,\ \al+\phi\ \text{invertible}
\implies \forall\al\in\C^*,\ \al+\Th\phi\ \text{invertible}$;
deduce that $F\defeq \Th t \in t\C[[t]]$;
then, for any $\phi\in\C[[t]]$ and $k\in\N$, show that $\val(\Th\phi-C_F\phi)\ge
k$ by writing $\phi = P + t^k \psi$ with a polynomial~$P$ and conclude.)
\end{exo}



\section{The formal Borel transform and the space of $1$-Gevrey formal series $\C[[z\ii]]_1$}


\parag
We now define a map on the space $z\ii\C[[z\ii]]$ of formal series
without constant term (recall the notation~\eqref{eqdefznuCzii}):

\begin{Def}
The \emph{formal Borel transform} is the linear map 
$\cB \col z\ii\C[[z\ii]] \to \C[[\ze]]$ 
defined by
\[
\cB \col \ti\ph = \sum_{n=0}^\infty a_n z^{-n-1} 
\mapsto
\hat\ph = \sum_{n=0}^\infty a_n \frac{\ze^n}{n!}.
\]
\end{Def}

In other words, we simply shift the powers by one unit and divide the $n$th
coefficient by~$n!$.
Changing the name of the indeterminate from~$z$ (or~$z\ii$) into~$\ze$ is only a
matter of convention, however we strongly advise against keeping the same symbol.

The motivation for introducing~$\cB$ will appear in
Sections~\ref{secLaplaceTrsf} and~\ref{secfineBL} with the use of the Laplace transform.

The map~$\cB$ is obviously a linear isomorphism between the spaces 
$z\ii\C[[z\ii]]$ and $\C[[\ze]]$. 
Let us see what happens with the \emph{convergent} formal series of the first of these
spaces.
We say that $\ti\ph(z) \in \C[[z\ii]]$ is convergent at~$\infty$ (or simply
`convergent') if the associated formal series $\ti\phi(t) = \ti\ph(1/z) \in
\C[[t]]$ has positive radius of convergence.
The set of convergent formal series at~$\infty$ is denoted $\C\{z\ii\}$; the ones
without constant term form a subspace denoted $z\ii\C\{z\ii\}$.

\begin{lemma}	\label{lemCVcase}
Let $\ti\ph\in z\ii\C[[z\ii]]$.
Then $\ti\ph\in z\ii\C\{z\ii\}$ if and only if its formal Borel transfom
$\hat\ph = \cB\ti\ph$
has infinite radius of convergence and defines an entire function of bounded
exponential type,
\ie there exist $A,c>0$ such that $|\hat\ph(\ze)|\le A\,\ee^{c|\ze|}$
for all $\ze\in\C$.
\end{lemma}
\begin{proof}
Let $\ti\ph(z) = \sum_{n\ge0} a_n z^{-n-1}$. 
This formal series is convergent if and only if there exist $A,c>0$ such that,
for all $n\in\N$, $|a_n| \le A c^n$.

If it is so, then $\left| a_n \ze^n / n! \right| \le A |c\ze|^n \ n!$, whence the
conclusion follows.

Conversely, suppose $\hat\ph = \cB\ti\ph$ sums to an entire function satisfying 
$|\hat\ph(\ze)|\le A\,\ee^{c|\ze|}$ for all $\ze\in\C$
and fix $n\in\N$.
We have $a_n = \hat\ph^{(n)}(0)$ and, applying the Cauchy inequality with a
circle $C(0,R) = \{\, \ze\in\C \mid |\ze|=R \,\}$, we get
\[
|a_n| \le \frac{n!}{R^n} \max_{C(0,R)}|\hat\ph|
\le \frac{n!}{R^n} A \, \ee^{cR}.
\]
Choosing $R=n$ and using $n! = 1\times2\times\cdots\times n\le n^n$,
%
%
%
we obtain $|a_n| \le A (\ee^c)^n$, which concludes the proof.
\end{proof}

The most basic example is the geometric series
\begin{equation}	\label{eqdefchic}
\ti\chi_c(z) = z\ii (1 - c z\ii)\ii = T_{-c}(z\ii)
\end{equation}
convergent for $|z|>|c|$, where $c\in\C$ is fixed.
Its formal Borel transform is the exponential series
\begin{equation}	\label{eqBorchic}
\hat\chi_c(\ze) = \ee^{c\ze}.
\end{equation}


\parag
In fact, we shall be more interested in formal series of~$\C[[\ze]]$ having positive but not
necessarily infinite radius of convergence.
They will correspond to power expansions in~$z\ii$ satisfying Gevrey estimates
similar to the ones encountered in Remark~\ref{remdefGevrey}:

\begin{Def}
We call \emph{$1$-Gevrey formal series} any formal series 
$\ti\ph(z) = \sum_{n\ge0} a_n z^{-n}\in \C[[z\ii]]$
for which there exist $A,B>0$ such that
$|a_n| \le A B^n n!$ for all $n\ge0$.
$1$-Gevrey formal series make up a vector space denoted by~$\C[[z\ii]]_1$.
\end{Def}

\begin{lemma}	\label{lemBorelGevreyCV}
Let $\ti\ph \in \zcz$ and $\hat\ph = \cB\ti\ph \in\C[[\ze]]$.
Then $\hat\ph\in\C\{\ze\}$
(\ie the formal series $\hat\ph(\ze)$ has positive radius of convergence)
if and only if $\ti\ph \in \C[[z\ii]]_1$.
\end{lemma}

\begin{proof}
Obvious.
\end{proof}

In other words, a formal series without constant term is $1$-Gevrey if and only
if its formal Borel transform is convergent.
The space of $1$-Gevrey formal series without constant term will be denoted
$\zcz_1 = \cB\ii \big(\C\{\ze\}\big)$, thus
\begin{equation}	\label{eqdecompCZiiGev}
\C[[z\ii]]_1 = \C \oplus \zcz_1.
\end{equation}


\parag
We leave it to the reader to check the following elementary properties:
\begin{lemma}	\label{lemelemptiescB}
If $\ti\ph\in z\ii\C[[z\ii]]$ and $\hat\ph = \cB\ti\ph \in \C[[\ze]]$, then
\begin{itemize}
\item
$\pa\ti\ph\in z^{-2}\C[[z\ii]]$ and $\cB(\pa\ti\ph) = -\ze\hat\ph(\ze)$,
\item
$T_c\ti\ph\in z\ii\C[[z\ii]]$ and $\cB(T_c\ti\ph) = \ee^{-c\ze}\hat\ph(\ze)$ for
any $c\in\C$,
\item
$\cB(z\ii\ti\ph) = \int_0^\ze \hat\ph(\ze_1)\,\dd\ze_1$,
\item
if $\ti\ph\in z^{-2}\C[[z\ii]]$ then $\cB(z\ti\ph) = \dfrac{\dd\hat\ph}{\dd\ze}$.
\end{itemize}
\end{lemma}

In the third property, the integration in the \rhs\ is to be interpreted
termwise.
The second property can be used to deduce~\eqref{eqBorchic} from the fact that,
according to~\eqref{eqdefchic}, $\ti\chi_c = T_{-c}(\ti\chi_0)$ and $\ti\chi_0 =
z\ii$ has Borel tranform~$= 1$.

Since $\frac{\ee^{-\ze}-1}{\ze}$ is invertible in $\C[[\ze]]$ and in
$\C\{\ze\}$, the second property implies 

\begin{cor}	\label{coreqdifflin}
Given $\ti\psi \in z^{-2}\C[[z\ii]]$, with Borel transform
$\hat\psi(\ze)\in\ze\C[[\ze]]$, the equation
\[
\ti\ph(z+1) - \ti\ph(z) = \ti\psi(z)
\]
admits a unique solution~$\ti\ph$ in $\zcz$, whose Borel transform is
given by
\[
\hat\ph(\ze) = \frac{1}{\ee^{-\ze}-1} \hat\psi(\ze). 
\]
If $\ti\psi(z)$ is $1$-Gevrey, then so is the solution~$\ti\ph(z)$.
\end{cor}


\section{The convolution in $\C[[\ze]]$ and in $\C\{\ze\}$}	\label{sec_convols}


\parag
The convolution product, denoted by the symbol~$*$, is defined as the push-forward
by~$\cB$ of the Cauchy product:
\begin{Def}
  Given two formal series $\hat\ph,\hat\psi\in\C[[\ze]]$, their
  \emph{convolution product} is
$\hat\ph*\hat\psi := \cB(\ti\ph \ti\psi)$,
where $\ti\ph = \cB\ii\hat\ph$, $\ti\psi = \cB\ii\hat\psi$.
\end{Def}

At the level of coefficients, we thus have
\begin{equation}	\label{eqdefcoeffconvol}
\hat\ph=\sum_{n\ge0} a_n \frac{\ze^n}{n!},\ens
\hat\psi=\sum_{n\ge0} b_n \frac{\ze^n}{n!} 
\quad\Longrightarrow\quad
\hat\ph*\hat\psi = \sum_{n\ge0} c_n \frac{\ze^{n+1}}{(n+1)!} 
\ens\textrm{with}\;
c_n = \sum_{p+q=n} a_p b_q.
\end{equation}

The convolution product is bilinear, commutative and
associative in $\C[[\ze]]$ (because the Cauchy product is bilinear, commutative and
associative in $z\ii\C[[z\ii]]$).
It has no unit in $\C[[\ze]]$ (since the Cauchy product, when restricted to
$\zcz$, has no unit).
One remedy consists in \emph{adjoining} a unit:
consider the vector space $\C\times\C[[\ze]]$, in which we denote the element
$(1,0)$ by~$\de$; we can write this space as $\C\de\oplus\C[[\ze]]$ if we
identify the subspace $\{0\}\times\C[[\ze]]$ with $\C[[\ze]]$. Defining the
product by
\[ (a\de+\hat\ph) * (b\de+\hat\psi) \defeq 
ab\de + a\hat\psi + b\hat\ph + \hat\ph*\hat\psi, \]
we extend the convolution law of~$\C[[\ze]]$ and
get a unital algebra $\C\de\oplus\C[[\ze]]$ in which $\C[[\ze]]$ is embedded;
by setting 
\[
\cB 1 \defeq \de,
\]
we extend~$\cB$ as an algebra isomorphism between $\C[[z\ii]]$ and
$\C\de\oplus\C[[\ze]]$.
The formula 
\beglabel{eqdefhatpa}
\hat\pa \col a\de+\hat\ph(\ze) \mapsto -\ze\hat\ph(\ze)
\elabel
defines a derivation of $\C\de\oplus\C[[\ze]]$ and the extended~$\cB$ appears as an
isomorphism of differential algebras
\[
\cB \col \big( \C[[z\ii]],\pa \big) \xrightarrow{\sim} \big( \C\de\oplus\C[[\ze]], \hat\pa \big)
\]
(simple consequence of the first property in Lemma~\ref{lemelemptiescB}).
It induces a linear isomorphism 
\beglabel{eqalgisGevdeCV}
\cB \col \C[[z\ii]]_1 \xrightarrow{\sim} \C\de\oplus\C\{\ze\}
\elabel
(in view of~\eqref{eqdecompCZiiGev} and Lemma~\ref{lemBorelGevreyCV}),
which is in fact an algebra isomorphism: we shall see in
Lemma~\ref{lemconvoldisc} that $\C\de\oplus\C\{\ze\}$ is a subalgebra
of $\C\de\oplus\C[[\ze]]$, and hence $\C[[z\ii]]_1$ is a subalgebra of
$\C[[z\ii]]$.

\begin{rem}
For $c \in \C$, the formula
\beglabel{eqdefcounterTc}
\hat T_c \col a\de+\hat\ph(\ze) \mapsto a\de + \ee^{-c\ze}\hat\ph(\ze)
\elabel
defines a differential algebra automorphism of $\big( \C\de\oplus\C[[\ze]], \hat\pa
\big)$, which is the counterpart of the operator~$T_c$ via the extended Borel transform. 
\end{rem}


\parag
When particularized to convergent formal series of the indeterminate~$\ze$, the
convolution can be given a more analytic description:
\begin{lemma}	\label{lemconvoldisc}
Consider two convergent formal series $\hat\ph, \hat\psi\in\C\{\ze\}$.
Let $R>0$ be smaller than the radius of convergence of each of them and
denote by $\Phi$ and~$\Psi$ the holomorphic functions defined by~$\hat\ph$
and~$\hat\psi$ in the disc $D(0,R) = \{\, \ze\in\C \mid |\ze| < R \,\}$.
Then the formula
\begin{equation}	\label{eqdefanalyticconvol}
\Phi*\Psi(\ze) = \int_0^\ze \Phi(\ze_1)\Psi(\ze-\ze_1) \,\dd\ze_1
\end{equation}
defines a function $\Phi*\Psi$ holomorphic in~$D(0,R)$ which is the sum of the
formal series $\hat\ph*\hat\psi$ (the radius of convergence of which is thus at
least~$R$). 
\end{lemma}

\begin{proof}
By assumption, the power series
\[
\hat\ph(\ze) = \sum_{n\ge0} a_n \frac{\ze^n}{n!}
\ens\text{and}\ens
\hat\psi(\ze) = \sum_{n\ge0} b_n \frac{\ze^n}{n!}
\]
sum to $\Phi(\ze)$ and~$\Psi(\ze)$ for any $\ze$ in $D(0,R)$.

Formula~\eqref{eqdefanalyticconvol} defines a function holomorphic
in~$D(0,R)$, since
$\Phi*\Psi(\ze) = \int_0^1 F(s,\ze) \,\dd s$ with 
\begin{equation}	\label{eqdefFconvol}
(s,\ze) \;\mapsto\; F(s,\ze) = \ze\Phi(s\ze)\Psi\big((1-s)\ze\big)
\end{equation}
continuous in~$s$, holomorphic in~$\ze$ and bounded in $[0,1]\times D(0,R')$ for any $R'<R$.

Now, manipulating $F(s,\ze)$ as a product of absolutely convergent series, we write
\[
F(s,\ze) = \sum_{p,q\ge0} a_p b_q \frac{(s\ze)^p}{p!} \frac{\big((1-s)\ze\big)^q}{q!} \ze
= \sum_{n\ge0} F_n(s) \ze^{n+1}
\]
with $F_n(s) = \sum_{p+q=n}  a_p b_q \frac{s^p}{p!} \frac{(1-s)^q}{q!}$;
the elementary identity
$\int_0^1 \frac{s^p}{p!} \frac{(1-s)^q}{q!}\,\dd s = \frac{1}{(p+q+1)!}$ yields
$\int_0^1 F_n(s)\,\dd s = \frac{c_n}{(n+1)!}$ with $c_n = \sum_{p+q=n} a_p b_q$,
hence 
\[\
\Phi*\Psi(\ze) = \sum_{n\ge0} c_n \frac{\ze^{n+1}}{(n+1)!}
\]
for any $\ze\in D(0,R)$;
recognizing in the \rhs\ the formal series $\hat\ph*\hat\psi$ (\cf
\eqref{eqdefcoeffconvol}), we conclude that this formal series has radius of
convergence~$\ge R$ and sums to~$\Phi*\Psi$.
\end{proof}


For instance, since $\cB z\ii = 1$, the \lhs\ in the third property of
Lemma~\ref{lemelemptiescB} can be written $(1*\hat\ph)(\ze)$ and, if
$\ti\ph(z)\in \zcz_1$, the integral $\int_0^\ze
\hat\ph(\ze_1)\,\dd\ze_1$ in the \rhs\ can now be given its usual analytical
meaning: it is the antiderivative of~$\hat\ph$ which vanishes at~$0$.


We usually make no difference between a convergent formal
series~$\hat\ph$ and the holomorphic function~$\Phi$ that it defines in a
neighbourhood of the origin; for instance we usually denote them by the same
symbol and consider that the convolution law defined by the
integral~\eqref{eqdefanalyticconvol} coincides with the restriction
to~$\C\{\ze\}$ of the convolution law of~$\C[[\ze]]$.
However, as we shall see from Section~\ref{sec_resurfunct} onward, things get more
complicated when we consider the analytic continuation in the large of such
holomorphic functions. 
Think for instance of a convergent~$\hat\ph(\ze)$ which is the Taylor expansion
at~$0$ of a function holomorphic in $\C\setminus\Om$, where $\Om$ is a discrete
subset of~$\C^*$ (\eg a function which is meromorphic in~$\C$ and regular
at~$0$): in this case $\hat\ph$ has an analytic continuation in $\C\setminus\Om$
whereas, as a rule, its antiderivative $1*\hat\ph$ has only a multiple-valued
continuation there\ldots


\parag
We end this section with an example which is simple (because it deals with
explicit entire functions of~$\ze$) but useful:
\begin{lemma}
Let $p,q\in\N$ and $c\in\C$. Then
\begin{equation}	\label{eqconvolexppq}
\Big( \frac{\ze^p}{p!} \ee^{c\ze} \Big) *
\Big( \frac{\ze^q}{q!} \ee^{c\ze} \Big) =
\frac{\ze^{p+q+1}}{(p+q+1)!} \ee^{c\ze}.
\end{equation}
\end{lemma}

\begin{proof}
One could compute the convolution integral \eg by induction on~$q$, but one can
also notice that
$\frac{\ze^p}{p!} \ee^{c\ze}$
is the formal Borel transform of $T_{-c} z^{-p-1}$ (by virtue of
the second property in Lemma~\ref{lemelemptiescB}),
therefore the \lhs\ of~\eqref{eqconvolexppq} is the Borel transform of
$(T_{-c} z^{-p-1})(T_{-c} z^{-q-1}) = T_{-c} z^{-p-q-2}$.
\end{proof}


\newpage

\centerline{\Large\sc The Borel-Laplace summation along $\R^+$}
\addcontentsline{toc}{part}{\sc The Borel-Laplace summation along $\R^+$}

\vspace{.3cm}


\section{The Laplace transform}	\label{secLaplaceTrsf}


The Laplace transform of a function $\hat\ph \col \R^+ \to \C$ is the
function~$\cL^0\hat\ph$ defined by the formula
\begin{equation}	\label{eqdefcLzero}
(\cL^0\hat\ph)(z) = \int_0^{+\infty} \ee^{-z\ze} \hat\ph(\ze)\,\dd\ze.
\end{equation}
Here we assume~$\hat\ph$ continuous (or at
least locally integrable on~$\R^{*+}$ and integrable on $[0,1]$) and
\begin{equation}	\label{ineqhatphAc}
|\hat\ph(\ze)| \le A \, \ee^{c_0\ze}, \qquad \ze\ge1,
\end{equation}
for some constants $A>0$ and $c_0\in\R$, so that the above integral makes sense for
any complex number~$z$ in the half-plane 
\[
\Pi_{c_0} \defeq \{\, z\in\C\mid \RE z > c_0 \,\}.
\]
Standard theorems ensure that~$\cL^0\hat\ph$ is holomorphic in~$\Pi_{c_0}$
(because $|\ee^{-z\ze}| = \ee^{-\ze\RE z} \le \ee^{-c_1\ze}$ for any
$z\in\Pi_{c_1}$, hence, for any $c_1>c_0$, we can find $\Phi\col\R^+\to\R^+$
integrable and independent of~$z$ such that $|\ee^{-z\ze} \hat\ph(\ze)| \le
\Phi(\ze)$ and deduce that $\cL^0\hat\ph$ is holomorphic on~$\Pi_{c_1}$).

\begin{lemma}
For any $n\in\N$, $\cL^0\big(\frac{\ze^n}{n!}\big)(z) = z^{-n-1}$ on~$\Pi_0$.
\end{lemma}

\begin{proof}
The function $\cL^0\big(\frac{\ze^n}{n!}\big)$ is holomorphic in $\Pi_{c_0}$ for any
$c_0>0$, thus in~$\Pi_0$.
The reader can check by induction on~$n$ that $\int_0^{+\infty} \ee^{-s}
s^n\,\dd s = n!$ and deduce the result for $z>0$ by the change of variable $\ze
= s/z$, and then for $z\in\Pi_0$ by analytic continuation.
\end{proof}

In fact, for any complex number $\nu$ such that $\RE\nu>0$,
$\cL^0\big(\frac{\ze^{\nu-1}}{\Ga(\nu)}\big) = z^{-\nu}$ for $z\in\Pi_0$,
where~$\Ga$ is Euler's gamma function (see Section~\ref{sec:Stirling}).

We leave it to the reader to check

\begin{lemma}	\label{lemelemLapl}
Let $\hat\ph$ as above, $\ph\defeq\cL^0\hat\ph$ and $c\in\C$. 
Then each of the functions 
$-\ze\hat\ph(\ze)$, $\ee^{-c\ze}\hat\ph(\ze)$ 
or $1*\hat\ph(\ze) = \int_0^\ze \hat\ph(\ze_1)\,\dd\ze_1$
satisfies estimates of the form~\eqref{ineqhatphAc} and
\begin{itemize}
\item
$\cL^0(-\ze\hat\ph) = \dfrac{\dd\ph}{\dd z}$,
\item
$\cL^0(\ee^{-c\ze}\hat\ph) = \ph(z+c)$,
\item
$\cL^0(1*\hat\ph) = z\ii\ph(z)$,
\item
if moreover $\hat\ph$ is continuously derivable on~$\R^+$ with
$\frac{\dd\hat\ph}{\dd\ze}$ satisfying estimates of the
form~\eqref{ineqhatphAc}, then
$\cL^0\bigg(\dfrac{\dd\hat\ph}{\dd\ze}\bigg) = z\ph(z) - \hat\ph(0)$.
\end{itemize}
\end{lemma}

\begin{rem}
Assume that $\hat\ph \col \R^+ \to \C$ is bounded and locally integrable. Then
$\cL^0\hat\ph$ is holomorphic in $\{\RE z>0\}$.
If one assumes moreover that $\cL^0\hat\ph$ extends holomorphically to a
neighbourhood of $\{\RE z\ge0\}$, then the limit of
$\int_0^T \hat\ph(\ze)\,\dd\ze$
as $T\to\infty$ exists and equals~$(\cL^0\hat\ph)(0)$;
see \cite{NewZag} for a proof of this statement and its application to a
remarkably short proof of the Prime Number Theorem (less than three pages!).
\end{rem}


\section{The fine Borel-Laplace summation}	\label{secfineBL}


\parag
We shall be particularly interested in the Laplace transforms of functions that
are analytic in a neighbourhood of~$\R^+$ and that we view as analytic
continuations of holomorphic germs at~$0$.

\begin{Def}	\label{defcNc}
We call \emph{half-strip} any set of the form
$S_\de = \{\, \ze\in\C \mid \dist(\ze,\R^+) < \de \,\}$
with a $\de>0$.
For $c_0\in\R$, we denote by~$\cN_{c_0}(\R^+)$ the set consisting of all convergent
formal series~$\hat\ph(\ze)$ defining a holomorphic function near~$0$ which extends
analytically to a half-strip~$S_\de$ with
\[
|\hat\ph(\ze)| \le A \, \ee^{c_0|\ze|}, \qquad \ze\in S_\de,
\]
where $A$ is a positive constant
(we use the same symbol~$\hat\ph$ to denote the function in~$S_\de$ and the
power series which is its Taylor expansion at~$0$).
We also set
\[
\cN(\R^+) = \bigcup_{c_0\in\R} \cN_{c_0}(\R^+)
\]
(increasing union).
\end{Def}


\begin{thm}	\label{thmAsympLapl}
Let $\hat\ph \in \cN_{c_0}(\R^+)$, $c_0 \ge 0$.
Set $a_n \defeq \hat\ph\pp n(0)$ for every $n\in\N$
and $\ph = \cL^0\hat\ph$.
Then for any $c_1>c_0$ there exist $L,M>0$ such that
\begin{equation}	\label{ineqGevasexp}
|\ph(z) - a_0 z\ii - a_1 z^{-2} - \cdots - a_{N-1} z^{-N}|
\le L M^N N! |z|^{-N-1},
\qquad z\in\Pi_{c_1}, \; N\in\N.
\end{equation}
\end{thm}

\begin{proof}
Without loss of generality we can assume $c_0>0$.
Let $\de,A>0$ be as in Definition~\ref{defcNc}.
We first apply the Cauchy inequalities in the discs $D(\ze,\de)$ of radius~$\de$
centred on the points $\ze\in\R^+$:
\begin{equation}	\label{ineqCauchyRde}
|\hat\ph\pp n(\ze)| \le \frac{n!}{\de^n} \sup_{D(\ze,\de)} |\hat\ph|
\le n! \de^{-n} A' \ee^{c_0\ze},
\qquad \ze\in\R^+,\; n\in\N,
\end{equation}
where $A' = A\,\ee^{c_0\de}$.
In particular, the coefficient $a_N = \hat\ph\pp N(0)$ satisfies
\begin{equation}	\label{ineqcoeffaN}
|a_N| \le N! \de^{-N} A'
\end{equation}
for any $N\in\N$.
Let us introduce the function
\[
R(\ze) \defeq \hat\ph(\ze) - a_0 - a_1 \ze - \cdots - a_N \frac{\ze^N}{N!},
\]
which belongs to $\cN_{c_0}(\R^+)$ (because $c_0>0$) and has Laplace transform
\[
\cL^0 R(z) = \ph(z) - a_0 z\ii - a_1 z^{-2} - \cdots - a_N z^{-N-1}.
\]
Since $0=R(0)=R'(0)=\cdots=R\pp N(0)$, the last property in
Lemma~\ref{lemelemLapl} implies
$\cL^0 R(z) = z\ii \cL^0 R'(z) 
= z^{-2} \cL^0 R''(z)
= \cdots = z^{-N-1} \cL^0 R\pp{N+1}(z)$
and, taking into account $R\pp{N+1} = \hat\ph\pp{N+1}$, we end up with
\[
\ph(z) - a_0 z\ii - \cdots - a_{N-1} z^{-N}
= a_N z^{-N-1} + z^{-N-1} \cL^0 \hat\ph\pp{N+1}(z).
\]
For $z\in\Pi_{c_1}$, $|\cL^0(\ee^{c_0\ze})(z)| \le \frac{1}{\RE z - c_0} \le \frac{1}{c_1-c_0}$,
thus inequality~\eqref{ineqCauchyRde} implies that
$|\cL^0\hat\ph\pp{N+1}(z)| \le (N+1)! \de^{-N-1} \frac{A'}{c_1-c_0} %
\le N! (2/\de)^N \frac{A'}{\de(c_1-c_0)}$.
Together with~\eqref{ineqcoeffaN}, this yields the conclusion with $M=2/\de$ and
$L = A'\big(1+\frac{1}{\de(c_1-c_0)}\big)$.
\end{proof}


\parag
Here we see the link between the Laplace transform of analytic functions and the
formal Borel transform: the Taylor series at~$0$ of~$\hat\ph(\ze)$ is
$\sum a_n \frac{\ze^n}{n!}$, thus the finite sum which is subtracted from~$\ph(z)$ in the \lhs\
of~\eqref{ineqGevasexp} is nothing but a partial sum of the formal series
$\ti\ph(z) \defeq \cB\ii\hat\ph = \sum a_n z^{-n-1} \in \zcz_1$.

The connection between the formal series~$\ti\ph$ and the
function~$\ph$ which is expressed by~\eqref{ineqGevasexp} is a
particular case of a kind of asymptotic expansion, called {$1$-Gevrey
  asymptotic expansion}. Let us make this more precise:
\begin{Def}   \label{DefUnifAsExp}
Given $\gD \subset \C^*$ unbounded, a
function $\phi \col \gD \to \C$ and a formal series $\ti\phi(z) =
\sum_{n\ge0} c_n z^{-n} \in\C[[z\ii]]$,
we say that~\emph{$\phi$ admits~$\ti\phi$ as uniform asymptotic expansion
in~$\gD$} if there exists a sequence of positive numbers
$(K_N)_{N\in\N}$ such that
\beglabel{inequnifasexpP}
\abs*{ \phi(z) - c_0 - c_1 z\ii - \cdots - c_{N-1}z^{-(N-1)} } 
\le K_N \abs{z}^{-N},
\qquad z \in \gD, \quad N \in\N.
\elabel
We then use the notation
\[
\phi(z) \sim \ti\phi(z) 
\quad \text{uniformly for $z\in\gD$.}
\]
If there exist $L,M>0$ such that~\eqref{inequnifasexpP} holds with the
sequence $K_N = L M^N N!$, then we say that~\emph{$\phi$ admits~$\ti\phi$ as
uniform $1$-Gevrey asymptotic expansion in~$\gD$} and we use the
notation 
\[
\phi(z) \sim_1 \ti\phi(z) 
\quad \text{uniformly for $z\in\gD$.}
\]
\end{Def}

The reader is referred to \cite{ML} for more on asymptotic expansions.
As for now, we content ourselves with observing that, given~$\phi$
and~$\gD$,
\begin{enumerate}[--]
\item there can be at most one formal series~$\ti\phi$ such that
$\phi(z)\sim\ti\phi(z)$ uniformly for $z\in\gD$;
\item
if $\phi(z)\sim_1\ti\phi(z)$ uniformly for $z\in\gD$, then $\ti\phi\in\C[[z\ii]]_1$.
\end{enumerate}
(Indeed, if~\eqref{inequnifasexpP} holds, then the coefficients
of~$\ti\phi$ are inductively determined by
\[
c_N = \lim_{\substack{ \abs{z}\to\infty \\ z \in \gD}} z^N \rho_N(z),
\qquad
\rho_N(z) \defeq \phi(z) - \sum_{n=0}^{N-1} c_n z^{-n}
\]
because $\abs{\rho_N(z) - c_N z^{-N}} \le K_{N+1} z^{-N-1}$,
and it follows that $\abs{c_N} \le K_N$.)

Theorem~\ref{thmAsympLapl} can be rephrased as:
\begin{quote} \em
If $\hat\ph \in \cN_{c_0}(\R^+)$ with $c_0\ge0$, then the function
$\ph \defeq \cL^0\hat\ph$ (which is holomorphic in~$\Pi_{c_0}$) and the formal
series $\ti\ph \defeq \cB\ii\hat\ph$ (which belong to $\zcz_1$)
satisfy
\beglabel{eqnotasimun}
\ph(z) \sim_1 \ti\ph(z)
\quad \text{uniformly for $z\in\Pi_{c_1}$}
\elabel
for any $c_1>c_0$.
\end{quote}


\parag
Theorem~\ref{thmAsympLapl} can be exploited as a tool for ``resummation'':
if it is the formal series $\ti\ph(z)\in\zcz_1$ which is given in the first
place, we may apply the formal Borel transform to get $\hat\ph(\ze)
\in\C\{\ze\}$; if it turns out that~$\hat\ph$ belongs to the
subspace~$\cN(\R^+)$ of~$\C\{\ze\}$, then we can apply the Laplace transform and
get a holomorphic function~$\ph(z)$ which admits~$\ti\ph(z)$ as $1$-Gevrey
asymptotic expansion.
This process, which allows us to go from the formal series~$\ti\ph(z)$ to the
function $\ph = \cL^0\cB\ti\ph$, is called \emph{fine Borel-Laplace summation} (in the
direction of~$\R^+$).

The above proof of Theorem~\ref{thmAsympLapl} is taken from \cite{Mal}, in which
the reader will also find a converse statement
(see also Nevanlinna's theorem in \cite{ML}):
given $\ti\ph\in \zcz$,
the mere existence of a holomorphic function~$\ph$ which admits~$\ti\ph$
as uniform $1$-Gevrey asymptotic expansion in a half-plane of the form
$\Pi_{c_1}$ entails that $\cB\ti\ph\in\cN(\R^+)$; moreover, such a
holomorphic function~$\ph$ is then unique (we skip the proof of these
facts).
In this situation, the holomorphic function~$\ph(z)$ can be viewed as a kind of sum
of~$\ti\ph(z)$, although this formal series may be divergent,
and the formal series~$\ti\ph$ itself is said to be \emph{fine-summable in the direction
of~$\R^+$}.

If we start with a convergent formal series, say $\ti\ph(z) \in
z\ii\C\{z\ii\}$ supposed to be convergent for $|z|>c_0$, then the reader can
check that $\cB\ti\ph \in \cN_{c_1}(\R^+)$ for any $c_1>c_0$, thus $\ti\ph(z)$
is fine-summable and $\cL^0\cB\ti\ph$ is holomorphic in the
half-plane~$\Pi_{c_0}$. 
We shall see in Section~\ref{secvarydir} that $\cL^0\cB\ti\ph$ is nothing but the
restriction to~$\Pi_{c_0}$ of the ordinary sum of~$\ti\ph(z)$.


\parag
The formal series without constant term that are fine-summable in the direction
of~$\R^+$ clearly form a linear subspace of $\zcz_1$.
To cover the case where there is a non-zero constant term, we make use of the
convolution unit $\de=\cB 1$ introduced in Section~\ref{sec_convols}.
We extend the Laplace transform by setting $\cL^0\de \defeq 1$ and, more generally,
\[
\cL^0 (a\,\de+\hat\ph) \defeq a + \cL^0\hat\ph
\]
for a complex number~$a$ and a function~$\hat\ph$.

\begin{Def}    \label{DeffinesummRp}
A formal series of $\C[[z\ii]]$ is said to be \emph{fine-summable in the direction
of~$\R^+$} if it can be written in the form $\ti\ph_0(z) = a+\ti\ph(z)$ with $a\in\C$ and
$\ti\ph \in \cB\ii\big(\cN(\R^+)\big)$, \ie if its formal Borel transform
$\cB\ti\ph_0 = a\,\de+\hat\ph(\ze)$ belongs to the subspace $\C\,\de \oplus \cN(\R^+)$ of
$\C\,\de \oplus \C[[\ze]]$.
Its Borel sum is then defined as the function $\cL^0 (a\,\de+\hat\ph)$,
which is holomorphic in a half-plane~$\Pi_c$ (choosing $c\in\R$ large enough).

The operator of \emph{Borel-Laplace summation in the direction of~$\R^+$} is defined as
the composition $\gS^0\defeq\cL^0\circ\cB$
acting on all such formal series~$\ti\ph_0(z)$.
\end{Def}

Clearly, as a consequence of Theorem~\ref{thmAsympLapl} and
Definition~\ref{DeffinesummRp}, we have 
\begin{cor}   \label{corFineSumRp}
If $\ti\ph_0 \in \C[[z\ii]]$ is fine-summable in the direction
of~$\R^+$, then there exists $c>0$ such that the function
$\gS^0\ti\ph_0$ is holomorphic in~$\Pi_c$ and satisfies
\[
\gS^0\ti\ph_0(z) \sim_1 \ti\ph_0(z)
\quad \text{uniformly for $z\in\Pi_c$.}
\]
\end{cor}


\begin{rem}	\label{remnotmaxdom}
Beware that $\Pi_c$ is usually not the maximal domain of holomorphy of the Borel
sum $\gS^0\ti\ph_0$: it often happens that this function admits analytic
continuation in a much larger domain and, in that case, $\Pi_c$ may or may not
be the maximal domain of validity of the uniform $1$-Gevrey asymptotic expansion
property.
\end{rem}


\parag
We now indicate a simple result of stability under convolution:

\begin{thm}	\label{thmcNstabconv}
The space $\cN(\R^+)$ is a subspace of $\C\{\ze\}$ stable by convolution.
Moreover, if $c_0\in\R$ and $\hat\ph,\hat\psi\in\cN_{c_0}(\R^+)$, then 
$\hat\ph*\hat\psi \in\cN_{c_1}(\R^+)$ for every $c_1>c_0$ and
\begin{equation}	\label{eqLaplconvolprod}
\cL^0(\hat\ph*\hat\psi) = (\cL^0\hat\ph) (\cL^0\hat\psi)
\end{equation}
in the half-plane $\Pi_{c_0}$.
\end{thm}


\begin{cor}	\label{corfinesummsubalg}
The space $\C \oplus \cB\ii\big(\cN(\R^+)\big)$ of all fine-summable formal
series in the direction of~$\R^+$ is a subalgebra of~$\C[[z\ii]]$ which
contains the convergent formal series.
The operator of Borel-Laplace summation~$\gS^0$ satisfies
\begin{gather}	
\label{eqelemptiesBL}
\gS^0\left( \frac{\dd\ti\ph_0}{\dd z} \right) = \frac{\dd\,}{\dd z}\big(
\gS^0\ti\ph_0 \big),
\qquad
\gS^0\big( \ti\ph_0(z+c) \big) = (\gS^0\ti\ph_0)(z+c) \\[1ex]
\label{eqBLhomomalg}
\gS^0(\ti\ph_0\ti\psi_0) = (\gS^0\ti\ph_0) (\gS^0\ti\psi_0)
\end{gather}
for any $c\in\C$ and fine-summable formal series $\ti\ph_0$, $\ti\psi_0$.
\end{cor}


Later, we shall see that Borel-Laplace summation is also compatible with
the non-linear operation of composition of formal series.

\begin{proof}[Proof of Theorem~\ref{thmcNstabconv}]
Suppose $\hat\ph,\hat\psi \in \cN(\R^+)$, 
with~$\hat\ph$ holomorphic in a half-strip~$S_{\de'}$ in which
$|\hat\ph(\ze)| \le A'\,\ee^{c_0'|\ze|}$,
and~$\hat\psi$ holomorphic in a half-strip~$S_{\de''}$ in which
$|\hat\psi(\ze)| \le A''\,\ee^{c_0''|\ze|}$.
Let $\de = \min\{\de',\de''\}$ and $c_0 = \max\{c_0',c_0''\}$.

We write $\hat\chi(\ze) = \int_0^1 F(s,\ze) \,\dd s$ with 
$F(s,\ze) = \ze\hat\ph(s\ze)\hat\psi\big((1-s)\ze\big)$ 
and argue as in the proof of Lemma~\ref{lemconvoldisc}:
$F$ is continuous in~$s$ and holomorphic in~$\ze$ for $(s,\ze) \in [0,1]\times S_\de$,
with 
\begin{equation}	\label{ineqFAA}
|F(s,\ze)| \le |\ze| A'A'' \ee^{c_0' s |\ze| + c_0'' (1-s) |\ze|}
\le A'A'' |\ze| \ee^{c_0|\ze|}.
\end{equation}
In particular $F$ is bounded in $[0,1]\times C$ for any compact subset~$C$
of~$S_\de$, thus $\hat\chi$ is holomorphic
in~$S_\de$. Inequality~\eqref{ineqFAA} implies $|\hat\chi(\ze)| \le A'A'' |\ze|
\ee^{c_0|\ze|} = O\big( \ee^{c_1|\ze|} \big)$ for any $c_1>c_0$,
hence $\hat\chi\in\cN_{c_1}(\R^+)$.
The identity~\eqref{eqLaplconvolprod} follows from Fubini's theorem.
\end{proof}


\begin{proof}[Proof of Corollary~\ref{corfinesummsubalg}]
Let $\ti\ph_0 = a\,\de + \ti\ph$ and $\ti\psi_0 = b + \ti\psi$ with $a,b\in\C$
and $\ti\ph,\ti\psi\in\zcz$.
We already mentioned the fact that if $\ti\ph\in z\ii\C\{z\ii\}$ then $\ti\ph$
is fine-summable, thus $\ti\ph_0$ is fine-summable in that case.

Suppose $\ti\ph,\ti\psi \in \cB\ii\big( \cN(\R^+) \big)$.
Property~\eqref{eqelemptiesBL} follows from Lemmas~\ref{lemelemptiescB}
and~\ref{lemelemLapl}, since the constant~$a$ is killed by~$\frac{\dd\,}{\dd z}$
and left invariant by~$T_c$.
Since $\ti\ph_0\ti\psi_0 = ab + a\ti\psi + b\ti\ph + \ti\ph\ti\psi$ has formal
Borel transform $ab\,\de + a\hat\psi + b\hat\ph + \hat\ph*\hat\psi$,
Theorem~\ref{thmcNstabconv} implies that
$\ti\ph_0\ti\psi_0 \in \C \oplus \cB\ii\big(\cN(\R^+)\big)$ and,
since $\gS^0(ab) = ab$,
property~\eqref{eqBLhomomalg} follows by linearity from Lemma~\ref{lemconvoldisc} and
Theorem~\ref{thmcNstabconv} applied to $\cB\ti\ph*\cB\ti\psi$.
\end{proof}


\section{The Euler series}


The Euler series 
$\ti\Phi\eul(t) = \sum_{n\ge0} (-1)^n n! t^{n+1}$
is a classical example of divergent formal series. We write it ``at~$\infty$'' as 
\begin{equation}	\label{eqdefEulerseriesatinfty}
\ti\ph\eul(z) = \sum_{n\ge0} (-1)^n n! z^{-n-1}.
\end{equation}
Clearly, its Borel transform is the geometric series
\begin{equation}	\label{eqBorEulseries}
\hat\ph\eul(\ze) = \sum_{n\ge0} (-1)^n \ze^n = \frac{1}{1+\ze},
\end{equation}
which is convergent in the unit disc and sums to a meromorphic function.
The divergence of~$\ti\ph\eul(z)$ is reflected in the non-entireness
of~$\hat\ph\eul$, which has a pole at~$-1$
(\cf Lemma~\ref{lemCVcase}).

Observe that $\ti\Phi\eul(t)$ can be obtained as the unique formal solution to
a differential equation, the so-called Euler equation:
\[
t^2 \frac{\dd\ti\Phi}{\dd t} + \ti\Phi = t.
\]
With our change of variable $z=1/t$, the Euler equation becomes $-\pa\ti\ph +
\ti\ph = z\ii$; 
applying the formal Borel transform to the equation itself is an
efficient way of checking the formula for~$\hat\ph\eul(\ze)$: a formal series without
constant term~$\ti\ph$ is solution if and only if its Borel transform~$\hat\ph$
satisfies $(\ze+1)\hat\ph(\ze)=1$ (\cf Lemma~\ref{lemelemptiescB}) and,
since $1+\ze$ is invertible in the ring
$\C[[\ze]]$, the only possibility is $\hat\ph\eul(\ze) = (1+\ze)\ii$.

Formula~\eqref{eqBorEulseries} shows that $\hat\ph\eul(\ze)$ is holomorphic and
bounded in a neighbourhood of~$\R^+$ in~$\C$, hence $\hat\ph\eul \in
\cN_0(\R^+)$.
The Euler series is thus fine-summable in the direction of~$\R^+$ and has a
Borel sum $\ph\eul = \cL^0\cB\ti\ph\eul$ holomorphic in the half-plane
$\Pi_0 = \{\, \RE z >0 \,\}$.
The first part of~\eqref{eqelemptiesBL} shows that this function~$\ph\eul$ is a
solution of the Euler equation in the variable~$z$.

%
\begin{rem}
The series~$\ti\Phi\eul(t)$ appears in Euler's famous 1760 article \emph{De seriebus
  divergentibus}, in which Euler introduces it as a tool in one of his
methods to study the divergent numerical series 
\[
1 - 1! + 2! - 3! + \cdots,
\]
which he calls Wallis' series---see \cite{Barb} and \cite{gazRamisA}.
Following Euler, we may adopt $\ph\eul(1) \simeq 0.59634736\ldots$ as the numerical
value to be assigned this divergent series.
\end{rem}
%

The discussion of this example continues in Section~\ref{secEulStokes}; in
particular, we shall see how Borel sums can be defined in other half-planes than
the ones bisected by~$\R^+$ and that $\ph\eul$ admits an analytic continuation
outside~$\Pi_0$ (\cf Remark~\ref{remnotmaxdom}).


\vspace{1.2cm}

\centerline{\Large\sc $1$-summable formal series in an arc of directions}
\addcontentsline{toc}{part}{\sc $1$-summable formal series in an arc of directions}

\vspace{.3cm}


\section{Varying the direction of summation}	\label{secvarydir}


\parag
Let $\th\in\R$. By $\eith\R^+$ we mean the oriented half-line which can be
parametrised as $\{\, \xi\,\eith, \; \xi\in\R^+ \,\}$.
Correspondingly, we define the Laplace transform of a function $\hat\ph \col
\eith\R^+ \to \C$ by the formula
\begin{equation}	\label{eqdefcLth}
(\cL^\th\hat\ph)(z) = \int_0^{+\infty} \ee^{-z\xi\,\eith} \hat\ph(\xi\,\eith)\eith\,\dd\xi,
\end{equation}
with obvious adaptations of the assumptions we had at the beginning of
Section~\ref{secLaplaceTrsf}, 
in particular $|\hat\ph(\ze)| \le A \, \ee^{c_0|\ze|}$ for $\ze\in\eith[1,+\infty)$,
so that $\cL^\th\hat\ph$ is a well-defined
function holomorphic in a half-plane
\[
\Pi^\th_{c_0} \defeq \{\, z\in\C\mid \RE (z\,\eith) > c_0 \,\}.
\]
Since $\lan z,w \ran \defeq \RE(z\bar w)$ defines the standard real scalar product on
$\C \simeq \R\oplus\I\R$, we see that $\Pi^\th_{c_0}$ is the half-plane bisected by
the half-line $\ee^{-\I\th}\R^+$ obtained from $\Pi_{c_0} = \Pi^0_{c_0}$ by the
rotation of angle $-\th$.

The operator~$\cL^\th$ is the Laplace transform in the direction~$\th$; the
reader can check that it satisfies properties analogous to those explained in
Sections~\ref{secLaplaceTrsf} and~\ref{secfineBL} for~$\cL^0$.


\begin{Def}	\label{Deffinesumth}
  A formal series $\ti\ph_0(z)\in\C[[z\ii]]$ is said to be
  \emph{fine-summable in the direction~$\th$} if it can be written
$\ti\ph_0 = a + \ti\ph$ with $a\in\C$ and
$\ti\ph \in \cB\ii\big(\cN(\eith\R^+)\big)$, 
where the space $\cN(\eith\R^+)$ is defined by replacing $S_\de$ with
$S^\th_\de \defeq  \{\, \ze\in\C \mid \dist(\ze,\eith\,\R^+) < \de \,\}$
in Definition~\ref{defcNc}
(see Figure~\ref{figLaplTrsfMaj} on p.~\pageref{figLaplTrsfMaj}).
\end{Def}


The Laplace transform~$\cL^\th$ is well-defined in $\cN(\eith\R^+)$; we extend
it as a linear map on $\C\,\de \oplus \cN(\eith\R^+)$ by setting $\cL^\th\de
\defeq 1$ and define the \emph{Borel-Laplace summation operator} as the composition
\begin{equation}	\label{eqdefgSth}
\gS^\th \defeq \cL^\th \circ \cB
\end{equation}
acting on all fine-summable formal series in the direction~$\th$.
There is an analogue of Corollary~\ref{corFineSumRp}:
\begin{quote} \em
  If $\ti\ph_0 \in \C[[z\ii]]$ is fine-summable in the
  direction~$\th$, then there exists $c>0$ such that the function
  $\gS^\th\ti\ph_0$ is holomorphic in~$\Pi^\th_c$ and satisfies
\[
\gS^\th\ti\ph_0(z) \sim_1 \ti\ph_0(z)
\quad \text{uniformly for $z\in\Pi^\th_c$.}
\]
\end{quote}
There is also an analogue of Corollary~\ref{corfinesummsubalg}:
the product of two fine-summable formal series is fine-summable
and $\gS^\th$ satisfies properties analogous to~\eqref{eqelemptiesBL}
and~\eqref{eqBLhomomalg}.


\parag	\label{paragcNIga}
The case of a function~$\hat\ph$ holomorphic in a sector is of particular interest,
we thus give a new definition in the spirit of Definitions~\ref{defcNc}
and~\ref{Deffinesumth}, replacing half-strips by sectors:

\begin{Def}	\label{defcNga}
Let $I$ be an open interval of~$\R$ and $\ga\col I\to\R$ a
locally bounded function.\footnote{%
A function $\ga\col I\to\R$ is said to be locally bounded if any point~$\th$
of~$I$ admits a neighbourhood on which~$\ga$ is bounded. Equivalently, the
function is bounded on any compact subinterval of~$I$.
}
For any locally bounded function $\al \col I \to \R^+$, we denote by
$\cN(I,\ga,\al)$ the set consisting of all convergent formal series~$\hat\ph(\ze)$
defining a holomorphic function near~$0$ which extends analytically to
the open sector $\{\, \xi\,\eith \mid \xi>0,\; \th \in I \,\}$ 
and satisfies
\[
|\hat\ph(\xi\,\eith)| \le \al(\th) \, \ee^{\ga(\th)\xi}, \qquad \xi>0,\; \th \in I.
\]
We denote by $\cN(I,\ga)$ the set of all~$\hat\ph(\ze)$ for which there exists
a locally bounded function~$\al$ such that $\hat\ph\in\cN(I,\ga,\al)$. 
We denote by $\cN(I)$ the set of all~$\hat\ph(\ze)$ for which there exists
a locally bounded function~$\ga$ such that $\hat\ph\in\cN(I,\ga)$. 

\end{Def}


For example, in view of~\eqref{eqBorEulseries}, the Borel transform
$\hat\ph\eul(\ze)$ of the Euler series belongs to $\cN(I,0,\al)$ with
$I=(-\pi,\pi)$ and 
\[
\al(\th) = \begin{cases} 
\hfil 1 \hfil & \text{if $|\th|\le\pi/2$,} \\
1/|\sin\th| & \text{else.}
\end{cases}
\]

Clearly, if $\hat\ph \in \cN(I,\ga)$ and $\th\in I$, then $z \mapsto (\cL^\th\hat\ph)(z)$ is
defined and holomorphic in $\Pi^\th_{\ga(\th)}$.

\begin{lemma}	\label{lemcNgaalfinesum}
Let $\ga$ and~$I$ be as in Definition~\ref{defcNga}.
Then, for every $\th\in I$, there exists a number $c=c(\th)$ such that 
$\cN(I,\ga) \subset \cN_c(\eith\,\R^+)$;
one can choose~$c$ to be the supremum of~$\ga$ on an arbitrary neighbourhood
of~$\th$.
\end{lemma}

The proof is left as an exercise.


Lemma~\ref{lemcNgaalfinesum} shows that a~$\hat\ph$ belonging to $\cN(I,\ga)$
is the Borel transform of a formal series~$\ti\ph(z)$ which is fine-summable in
any direction $\th\in I$;
for each $\th\in I$, we get a function $\cL^\th\hat\ph$ holomorphic in the
half-plane $\Pi^\th_{\ga(\th)}$, with the property of uniform $1$-Gevrey asymptotic
expansion
\[
\cL^\th\hat\ph(z) \sim_1 \ti\ph(z)
\quad \text{uniformly for $z\in\Pi^\th_{\ga'(\th)}$,}
\]
where $\ga'(\th)>0$ is large enough to be larger than a local bound of~$\ga$.
We now show that these various functions match, at least if the length of~$I$ is
less than~$\pi$, so that we can glue them and define a Borel sum
of~$\ti\ph(z)$ holomorphic in the union of all the half-planes
$\Pi^\th_{\ga(\th)}$.


\begin{lemma}	\label{lemLaplmatch}
Suppose $\hat\ph \in \cN(I,\ga)$ with $\ga$ and~$I$ as in 
Definition~\ref{defcNga} and suppose 
\[
\th_1,\th_2\in I, \qquad 0 < \th_2-\th_1 < \pi.
\]
Then $\Pi^{\th_1}_{\ga(\th_1)} \cap \Pi^{\th_2}_{\ga(\th_2)}$ is a non-empty sector in
restriction to which the functions $\cL^{\th_1}\hat\ph$ and $\cL^{\th_2}\hat\ph$
coincide.
\end{lemma}


\begin{proof}
The non-emptiness of the intersection of the half-planes 
$\Pi^{\th_1}_{\ga(\th_1)}$ and $\Pi^{\th_2}_{\ga(\th_2)}$ is an elementary geometric
fact which follows from the assumption $0 < \th_2-\th_1 < \pi$: 
this set is the sector $\gD =
\{\, z_* + r\,\eith \mid r>0, \; \th\in(-\th_1-\frac{\pi}{2},-\th_2+\frac{\pi}{2}) \,\}$,
where $\{z_*\}$ is the intersection of the lines 
$\ee^{-\I\th_1}\big(\ga(\th_1)+\I\R\big)$ and 
$\ee^{-\I\th_2}\big(\ga(\th_2)+\I\R\big)$.

Let $\al\col I\to\R^+$ be a locally bounded function such that
$\hat\ph\in\cN(I,\ga,\al)$.
Let $c = \sup_{[\th_1,\th_2]}\ga$ and $A =\sup_{[\th_1,\th_2]}\al$ (both $c$
and~$A$ are finite by the local boundedness assumption).
By the identity theorem for holomorphic functions,
it is sufficient to check that $\cL^{\th_1}\hat\ph$ and $\cL^{\th_2}\hat\ph$
coincide on the set $\gD_1 = \Pi^{\th_1}_{c+1} \cap \Pi^{\th_2}_{c+1}$, since
$\gD_1$ is a non-empty sector contained in~$\gD$.

Let $z\in\gD_1$. We have $\RE(z\,\eith) > c+1$ for all $\th\in[\th_1,\th_2]$
(simple geometric property, or property of the superlevel sets of the cosine
function) thus, for any $\ze\in\C^*$,
\beglabel{argineqmod}
\arg\ze \in [\th_1,\th_2] \ens\Longrightarrow\ens
| \ee^{-z\ze} \hat\ph(\ze)| \le A \, \ee^{-|\ze|}.
\elabel
The two Laplace transforms can be written
\[
\cL^{\th_j}\hat\ph(z) = 
\int_0^{\ee^{\I\th_j}\infty} \ee^{-z\ze} \hat\ph(\ze)\,\dd\ze =
\lim_{R\to\infty} \int_0^{R\,\ee^{\I\th_j}} \ee^{-z\ze} \hat\ph(\ze)\,\dd\ze,
\qquad j=1,2,
\]
but, for each $R>0$, the Cauchy theorem implies 
\[
\left( \int_0^{R\,\ee^{\I\th_2}} - \int_0^{R\,\ee^{\I\th_1}} \right)
\ee^{-z\ze} \hat\ph(\ze)\,\dd\ze =
\int_C \ee^{-z\ze} \hat\ph(\ze)\,\dd\ze, 
\qquad C = \{\, R\,\eith \mid \th \in [\th_1,\th_2] \,\}
\]
and, by~\eqref{argineqmod}, this difference has a modulus $\le A R (\th_2-\th_1)
\ee^{-R}$, hence it tends to~$0$ as $R\to\infty$.
\end{proof}

\begin{figure}
\begin{center}

\includegraphics[scale=1]{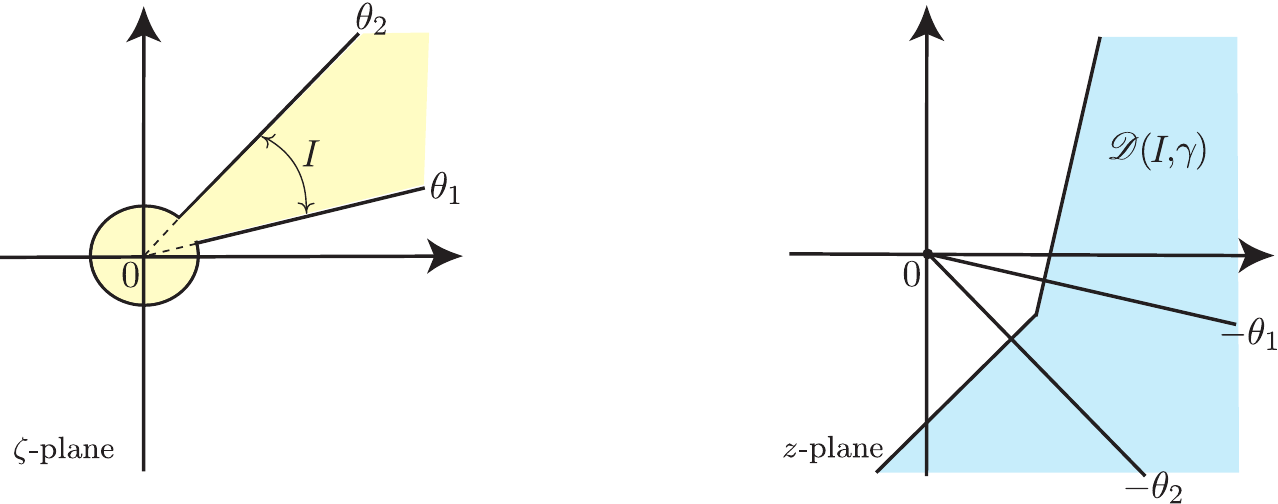} 

\bigskip 

\caption{%
\emph{$1$-summability in an arc of directions.}
Left: $\hat\ph(\ze)\in\cN(I,\ga)$ is holomorphic in the union of a disc and a
sector.
Right: the domain $\gD(I,\ga)$ where $\cL^I\hat\ph(z)$ is holomorphic.}
\label{fig:sectSig}

\end{center}
\end{figure}

\parag	\label{paraggDIga}
Lemma~\ref{lemLaplmatch} allows us to glue toghether the various Laplace
transforms:

\begin{Def}	\label{defglueLapl}
For $I$ open interval of~$\R$ of length $|I|\le\pi$ and $\ga\col I\to\R$ locally
bounded, we define
\[
\gD(I,\ga) = \bigcup_{\th\in I} \Pi^\th_{\ga(\th)},
\]
which is an open subset of~$\C$ (see Figure~\ref{fig:sectSig}),
and, for any $\hat\ph \in \cN(I,\ga)$, we define a function $\cL^I\hat\ph$
holomorphic in $\gD(I,\ga)$ by
\[
\cL^I\hat\ph(z) = \cL^\th\hat\ph(z) 
\quad\text{with $\th\in I$ such that $z\in\Pi^\th_{\ga(\th)}$}
\]
for any $z\in\gD(I,\ga)$.
\end{Def}

Observe that, for a given $z\in\gD(I,\ga)$, there are infinitely many possible
choices for~$\th$, which all give the same result by virtue of
Lemma~\ref{lemLaplmatch};
$\gD(I,\ga)$ is a ``sectorial neighbourhood of~$\infty$'' centred on the
ray $\arg z = -\th^*$ with aperture $\pi+|I|$,
where $\th^*$ denotes the midpoint of~$I$,
in the sense that, for every $\eps>0$, it contains a sector bisected by the
half-line of direction~$-\th^*$ with opening $\pi+|I|-\eps$
(see \cite{CNP}).


We extend the definition of the linear map~$\cL^I$ to $\C\,\de \oplus
\cN(I,\ga)$ by setting $\cL^I\de \defeq 1$.

\begin{Def}	\label{defgSI}
Given an open interval~$I$, we say that a formal series $\ti\ph_0(z) \in \C[[z\ii]]$ is
\emph{$1$-summable in the directions of~$I$} if 
$\cB\ti\ph_0 \in \C\,\de \oplus \cN(I)$.
%
%
The \emph{Borel-Laplace summation operator} is defined as the composition
\begin{equation}	\label{eqdefgsI}
\gS^I \defeq \cL^I\circ\cB
\end{equation}
acting on all such formal series, which produces functions holomorphic in
sectorial neighbourhoods of~$\infty$ of the form~$\gD(I,\ga)$, with locally
bounded functions $\ga\col I\to \R$.
\end{Def}

There is an analogue of Corollary~\ref{corfinesummsubalg}:
the product of two formal series which are $1$-summable in the directions of~$I$
is itself $1$-summable in these directions, as a consequence of
Lemma~\ref{lemcNgaalfinesum} and of the stability under multiplication of
fine-summable series, and the properties~\eqref{eqelemptiesBL}
and~\eqref{eqBLhomomalg} hold for the summation operator~$\gS^I$ too.

As for the property of asymptotic expansion, it takes the following
form: if $\ti\ph_0(z)$ is $1$-summable in the directions of~$I$, then there exists
$\ga\col I\to\R$ locally bounded such that 
\[
\text{$J$ relatively compact subinterval of~$I$}
\Imp
\gS^I\ti\ph_0(z) \sim_1 \ti\ph_0(z)
\quad \text{uniformly for $z \in \gD(J,\ga_{|J})$}
\]
(use Theorem~\ref{thmAsympLapl} and Lemma~\ref{lemcNgaalfinesum}).
We introduce the notation
\beglabel{ineqsimI}
\gS^I\ti\ph_0(z) \sim_1 \ti\ph_0(z)
\quad \text{for $z \in \gD(I,\ga)$}
\elabel
for this property, thus dropping the adverb ``uniformly''. 
Indeed we cannot claim that $\gS^I\ti\ph_0$ admits~$\ti\ph_0$ as
uniform $1$-Gevrey asymptotic expansion for $z \in \gD(I,\ga)$ (this
might simply be wrong for any locally bounded function $\ga \col I\to\R$):
uniform estimates are guaranteed only when restricting to relatively
compact subintervals.

The reader may check that the above definition of $1$-summability in an arc of
directions~$I$ coincides with the definition of $k$-summability in the
directions of~$I$ given in \cite{ML} when $k=1$.


\begin{rem}	\label{remStokesPhen}
  Suppose that $\ti\ph_0(z) \in \cB\ii\big(\C\,\de \oplus
  \cN(I,\ga)\big)$, so that the Borel sum $\ph_0(z) =
  \gS^I\ti\ph_0(z)$ is holomorphic in~$\gD(I,\ga)$ with the asymptotic
  property~\eqref{ineqsimI}.
Of course it may happen that~$\ti\ph_0$ is $1$-summable in the directions of an
interval which is larger than~$I$, in which case there will be an analytic
continuation for~$\ph_0$ with $1$-Gevrey asymptotic expansion in a sectorial
neighbourhood of~$\infty$ of aperture larger than $\pi+|I|$.
But even if it is not so it may happen that~$\ph_0$ admits analytic continuation
outside~$\gD(I,\ga)$. 

An interesting phenomenon which may occur in that case is the so-called Stokes
phenomenon: the asymptotic behaviour at~$\infty$ of the analytic continuation
of~$\ph_0$ may be totally different of what it was in the directions
of~$\gD(I,\ga)$, typically one may encounter oscillatory behaviour along the
limiting directions $-\th^* \pm \demi\big({\pi}+|I|\big)$ (where $\th^*$ is the
midpoint of~$I$) and exponential growth beyond these directions.
Examples can be found in Section~\ref{secEulStokes} (Euler series:
Remark~\ref{remEulStokes} and Exercise~\ref{exoEulStokes}) and
\S~\ref{parStirlStokes} (exponential of the Stirling series).
\end{rem}


\parag	\label{paragGrdSect}
What if $|I|>\pi$? 
First observe that, if $|I|\ge2\pi$, then $\cN(I)$ coincides with the set of
entire functions of bounded exponential type and the corresponding formal series
in~$z$ are precisely the convergent ones by Lemma~\ref{lemCVcase}:
\[
|I|\ge2\pi \quad\Longrightarrow\quad
\cB\ii\big( \C\,\de \oplus \cN(I) \big) = \C\{z\ii\}.
\]
This case will be dealt with in \S~\ref{parBorSumCV}. We thus suppose $\pi < |I|
< 2\pi$.

For $\hat\ph\in\cN(I,\ga)$, we can still define a family of holomorphic functions 
$\ph_\th \defeq \cL^\th\hat\ph$ holomorphic on $\pi_\th\defeq
\Pi^\th_{\ga(\th)}$ ($\th\in I$), 
with the property that 
$0<\th_2-\th_1<\pi \;\Longrightarrow\; \pi_{\th_1}\cap\pi_{\th_2}\neq\emptyset$
and $\ph_{\th_1}\equiv\ph_{\th_2}$ on $\pi_{\th_1}\cap\pi_{\th_2}$,
but the trouble is that also for $\pi < \th_2 - \th_1 < 2\pi$ is the
intersection of half-planes $\pi_{\th_1}\cap\pi_{\th_2}$ non-empty and then
nothing guarantees that $\ph_{\th_1}$ and $\ph_{\th_2}$ match on
$\pi_{\th_1}\cap\pi_{\th_2}$.

The remedy consists in lifting the half-planes $\pi_\th$ and their union
$\gD(I,\ga)$ to the Riemann surface of the logarithm $\Clog = \{\,
r\,\eel^{\I t} \mid r>0,\; t\in\R \,\}$
(see Section~\ref{secClog} for the definition of~$\Clog$ and the notation
$\eel^{\I t}$ which represents a point ``above'' the complex number~$\ee^{\I t}$).
For this, we suppose $\ga(\th)>0$, so that $\pi_\th$ is the set of
all complex numbers $z = r\,\ee^{\I t}$ with
$r>\ga(\th)$ and
$t\in\big(-\th-\arccos\frac{\ga(\th)}{r},-\th+\arccos\frac{\ga(\th)}{r}\big)$ 
(and adding any integer multiple of~$2\pi$ to~$t$ yields the same~$z$).
We set
\[
\ti\pi_\th \defeq \{\, z = r\,\eel^{\I t} \in \Clog \mid
r>\ga(\th), \;
t\in\big(-\th-\arccos\tfrac{\ga(\th)}{r},-\th+\arccos\tfrac{\ga(\th)}{r}\big)
\,\}, \qquad
\ti\gD(I,\ga) \defeq \bigcup_{\th\in I} \ti\pi_\th
\]
(this time $r\,\eel^{\I t}$ and $r\,\eel^{\I (t+2\pi m)}$ are regarded as
different points of~$\Clog$)
and consider $\ph_\th = \cL^\th\hat\ph$ as holomorphic in~$\ti\pi_\th$.
By gluing the various $\ph_\th$'s we now get a function which is holomorphic in
$\ti\gD(I,\ga) \subset \Clog$ and which we denote by~$\cL^I\hat\ph$.

The overlap between the half-planes $\pi_{\th_1}$ and~$\pi_{\th_2}$ for $\th_2 -
\th_1 > \pi$ is no longer a problem since their lifts $\ti\pi_{\th_1}$
and~$\ti\pi_{\th_2}$ do not intersect (they do not lie in the same sheet
of~$\Clog$) and~$\cL^I\hat\ph$ may behave differently on them.%
\footnote{Notice that $\cN(I,\ga) = \cN(2\pi+I,\ga)$, but the
  functions $\cL^\th\hat\ph$ and $\cL^{\th+2\pi}\hat\ph$ must now be
  considered as different: they are a priori defined in domains
  $\ti\pi_\th$ and $\ti\pi_{\th+2\pi}$ which do not intersect
  in~$\Clog$. Besides, it may happen that~$\cL^\th\hat\ph$ admit an
  analytic continuation in a part of $\ti\pi_{\th+2\pi}$ which does
  not coincide with $\cL^{\th+2\pi}\hat\ph$.}

\emph{Therefore, one can extend Definition~\ref{defgSI} to the case of an interval~$I$ of
length $>\pi$ and define $1$-summability in the directions of~$I$ and the
summation operator $\gS^I$ the same way, except that the Borel
sum~$\gS^I\ti\ph_0$ of a $1$-summable formal series~$\ti\ph_0$ is now a function
holomorphic in an open subset of the Riemann surface of the logarithm~$\Clog$.}


\parag	\label{parBorSumCV}
As already announced, the Borel sum of a convergent formal series coincides
with its ordinary sum:
\begin{lemma}   \label{lemSUMCVcase}
Suppose $\ti\ph_0\in \C\{z\ii\}$ and call~$\ph_0(z)$ the holomorphic function it
defines for $|z|$ large enough. Then $\ti\ph_0$ is $1$-summable in the directions
of any interval~$I$ and $\gS^I\ti\ph_0$ coincides with~$\ph_0$.
\end{lemma}

\begin{proof}
Let $\ti\ph_0 = a + \ti\ph$ with $a\in\C$ and $\ti\ph(z) = \sum a_n z^{-n-1}$,
so $\ph(z) = a + \sum a_n z^{-n-1}$ for $|z|$ large enough.
By Lemma~\ref{lemCVcase}, $\hat\ph = \cB\ti\ph$ is a convergent formal
series summing to an entire function and there exists $c>0$ such that
$\hat\ph\in\cN_{c}(\eith\,\R^+)$ for all $\th\in\R$.
Lemma~\ref{lemLaplmatch} allows us to glue together the Laplace transforms
$\cL^\th\hat\ph$: we get one function $\ph_*$ holomorphic in
$\bigcup_{\th\in\R} \Pi^\th_c = \{\,|z|>c\,\}$, with the asymptotic expansion
property $\ph_*(z) \sim_1 \ti\ph(z)$ uniformly for $\{\,|z|>c_1\,\}$
for any $c_1>c$.

The function $\Phi_*\col t\mapsto \ph_*(1/t)$ is thus holomorphic in the punctured disc 
$\{\, 0 < |t| < 1/c \,\}$. Inequality~\eqref{inequnifasexpP} with $N=0$ shows
that~$\Phi_*$ is bounded, thus the origin is a removable singularity and
$\Phi_*$ is holomorphic at $t=0$.
Now inequality~\eqref{inequnifasexpP} with $N=1,2,\ldots$ shows that $\sum a_n
t^{n+1}$ is the Taylor expansion at~$0$ of $\Phi_*(t)$, hence $a+\ph_*(1/t) \equiv
\ph_0(1/t)$. 
\end{proof}


\section{Return to the Euler series}
\label{secEulStokes}

%
As already mentioned (right after Definition~\ref{defcNga}), 
$\hat\ph\eul \in \cN(I,0)$ with $I=(-\pi,\pi)$.
We can thus extend the domain of analyticity of $\ph\eul = \cL^0\hat\ph\eul$, a
priori holomorphic in $\pi_0 = \{\, \RE z >0 \,\}$, by gluing the Laplace
transforms $\cL^\th\hat\ph\eul$, $-\pi<\th<\pi$, each of which is holomorphic in
the open half-plane $\pi_\th$ bisected by the ray of direction~$-\th$ and having
the origin on its boundary.
But if we take no precaution this yields a multiple-valued function:
there are two possible values for $\RE z < 0$, according as one uses $\th$ close
to~$\pi$ or to~$-\pi$.

\begin{figure}
\begin{center}
\includegraphics[scale=1]{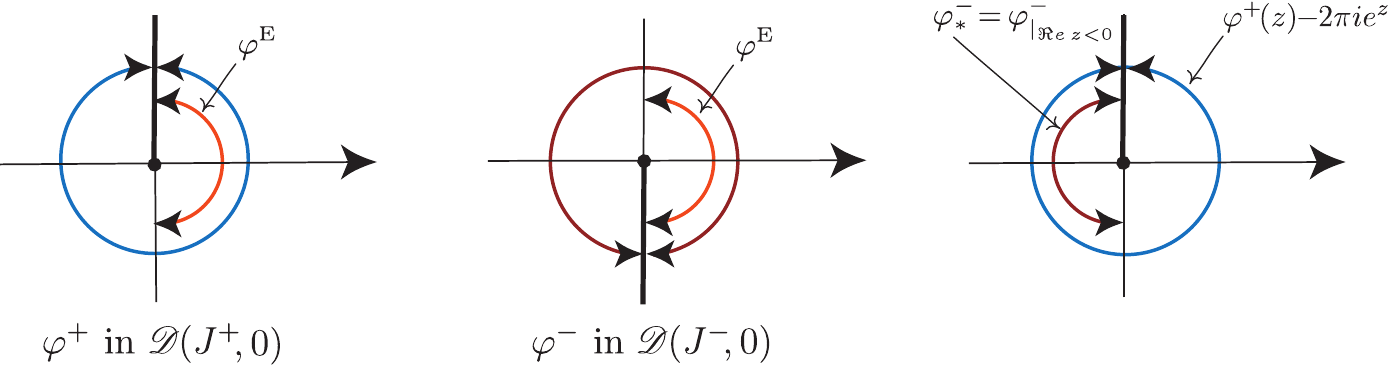} 

\bigskip 

\caption{%
\emph{Borel sums of the Euler series.}
Left and middle: $\ph^\pm$ extends $\ph\protect\eul$ in the cut plane $\gD(J^\pm,0)$.
Right: $\ph^+(z)-2\pi\I\,\ee^z$ extends $\ph^-_* = \ph^-_{|\{\RE z < 0\}}$
in the cut plane $\gD(J^+,0)$.}

\label{figEulStokes}

\end{center}
\end{figure}

A first way of presenting the situation consists in considering 
the subinterval $J^+ = (0,\pi)$, the Borel sum $\ph^+=\gS^{J^+}\ti\ph\eul$
holomorphic in $\gD(J^+,0) = \C\setminus\I\R^+$ which extends
analytically~$\ph\eul$ there,
and $J^- = (-\pi,0)$, $\ph^-=\gS^{J^-}\ti\ph\eul$ analytic continuation of~$\ph\eul$
in $\gD(J^-,0) = \C\setminus\I\R^-$.
See the first two parts of Figure~\ref{figEulStokes}.

The intersection of the domains $\C\setminus\I\R^+$ and $\C\setminus\I\R^-$ has
two connected components, the half-planes $\{\, \RE z >0 \,\}$ and $\{\, \RE z <0
\,\}$; both functions $\ph^+$ and $\ph^-$ coincide with~$\ph\eul$ on the former, whereas a
simple adaptation of the proof of Lemma~\ref{lemLaplmatch} involving Cauchy's
residue theorem yields
\begin{equation}	\label{eqmonodEul}
\RE z < 0 \ens\Longrightarrow\ens
\ph^+(z) - \ph^-(z) = 2\pi\I \,\ee^z.
\end{equation}
(This corresponds to the cohomological viewpoint presented in~\cite{ML}:
$(\ph^+,\ph^-)$ defines a $0$-cochain.)


Another way of putting it is to declare that $\ph\eul=\gS^I\ti\ph\eul$ is a holomorphic
function on 
\[
\ti\gD(I,0) = \{\, z \in \Clog \mid -\tfrac{3\pi}{2} < \arg z < \tfrac{3\pi}{2} \,\}
\]
(\cf Section~\ref{paragGrdSect})
and to rewrite~\eqref{eqmonodEul} as
\begin{equation}	\label{eqmonodEulbis}
\frac{\pi}{2} < \arg z < \frac{3\pi}{2}  \ens\Longrightarrow\ens
\ph\eul(z\,\eel^{-2\pi\I}) - \ph\eul(z) = 2\pi\I \,\ee^z.
\end{equation}


\begin{rem}	\label{remEulStokes}
[\textbf{Stokes phenomenon for~$\ph\eul$.}]
Let us consider the restriction~$\ph^-_*$ of the above function~$\ph^-$ to the
left half-plane $\{\, \RE z < 0 \,\}$.
Using~\eqref{eqmonodEul} we can write it as $\ph^+(z) - 2\pi\I\,\ee^z$, where
$\ph^+$ is holomorphic in an open sector bisected by~$\I\R^-$, namely the cut
plane $\gD^+=\C\setminus\I\R^+$, and the other term is an entire function: this
provides the analytic continuation of~$\ph^-_*$ through the cut~$\I\R^-$
to the whole of~$\gD^+$. See the third part of Figure~\ref{figEulStokes}.

Observe that $\ph^+\sim_1\ti\ph\eul$ in~$\gD^+$, in particular it tends to~$0$
at~$\infty$ along the directions contained in~$\gD^+$, while the
exponential~$\ee^z$ oscillates along~$\I\R^-$ and is exponentially growing in
the right half-plane: we see that, for~$\ph^-$, the asymptotic behaviour encoded
by~$\ti\ph\eul$ in the left half-plane breaks when we cross the limiting
direction~$\I\R^-$; the asymptotic behaviour of the analytic continuation is
oscillatory on~$\I\R^-$ (up to a correction which tends to~$0$) and after the
crossing we find exponential growth.

A similar analysis can be performed with $\ph^+_* = \ph^+_{|\{\RE z<0\}}$ when
one crosses $\I\R^+$, writing it as $\ph^-(z)+2\pi\I\,\ee^z$.
This is a manifestation of the Stokes phenomenon evoked in Remark~\ref{remStokesPhen}.
\end{rem}


\begin{exo}	\label{exoEulStokes}
Use~\eqref{eqmonodEulbis} to prove that $\ph\eul$ is the restriction to
$\ti\gD(I,0)$ of a function which is holomorphic in the whole of~$\Clog$.
(Hint: Show that the formula
$z\in\Clog \mapsto \ph(z) \defeq \ph\eul(z\,\eel^{-2\pi\I m}) - 2\pi\I m \,\ee^z$
if $m\in\Z$ and $\arg z \in (2\pi m - \frac{3\pi}{2},2\pi m + \frac{3\pi}{2})$
makes sense.)
In which sectors of~$\Clog$ is the Euler series asymptotic to this function?
\end{exo}


\begin{exo}
What kind of singularity has~$\ph\eul(z)$ when $|z|\to 0$?
(Hint: Find an elementary function $L(z)$ such that $L(z\,\eel^{-2\pi\I}) - L(z)
= -2\pi\I \,\ee^z$ and consider $\ph\eul+L$.)
\end{exo}


Observe that the Euler equation $-\frac{\dd\ph}{\dd z} + \ph = z\ii$ is a
non-homogeneous linear differential equation; the solutions of the associated
homogeneous equation are the functions $\la\,\ee^z$, $\la\in\C$.
By virtue of the general properties of the summation operator~$\gS^\th$, any
Borel sum of~$\ti\ph\eul$ is an analytic solution of the Euler equation.
In particular, the Borel sums~$\ph^+$ and~$\ph^-$ are solutions each in its own
domain of definition; 
on formula~\eqref{eqmonodEul} we can check that their restrictions to $\{\RE
z<0\}$ differ by a solution of the homogeneous equation, as should be. 
In fact, any two branches of the analytic continuation of~$\ph\eul$ differ by an
integer multiple of $2\pi\I\,\ee^z$.
Among all the solutions of the Euler equation, $\ph\eul$ can be characterised as
the only one which tends to~$0$ when $z\to\infty$ along a ray of direction
$\in (-\frac{\pi}{2},\frac{\pi}{2})$
(whereas, in the directions of $(\frac{\pi}{2},\frac{3\pi}{2})$, this is no
longer a distinctive property of~$\ph\eul$: all the solutions tend to~$0$ in
those directions!).


\begin{exo}
  How can one use the so-called method of ``variation of constants''
  to find directly an integral formula for the solution~$\ph\eul$ of
  the Euler equation?
\end{exo}


\section{The Stirling series} \label{sec:Stirling}


The Stirling series is a classical example of divergent formal series, which is
connected to Euler's gamma function.
The latter is the holomorphic function defined by the formula
\begin{equation}	\label{eqdefgamma}
\Ga(z) = \int_0^{+\infty} t^{z-1} \ee^{-t} \, \dd t
\end{equation}
for any $z\in\C$ with $\RE z > 0$ (so as to ensure the convergence of the
integral).
Integrating by parts, one gets the functional equation
\begin{equation}	\label{eqfunctgamma}
\Ga(z+1) = z \Ga(z).
\end{equation}
This equation provides the analytic continuation of~$\Ga$ for $z\in\C\setminus(-\N)$ in the form
\beglabel{eqmeromcontGa}
\Ga(z) = \frac{\Ga(z+n)}{z(z+1)\cdots(z+n-1)}
\elabel 
with any non-negative integer $n > -\RE z$; thus $\Ga$ is meromorphic in~$\C$
with simple poles at the non-positive integers.
Since $\Ga(1) = 1$, the functional equation also shows that 
\beglabel{eqGammaInterpolFact}
\Ga(n+1) = n!, \qquad n\in\N.
\elabel

Our starting point will be Stirling's formula for the restriction of~$\Ga$ to
the positive real axis:


\begin{lemma}
\begin{equation}	\label{eqStirlingf}
\left( \frac{x}{2\pi} \right)^{\demi} x^{-x} \ee^x \, \Ga(x) \xrightarrow[x\to+\infty]{} 1.
\end{equation}
\end{lemma}


\begin{proof}
This is an exercise in real analysis (and, as such, the following proof has nothing to do
with the rest of the text!). In view of the functional equation, it is sufficient to
prove that the function
\[
f(x) \defeq \frac{\Ga(x+1)}{x^{x+\demi} \ee^{-x}}
= \int_0^{+\infty} \frac{t^{x} \ee^{-t}}{x^{x} \ee^{-x}}  \frac{\dd t}{x^{1/2}}
\]
tends to $\sqrt{2\pi}$ as $x\to+\infty$.
The idea is that the main contribution in this integral arises for~$t$ close
to~$x$ and that, for $t = x + s$ with $s\to0$,
$\frac{t^{x} \ee^{-t}}{x^{x} \ee^{-x}} \sim \exp(-\frac{s^2}{2x})$
and 
$\int_{-x}^{+\infty} \exp(-\frac{s^2}{2x}) \,\frac{\dd s}{x^{1/2}} =
\int_{-\sqrt{x}}^{+\infty} \exp(-\frac{\xi^2}{2})\,\dd\xi$,
which converges to 
\begin{equation}	\label{eqGaussInt}
\int_{-\infty}^{+\infty} \ee^{-\xi^2/2}\,\dd\xi = \sqrt{2\pi}
\end{equation}
as $x\to+\infty$.
We now provide estimates to convert this into rigorous arguments.

We shall always assume $x\ge1$. 
The change of variable $t = x + \xi\sqrt{x}$ yields
\begin{equation}	\label{eqfexpg}
f(x) = \int_{-\infty}^{+\infty} \ee^{g(x,\xi)}\,\dd\xi,
\quad\text{with}\ens
g(x,\xi) \defeq \bigg( x \log\big(1 + \frac{\xi}{\sqrt{x}} \big) - \xi\sqrt{x} \bigg)
\mathds{1}_{\{\xi > -\sqrt{x}\}}.
\end{equation}
Integrating 
$\frac{1}{1+\sig} = 1  -\frac{\sig}{1+\sig} = 1 - \sig + \frac{\sig^2}{1+\sig}$,
we get
$\log(1+\tau) = \tau - \int_0^\tau \frac{\sig\,\dd\sig}{1+\sig} 
= \tau - \tau^2/2 + \int_0^\tau \frac{\sig^2\,\dd\sig}{1+\sig}$
for any $\tau>-1$,
whence
\begin{equation}	\label{eqggeomgeom}
g(x,\xi) = - x \int_0^{\xi/\sqrt{x}} \frac{\sig\,\dd\sig}{1+\sig} 
= - \frac{\xi^2}{2} + x \int_0^{\xi/\sqrt{x}} \frac{\sig^2\,\dd\sig}{1+\sig} 
\end{equation}
for any $\xi > -\sqrt{x}$.
Since $\int_0^\tau \frac{\sig^2\,\dd\sig}{1+\sig} = O(\tau^3)$ as $\tau\to 0$, the last
part of~\eqref{eqggeomgeom} shows that  
\[
g(x,\xi) \xrightarrow[x\to+\infty]{} -\xi^2/2 
\quad\text{for each $\xi\in\R$}.
\]
We shall use the first part of~\eqref{eqggeomgeom} to show that  
\begin{enumerate}[(i)]
\item \label{itemgche}
for $-\sqrt{x} < \xi \le0$, $g(x,\xi) \le -\xi^2/2$, 
whence $e^{g(x,\xi)} \le G_1(\xi) \defeq \ee^{-\xi^2/2}$;
\item \label{itemdte}
for $0\le \xi \le \sqrt{x}$, $g(x,\xi) \le -\xi^2/4$, 
whence $e^{g(x,\xi)} \le G_2(\xi) \defeq \ee^{-\xi^2/4}$;
\item \label{itemextrdte}
for $\xi \ge \sqrt{x}$, $g(x,\xi) \le -\xi/2$, 
whence $e^{g(x,\xi)} \le G_3(\xi) \defeq \ee^{-|\xi|/2}$.
\end{enumerate}
This is sufficient to conclude by means of Lebesgue's dominated convergence
theorem, since this will yield
$e^{g(x,\xi)} \le G_1(\xi) + G_2(\xi) + G_3(\xi)$ for all $x\ge1$ and $\xi\in\R$
and the function $G_1+G_2+G_3$ is independent of~$x$ and integrable on~$\R$,
thus \eqref{eqfexpg} implies 
$\dst f(x) \xrightarrow[x\to+\infty]{} \int_{-\infty}^{+\infty}
\lim_{x\to+\infty} \ee^{g(x,\xi)} \,\dd\xi$
and \eqref{eqGaussInt} yields the final result.

\medskip

\noindent-- Proof of~\eqref{itemgche}:
Assume $-\sqrt{x} < \xi \le0$. Changing $\sig$ into $-\sig$ and
integrating the inequality
$\frac{\sig}{1-\sig} \ge \sig$ over $\sig \in \big[0,|\xi|/\sqrt{x}\big]$,
we get
$g(x,\xi) = -x \int_0^{|\xi|/\sqrt{x}}\frac{\sig\,\dd\sig}{1-\sig} \le
-|\xi|^2/2$.

\medskip

\noindent-- Proof of~\eqref{itemdte}:
Assume $0\le \xi \le \sqrt{x}$, observe that
$\frac{\sig}{1+\sig} \ge \frac{\sig}{2}$ for $0\le \sig \le \xi/\sqrt{x}$
and integrate.

\medskip

\noindent-- Proof of~\eqref{itemextrdte}:
Assume $\xi \ge \sqrt{x} \ge 1$. Noticing that 
$\frac{\sig}{1+\sig} \ge \demi$ for $\sig\ge1$, we get 
$\int_0^{\xi/\sqrt{x}} \frac{\sig\,\dd\sig}{1+\sig} \ge 
\int_1^{\xi/\sqrt{x}} \frac{\sig\,\dd\sig}{1+\sig} \ge \frac{\xi}{2\sqrt{x}}$,
hence $g(x,\xi) \le - \demi\xi\sqrt{x} \le - \frac{\xi}{2}$.

\end{proof}


Observe that the \lhs\ of~\eqref{eqStirlingf} extends to a holomorphic function
in a cut plane:
\begin{equation}	\label{eqleftStirlingf}
\la(z) \defeq \frac{1}{\sqrt{2\pi}} z^{\demi-z} \ee^z \, \Ga(z),
\qquad z\in\C\setminus\R^-
\end{equation}
(using the principal branch of the logarithm~\eqref{eqdefLogppalbr} to define
$z^{\demi-z} \defeq \ee^{(\demi-z)\Log z}$;
in fact, $\la$ has a meromorphic continuation to the Riemann surface of the
logarithm~$\Clog$ defined in Section~\ref{secClog}).


\begin{thm}	\label{thmStirlings}
Let $I=(-\frac{\pi}{2}, \frac{\pi}{2})$. 
The above function~$\la$ can be written $\ee^{\gS^I\ti\mu}$,
%
%
where $\ti\mu(z) \in \zcz$ is a divergent odd formal series which is
$1$-summable in the directions of~$I$, whose formal Borel transform belongs to
$\cN(I,0)$ and is explicitly given by
\begin{equation}	\label{eqBorelStirling}
\hat\mu(\ze) = \ze^{-2}\left( \frac{\ze}{2}\coth\frac{\ze}{2} - 1 \right),
\qquad \ze\in \C\setminus(\De^+\cup\De^-)
\end{equation}
where $\De^\pm$ is the half-line $\pm2\pi\I [1,+\infty)$,
and whose Borel sum $\gS^I\ti\mu$ is holomorphic in the cut plane $\gD(I,0) = \C\setminus\R^-$.
\end{thm}

It is the formal series~$\ti\mu(z)$, the asymptotic expansion of $\log\la(z)$,
that is usually called the Stirling series.


\begin{exo}	\label{exoBernoulli}
Compute the Taylor expansion of the \rhs\ of~\eqref{eqBorelStirling}
in terms of the Bernoulli numbers~$B_{2k}$ defined by 
$\dst \frac{\ze}{\ee^\ze - 1} = 1 - \demi \ze + 
\sum_{k\ge1} \frac{B_{2k}}{(2k)!} \ze^{2k}$
(so $B_2 = 1/6$, $B_4 = -1/30$, $B_6 = 1/42$, etc.).
Deduce that
\beglabel{eqexplicitStirl}
\ti\mu(z) = \sum_{k\ge1} \frac{B_{2k}}{2k(2k-1)} z^{-2k+1}
= \frac{1}{12} z\ii - \frac{1}{360} z^{-3} + \frac{1}{1260} z^{-5} + \cdots.
\elabel
\end{exo}


We shall see in \S~\ref{parStirlStokes} that one can pass from~$\ti\mu$ to its
exponential and get an improvement of~\eqref{eqStirlingf} in the form of
\begin{cor}[Refined Stirling formula]	\label{corRefStirlF}
The formal series $\ti\la(z) \defeq \ee^{\ti\mu(z)}$ is $1$-summable in the
directions of $(-\frac{\pi}{2},\frac{\pi}{2})$ and its Borel sum is the
function~$\la$, with
\beglabel{eqrefineStirlingf}
\la(z) = \tfrac{1}{\sqrt{2\pi}} z^{\demi-z} \ee^z \, \Ga(z) 
\sim_1 \ti\la(z) = 1 + \sum_{n\ge0} g_n z^{-n-1}
\ens \text{uniformly for $|z|>c$ and $\arg z \in (-\be,\be)$}
\elabel
for any $c>0$ and $\be\in (0,{\pi})$,
with rationals $g_0,g_1,g_2,\ldots$ computable in terms of the
Bernoulli numbers:
\begin{align*}
g_0 & = \dem B_2 \\
g_1 & = \tfrac{1}{8} B_2^2 \\
g_2 & = \tfrac{1}{48} B_2^3 + \tfrac{1}{12} B_4 \\
g_3 & = \tfrac{1}{384} B_2^4 + \tfrac{1}{24} B_2 B_4 \\
g_4 & = \tfrac{1}{3840} B_2^5 + \tfrac{1}{96} B_2^2 B_4 + \tfrac{1}{30} B_6 \\
      & \hphantom{=} \vdots
\end{align*}
\end{cor}


Inserting the numerical values of the Bernoulli numbers,\footnote{%
\label{footextsimun}
and extending the notation ``$\sim_1$'' used in~\eqref{eqnotasimun}
or~\eqref{ineqsimI} by writing 
$F(z) \sim_1 G(z) \ti\ph_0(z)$ whenever 
$F(z)/G(z) \sim_1 \ti\ph_0(z)$
} 
we get
\beglabel{eqasymptGa}
\Ga(z) \sim_1 \ee^{-z} z^{z-\demi} \sqrt{2\pi} \Big( 1
+ \tfrac{1}{12} z\ii
+ \tfrac{1}{288} z^{-2} 
- \tfrac{139}{51840} z^{-3} 
- \tfrac{571}{2488320} z^{-4}
+ \tfrac{163879}{209018880} z^{-5}
+ \cdots \Big)
\elabel
uniformly in the domain specified in~\eqref{eqrefineStirlingf}.


\begin{proof}[Proof of Theorem~\ref{thmStirlings}]
\textbf{a)}
We first consider $\la(x) = \frac{1}{\sqrt{2\pi}} x^{\demi-x} \ee^x \, \Ga(x)$
for $x>0$. The functional equation~\eqref{eqfunctgamma} yields
\[
\la(x+1) = (1+x\ii)^{-\demi-x} \ee\,\la(x).
\]
Formula~\eqref{eqdefgamma} shows that, for $x>0$, $\Ga(x)>0$ thus also $\la(x)>0$
and we can define 
\begin{equation}	\label{eqdefmux}
\mu(x) \defeq \log \la(x), \qquad x>0.
\end{equation}
This function is a particular solution of the linear difference equation
\begin{equation}	\label{eqlindiffmu}
\mu(x+1) - \mu(x) = \psi(x),
\end{equation}
where $\psi(x) \defeq \log\big( (1+x\ii)^{-\demi-x} \ee \big) = 
1 - (\demi+x)\log(1+x\ii)$.
\medskip

\noindent\textbf{b)}
Using the principal branch of the logarithm~\eqref{eqdefLogppalbr}, holomorphic
in $\C\setminus\R^-$, we see that $\psi$ is the restriction to $(0,+\infty)$ of
a function which is holomorphic in $\C\setminus [-1,0]$:
\[
\psi(z) = -\demi \Log(1+z\ii) + z \big(z\ii - \Log(1+z\ii)\big).
\]
We observe that $\psi$ is holomorphic at~$\infty$ (\ie $t\mapsto\psi(1/t)$ is
holomorphic at the origin); moreover $\psi(z) = O(z^{-2})$ and its Taylor series
at~$\infty$ is
\[
\ti\psi(z) = \demi \ti L(z) + z \big(z\ii + \ti L(z)\big) \in z^{-2}\C\{z\ii\},
\qquad \ti L(z) \defeq - \sum_{n\ge1} \frac{(-1)^{n-1}}{n} z^{-n}.
\]
With a view to applying Corollary~\ref{coreqdifflin}, we compute the Borel
transform $\hat\psi = \cB\ti\psi$:
using $\hat L(\ze) = - \sum_{n\ge1} (-\ze)^{n-1}/n! = \ze\ii(\ee^{-\ze}-1)$
and the last property in Lemma~\ref{lemelemptiescB}, we get
\[
\hat\psi(\ze) = \demi \hat L(\ze) + \frac{\dd\,}{\dd\ze}(1+\hat L)
= \demi \ze\ii (\ee^{-\ze}-1) - \ze^{-2} (\ee^{-\ze}-1) - \ze\ii \ee^{-\ze}.
\]
\medskip

\noindent\textbf{c)}
Corollary~\ref{coreqdifflin} shows that the difference equation 
$\ti\ph(z+1) - \ti\ph(z) = \ti\psi(z)$
has a unique solution in $\zcz$, whose Borel transform is
\[
- \ze^{-2} + \demi \ze\ii - \ze\ii \frac{\ee^{-\ze}}{\ee^{-\ze}-1}
= \ze^{-2} \left( - 1 + \ze \bigg( \demi + \frac{1}{\ee^\ze-1} \bigg) \right)
= \hat\mu(\ze),
\]
where $\hat\mu(\ze)$ is defined by~\eqref{eqBorelStirling}. 
The formal series~$\hat\mu(\ze)$ is convergent and defines an even holomorphic
function which extends analytically to $\C\setminus(\De^+\cup\De^-)$ (in fact,
it even extends meromorphically to~$\C$, with simple poles on $2\pi\I\,\Z^*$).

\medskip

\noindent\textbf{d)}
Let us check that $\hat\mu \in \cN(I,0)$ with $I=(-\frac{\pi}{2},
\frac{\pi}{2})$.  
For $\th_0\in (0,\frac{\pi}{2})$, we shall bound~$|\hat\mu|$ in the sector
$\Sig = \{\, \xi\,\eith\mid \xi\ge0, \; \th\in[-\th_0,\th_0]\,\}$.
Let $\eps \defeq \min\{\pi,2\pi\cos\th_0\}$, so that $\Sig$ does not intersect
the discs $D(\pm 2\pi\I,\eps)$.
Since $\eps>0$, the number
\[
A(\eps) \defeq \sup \Big\{ \Big|\coth\frac{\ze}{2}\Big|,\;
\ze \in \C \setminus \bigcup_{m\in\Z} D(2\pi\I\, m,\eps) \Big\}
\]
is finite, because $\ze\mapsto\coth\frac{\ze}{2}$ is
$2\pi\I$-periodic, continuous in the closed set $\{\,|\IM\ze|\le\pi\,\}\setminus
D(0,\eps)$ and tends to~$\pm1$ as $\RE \ze \to \pm\infty$;~$A$ is in fact a
decreasing function of~$\eps$.
For $\ze\in\Sig\setminus D(0,1)$, we have 
$|\hat\mu(\ze)| \le \demi|\ze|\ii A(\eps) + |\ze|^{-2} \le A(\eps)+1$.
Since $\hat\mu$ is holomorphic in the disc $D(0,2\pi)$, the number
$B \defeq \sup \{ |\hat\mu(\ze)|,\; \ze\in D(0,1) \}$ is finite too, and we end
up with
\[
|\hat\mu(\ze)| \le \max\{A(\eps)+1,B\}, \qquad \ze\in\Sig,
\]
whence we can conclude $\hat\mu \in \cN(I,0,\al)$ with
$\al(\th) = \max\big\{ A\big(\eps(\th) \big)+1, B \big\}$, 
$\eps(\th) = \min\{\pi,2\pi|\cos\th|\}$.
\medskip 

\noindent\textbf{e)}
On the one hand, we have a solution $x\mapsto\mu(x)$ of
equation~\eqref{eqlindiffmu}: $\mu(x+1)-\mu(x)=\psi(x)$; this solution
is defined for $x>0$ and Stirling's formula~\eqref{eqStirlingf} implies that
$\mu(x)$ tends to~$0$ as $x\to+\infty$.

On the other hand, we have a formal solution $\ti\mu(z)$ to the equation 
$\ti\mu(z+1)-\ti\mu(z) = \ti\psi(z)$, which is $1$-summable, with a Borel sum
$\mu^+(z) \defeq \gS^I\ti\mu(z)$ holomorphic in $\gD(I,0)= \C\setminus\R^-$.
The property~\eqref{eqelemptiesBL} for the summation operator~$\gS^I$ implies
that
\[
\mu^+(z+1)-\mu^+(z) = \gS^I\ti\psi(z), 
\qquad z\in \C\setminus\R^-.
\]
But $\ti\psi$ is the convergent Taylor expansion of~$\psi$ at~$\infty$,
$\gS^I\ti\psi$ is nothing but the analytic continuation of~$\psi_{|(0,+\infty)}$.
The restriction of~$\mu^+$ to $(0,+\infty)$ is thus a solution to the same difference
equation~\eqref{eqlindiffmu}. 
Moreover, the $1$-Gevrey asymptotic property implies that 
$\mu^+(x)$ tends to~$0$ as $x\to+\infty$.

The difference $x\mapsto \De(x)\defeq\mu^+(x)-\mu(x)$ thus satisfies 
$\De(x+1)-\De(x)=0$ and it tends to~$0$ as $x\to+\infty$, hence $\De\equiv0$.
\end{proof}


\begin{rem}	\label{remnonodd}
Our chain of reasoning consisted in considering $\log\la_{|(0,+\infty)}$ and
obtaining its analytic continuation to $\C\setminus\R^-$ in the form
$\gS^I\ti\mu$. 
As a by-product, we deduce that the holomorphic function~$\la$ does not vanish on
$\C\setminus\R^-$ (being the exponential of a holomorphic function), hence the
function~$\Ga$ itself does not vanish on $\C\setminus\R^-$, nor does its
meromorphic continuation anywhere in the complex plane in view
of~\eqref{eqmeromcontGa}.

The formal series~$\ti\mu(z)$ is odd because $\hat\mu(\ze)$ is even and the
Borel transform~$\cB$ shifts the powers by one unit.
This does not imply that $\gS^I\ti\mu$ is odd! 
The direct consequence of the oddness of~$\ti\mu$ is rather the following:
$\ti\mu$ is $1$-summable in the directions of $J=(\frac{\pi}{2},
\frac{3\pi}{2})$ and the Borel sums $\mu^+=\gS^I\ti\mu$ and
$\mu^-=\gS^J\ti\mu$ are related by
\[
\mu^-(z) = - \mu^+(-z), \qquad z\in\C\setminus\R^+,
\]
because a change of variable in the Laplace integral yields $\cL^\th\hat\mu(z) =
- \cL^{\th+\pi}\hat\mu(-z)$.
The function~$\mu^-$ is in fact another solution of the difference
equation~\eqref{eqlindiffmu}.
\end{rem}


\begin{exo}	\label{exoStirlprepStokes}
\begin{enumerate}[$-$]
\item
With the notations of Remark~\ref{remnonodd}, prove that
\[
\mu^+(z) - \mu^-(z) = \sum_{m\ge1} \frac{1}{m} \ee^{-2\pi\I\,mz},
\qquad \IM z < 0
\]
by means of a residue computation (taking advantage of the existence of a
meromorphic continuation to~$\C$ for~$\hat\mu(\ze)$, with simple poles on
$2\pi\I\Z^*$, according to~\eqref{eqBorelStirling}).
\item
Deduce that, when we increase $\arg z$ above~$\pi$ or diminish it below~$-\pi$,
the function~$\mu^+(z)$ has a multiple-valued analytic continuation with
logarithmic singularities at negative integers.
\item
Deduce that $\la(z) = \frac{1}{(1-\ee^{-2\pi\I z})\la(-z)}$ for $\IM z < 0$, 
thus the restriction $\la_{|\{\IM z < 0\}}$
extends meromorphically to $\C\setminus\R^+$ with simple poles at the negative
integers.
\item
Compute the residue of this meromorphic continuation at a negative integer~$-k$
and check that the result is consistent with formula~\eqref{eqleftStirlingf} and
the fact that the residue of the simple pole of~$\Ga$ at~$-k$ is $(-1)^k/k!$.
(Answer: $-\frac{\I k^{k+\demi}\ee^{-k}}{k!\sqrt{2\pi}}$.)
\item
Repeat the previous computations with $\IM z>0$. Does one obtain the same
meromorphic continuation to $\C\setminus\R^+$ for $\la_{|\{\IM z > 0\}}$?
(Answer: no! But why?)
\item
Prove the reflection formula
\beglabel{eqreflectionGa}
\Ga(z) \Ga(1-z) = \frac{\pi}{\sin(\pi z)}.
\elabel
\end{enumerate}
\end{exo}


\begin{exo}
Using~\eqref{eqfunctgamma}, write a functional equation for the
logarithmic derivative $\psi(z) \defeq \Ga'(z)/\Ga(z)$.
Is there any solution of this equation in $\C[[z\ii]]$?
Using the principal branch of the logarithm~\eqref{eqdefLogppalbr} and
taking for granted that $\chi(z) \defeq \psi(z) - \Log z$ tends to~$0$
as~$z$ tends to~$+\infty$ along the real axis, show that $\chi(z)$ is the
Borel sum of a $1$-summable formal series (to be computed explicitly).
\end{exo}


\section{Return to Poincar\'e's example} \label{sec:ReturnPoinc}


In Section~\ref{secExPoin}, we saw Poincar\'e's example of a meromorphic
function~$\phi(t)$ of~$\C^*$ giving rise to a divergent formal
series~$\ti\phi(t)$ (formulas~\eqref{eqdefphiP} and~\eqref{eqdefformalphiP}).
There, $w = \ee^s$ was a parameter, with $|w|<1$, \ie $\RE s < 0$, and we had
\[
\phi(t) = \sum_{k\ge0} \frac{w^k}{1+kt}, \qquad
\ti\phi(t) = \sum_{n\ge0} a_n t^n
\]
with well-defined coefficients $a_n = (-1)^n b_n$ depending on~$s$.

To investigate the relationship between $\phi(t)$ and~$\ti\phi(t)$, we now set
\beglabel{eqdefphpoin}
\ph\poin(z) = z\ii\phi(z\ii) = \sum_{k\ge0} \frac{w^k}{z+k}, \qquad
\ti\ph\poin(z) = z\ii\ti\phi(z\ii) = \sum_{n\ge0} a_n z^{-n-1}
\elabel
(to place ourselves at~$\infty$ and get rid of the constant term)
so that $\ph\poin$ is a meromorphic function of~$\C$ with simple poles at
non-positive integers and $\ti\ph\poin(z)\in\zcz$.
The formal Borel transform $\hat\ph\poin(\ze)$ of $\ti\ph\poin(z)$ was already computed
under the name~$F(\ze)$ (\cf formula~\eqref{eqdefBorphiP} and the paragraph
which contains it):
\beglabel{eqhatphpoin}
\hat\ph\poin(\ze) = \frac{1}{1-\ee^{s-\ze}}.
\elabel
The natural questions are now: Is~$\ti\ph\poin$ $1$-summable in any arc of
directions and is~$\ph\poin$ its Borel sum?
We shall see that the answers are affirmative, with the help of a difference
equation:


\begin{lemma}	\label{lempoinuniqsol}
The function~$\ph\poin$ of~\eqref{eqdefphpoin} satisfies the functional equation
\beglabel{eqdiffpoin}
\ph(z) - w \ph(z+1) = z\ii.
\elabel
For any $z_0\in\C\setminus\R^-$, the restriction of~$\ph\poin$ to the half-line
$z_0+\R^+$ is the only bounded solution of~\eqref{eqdiffpoin} on this half-line.
\end{lemma}

\begin{proof}
We easily see that $w\ph\poin(z+1) = \sum \frac{w^{k+1}}{z+1+k} =
\ph\poin(z)-z\ii$ for any $z\in\C\setminus(-\N)$.
The boundedness of~$\ph\poin$ on the half-lines stems from the fact
that, for $z\in z_0+\R^+$ and $k\in\N$, 
$|z+k| \ge |\IM(z+k)| = |\IM z_0|$ and, if $\IM z_0=0$,
$|z+k| \ge z_0 > 0$,
hence, in all cases, $\big| \frac{w^k}{z+k} \big| \le A(z_0) |w|^k$ with $A(z_0)>0$
independent of~$z$.

As for the uniqueness: suppose $\ph_1$ and~$\ph_2$ are bounded functions on
$z_0+\R^+$ which solve~\eqref{eqdiffpoin}, then
$\psi\defeq\ph_2-\ph_1$ is a bounded solution of the equation 
$\psi(z) - w \psi(z+1) = 0$, 
which implies
$\psi(z) = w^n \psi(z+n)$ for any $z\in z_0+\R^+$ and $n\in\N$; we
get $\psi(z)=0$ by taking the limit as $n\to\infty$.
\end{proof}


But equation~\eqref{eqdiffpoin}, written in the form $\ph - w T_1\ph = z\ii$,
can also be considered in $\C[[z\ii]]$.

\begin{lemma}	\label{lempoinformsol}
The formal series~$\ti\ph\poin$ of~\eqref{eqdefphpoin} is the unique solution
of~\eqref{eqdiffpoin} in $\C[[z\ii]]$.
\end{lemma}

\begin{proof}
It is clear that the constant term of any formal solution of~\eqref{eqdiffpoin}
must vanish. We thus consider a formal series $\ti\ph(z) \in \zcz$.
Let us denote its formal Borel transform by $\hat\ph(\ze) \in \C[[\ze]]$;
in view of the second property of Lemma~\ref{lemelemptiescB}, $\ti\ph$ is
solution of~\eqref{eqdiffpoin} if and only if $(1-w\,\ee^{-\ze}) \hat\ph(\ze)
= 1$.
There is a unique solution because $1-w\,\ee^{-\ze}$ is invertible in
$\C[[\ze]]$ (recall that $w\neq1$ by assumption) and its Borel transform is
$(1-w\,\ee^{-\ze})\ii$, which according to~\eqref{eqhatphpoin} coincides
with~$\hat\ph^P(\ze)$ (recall that $w=\ee^s$).
\end{proof}

\begin{figure}
\begin{center}

\includegraphics[scale=1]{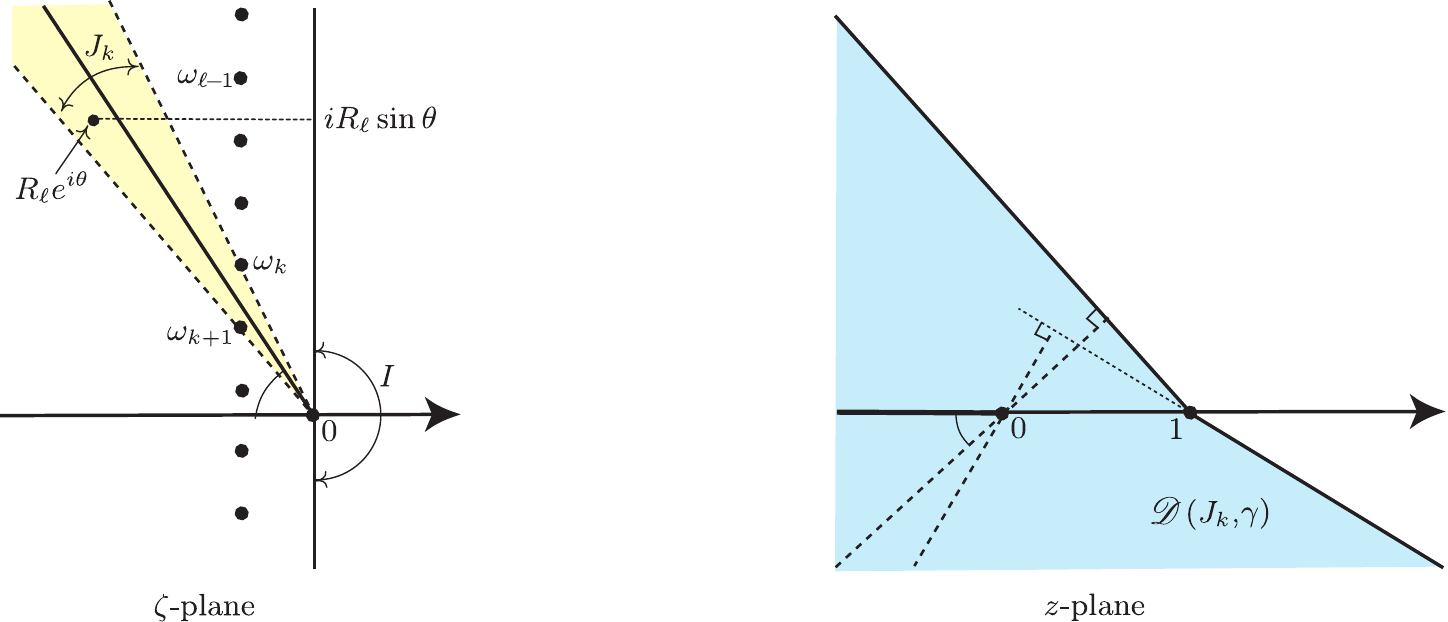}

\bigskip 

\caption{%
\emph{Borel-Laplace summation for Poincar\'e's example.}
%
%
}
\label{fig:figpoin}

\end{center}
\end{figure}

\begin{thm}	\label{thmSommaPoin}
The formal series~$\ti\ph\poin$ is $1$-summable in the directions of $I =
(-\frac{\pi}{2}, \frac{\pi}{2})$ and fine-summable in the directions
$\pm\frac{\pi}{2}$,
with $\hat\ph\poin \in \cN(I,0) \cap \cN_0(\I\R^+) \cap \cN_0(-\I\R^-)$.
Its Borel sum $\gS^I\ti\ph^P$ coincides with the function~$\ph\poin$ in
$\gD(I,0) = \C\setminus\R^-$.

Let $\om_k = s - 2\pi\I k$ for $k\in\Z$.
Then, for each $k\in\Z$, the formal series~$\ti\ph\poin$ is $1$-summable in the
directions of $J_k = (\arg\om_k,\arg\om_{k+1}) \subset (\frac{\pi}{2},
\frac{3\pi}{2})$,
with $\hat\ph\poin \in \cN(J_k,\ga)$, $\ga(\th) \defeq \cos\th$,
thus $\gD(J_k,\ga)$ is a sectorial neighbourhood of~$\infty$ containing the
real half-line $(-\infty,1)$ (see Figure~\ref{fig:figpoin}).
The Borel sum of~$\ti\ph\poin$ in the directions of~$J_k$ is a solution
of equation~\eqref{eqdiffpoin} which differs from~$\ph\poin$ by
\beglabel{eqdiffsumpoin}
\ph\poin(z) - \gS^{J_k}\ti\ph\poin(z) = 
2\pi\I \frac{\ee^{-\om_k z}}{1 - \ee^{-2\pi\I z}}
= - 2\pi\I \frac{\ee^{-\om_{k+1} z}}{1 - \ee^{2\pi\I z}}.
\elabel
\end{thm}


\begin{rem}
As a consequence of~\eqref{eqdiffsumpoin},
we rediscover the fact that~$\ph\poin$ not only is holomorphic in
$\C\setminus\R^-$ but also extends to a meromorphic function of~$\C$,
with simple poles at non-positive integers (because we can express it as the sum of 
$2\pi\I \frac{\ee^{-s z}}{1 - \ee^{-2\pi\I z}}$,
meromorphic on~$\C$, and $\gS^{J_0}\ti\ph\poin$, holomorphic in a sectorial
neighbourhood of~$\infty$ which contains~$\R^-$).
Similarly, each function $\gS^{J_k}\ti\ph\poin$ is meromorphic in~$\C$, with
simple poles at the positive integers.

In the course of the proof of formula~\eqref{eqdiffsumpoin}, it will be clear
that its \rhs\ is exponentially flat at~$\infty$ in the appropriate directions,
as one might expect since it has $1$-Gevrey asymptotic expansion reduced to~$0$.
This \rhs\ is of the form $\psi(z) = \ee^{-s z} \chi(z)$ with a $1$-periodic
function~$\chi$; it is easy to check that this is the general form of the
solution of the homogeneous difference equation $\psi(z) - w \psi(z+1) = 0$.
\end{rem}

The proof of Theorem~\ref{thmSommaPoin} makes use of

\begin{lemma}	\label{lemestimpoin}
Let $\sig \in (0,-\RE s)$ and $\de>0$. Then there exist $A = A(\sig)>0$ and $B =
B(\de)>0$ such that, for
any $\ze\in\C$,
\begin{gather}
\label{ineqhatpoindroit}
\RE\ze \ge -\sig \quad\Longrightarrow\quad
|\hat\ph\poin(\ze)| \le A, \\[1ex]
\label{ineqhatpoindist}
\dist(\ze,s+2\pi\I\Z) \ge \de \quad\Longrightarrow\quad
|\hat\ph\poin(\ze)| \le B \,\ee^{\RE\ze}.
\end{gather}
\end{lemma}


\begin{proof}[Lemma~\ref{lemestimpoin} implies Theorem~\ref{thmSommaPoin}]
Inequality~\eqref{ineqhatpoindroit} implies that
\[
\hat\ph\poin \in \cN(I,0) \cap \cN_0(\I\R^+) \cap \cN_0(-\I\R^-),
\]
whence the first summability statements follow.
Lemma~\ref{lempoinformsol} and the property~\eqref{eqelemptiesBL} for the
summation operator~$\gS^I$ imply that $\gS^I\ti\ph\poin$ is a solution
of~\eqref{eqdiffpoin};
this solution is bounded on the half-line $[1,+\infty)$, because of the
property~\eqref{ineqsimI} (in fact it tends to~$0$ on any half-line of the form
$z_0+\R^+$), thus it coincides with~$\ph\poin$ by virtue of
Lemma~\ref{lempoinuniqsol}.

Since $\RE\ze = \ga(\arg\ze)|\ze|$, 
inequality~\eqref{ineqhatpoindist} implies that
\[
\hat\ph\poin \in \cN(J_k,\ga,\al_k)
\]
with $\al_k \col \th \in J_k \mapsto B\big(\de_k(\th)\big)$,
%
%
$\de_k(\th) = \min\big\{ \dist(\om_k,\eith\R^+), 
\dist(\om_{k+1},\eith\R^+) \big\}$,
%
%
whence the $1$-summability in the directions of~$J_k$ follows.
Again, the Borel sum is a solution of the difference
equation~\eqref{eqdiffpoin}, a priori defined and holomorphic in
$\gD(J_k,\ga)$, which is the union of the half-planes $\Pi^\th_{\ga(\th)}$ for
$\th\in J_k$; one can check that each of these half-planes has the point~$1$ on
its boundary and that the intersection~$\gD$ of $\gD(J_k,\ga)$ with
$\C\setminus\R^-$ is connected. Thus, to conclude, it is sufficient to
prove~\eqref{eqdiffsumpoin} for $z$ belonging to one of the open subdomains
$\gD_1^+ \defeq \Pi^\th_{\ga(\th)+1}\cap \Pi^{\pi/2}_1$ or
$\gD_1^- \defeq \Pi^\th_{\ga(\th)+1}\cap \Pi^{-\pi/2}_1$,
with an arbitrary $\th\in J_k$ (none of them is empty).

Without loss of generality we can suppose $\th\neq\pi$. 
If $\th \in ( \frac{\pi}{2}, \pi)$, we proceed as follows:
for any integer $\ell\le k$, the horizontal line through the midpoint of
$(\om_\ell,\om_{\ell-1})$ cuts the half-lines $\eith\R^+$ and $\I\R^+$ in the
points $R_\ell\,\eith$ and $\I R_\ell\sin\th$, where $R_\ell$ is a positive real
number which tends to $+\infty$ as $\ell\to\infty$ (see Figure~\ref{fig:figpoin}).
Thus, for $z\in \gD_1^+$, we have 
\begin{align*}
\ph\poin(z) = \cL^{\pi/2}\hat\ph\poin(z) &= 
\lim_{\ell\to\infty} \int_0^{\I R_\ell\sin\th} \ee^{-z\ze} \hat\ph\poin(\ze)\,\dd\ze,
\\[1ex]
\gS^{J_k}\ti\ph\poin(z) = \cL^{\th}\hat\ph\poin(z) &= 
\lim_{\ell\to\infty} \int_0^{R_\ell\,\eith} \ee^{-z\ze} \hat\ph\poin(\ze)\,\dd\ze.
\end{align*}
Formula~\eqref{eqhatphpoin} shows that~$\hat\ph\poin$ is meromorphic, with
simple poles at the points~$\om_m$, $m\in\Z$, and residue $=1$ at each of these poles.
Cauchy's Residue Theorem implies that, for each $\ell\le k$, 
\beglabel{eqrespoin}
\left( \int_0^{\I R_\ell\sin\th} - \int_0^{R_\ell\,\eith} \right)
\ee^{-z\ze} \hat\ph\poin(\ze)\,\dd\ze =
2\pi\I \sum_{m=\ell}^k \ee^{-\om_m z}
+ \int_{L_\ell} \ee^{-z\ze} \hat\ph\poin(\ze)\,\dd\ze,
\elabel
where $L_\ell$ is the line segment $[R_\ell\,\eith,\I R_\ell\sin\th]$.
As in the proof of Lemma~\ref{lemLaplmatch}, we have 
\[
\arg\ze \in \left[ \tfrac{\pi}{2},\th \right] 
\quad\Longrightarrow\quad
|\ee^{-z\ze}| \le \ee^{-|\ze|(\ga(\th)+1)} = \ee^{-\RE\ze -|\ze|}
\]
(we have used $1 \ge \ga(\th)+1$), thus
$\ze\in L_\ell \,\Longrightarrow\,
|\ee^{-z\ze}\hat\ph\poin(\ze)| \le B(\pi) \, \ee^{-|\ze|} 
\le B(\pi) \, \ee^{-R_\ell\sin\th}$.
Hence the integral in the \rhs\ of~\eqref{eqrespoin} tends to~$0$ and we are
left with the geometric series $\ee^{-\om_k z} + \ee^{-\om_{k-1} z} + \cdots
= \ee^{-\om_k z} \sum_{n\ge0} \ee^{-2\pi\I n z}$
(since $-\om_m z = -\om_k z - 2\pi\I(k-m)z$), which
yields~\eqref{eqdiffsumpoin}.

If $\th \in (\pi, \frac{3\pi}{2})$, we rather take $\ell\ge k+1$ and
$z\in\gD_1^-$ and end up with
\begin{multline*}
\ph\poin(z) - \gS^{J_k}\ti\ph\poin(z) = 
\left( \int_0^{-\I\infty} - \int_0^{\eith\infty} \right)
\ee^{-z\ze} \hat\ph\poin(\ze)\,\dd\ze =\\
- 2\pi\I \sum_{m=k+1}^\infty \ee^{-\om_m z}
= - 2\pi\I \ee^{-\om_k z} \sum_{n\ge1} \ee^{2\pi\I n z},
\end{multline*}
which yields the same formula.
\end{proof}


\begin{proof}[Proof of Lemma~\ref{lemestimpoin}]
In view of~\eqref{eqhatphpoin}, for $\RE\ze\ge-\sig$ we have
$|\ee^{s-\ze}| \le \ee^{\sig+\RE s} <1$ and inequality~\eqref{ineqhatpoindroit}
thus holds with $A = (1-\ee^{\sig+\RE s})\ii$.

Formula~\eqref{eqhatphpoin} can be rewritten as $\hat\ph\poin(\ze) =
\frac{\ee^\ze}{\ee^\ze-\ee^s}$.
Let $C_\de \defeq \{\, \ze\in\C \mid \dist(\ze,s+2\pi\I\Z) \ge \de \,\}$ and
$F(\ze) \defeq |\ee^\ze-\ee^s|$.
The function~$F$ is $2\pi\I$-periodic and does not vanish on~$C_\de$;
since $F(\ze)$ tends to~$+\infty$ as $\RE\ze \to +\infty$
and to~$|w|$ as $\RE\ze \to -\infty$, we can find $R>0$ such that
$F(\ze) \ge |w|/2$ for $|\RE\ze| \ge R$, 
while $M\defeq \min\{\, F(\ze) \mid \ze\in C_\de,\, |\RE\ze| \le R, \, |\IM\ze|\le\pi \,\}$
is a well-defined positive number by compactness;
\eqref{ineqhatpoindist} follows with $B = \max\{ 2/|w|, 1/M \}$.
\end{proof}


\section{Non-linear operations with $1$-summable formal series}	\label{secnonlinsumm}


\parag
The stability under multiplication of the space of $1$-summable formal series
associated with an interval~$I$ was already mentioned (right after
Definition~\ref{defgSI}), but it is often useful to have more quantitative
information on what happens in the variable~$\ze$, which amounts to controlling
better the convolution products.

\begin{lemma}	\label{lemestimconv}
Suppose that $\th\in\R$ and  we are given locally integrable functions
$\hat\ph_1, \hat\ph_2 \col \eith\,\R^+\to\C$ and
$\Phi_1, \Phi_2 \col \eith\,\R^+\to\R^+$ 
such that 
\[
|\hat\ph_j(\ze)| \le \Phi_j\big(|\ze|\big), \qquad \ze\in\eith\,\R^+
\]
for $j=1,2$ and $\Phi_1$, $\Phi_2$ are integrable on~$[0,1]$.
Then the convolution products 
$\hat\ph_3 = \hat\ph_1*\hat\ph_2$ 
and $\Phi_3 = \Phi_1*\Phi_2$ 
defined by formula~\eqref{eqdefanalyticconvol} satisfy
\[
|\hat\ph_3(\ze)| \le \Phi_3\big(|\ze|\big), \qquad \ze\in\eith\,\R^+.
\]
\end{lemma}

\begin{proof}
Write~$\hat\ph_3(\ze)$ as $\int_0^1
\hat\ph_1(s\ze)\hat\ph_2\big((1-s)\ze\big)\ze\,\dd s$
and~$\Phi_3(\xi)$ as $\int_0^1
\Phi_1(s\xi)\Phi_2\big((1-s)\xi\big)\xi\,\dd s$.
\end{proof}


\begin{lemma}	\label{lemstarshapedstb}
Suppose $\De$ is an open subset of~$\C$ which is star-shaped \wrt~$0$ (\ie it is
non-empty and, for every $\ze\in\De$, the line segment $[0,\ze]$ is included
in~$\De$).
Suppose $\hat\ph_1$ and~$\hat\ph_2$ are holomorphic in~$\De$.
Then their convolution product (which is well defined since $0\in\De$) is also
holomorphic in~$\De$.
\end{lemma}

\begin{proof}
The function $(s,\ze)\mapsto \hat\ph_1(s\ze)\hat\ph_2\big((1-s)\ze\big)$
is continuous in~$s$, holomorphic in~$\ze$ and bounded in $[0,1]\times K$ for
any compact subset~$K$ of~$\De$.
\end{proof}


\parag
As an application, we show that $1$-summability is compatible with the
composition operator associated with a $1$-summable formal series and with
substitution into a convergent power expansion:

\begin{thm}	\label{thmcompatsumcirc}
Suppose $I$ is an open interval of~$\R$, $\ti\ph_0(z) = a + \ti\ph(z)$
and~$\ti\psi_0(z)$ are $1$-summable formal series in the directions of~$I$, with
$a\in\C$ and $\ti\ph(z)\in\zcz$, and $H(t)\in \C\{t\}$.
Then the formal series $\ti\psi_0\circ(\id+\ti\ph_0)$ and
$H\circ\ti\ph$ are $1$-summable in the directions of~$I$ and
\begin{equation}	\label{eqcompatsumcirc}
\gS^I\big( \ti\psi_0\circ(\id+\ti\ph_0) \big) = 
(\gS^I\ti\psi_0)\circ( \id + \gS^I\ti\ph_0 ),
\qquad
\gS^I( H\circ\ti\ph ) = H \circ \gS^I\ti\ph.
\end{equation}
More precisely, if $\cB\ti\ph\in\cN(I,\ga,\al)$ and
$\cB\ti\psi_0\in\C\,\de\oplus\cN(I,\ga)$
with $\al,\ga\col I\to\R$ locally bounded, $\al\ge0$,
and $\rho$ is a positive number smaller than the radius of convergence of~$H$,
then
\begin{align}	
\label{ineqestimtypecirc}
\cB\big( \ti\psi_0\circ(\id+\ti\ph_0) \big) &\in \C\,\de\oplus\cN(I,\ga_1),
& \hspace{-3em} \ga_1 &\defeq \ga+|a|+\sqrt{\al},\\
\cB( H\circ\ti\ph ) &\in \C\,\de\oplus\cN(I,\ga_2),
& \hspace{-3em} \ga_2 &\defeq \ga+\rho\ii\al,
\end{align}
\beglabel{ineqinclusions}
z \in \gD(I,\ga_1) \ens\Longrightarrow\ens z + \gS^I\ti\ph_0(z) \in \gD(I,\ga),
\qquad
z \in \gD(I,\ga_2) \ens\Longrightarrow\ens |\gS^I\ti\ph(z)| \le \rho
\elabel
and the identities in~\eqref{eqcompatsumcirc} hold in $\gD(\ga_1,I)$ and
$\gD(\ga_2,I)$ respectively.
\end{thm}

\begin{proof}
By assumption, $\hat\ph = \cB\ti\ph \in \cN(I,\ga,\al)$. The
properties~\eqref{ineqinclusions} are easily obtained as a consequence of
\beglabel{ineqgSImapPi}
z \in \Pi^\th_{\ga_j(\th)} \quad\Longrightarrow\quad
\gS^I\ti\ph(z) = \cL^\th\hat\ph(z)
\ens\text{and}\ens
|\cL^\th\hat\ph(z)| \le \frac{\al(\th)}{\ga_j(\th) - \ga(\th)}
\elabel
for any $\th\in I$ and $j=1,2$.

Let $\ti\psi_0(z) = b + \ti\psi(z)$ and $H = c + h(t)$ with $b,c\in\C$ and
$\ti\psi(z) \in \zcz$, $h(t) \in t\C\{t\}$, so that
\begin{align}
\ti\psi_0\circ(\id+\ti\ph_0) &= b + \ti\la, \qquad 
& \hspace{-5em} \ti\la &\defeq \ti\psi\circ(\id+\ti\ph_0), \\
H\circ\ti\ph &= c + \ti\mu, \qquad 
& \hspace{-5em} \ti\mu &\defeq h\circ\ti\ph.
\end{align}
We recall that $\ti\la$ and~$\ti\mu$ are defined by the formally convergent series of formal series
\begin{equation}	\label{eqdeftilactimu}
\ti\la = \sum_{k\ge0} \frac{1}{k!} (\pa^k\ti\psi) (\ti\ph_0)^k,
\qquad
\ti\mu = \sum_{k\ge1} h_k \ti\ph^k,
\end{equation}
where we use the notation $h(t) = \sum_{k\ge1} h_k t^k$.

Correspondingly, in we have formally convergent series of formal series in
$\C[[\ze]]$:
for instance, the Borel transform of~$\ti\mu$ is
\begin{equation}	\label{eqserhatmu}
\hat\mu = \sum_{k\ge1} h_k \hat\ph^{*k},
\quad
\text{where} \ens
\hat\ph^{*k} = \underbrace{\hat\ph * \cdots * \hat\ph}_\text{$k$ factors} 
\in \ze^{k-1}\C[[\ze]].
\end{equation}
But the series in the \rhs\ of~\eqref{eqserhatmu} can be viewed as a series of
holomorphic functions, since $\hat\ph$ is holomorphic in the union of a disc
$D(0,R)$ and of the sector $\Sig = \{\, \xi\,\eith \mid \xi>0,\; \th \in I \,\}$:
the open set $D(0,R)\cup\Sig$ is star-shaped \wrt~$0$, thus
Lemma~\ref{lemstarshapedstb} applies and each $\hat\ph^{*k}$ is holomorphic in 
$D(0,R)\cup\Sig$.
We shall prove the normal convergence of this series of functions in each
compact subset of $D(0,R)\cup\Sig$ and provide appropriate bounds.

Choosing $R>0$ smaller than the radius of convergence of~$\hat\ph$, we have
\begin{align*}
|\hat\ph(\ze)| &\le A, & \hspace{-5em} \ze &\in D(0,R), \\[1ex]
|\hat\ph(\ze)| &\le \Phi_\th(\xi)\defeq \al(\th)\,\ee^{\ga(\th)\xi}, 
& \hspace{-5em} \ze &\in \Sig,
\end{align*}
with a positive number~$A$, using the notations $\xi=|\ze|$ and $\th=\arg\ze$ in
the second case.
The computation of $\Phi_\th^{*k}(\xi)$ is easy, since $\Phi_\th$ can be viewed
as the restriction to~$\R^+$ of the Borel transform of 
$\al(\th) T_{-\ga(\th)}(z\ii)$; Lemma~\ref{lemestimconv} thus yields
\begin{align}
\label{ineqhatphkD}
|\hat\ph^{*k}(\ze)| &\le A^k \frac{\xi^{k-1}}{(k-1)!},
& \hspace{-4em} \ze &\in D(0,R), \\[1ex]
\label{ineqhatphkS}
|\hat\ph^{*k}(\ze)| &\le \Phi_\th^{*k}(\xi) = \al(\th)^k
\frac{\xi^{k-1}}{(k-1)!} \ee^{\ga(\th)\xi},
& \hspace{-4em} \ze &\in \Sig.
\end{align}
These inequalities, together with the fact that there exists $B>0$ such that
$|h_k| \le B \rho^{-k}$ for all $k\ge1$ (because $\rho$ is smaller than the
radius of convergence of~$H$), imply that the series of functions $\sum h_k
\hat\ph^{*k}$ is uniformly convergent in every compact subset of $D(0,R)\cup\Sig$;
the sum of this series is a holomorphic function whose Taylor coefficients
at~$0$ coincide with those of~$\hat\mu$, hence $\hat\mu(\ze) \in \C\{\ze\}$ and
$\hat\mu$ extends analytically to $D(0,R)\cup\Sig$.

Inequalities~\eqref{ineqhatphkS} also show that, for $\ze\in\Sig$,
\[ | h_k \hat\ph^{*k}(\ze) | \le \al(\th)B\rho\ii 
\frac{\big( \rho\ii\al(\th)\xi \big)^{k-1}}{(k-1)!} \ee^{\ga(\th)\xi},\]
hence $|\hat\mu(\ze)| \le
\al(\th)B\rho\ii \exp{\big( (\ga(\th) + \rho\ii\al(\th))\xi \big)}$,
\ie $\hat\mu \in \cN(I,\ga + \rho\ii\al)$. The dominated
convergence theorem shows that, for each $\th\in I $ and
$z\in\Pi^\th_{\ga_2(\th)}$, $\cL^\th\hat\mu(z)$ coincides with the convergent sum of the
series
$\sum h_k (\cL^\th\hat\ph^{*k})(z) = \sum h_k \big(\cL^\th\hat\ph(z)\big)^k$,
which is 
$h\big(\cL^\th\hat\ph(z)\big)$,
whence $\gS^I\ti\mu(z) \equiv h\big( \gS^I\ti\ph(z) \big)$.

We now move on to the case of~$\ti\la$.
Without loss of generality we can suppose that $a=0$, \ie that there is no
translation term in~$\ti\ph_0$, since 
$\ti\la = (T_a\ti\psi)\circ (\id+\ti\ph)$, 
thus it will be sufficient to apply the translationless case
of~\eqref{eqcompatsumcirc} and~\eqref{ineqestimtypecirc} to $T_a\ti\psi \in
\cB\ii\big(\cN(I,\ga+|a|)\big)$:
the identity~\eqref{eqelemptiesBL} for~$\gS^I$ will yield 
$\gS^I \Big( (T_a\ti\psi)\circ (\id+\ti\ph) \Big) =
(\gS^I T_a\ti\psi)\circ (\id+\gS^I\ti\ph) =
(\gS^I\ti\psi)\circ (\id+a)\circ(\id+\gS^I\ti\ph) =
(\gS^I\ti\psi)\circ (\id+a+\gS^I\ti\ph)$.

When $a=0$,
in view of~\eqref{eqdeftilactimu} and the first property in Lemma~\ref{lemelemptiescB},
the formal series $\hat\la \defeq \cB\ti\la \in \C[[\ze]]$ is given by the
formally convergent series of formal series
\[
\hat\la = \sum_{k\ge0} \hat\chi_k,
\qquad
\hat\chi_k \defeq \frac{1}{k!} \big((-\ze)^k\hat\psi\big) * \hat\ph^{*k}.
\]
We now view the \rhs\ as a series of holomorphic functions. Diminishing~$R$ if
necessary so as to make it smaller than the radius of convergence of~$\hat\psi$
and taking $\al'\col I\to\R^+$ locally bounded such that
$\hat\psi\in\cN(I,\ga,\al')$, we can find $A'>0$ such that
\begin{align*}
|\hat\psi(\ze)| &\le A', & \hspace{-5em} \ze &\in D(0,R), \\[1ex]
|\hat\psi(\ze)| &\le \Psi_\th(\xi)\defeq \al'(\th)\,\ee^{\ga(\th)\xi}, 
& \hspace{-5em} \ze &\in \Sig.
\end{align*}
Lemma~\ref{lemestimconv} and~\ref{lemstarshapedstb} show that the $\hat\chi_k$'s
are holomorphic in $D(0,R)\cup\Sig$ and satisfy
\begin{align}
\label{ineqhatchikD}
|\hat\chi_k(\ze)| &\le A' \frac{\xi^{k}}{k!} * A^k \frac{\xi^{k-1}}{(k-1)!}
= A' A^k \frac{\xi^{2k}}{(2k)!},
& \hspace{-3em} \ze &\in D(0,R), \\[1ex]
%
%
|\hat\chi_k(\ze)| &\le \Big(\frac{\xi^{k}}{k!}\Psi_\th\Big)*\Phi_\th^{*k}(\xi)
= \al'(\th)\al^k(\th) \frac{\xi^{2k}}{(2k)!} \ee^{\ga(\th)\xi},
& \hspace{-3em} \ze &\in \Sig
\end{align}
(we used~\eqref{ineqhatphkD},~\eqref{ineqhatphkS} and~\eqref{eqconvolexppq}).
The series $\sum\hat\chi_k$ is thus uniformly convergent in the compact subsets
of $D(0,R)\cup\Sig$ and sums to a holomorphic function, whose Taylor series
at~$0$ is~$\hat\la$.
Hence we can view $\hat\la$ as a holomorphic function and the last
inequalities imply that
$|\hat\la(\ze)| \le \al'(\th) \cosh\big(\sqrt{\al(\th)}\xi\big)\ee^{\ga(\th)\xi}
\le \al'(\th) \ee^{(\sqrt{\al(\th)}+\ga(\th))\xi}$
for $\ze\in\Sig$.
This yields $\hat\la \in \cN(I,\ga+\sqrt{\al})$ and, since  
$\cL^\th\hat\chi_k = \frac{1}{k!} 
\big(\big(\frac{\dd\,}{\dd z}\big)^k \cL^\th\hat\psi\big)
( \cL^\th\hat\ph )^k$
(use the first property in Lemma~\ref{lemelemLapl} and the
identity~\eqref{eqLaplconvolprod} for~$\cL^\th$), 
the dominated convergence theorem yields $\gS^I\ti\la =
(\gS^I\ti\psi)\circ( \id + \gS^I\ti\ph )$.
\end{proof}


\begin{exo}
Prove the following multivariate version of the result on substitution in a
convergent series:
suppose that $r\ge1$, $H(t_1,\ldots,t_r) \in \C\{t_1,\ldots,t_r\}$, $I$ is an
open interval of~$\R$ and $\ti\ph_1(z),\ldots,\ti\ph_r(z)\in\zcz$ are
$1$-summable in the directions of~$I$;
then the formal series 
\[ 
\ti\chi(z) \defeq H\big(\ti\ph_1(z),\ldots,\ti\ph_r(z)\big)
\] 
is $1$-summable in the directions
of~$I$ and $\gS^I\ti\chi = H \circ (\gS^I\ti\ph_1,\ldots,\gS^I\ti\ph_r)$.
\end{exo}


\parag	\label{parStirlStokes}
\emph{Proof of Corollary~\ref{corRefStirlF}.}
As a consequence of Theorem~\ref{thmcompatsumcirc}, using $H(t)=\ee^t$, we obtain the
$1$-summability in the directions of $I=(-\frac{\pi}{2}, \frac{\pi}{2})$ of the
exponential~$\ti\la$ of the Stirling series~$\ti\mu$, whence the refined Stirling
formula~\eqref{eqrefineStirlingf} for $\la = \ee^{\gS^I\ti\mu} = \gS^I\ti\la$.


\begin{exo}
We just obtained that
\[
\Ga(z) \sim_1 \ee^{-z} z^{z-\demi} \sqrt{2\pi} \Big( 1
+ \sum_{k\ge0} g_k z^{-k-1} \Big)
\quad \text{uniformly for $|z|>c$ and $\arg z \in (-\be,\be)$}
\]
for any $c>0$ and $\be\in (0,{\pi})$ (with the extended notation of
footnote~\ref{footextsimun}).
Show that
\[
\frac{1}{\Ga(z)} \sim_1 \tfrac{1}{\sqrt{2\pi}} \ee^{z} z^{-z+\demi} \Big( 1
+ \sum_{k\ge0} (-1)^{k+1} g_k z^{-k-1} \Big)
\quad \text{uniformly for $|z|>c$ and $\arg z \in (-\be,\be)$}
\]
for the same values of~$c$ and~$\be$.
%
%
\end{exo}


\begin{rem}
Since $n! = n\Ga(n)$ by~\eqref{eqGammaInterpolFact}
and~\eqref{eqfunctgamma}, we get
\begin{align*}
n! & \sim \frac{ n^n \sqrt{2\pi n} }{\ee^n} \Big( 1
+ \frac{1}{12 n}
+ \frac{1}{288 n^2}
- \frac{139}{51840 n^3}
- \frac{571}{2488320 n^4}
+ \frac{163879}{209018880 n^5}
+ \cdots \Big) \\[1ex]
\frac{1}{n!}  & \sim \frac{\ee^n}{ n^n \sqrt{2\pi n} } \Big( 1
- \frac{1}{12 n}
+ \frac{1}{288 n^2}
+ \frac{139}{51840 n^3}
- \frac{571}{2488320 n^4}
- \frac{163879}{209018880 n^5}
+ \cdots \Big)
\end{align*}
uniformly for $n\in\N^*$.
See \cite{VDA} for a direct proof.
\end{rem}


\begin{rem}   \label{remStokeslaStirl}
In accordance with Remark~\ref{remStokesPhen}, we observe a kind of Stokes
phenomenon for the function~$\la$: it is a priori holomorphic in the cut plane
$\C\setminus\R^-$, or equivalently in the sector $\{\, -\pi < \arg z < \pi \,\}$
of the Riemann surface of the logarithm~$\Clog$, but
Exercise~\ref{exoStirlprepStokes} gives the `reflection formula'
$\la(z) = \frac{1}{(1-\ee^{-2\pi\I z})\la(\ee^{\I\pi} z)}$ for $-\pi < \arg z < 0$, 
which yields a meromorphic continuation for~$\la$ in the larger sector
$\{\, -2\pi < \arg z < \pi \,\}$ 
(with the points $k\,\ee^{-\I\pi}$, $k\in\N^*$, as only poles);
the asymptotic property $\la(z) \sim_1 \ti\la(z)$ is valid in the directions
of $(-\pi,\pi)$ but not in those of $(-2\pi,-\pi]$:
the ray $\ee^{-\I\pi}\R^+$ is singular and the reflection formula implies that,
in the directions of $(-2\pi,-\pi)$, 
$\la(z) \sim -\ee^{2\pi\I z}$, which is exponentially small 
(and $\ee^{-2\pi\I z} \la(z) \sim_1 -\ti\la(z)$ there).

In fact, iterating the reflection formula we find a meromorphic continuation to
the whole of~$\Clog$, with a `monodromy relation' $\la(z) = - \ee^{2\pi\I z}
\la(z\,\eel^{2\pi\I})$ (with the notations of Section~\ref{secClog}).
Outside the singular rays, the asymptotic behaviour is given by
\[
\la(z) = (-1)^n \ee^{-2\pi\I nz}\, \la(z\,\eel^{-2\pi\I n})
\sim_1 (-1)^n \ee^{-2\pi\I nz}\, \ti\la(z)
\]
uniformly for {$|z|$ large enough} and $2\pi n-\be < \arg z < 2\pi n+\be$,
with arbitrary $n\in\Z$ and $\be\in(0,\pi)$.
Except in the initial sector of definition ($n=0$), we thus find exponential
decay and growth alternating at each crossing of a singular ray
$\eel^{(2n-1)\I\pi}\R^+$ or of a ray $\eel^{2n\I\pi}\R^+$ on which the behaviour
is oscillatory,
according to the sign of
$n\IM z$ (since $|\ee^{-2\pi\I nz}| = \ee^{2\pi n\IM z}$).

The last properties can also be deduced from formula~\eqref{eqleftStirlingf}.
\end{rem}


\parag
We leave it to the reader to adapt the results of this section to fine-summable
formal series in a direction~$\th$.

\newpage



\vspace{1.2cm}

\centerline{\Large\sc Formal tangent-to-identity diffeomorphisms}
\addcontentsline{toc}{part}{\sc Formal tangent-to-identity diffeomorphisms}

\vspace{.3cm}


\section{Germs of holomorphic diffeomorphisms}	\label{secGermsholdiffeos}

A \emph{holomorphic local diffeomorphism around~$0$} is a holomorphic map $F\col U \to \C$,
where $U$ is an open neighbourhood of~$0$ in~$\C$, such that $F(0)=0$ and
$F'(0)\neq0$.
The local inversion theorem shows that there is an open neighbourhood~$V$ of~$0$
contained in~$U$ such that $F(V)$ is open and~$F$ induces a biholomorphism
from~$V$ to~$F(V)$.
When we are not too much interested in the precise domains~$U$ or~$V$ but are
ready to replace them by smaller neighrbouhoods of~$0$, we may consider the \emph{germ
of~$F$ at~$0$}. This means that we consider the equivalence class of~$F$ for the
following equivalence relation:
two holomorphic local diffeomorphisms are equivalent if there exists
an open neighbourhood of~$0$ on which their restrictions coincide.

It is easy to see that a germ of holomorphic diffeomorphism at~$0$ can be
identified with the Taylor series at~$0$ of any of its representatives.
Moreover, our equivalence relation is compatible with the composition and the inversion
of holomorphic local diffeomorphisms. 
Consequently, the germs of holomorphic diffeomorphisms at~$0$ make up a
(nonabelian) group, isomorphic to 
\[
\{\, F(t)\in t\C\{t\} \mid F'(0)\neq 0 \,\} =
\Big\{ F(t) = \sum_{n\ge1} c_n t^n \in \C\{t\} \mid c_1 \neq 0 \Big\}.
\]
The coefficient $c_1 = F'(0)$ is called the ``multiplier'' of~$F$.
Obviously, for two germs of holomorphic diffeomorphisms~$F$ and~$G$,
$(F\circ G)'(0) = F'(0) G'(0)$.
Therefore, the germs~$F$ of holomorphic diffeomorphisms at~$0$ such
that $F'(0)=1$ make up a subgroup; such germs are said to be ``tangent-to-identity''.

Germs of holomorphic diffeomorphisms can also be considered at~$\infty$:
via the inversion $t\mapsto z=1/t$, a germ $F(t)$ at~$0$ is conjugate to $f(z) =
1/F(1/z)$.
From now on, we focus on the tangent-to-identity case
\beglabel{eqdefFtsigtau}
F(t) = t - \sig t^2 - \tau t^3 + \cdots 
= t(1-\sig t-\tau t^2 + \cdots)
\in\C\{t\}
\qquad (\sig,\tau\in\C).
\elabel
This amounts to considering germs of holomorphic diffeomorphisms at~$\infty$ of the form
\beglabel{eqdeffzFt}
f(z) = z(1-\sig z\ii-\tau z^{-2} + \cdots)\ii
= z + \sig + (\tau+\sig^2) z\ii + \cdots
\in \id + \C\{z\ii\}.
\elabel
For such a germ~$f$, there exists $c>0$ large enough and a representative which
is an injective holomorphic function in $\{\,|z|>c\,\}$.
We use the notations
\[
\gG \defeq \id + \C\{z\ii\}
\]
for the group of tangent-to-identity germs of holomorphic
diffeomorphisms at~$\infty$, and
\[
\gG_\sig \defeq \id + \sig + z\ii\C\{z\ii\}
\]
when we want to keep track of the coefficient~$\sig$ in~\eqref{eqdeffzFt}.
Notice that, if $f_1 \in \gG_{\sig_1}$ and $f_2 \in \gG_{\sig_2}$,
then
$f_1\circ f_2 \in \gG_{\sig_1+\sig_2}$.


\section{Formal diffeomorphisms}	\label{secFormalDiffeos}

Even if we are interested in properties of the group~$\gG$, or even of a
single element of~$\gG$, it is useful 
(as we shall see in Sections \ref{sec:GermsDeb}--\ref{sec:Bridge})
to drop the convergence requirement
and consider the larger set
\[
\ti\gG = \id + \C[[z\ii]].
\]
This is the set of \emph{formal tangent-to-identity diffeomorphisms at~$\infty$},
which we view as a complete metric space by means of the distance
\[
d\big(\ti f, \ti h\big) \defeq 2^{-\val(\ti\chi-\ti\ph)}, 
\qquad \ti f = \id+\ti\ph, \; \ti h = \id+\ti\chi,
\qquad \ti\ph,\ti\chi\in\C[[z\ii]],
\]
as we did for~$\C[[z\ii]]$ in \S~\ref{parKrull}.
Notice that $\gG$ appears as a dense subset of~$\ti\gG$.
We also use the notation
\[
\ti\gG_\sig = \id+\sig+\zcz = \big\{\, \ti f(z) = z + \sig + \ti\ph(z) \mid
\ti\ph \in \zcz \,\big\} \subset \ti\gG
\]
for any $\sig\in\C$. 
Via the inversion $z\mapsto 1/z$, the elements of~$\ti\gG$ are conjugate to
formal tangent-to-identity diffeomorphisms at~$0$, \ie formal series of the
form~\eqref{eqdefFtsigtau} but without the convergence condition (the
corresponding~$F(t)$  is in~$\C[[t]]$ but not necessarily in~$\C\{t\}$); 
the elements of~$\ti\gG_\sig$ are conjugate to formal series of the form
$F(t) = t - \sig t^2 + \cdots \in \C[[t]]$, by the formal analogue of~\eqref{eqdeffzFt}.


\begin{thm}	\label{thmgroupformdiffinfty}
The set~$\ti\gG$ is a nonabelian topological group for the composition law
\beglabel{eqdefcompos}
\ti f\circ \ti h \defeq \id + \ti\chi + \ti\ph\circ(\id+\ti\chi),
\qquad \ti f = \id+\ti\ph, \; \ti h = \id+\ti\chi,
\qquad \ti\ph,\ti\chi\in\C[[z\ii]],
\elabel
with $\ti\ph\circ(\id+\ti\chi)$ defined by~\eqref{eqdefcircIDchi}.
The subset 
\[
\ti\gG_0 = \id + \zcz
\]
is a subgroup of~$\ti\gG$.
\end{thm}


Notice that the definition~\eqref{eqdefcompos} of the composition law
in~$\ti\gG$ can also be written
\beglabel{eqdefcomposvar}
\ti f \circ \ti h = \sum_{k\ge0} \frac{1}{k!} \ti\chi^k \, \pa^k\ti f,
\qquad \ti h = \id+\ti\chi,
\elabel
with the convention $\pa^0\ti f = \ti f = \id+\ti\ph$, $\pa\ti f = 1+\pa\ti\ph$
and $\pa^k\ti f = \pa^k \ti\ph$ for $k\ge2$.


\begin{proof}[Proof of Theorem~\ref{thmgroupformdiffinfty}]
The composition~\eqref{eqdefcompos} is a continuous map $\ti\gG\times\ti\gG \to
\ti\gG$ because, 
for $\ti f, \ti f^*, \ti h, \ti h^* \in \ti\gG$, formula~\eqref{eqdefcomposvar} implies
\beglabel{eqdiffcirchsth}
\ti f\circ \ti h^* - \ti f\circ \ti h = (\ti h^*-\ti h) \int_0^1 \pa\ti f\circ
\big( (1-t)\ti h+t \ti h^* \big) \,\dd t
\elabel
(where
$\pa\ti f\circ \big( (1-t)\ti h+t \ti h^* \big)$ is a formal series whose
coefficients depend polynomially on~$t$ and integration is meant coefficient-wise);
this is a formal series of valuation $\ge \val(\ti h^*-\ti h)$,
by virtue of~\eqref{ineqvalcirc},
hence the difference
\[
\ti f^*\circ \ti h^* - \ti f\circ \ti h = 
(\ti f^*-\ti f)\circ\ti h^* + \ti f\circ \ti h^* - \ti f\circ \ti h
\]
is a formal series of valuation $\ge \min\big\{ \val(\ti f^*-\ti f), \val(\ti
h^*-\ti h) \big\}$ (using again~\eqref{ineqvalcirc}), \ie 
\[
d(\ti f\circ \ti h,\ti f^*\circ \ti h^*)
\le \max \big\{ d(\ti f,\ti f^*), d(\ti h,\ti h^*) \big\}.
\]
The subset~$\ti\gG_0$ is clearly stable by composition.


The composition law of~$\ti\gG$, when restricted to~$\gG$, boils down to the
composition of holomorphic germs which is associative ($\gG$ is a group) and
$\gG$ is a dense subset of~$\ti\gG$, thus composition is associative in~$\ti\gG$
too.
It is not commutative in~$\ti\gG$ since it is not commutative in~$\gG$.
The element~$\id$ is clearly a unit for composition in~$\ti\gG$ thus we only
need to show that there is a well-defined continuous inverse map
$\ti h \in \ti\gG \mapsto \ti h\ic\in \ti\gG$
and that this map leaves~$\ti\gG_0$ invariant.


We first show that every element $\ti h \in \ti\gG$ has a unique left inverse
$\gL(\ti h)$.
Given $\ti h = \id+\ti\chi$, the equation $\ti f\circ\ti h = \id$ is
equivalent to the fixed-point equation
\beglabel{eqdefgC}
\ti f = \gC(\ti f), \qquad
\gC(\ti f) \defeq \id - (\ti f\circ \ti h - \ti f) =
\id - \ti\chi \int_0^1 \pa\ti f \circ (\id+t\ti\chi) \,\dd t
\elabel
(we have used~\eqref{eqdiffcirchsth} to get the last expression of~$\gC$).
The map $\gC \col \ti\gG \to \ti\gG$ is a contraction of our complete metric
space, because the difference 
\beglabel{eqdiffCfCfstar}
\gC(\ti f^*) - \gC(\ti f)  = - \ti\chi \int_0^1 
\pa(\ti f^*-\ti f) \circ(\id+t\ti\chi) \,\dd t
\elabel
has valuation $\ge \val(\ti f^*-\ti f) + 1$
(because of~\eqref{ineqvalcirc}: 
$\val\big( \pa(\ti f^*-\ti f) \circ(\id+t\ti\chi) \big) = 
\val\big( \pa(\ti f^*-\ti f) \big) \ge 
\val(\ti f^*-\ti f) + 1$ for each $t$),
hence
$d\big( \gC(\ti f), \gC(\ti f^*) \big) \le \demi d(\ti f,\ti f^*)$.
The Banach fixed-point theorem implies that there is a unique solution $\ti f =
\gL(\ti h)$, obtained as the limit of the Cauchy sequence 
$\gL_n(\ti h) \defeq \underbrace{\gC \circ \cdots \circ \gC}_\text{$n$ times}(0)$
as $n\to\infty$.

We observe that, if $\ti h\in\ti\gG_0$, then $\gC(\ti\gG) \subset \ti\gG_0$,
thus $\gL_n(\ti h) \in \ti\gG_0$ for each $n\ge0$ and clearly $\gL(\ti
h)\in\ti\gG_0$ in that case.


The fact that each element has a unique left inverse implies that each element
is invertible: 
given $\ti h\in\ti\gG$, its left inverse $\ti f \defeq \gL(\ti h)$ is also a right inverse because
$\ti h^*\defeq\gL(\ti f)$ satisfies
$\ti h^* = \ti h^* \circ(\ti f\circ\ti h)= (\ti h^* \circ\ti f)\circ\ti h = \ti h$,
\ie $\ti h\circ \ti f = \id$.


Finally, we check that~$\gL$ is continuous.
For $\ti h,\ti h^*\in\ti\gG$, we denote by $\gC, \gC^*$ the corresponding maps
defined by~\eqref{eqdefgC}.
For any $\ti f, \ti f^*$, the difference
$\gC^*(\ti f) - \gC(\ti f) = \ti f\circ\ti h - \ti f\circ\ti h^*$ has valuation 
$\ge \val(\ti h^*-\ti h)$ (as already deduced from~\eqref{eqdiffcirchsth}),
while $\val\big(\gC^*(\ti f^*) - \gC^*(\ti f)\big) \ge  \val(\ti f^*-\ti f) +
1$ (as already deduced from~\eqref{eqdiffCfCfstar}),
hence $d\big(\gC(\ti f),\gC^*(\ti f^*)\big) \le \max\big\{ 
d(\ti h, \ti h^*), \demi d(\ti f, \ti f^*) \big\}$.
It follows by induction that
$d\big( \gL_n(\ti h),\gL_n(\ti h^*) \big) =
d\big( \gC(\gL_{n-1}(\ti h)),\gC^*(\gL_{n-1}(\ti h^*)) \big)
\le d(\ti h, \ti h^*)$ for every $n\ge1$,
hence $d\big( \gL(\ti h),\gL(\ti h^*) \big) \le d(\ti h, \ti h^*)$.
\end{proof}

Notice that $\ti\gG_0 = \{\, \ti f \in \ti\gG \mid d(\id,\ti f) \le \demi \,\}
= \{\, \ti f \in \ti\gG \mid d(\id,\ti f) < 1 \,\}$
is a closed ball as well as an open ball, thus it is both closed and open for
the Krull topology of~$\ti\gG$.

\section{Inversion in the group~$\protect\ti{\protect\gG}$} \label{secinversgrpgG}

There is an explicit formula for the inverse of an element of~$\ti\gG$, which is
a particular case of the Lagrange reversion formula (adapted to our framework):

\begin{thm}	\label{thmLagrInvformdiffinfty}
For any $\ti\chi\in\C[[z\ii]]$, the inverse of $\ti h = \id+\ti\chi$ can be written
as the formally convergent series of formal series
\beglabel{eqnLagrangeRev}
(\id+\ti\chi)\ic = \id + \sum_{k\ge1} \frac{(-1)^k}{k!} \pa^{k-1} (\ti\chi^k).
\elabel
\end{thm}


The proof of Theorem~\ref{thmLagrInvformdiffinfty} will make use of

\begin{lemma}	\label{lemHnzero}
Let $\ti\chi\in\C[[z\ii]]$ and $n\ge1$. Then, for any $\ti\psi\in\C[[z\ii]]$,
\beglabel{eqHnzero}
\sum_{k=0}^n (-1)^k \binom{n}{k} \ti\chi^{n-k} \, \pa^{n-1}(\ti\chi^k \ti\psi)
= 0.
\elabel
\end{lemma}

\begin{proof}[Proof of Lemma~\ref{lemHnzero}]
Let us call $H_n\ti\psi$ the \lhs\ of~\eqref{eqHnzero}.
We have $H_1\ti\psi = \ti\chi \, \pa^0\ti\psi - \pa^0(\ti\chi \ti\psi) = 0$. It is
thus sufficient to prove the recursive formula
\[
H_{n+1}\ti\psi = - \pa H_n(\ti\chi\ti\psi) + \ti\chi \, \pa H_n\ti\psi 
- n (\pa\ti\chi) H_n\ti\psi.
\]
To this end, we use the convention $\binom{n}{-1} = \binom{n}{n+1} = 0$ and
compute
\begin{align*}
- \pa H_n(\ti\chi\ti\psi) &= \sum_{k=-1}^n (-1)^{k+1} \binom{n}{k} \pa\big[
\ti\chi^{n-k}\, \pa^{n-1}(\ti\chi^{k+1}\ti\psi)
\big] \\
&= \sum_{k=0}^{n+1} (-1)^k \binom{n}{k-1} \pa\big[
\ti\chi^{n+1-k}\, \pa^{n-1}(\ti\chi^{k}\ti\psi)
\big]\\
\intertext{(shifting the summation index to get the last expression), while}
\ti\chi \, \pa H_n\ti\psi &= \sum_{k=0}^{n+1} (-1)^k \binom{n}{k} \ti\chi \, \pa\big[
\ti\chi^{n-k} \, \pa^{n-1}(\ti\chi^k \ti\psi)
\big].
\end{align*}
The Leibniz rule yields
\begin{multline*}
- \pa H_n(\ti\chi\ti\psi) + \ti\chi \, \pa H_n\ti\psi =
\sum_{k=0}^{n+1} (-1)^k
\Big[ \binom{n}{k-1} + \binom{n}{k} \Big]
\ti\chi^{n+1-k} \, \pa^{n}(\ti\chi^k \ti\psi) \\[1ex]
+ \sum_{k=0}^{n+1} (-1)^k
\Big[ (n+1-k)\binom{n}{k-1} + (n-k)\binom{n}{k} \Big]
\ti\chi^{n-k} (\pa\ti\chi) \pa^{n-1}(\ti\chi^k \ti\psi).
\end{multline*}
The expression in the former bracket is $\binom{n+1}{k}$, hence the first sum is
nothing but $H_{n+1}\ti\psi$;
the expression in the latter bracket is $n$ times
$\binom{n-1}{k}+\binom{n-1}{k-1} = \binom{n}{k}$, hence the second sum is
$n(\pa\ti\chi)H_n\ti\psi$. 
\end{proof}


\begin{proof}[Proof of Theorem~\ref{thmLagrInvformdiffinfty}]
Let $\ti h = \id + \ti\chi \in \ti\gG$.
Lemma~\ref{lemHnzero} shows that the \rhs\ of~\eqref{eqnLagrangeRev} defines a
left inverse for~$\ti h$. Indeed, denoting by $\ti f = \id + \ti\ph$ this \rhs,
we have
\[
\ti f \circ \ti h - \id = \ti\chi + \ti\ph\circ(\id+\ti\chi)
= \ti\chi + \sum_{\ell\ge0,\,k\ge1}
\frac{(-1)^k}{k!\ell!} \ti\chi^\ell \, \pa^{k+\ell-1}(\ti\chi^k)
= \sum_{n\ge1} \frac{1}{n!} \ti H_n
\]
with $\ti H_n = \sum (-1)^k \binom{n}{k} \ti\chi^\ell \, \pa^{n-1}(\ti\chi^k)$,
the last sum running over all pairs of non-negative integers $(k,\ell)$ such
that $k+\ell = n$
(absorbing the first~$\ti\chi$ in~$\ti H_1$ and taking care of $k=0$
according as $n=1$ or $n\ge2$; formal summability legitimates our Fubini-like
manipulation), then Lemma~\ref{lemHnzero} with $\ti\psi=1$ says
that $\ti H_n=0$ for every $n\ge1$.
\end{proof}


\begin{exo}
[\textbf{Lagrange reversion formula}]
Prove that, with the same convention as in~\eqref{eqdefcomposvar},
\[
\ti f \circ \ti h\ic = \ti f + 
\sum_{k\ge1} \frac{(-1)^k}{k!} \pa^{k-1}(\ti\chi^k\pa\ti f),
\qquad \ti h = \id + \ti\chi.
\]
(Hint: Use Lemma~\ref{lemHnzero} with $\ti\psi = \pa(\ti f-\id) = -1 + \pa\ti f$.)
\end{exo}


\begin{exo}
Let $h = \id+\chi \in \gG$, \ie with $\chi\in\C\{z\ii\}$. We can thus choose
$c_0,M>0$ such that $|\chi(z)| \le M$ for $|z|\ge c_0$.
Show that $h\ic(z)$ is convergent for $|z|\ge c_0+M$.
(Hint: Given $\de>M$, use the Cauchy inequalities to bound $|\pa^{k-1}(\chi^k)(z)|$ for
$|z|>c_0+\de$.)
\end{exo}

\section{The group of $1$-summable formal diffeomorphisms in an arc of directions}	\label{secgrponesummdiffeos}

Among all formal tangent-to-identity diffeomorphisms, we now distinguish those which are
$1$-summable in an arc of directions.

\begin{Def}   \label{DefSummaDiffeo}
Let $I$ be an open interval of~$\R$. Let $\ga,\al\col I\to\R$ be locally bounded
functions with $\al\ge0$.
For any $\sig\in\C$ we define
\begin{align*}
\ti\gG(I,\ga,\al) &\defeq \big\{\, \ti f = \id+\ti\ph_0 \mid 
\ti\ph_0 \in \cB\ii\big(\C\,\de \oplus \cN(I,\ga,\al) \big) \,\big\},
&\quad \ti\gG_\sig(I,\ga,\al) &\defeq \ti\gG(I,\ga,\al) \cap \ti\gG_\sig,\\
\ti\gG(I,\ga) &\defeq \big\{\, \ti f = \id+\ti\ph_0 \mid 
\ti\ph_0 \in \cB\ii\big(\C\,\de \oplus \cN(I,\ga) \big) \,\big\},
&\quad \ti\gG_\sig(I,\ga) &\defeq \ti\gG(I,\ga) \cap \ti\gG_\sig,\\
\ti\gG(I) &\defeq \big\{\, \ti f = \id+\ti\ph_0 \mid 
\ti\ph_0 \in \cB\ii\big(\C\,\de \oplus \cN(I) \big) \,\big\},
&\quad \ti\gG_\sig(I) &\defeq \ti\gG(I) \cap \ti\gG_\sig.
\end{align*}
We extend the definition of the Borel summation operator~$\gS^I$ to~$\ti\gG(I)$
by setting
\[
\ti f = \id+\ti\ph_0 \in \ti\gG(I,\ga) 
\quad\Longrightarrow\quad
\gS^I \ti f(z) = z + \gS^I\ti\ph^0(z), \qquad z\in\gD(I,\ga).
\]
\end{Def}

For $|I|\ge2\pi$, $\ti\gG(I)$ coincides with the group~$\gG$ of holomorphic
tangent-to-identity diffeomorphisms and $\gS^I$ is the ordinary summation operator for
Taylor series at~$\infty$, but 
\[
|I|<2\pi \quad\Longrightarrow\quad
\gG \,\subsetneqq\, \ti\gG(I) \,\subsetneqq\, \ti\gG.
\]
For $\ti f\in\ti\gG(I)$, the function $\gS^I\ti f$ is holomorphic in a sectorial
neighbourhood of~$\infty$ (but not in a full neighbourhood of~$\infty$ if $\ti
f\notin \gG$); we shall see that it defines an injective transformation in a
domain of the form~$\gD(I,\ga)$. We first study composition and inversion in~$\ti\gG(I)$.


\begin{thm}	\label{thmcompossummdiffeo}
Let $I$ be an open interval of~$\R$ and $\ga,\al\col I\to\R$ be locally bounded
functions with $\al\ge0$.
Let $\sig,\tau\in\C$ and $\ti f\in\ti\gG_\sig(I,\ga,\al)$, $\ti g\in\ti\gG_\tau(I,\ga)$.
Then $\ti g\circ\ti f \in \ti\gG_{\sig+\tau}(I,\ga_1)$
with $\ga_1 = \ga + |\sig|+\sqrt{\al}$,
the function~$\gS^I\ti f$ maps $\gD(I,\ga_1)$ in $\gD(I,\ga)$ and
\[
\gS^I(\ti g\circ\ti f) = (\gS^I\ti g) \circ (\gS^I\ti f)
\ens\text{on $\gD(I,\ga_1)$}.
\]
\end{thm}

\begin{proof}
Apply Theorem~\ref{thmcompatsumcirc} to $\ti\ph_0 \defeq \ti f - \id$ and
$\ti\psi_0 \defeq \ti g - \id$.
\end{proof}


\begin{thm}   \label{thmSummaDiffeo}
Let $\ti f \in \ti\gG_\sig(I,\ga,\al)$.
Then $\ti h \defeq \ti f\ic \in \ti\gG_{-\sig}(I,\ga^*,\al)$ with $\ga^* \defeq
\ga+|\sig|+2\sqrt{\al}$ and
\begin{align}
\gS^I\ti f\big(\gD(I,\ga_1)\big) &\subset \gD(I,\ga^*),&\ens
(\gS^I\ti h) \circ \gS^I\ti f &= \id \ens\text{on $\gD(I,\ga_1)$,} \\
\gS^I\ti h\big(\gD(I,\ga_2)\big) &\subset \gD(I,\ga), &\ens
(\gS^I\ti f) \circ \gS^I\ti h &= \id \ens\text{on $\gD(I,\ga_2)$,}
\end{align}
with $\ga_1 \defeq \ga + 2|\sig| + (1+\sqrt{2})\sqrt{\al}$
and $\ga_2 \defeq \ga + |\sig| + (1+\sqrt{2})\sqrt{\al}$.

Moreover, $\gS^I\ti f$ is injective on $\gD\big(I,\ga+(1+\sqrt{2})\sqrt{\al}\big)$.
\end{thm}

\begin{proof}
We first assume $\ti f\in \ti\gG_0(I,\ga,\al)$.
By~\eqref{eqnLagrangeRev}, we have $\ti h = \id + \ti\chi$ with $\ti\chi$ given
by a formally convergent series in $\zcz$:
\[
\ti\chi = \sum_{k\ge1} \ti\chi_k, \qquad
\ti\chi_k = \frac{(-1)^k}{k!} \pa^{k-1} (\ti\ph^k).
\]
Correspondingly, $\cB\ti\chi$ is given by a formally convergent series in
$\C[[\ze]]$:
\[
\hat\chi = \sum_{k\ge1} \hat\chi_k, \qquad
\hat\chi_k = -\frac{\ze^{k-1}}{k!} \hat\ph^{*k}
\]
(beware that the last expression involves multiplication by $-\frac{\ze^{k-1}}{k!}$, not
convolution!).
We argue as in the proof of Theorem~\ref{thmcompatsumcirc} and view~$\hat\chi$
as a series of holomorphic functions in the union of a disc $D(0,R)$ and a
sector~$\Sig$ in which~$\hat\ph$ itself is holomorphic; inequalities~\eqref{ineqhatphkD}
and~\eqref{ineqhatphkS} yield
\begin{align}
\label{ineqhatchikDR}
|\hat\chi_k(\ze)| &\le A^k \frac{\xi^{2(k-1)}}{k!(k-1)!},
& \hspace{-4em} \ze &\in D(0,R), \\[1ex]
\label{ineqhatchikfactkfactkm}
|\hat\chi_k(\ze)| &\le \al(\th)^k
\frac{\xi^{2(k-1)}}{k!(k-1)!} \ee^{\ga(\th)\xi},
& \hspace{-4em} \ze &\in \Sig,
\end{align}
where $\xi=|\ze|$ and $\th=\arg\ze$.
The series of holomorphic functions $\sum\hat\chi_k$ is thus uniformly
convergent in every compact subset of $D(0,R)\cup\Sig$ and its sum is a
holomorphic function whose Taylor series at~$0$ is~$\hat\chi$. 
Therefore $\hat\chi\in\C\{\ze\}$ extends analytically to $D(0,R)\cup\Sig$;
moreover, since $\frac{1}{k!(k-1)!} \le \frac{1}{k} \frac{2^{2(k-1)}}{(2(k-1))!}$,
\eqref{ineqhatchikfactkfactkm} yields 
\[
|\hat\chi(\ze)| \le
\sum_{k\ge1} \frac{\al(\th)^k}{k} \frac{(2\xi)^{2(k-1)}}{\big(2(k-1)\big)!} \ee^{\ga(\th)\xi}
\le \al(\th)\, \ee^{(\ga(\th)+2\sqrt{\al(\th)})\xi}
\]
for $\ze\in\Sig$.
Hence $\ti h \in \ti\gG_0(I,\ga+2\sqrt{\al},\al)$ when $\sig=0$.


In the general case, we observe that $\ti f = (\id+\sig)\circ \ti g$ with 
$\ti g \defeq (\id-\sig)\circ\ti f \in \ti\gG_0(I,\ga,\al)$,
thus $\ti g\ic = \id+\ti\chi \in \ti\gG_0(I,\ga+2\sqrt{\al},\al)$
and $\ti h = \ti f\ic = \ti g\ic \circ (\id-\sig) = \id - \sig + T_{-\sig}\ti\chi$,
which implies $\ti h\in \ti\gG(I,\ga+2\sqrt{\al}+|\sig|,\al)$ by the second property
in Lemma~\ref{lemelemptiescB}.


Since $\ti h\circ \ti f = \ti f\circ \ti h = \id$, we can apply
Theorem~\ref{thmcompossummdiffeo} and get
$(\gS^I\ti h) \circ \gS^I\ti f = \id$ and
$(\gS^I\ti f) \circ \gS^I\ti h = \id$ 
in appropriate domains; in fact, by analytic continuation, these identities will
hold in any domain $\gD(I,\ga+\de_1)$, \resp $\gD(I,\ga+\de_2)$, such that
\[
\gS^I\ti f\big(\gD(I,\ga+\de_1)\big) \subset \gD(I,\ga^*), \qquad
\gS^I\ti h\big(\gD(I,\ga+\de_2)\big) \subset \gD(I,\ga).
\]
Writing $\ti f = \id+\sig+\ti\ph$ with $\cB\ti\ph\in\cN(I,\ga,\al)$, with the
help of~\eqref{ineqgSImapPi} one can easily show that $\de_1=\ga_1-\ga$ and
$\de_2=\ga_2-\ga$ satisfy this.


For the injectivity statement, we write again $\ti f = (\id+\sig)\circ\ti g$ and
apply the previous result to $\ti g\in\gG_0(I,\ga,\al)$.
The function $\gS^I\ti g$ maps $\gD\defeq
\gD\big(I,\ga+(1+\sqrt{2})\sqrt{\al}\big)$ in the domain
$\gD\big(I,\ga+2\sqrt{\al}\big)$, on which $\gS^I(\ti g\ic)$ is well-defined,
and
$\gS^I(\ti g \ic) \circ \gS^I\ti g = \id$ on~$\gD$,
therefore $\gS^I\ti g$ is injective on~$\gD$, and so is the function $\gS^I\ti f = \sig+\gS^I\ti g$.
\end{proof}


\begin{cor}
For any open interval~$I$,
$\ti\gG(I)$ and $\ti\gG_0(I)$ are subgroups of~$\ti\gG$.
\end{cor}


\begin{exo} 
Consider the set $\id+\C[[z\ii]]_1$ of $1$-Gevrey tangent-to-identity formal
diffeomorphisms, so that
\[
\ti\gG(I) \;\subsetneqq\; \id+\C[[z\ii]]_1 \;\subsetneqq\; \ti\gG.
\]
Show that $\id+\C[[z\ii]]_1$ is a subgroup of~$\ti\gG$.
(Hint: View $\C[[z\ii]]_1$ as $\cB\ii\big(\C\,\de \oplus \C\{\ze\} \big)$ and
imitate the previous chain of reasoning.)
\end{exo}


We shall see in Section~\ref{sec:SummaIter} how $1$-summable formal diffeomorphisms
occur in the study of a holomorphic germ $f\in\gG_1$.


\vspace{.9cm}

\centerline{\Large\sc The algebra of resurgent functions}
\addcontentsline{toc}{part}{\sc The algebra of resurgent functions}

\vspace{.6cm}


\section{Resurgent functions, resurgent formal series}	
\label{sec_resurfunct}


Among $1$-Gevrey formal series, we have distinguished the subspace of those which
are $1$-summable in a given arc of directions and studied it in
Sections~\ref{secvarydir}--\ref{secgrponesummdiffeos}.
We shall now study another subspace of $\C[[z\ii]]_1$, which consists of
``resurgent formal series''.
As in the case of $1$-summability, 
we make use of the algebra isomorphism~\eqref{eqalgisGevdeCV} 
\[
\cB \col \C[[z\ii]]_1 \xrightarrow{\sim} \C\de\oplus\C\{\ze\}
\]
and give the definition not directly in terms of the formal series
themselves\footnote{See \cite{LodayRemy} for another point of view.},
but rather in terms of their formal Borel transforms, for which,
beyond convergence near the origin, we shall require a certain
property of analytic continuation.

For any $R>0$ and $\ze_0\in\C$ we use the notations 
\begin{gather}
\label{eqdefDst}
D(\ze_0,R) \defeq \{\, \ze\in\C \mid |\ze-\ze_0|< R \,\}, \\
\label{eqdefDDst}
\D_R \defeq D(0,R), \qquad
\D^*_R \defeq \D_R \setminus \{0\}.
\end{gather}

\begin{figure}
\begin{center}

\includegraphics[scale=1]{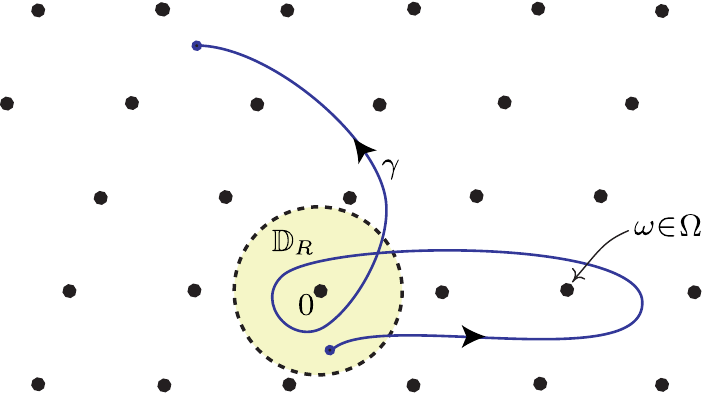} 

\bigskip 

\caption{%
\emph{$\Om$-continuability.}
Any path~$\ga$ starting in $\D_R^*$ and contained in $\C\setminus\Om$ must be a
path of analytic continuation for $\hat\ph \in \hat{\protect\gR}_\Om$.}
\label{fig:Om_cont}

\end{center}
\end{figure}

\begin{Def}	\label{DefOmCont} \label{DefOmResur}
Let $\Om$ be a non-empty closed discrete subset of~$\C$,
let $\hat\ph(\ze)\in\C\{\ze\}$ be a holomorphic germ at the origin.
We say that $\hat\ph$ is an \emph{$\Om$-continuable germ} if there exists $R>0$ not
larger than the radius of convergence of~$\hat\ph$ such that $\D^*_R\cap\Om=\emptyset$
and $\hat\ph$ admits analytic continuation along any path of $\C\setminus\Om$
originating from any point of~$\D^*_R$.
See Figure~\ref{fig:Om_cont}.
We use the notation
\[
\hat\gR_\Om \defeq \{\, \text{all $\Om$-continuable germs} \,\}
\subset \C\{\ze\}.
\]
We call \emph{$\Om$-resurgent function} any element of $\C\,\de \oplus \hat\gR_\Om$,
\ie any element of $\C\,\de \oplus \C\{\ze\}$ of the form $a\,\de + \hat\ph$
with $a =$ a complex number and $\hat\ph = $ an $\Om$-continuable germ.

We call \emph{$\Om$-resurgent formal series} any $\ti\ph_0(z) \in
\C[[z\ii]]_1$ whose formal Borel transform is an $\Om$-resurgent
function,
%
%
\ie any $\ti\ph_0$ belonging to
\[
\ti\gR_\Om \defeq \cB\ii \big( \C\,\de \oplus \hat\gR_\Om\big) \subset \C[[z\ii]]_1.
\]
\end{Def}


\begin{rem}	\label{rempathscont}
%
%
In the above definition, ``path'' means a continuous function $\ga
\col J\to\C\setminus\Om$, where~$J$ is a compact interval of~$\R$;
without loss of generality, all our paths will be assumed piecewise
continuously differentiable.
As is often the case with analytic continuation and Cauchy integrals, the
precise parametrisation of~$\ga$ will usually not matter, in the sense that we
shall get the same result from two paths
$\ga\col [a,b]\to \C\setminus\Om$ and 
$\ga'\col [a',b']\to \C\setminus\Om$
which only differ by a change of parametrisation ($\ga = \ga'\circ\sig$ with
$\sig\col [a,b]\to [a',b']$ piecewise continuously differentiable, increasing
and mapping $a$ to~$a'$ and $b$ to~$b'$).
\end{rem}

Our definitions are particular cases of \'Ecalle's definition of
``continuability without a cut'' (or ``endless continuability'') for germs, and
``resurgence'' for formal series
(we prescribe in advance the possible location of the singularities of the
analytic continuation of~$\hat\ph$, whereas the theory is developed in
Vol.~3 of \cite{Eca81} without this restriction).
Here we stick to the simplest cases;
typical examples with which we shall work are $\Om=\Z$ or $2\pi\I\,\Z$.


\begin{rem}	\label{reminitialpoint} 
Let ${\rho(\Om)} \defeq \min\big\{ |\om|, \; \om\in\Om\setminus\{0\} \big\}$.
Any $\hat\ph\in\hat\gR_\Om$ is a holomorphic germ at~$0$ with radius of
convergence $\ge{\rho(\Om)}$ and one can always take $R={\rho(\Om)}$ in
Definition~\ref{DefOmCont}. 
In fact, given an arbitrary $\ze_0 \in \D_{\rho(\Om)}$, we have
\[
\hat\ph\in\hat\gR_\Om \quad\Longleftrightarrow\quad
\left| \begin{aligned}
&\text{$\hat\ph$ germ of holomorphic function of~$\D_{\rho(\Om)}$ admitting analytic continuation}\\
&\text{along any path $\ga\col [0,1] \to \C$ such that $\ga(0) = \ze_0$ and $\ga\big( (0,1] \big)
\subset \C\setminus\Om$}
\end{aligned} \right.
\]
(even if $\ze_0 = 0$ and $0\in\Om$: there is no need to avoid~$0$ at the beginning of the
path, when we still are in the disc of convergence of~$\hat\ph$).
\end{rem}
%

\begin{exa}
Trivially, any entire function of~$\C$ defines an $\Om$-continuable germ; as a
consequence,
\[
\C\{z\ii\} \subset \ti\gR_\Om.
\]
Other elementary examples of $\Om$-continuable germs are the functions which are
holomorphic in $\C\setminus\Om$ and regular at~$0$, like $\frac{1}{(\ze-\om)^m}$
with $m\in\N^*$ and $\om\in\Om\setminus\{0\}$.
\end{exa}


\begin{lemma}	\label{lemEulPoinStrilRES}
\begin{enumerate}[--]
\item
The Euler series $\ti\ph\eul(z)$ defined by~\eqref{eqdefEulerseriesatinfty}
belongs to $\ti\gR_{\{-1\}}$.
\item
Given $w = \ee^s$ with $\RE s < 0$, the series $\ti\ph\poin(z)$ of Poincar\'e's
example~\eqref{eqdefphpoin} belongs to $\ti\gR_\Om$ with
$\Om \defeq s + 2\pi\I\,\Z$.
\item
The Stirling series $\ti\mu(z)$ of Theorem~\ref{thmStirlings} (explicitly given
by~\eqref{eqexplicitStirl}) belongs to $\ti\gR_{2\pi\I\,\Z}$.
\end{enumerate}
\end{lemma}

\begin{proof}
The Borel transforms of all these series have a meromorphic continuation:
\begin{enumerate}[--]
\item
Euler: $\hat\ph\eul(\ze) = (1+\ze)\ii$ by~\eqref{eqBorEulseries}.
\item
Poincar\'e: $\hat\ph\poin(\ze) = \frac{1}{1-\ee^{s-\ze}}$
by~\eqref{eqhatphpoin}.
\item
Stirling: $\hat\mu(\ze) = \ze^{-2}\left( \frac{\ze}{2}\coth\frac{\ze}{2} - 1
\right)$ by~\eqref{eqBorelStirling}.
\end{enumerate}
\end{proof}


\begin{exo}
Any $\{0\}$-continuable germ defines an entire function of~$\C$.
(Hint: view $\C$ as the union of a disc and two cut planes.)
\end{exo}


\begin{exo}
Give an example of a holomorphic germ at~$0$ which is not $\Om$-continuable for any
non-empty closed discrete subset~$\Om$ of~$\C$.
\end{exo}


But in all the previous examples the Borel transform was single-valued, whereas the
interest of Definition~\ref{DefOmCont} is to authorize multiple-valuedness when
following the analytic continuation.
For instance, the exponential of the Stirling series $\ti\la = \ee^{\ti\mu}$,
which gives rise to the refined Stirling formula~\eqref{eqrefineStirlingf}, has
a Borel transform with a multiple-valued analytic continuation and belongs to
$\ti\gR_{2\pi\I\,\Z}$, although this is more difficult to check (see
Sections~\ref{secSubstResCV} and~\ref{secAlienOpNL}).
We now give elementary examples which illustrate multiple-valued analytic
continuation.

\begin{nota}
If~$\hat\ph$ is a holomorphic germ at~$\ga(a)$ which admits an analytic
continuation along~$\ga$, we denote by
$\cont_\ga\hat\ph$
the resulting holomorphic germ at the endpoint~$\ga(b)$.
\end{nota}


\begin{exa}	\label{exalog}
Consider $\hat\ph(\ze) = \sum_{n\ge1} \frac{\ze^n}{n}$: this is a
holomorphic germ belonging to~$\hat\gR_{\{1\}}$ but its analytic continuation is
not single-valued. 
Indeed, the disc of convergence of~$\hat\ph$ is $\D_1$ and, for any
$\ze\in \D_1$, $\hat\ph(\ze) = \int_0^\ze \frac{\dd\xi}{1-\xi} = - \Log(1-\ze)$
with the notation~\eqref{eqdefLogppalbr} for the principal branch of the logarithm, 
hence the analytic continuation
of~$\hat\ph$ along a path~$\ga$ originating from~$0$,
avoiding~$1$ and ending at a point~$\ze_1$ is the holomorphic germ at~$\ze_1$
explicitly given by
\[
\cont_\ga\hat\ph(\ze) = \int_\ga \frac{\dd\xi}{1-\xi} + \int_{\ze_1}^\ze \frac{\dd\xi}{1-\xi}
\qquad \text{($\ze$ close enough to~$\ze_1$)},
\]
which yields a multiple-valued function in $\C\setminus\{1\}$
(two paths from~$0$ to~$\ze_1$ do not give rise to the same analytic
continuation near~$\ze_1$ unless they are homotopic in $\C\setminus\{1\}$).
The germ~$\hat\ph$ is $\Om$-continuable if and only if $1\in\Om$.
\end{exa}


\begin{exa}	\label{exaLogPole}
A related example of $\{0,1\}$-continuable germ with mutivalued analytic
continuation is given by
$\sum_{n\ge0} \frac{\ze^n}{n+1} = - \frac{1}{\ze}\Log(1-\ze)$, for which there
is a principal branch holomorphic in the cut plane $\C\setminus[1,+\infty)$ and
all the other branches have a simple pole at~$0$.
This germ is $\Om$-continuable if and only if $\{0,1\}\subset\Om$.
\end{exa}


\begin{exa}
If $\om\in\Om\setminus\{0\}$ and $\hat\psi \in\C\{\ze\}$ extends analytically to
$\C\setminus\Om$, then,
for any branch of the logarithm $\Llog$, 
the formula $\hat\ph(\ze) = \hat\psi(\ze)\Llog(\ze-\om)$ defines a
germ of $\hat\gR_\Om$ with non-trivial monodromy around~$\om$:
the branches of the analytic continuation of~$\hat\ph$ differ by integer
multiples of $2\pi\I\,\hat\psi$.
\end{exa}


\begin{exa}
If $\om\in\C^*$ and $m\in\N^*$, then
$\big(\Llog(\ze-\om)\big)^m \in \hat\gR_{\{\om\}}$
for any branch of the logarithm;
if moreover $\om\neq-1$, then
$\big(\Llog(\ze-\om)\big)^{-m}\in \hat\gR_{\{\om,\om+1\}}$.
\end{exa}


\begin{exa}	\label{exaIncGa}
Given $\al\in\C$, the incomplete Gamma function is defined for $z>0$ by
\[
%
\Ga(\al,z) \defeq \int_z^{+\infty} \ee^{-t} t^{\al-1} \, \dd t
%
\]
and it extends to a holomorphic function in $\C\setminus\R^-$
(notice that $\Ga(\al,z) \xrightarrow[z\to0]{} \Ga(\al)$ if $\RE\al>0$).
The change of variable $t = z(\ze+1)$ in the integral yields the formula
\beglabel{eqincGasumma}
\Ga(\al,z) = \ee^{-z} z^\al (\gS^I \hat\ph_\al)(z),
\qquad
\hat\ph_\al(\ze) \defeq (1+\ze)^{\al-1},
\elabel
where $I = (-\frac{\pi}{2},\frac{\pi}{2})$ and we use the principal
branch of the logarithm~\eqref{eqdefLogppalbr} to define the holomorphic
function $(1+\ze)^{\al-1}$ as
$\ee^{(\al-1)\Log(1+\ze)}$.
The germ $\hat\ph_\al$ is always $\{-1\}$-resurgent; it has multiple-valued
analytic continuation if $\al\not\in\Z$.
Hence
\beglabel{eqasymptincGa}
z^{-\al} \ee^z \Ga(\al,z) \sim_1
\ti\ph_\al(z) = \sum_{n\ge0} (\al-1)(\al-2)\cdots(\al-n) z^{-n-1},
\elabel
which is always a $1$-summable and $\{-1\}$-resurgent formal series (a polynomial in~$z\ii$ if
$\al\in\N^*$, a divergent formal series otherwise).
\end{exa}


$\hat\gR_\Om$ and $\ti\gR_\Om$ clearly are linear subspaces of $\C\{\ze\}$ and
$\C[[z\ii]]_1$. We end this section with elementary stability properties:

\begin{lemma}	\label{lemelemstab}
Let $\Om$ be any non-empty closed discrete subset of~$\C$.
Let $\hat B \in \hat\gR_\Om$. Then multiplication by~$\hat B$ leaves
$\hat\gR_\Om$ invariant.
In particular, for any $c\in \C$,
\[
\hat\ph(\ze)\in \hat\gR_\Om 
\ens\Longrightarrow\ens
-\ze\hat\ph(\ze)\in  \hat\gR_\Om \ens\text{and}\ens \ee^{-c\ze}\hat\ph(\ze) \in  \hat\gR_\Om.
\]
The operator $\frac{\dd\,}{\dd\ze}$ too leaves $\hat\gR_\Om$ invariant.

As a consequence, $\ti\gR_\Om$ is stable by $\pa = \frac{\dd\,}{\dd z}$
and~$T_c$.
Moreover, if $\ti\psi\in\ti\gR_\Om \cap z^{-2}\C[[z\ii]]$, then 
$z\ti\psi\in\ti\gR_\Om$
and the solution in $z\ii\C[[z\ii]]$ of the difference equation
\[
\ti\ph(z+1) - \ti\ph(z) = \ti\psi(z)
\]
belongs to $\ti\gR_{\Om\cup 2\pi\I\,\Z^*}$.
\end{lemma}

\begin{proof}
  Exercise (use the fact that multiplication by~$\hat B$ commutes with analytic
  continuation: the analytic continuation of $\hat B\hat\ph$ along a path~$\ga$ of
  $\C\setminus\Om$ starting in $\D_{\rho(\Om)}^*$ exists and equals $\hat B(\ze)
  \cont_\ga\hat\ph(\ze)$; then use Lemma~\ref{lemelemptiescB}, \eqref{eqdefhatpa},
  \eqref{eqdefcounterTc} and Corollary~\ref{coreqdifflin}).
\end{proof}


\section{Analytic continuation of a convolution product: the easy case}	\label{sec_contconvoleasy}


Lemma~\ref{lemelemstab} was dealing with the multiplication of two germs of
$\C\{\ze\}$, however we saw in Section~\ref{sec_convols} that the natural
product in this space is convolution.
The question of the stability of~$\hat\gR_\Om$ under convolution is much
subtler. Let us begin with an easy case, which is already of interest:
\begin{lemma}	\label{lemeasyconvol}
Let $\Om$ be any non-empty closed discrete subset of~$\C$ and suppose $\hat B$ is an entire
function of~$\C$.
Then, for any $\hat\ph\in\hat\gR_\Om$, the convolution product $\hat B*\hat\ph$
belongs to~$\hat\gR_\Om$;
its analytic continuation along a path~$\ga$ of $\C\setminus\Om$ starting
from a point $\ze_0 \in \D_{\rho(\Om)}$ and ending at a point~$\ze_1$ is the holomorphic germ
at~$\ze_1$ explicitly given by
\beglabel{eqcontAph}
\cont_\ga(\hat B*\hat\ph)(\ze) = 
\int_0^{\ze_0} \hat B(\ze-\xi) \hat\ph(\xi) \,\dd\xi 
+ \int_\ga \hat B(\ze-\xi) \hat\ph(\xi) \,\dd\xi 
+ \int_{\ze_1}^\ze \hat B(\ze-\xi) \hat\ph(\xi) \,\dd\xi
\elabel
for $\ze$ close enough to~$\ze_1$.
As a consequence, 
\beglabel{easycasez}
\ti B_0\in\C\{z\ii\}, \; \ti\ph_0 \in \ti\gR_\Om 
\ens\Longrightarrow\ens
\ti B_0 \ti\ph_0 \in \ti\gR_\Om.
\elabel
\end{lemma}


\begin{rem}	\label{remcautionbranch}
Formulas such as~\eqref{eqcontAph} require a word of caution: 
the value of $\hat B(\ze-\xi)$ is unambiguously defined whatever~$\ze$ and~$\xi$
are, but in the notation ``$\hat\ph(\xi)$'' it is understood that we are using the
appropriate branch of the possibily multiple-valued function~$\hat\ph$;
in such a formula, what branch we are using is clear from the context: 
\begin{enumerate}[$-$]
\item $\hat\ph$ is unambiguously defined in its disc of convergence~$D_0$ (centred
at~$0$) and the first integral thus makes sense for $\ze_0\in D_0$;
\item in the second integral
$\xi$ is moving along~$\ga$ which is a path of analytic continuation for~$\hat\ph$,
we thus consider the analytic continuation of~$\hat\ph$ along the piece of~$\ga$
between its origin and~$\xi$;
\item in the third integral, ``$\hat\ph$'' is to be understood as $\cont_\ga\hat\ph$, the germ at~$\ze_1$
resulting form the analytic continuation of~$\hat\ph$ along~$\ga$, this
integral then makes sense for any $\ze$ at a distance from~$\ze_1$ less than the radius of
convergence of $\cont_\ga\hat\ph$.
\end{enumerate}

Using a parametrisation $\ga\col[0,1]\to \C\setminus\Om$, with $\ga(0)=\ze_0$
and $\ga(1)=\ze_1$, and introducing the truncated paths
$\ga_s \defeq \ga_{|[0,s]}$ for any $s\in[0,1]$, the interpretation of the last two
integrals in~\eqref{eqcontAph} is
\begin{align}
\label{eqinterpretsec}
\int_\ga \hat B(\ze-\xi) \hat\ph(\xi) \,\dd\xi &\defeq
\int_0^1 \hat B(\ze-\ga(s)) (\cont_{\ga_s} \hat\ph)(\ga(s)) \ga'(s)\,\dd s,\\[1ex]
\label{eqinterpretthird}
\int_{\ze_1}^\ze \hat B(\ze-\xi) \hat\ph(\xi) \,\dd\xi &\defeq
\int_{\ze_1}^\ze \hat B(\ze-\xi) (\cont_\ga \hat\ph)(\xi) \,\dd\xi.
\end{align}
\end{rem}


\begin{proof}[Proof of Lemma~\ref{lemeasyconvol}]
The property~\eqref{easycasez} directly follows from the first statement:
write $\ti B_0 = a + \ti B$ and $\ti\ph_0 = b + \ti\ph$ with $a,b\in\C$ and $\ti
A,\ti\ph\in \zcz$ and apply Lemma~\ref{lemCVcase} to~$\ti B$.

To prove the first statement, we use a parametrisation
$\ga\col[0,1]\to\C\setminus\Om$ and the truncated paths $\ga_s \defeq
\ga_{|[0,s]}$:
we shall check that, for each $t\in[0,1]$, the formula
\beglabel{eqdefhatchit}
\hat\chi_t(\ze) \defeq 
\int_0^{\ze_0} \hat B(\ze-\xi) \hat\ph(\xi) \,\dd\xi 
+ \int_{\ga_t} \hat B(\ze-\xi) \hat\ph(\xi) \,\dd\xi
+ \int_{\ga(t)}^\ze \hat B(\ze-\xi) \hat\ph(\xi) \,\dd\xi
\elabel
(with the above conventions for the interpretation of ``$\hat\ph(\xi)$'' in the integrals)
defines a holomorphic germ at~$\ga(t)$ which is the analytic continuation of
$\hat B*\hat\ph$ along~$\ga_t$.

The holomorphic dependence of the integrals upon the parameter~$\ze$ is such that
$\ze \mapsto 
\int_0^{\ze_0} \hat B(\ze-\xi) \hat\ph(\xi) \,\dd\xi 
+
\int_{\ga_t} \hat B(\ze-\xi) \hat\ph(\xi) \,\dd\xi$
is an entire function of~$\ze$ and 
$\ze \mapsto \int_{\ga(t)}^\ze \hat B(\ze-\xi) \hat\ph(\xi) \,\dd\xi$
is holomorphic for~$\ze$ in the disc of convergence~$D_t$ of
$\cont_{\ga_t}\hat\ph$ (centred at~$\ga(t)$), 
we thus have a family of analytic elements $(\hat\chi_t,D_t)$, $t\in[0,1]$, along the path~$\ga$.

For $t$ small enough, the truncated path~$\ga_t$ is contained in~$D_0$; 
then, for $\ze\in D_0$,
the Cauchy theorem implies that $\hat\chi_t(\ze)$ coincides with $\hat
A*\hat\ph(\ze) = \int_0^\ze \hat B(\ze-\xi) \hat\ph(\xi) \,\dd\xi$
(since the rectilinear path $[0,\ze]$ is homotopic in~$D_0$ to the concatenation of
$[0,\ze_0]$, $\ga_t$ and $[\ga(t),\ze]$).


\begin{figure}
\begin{center}

\includegraphics[scale=1]{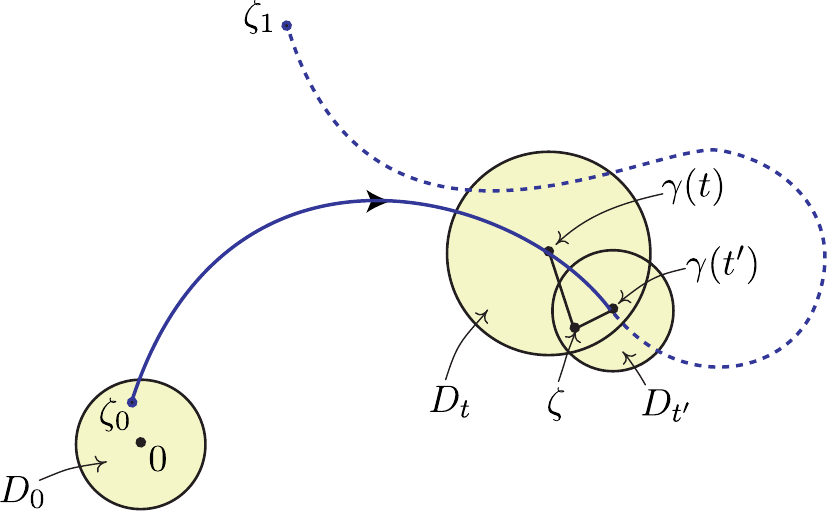} 

\bigskip 

\caption{%
Integration paths for the convolution in the easy case.}
\label{fig:Easyconv}

\end{center}
\end{figure}


For every $t\in[0,1]$, there exists $\eps>0$ such that
$\ga\big( (t-\eps,t+\eps)\cap[0,1] \big) \subset D_{t}$;
%
by compactness, we can thus find $N\in\N^*$ and $0=t_0 < t_1 < \cdots < t_N = 1$
so that $\ga\big( [t_j,t_{j+1}] \big) \subset D_{t_j}$ for every~$j$.
The proof will thus be complete if we check that, for any $t<t'$ in $[0,1]$,
\[
\ga\big( [t,t'] \big) \subset D_{t} 
\ens\Longrightarrow\ens
\hat\chi_{t} \equiv \hat\chi_{t'} \;\text{in $D_{t}\cap D_{t'}$}.
\]
This follows from the observation that, under the hypothesis
$\ga\big( [t,t'] \big) \subset D_{t}$,
\[
s\in[t,t'] \;\text{and}\; \xi\in D_{t}\cap D_s 
\ens\Longrightarrow\ens
\cont_{\ga_s}\hat\ph(\xi) = \cont_{\ga_{t}}\hat\ph(\xi),
\]
thus, when computing $\hat\chi_{t'}(\ze)$ with $\ze\in D_{t}\cap D_{t'}$, 
the third integral in~\eqref{eqdefhatchit} is
\[
\int_{\ga(t')}^\ze  \hat B(\ze-\xi) \cont_{\ga_{t'}}\hat\ph(\xi) \,\dd\xi
= \int_{\ga(t')}^\ze  \hat B(\ze-\xi) \cont_{\ga_{t}}\hat\ph(\xi) \,\dd\xi
\]
and, interpreting the second integral of~\eqref{eqdefhatchit} as
in~\eqref{eqinterpretsec}, we get
\begin{align*}
\hat\chi_{t'}(\ze) - \hat\chi_{t}(\ze) &= 
\int_{t}^{{t'}} \hat B(\ze-\ga(s)) \big( \cont_{\ga_s}\hat\ph \big)(\ga(s)) \ga'(s)\,\dd s  
+ \int_{\ga({t'})}^{\ga(t)}  \hat B(\ze-\xi)\big( \cont_{\ga_{t}}\hat\ph \big)(\xi) \,\dd\xi \\
&= \int_{t}^{{t'}} \hat B(\ze-\ga(s))\big( \cont_{\ga_{t}}\hat\ph \big)(\ga(s)) \ga'(s)\,\dd s 
+ \int_{\ga({t'})}^{\ga(t)}  \hat B(\ze-\xi)\big( \cont_{\ga_{t}}\hat\ph \big)(\xi) \,\dd\xi 
= 0
\end{align*}
(see Figure~\ref{fig:Easyconv}).
\end{proof}


\begin{rem}
Lemma~\ref{lemeasyconvol} can be used to prove the $\Om$-resurgence of certain
formal series solutions of non-linear functional equations---see
Section~\ref{sec:ResurIter} (with $\Om=2\pi\I\Z$) 
and \cite[\S8]{mouldSN} (with $\Om = \Z$).
\end{rem}
%


\section{Analytic continuation of a convolution product: an example}	\label{sec_contconvexa}


We now wish to consider the convolution of two $\Om$-continuable
holomorphic germs at~$0$ without assuming that any of
them extends to an entire function.
A first example will convince us that there is no hope to get stability under
convolution if we do not impose that $\Om$ be stable under addition.

Let $\om_1,\om_2\in\C^*$ and 
\[
\hat\ph(\ze) \defeq \frac{1}{\ze-\om_1}, \quad
\hat\psi(\ze) \defeq \frac{1}{\ze-\om_2}.
\]
Their convolution product is
\[
\hat\chi(\ze) \defeq \hat\ph*\hat\psi(\ze) = 
\int_0^\ze \frac{1}{(\xi-\om_1)(\ze-\xi-\om_2)} \,\dd\xi,
\qquad |\ze| < \min\big\{ |\om_1|,|\om_2| \big\}.
\]
The formula
\[
\frac{1}{(\xi-\om_1)(\ze-\xi-\om_2)}= \frac{1}{\ze-\om_1-\om_2} \left(
\frac{1}{\xi-\om_1} + \frac{1}{\ze-\xi-\om_2} 
\right)
\]
shows that, for any $\ze\neq\om_1+\om_2$ of modulus $< \min\big\{
|\om_1|,|\om_2| \big\}$, one can write
\begin{equation}
\hat\chi(\ze) = \frac{1}{\ze-\om_1-\om_2} \big( L_1(\ze) + L_2(\ze) \big),
\qquad L_j(\ze) \defeq \int_0^\ze \frac{\dd\xi}{\xi-\om_j}
\end{equation}
(with the help of the change of variable $\xi \mapsto \ze-\xi$ in the case of~$L_2$).

Removing the half-lines $\om_j[1,+\infty)$ from~$\C$, we obtain a cut
plane~$\De$ in which~$\hat\chi$ has a meromorphic continuation (since $[0,\ze]$ avoids the
points~$\om_1$ and~$\om_2$ for all $\ze\in\De$). 
We can in fact follow the meromorphic continuation of~$\hat\chi$ along any path
which avoids~$\om_1$ and~$\om_2$, because 
\[
L_j(\ze) = -\int_0^{\ze/\om_j} \frac{\dd\xi}{1-\xi} =
\Log\Big( 1-\frac{\ze}{\om_j} \Big) \in \hat\gR_{\{\om_j\}}
\]
(\cf example~\ref{exalog}).
We used the words ``meromorphic continuation'' and not ``analytic continuation''
because of the factor $\frac{1}{\ze-\om_1-\om_2}$.
The conclusion is thus only $\hat\chi \in \hat\gR_\Om$, with $\Om \defeq \{\om_1,\om_2,\om_1+\om_2\}$.

\smallskip

-- If $\om\defeq\om_1+\om_2\in \De$, the principal branch of~$\hat\chi$ (\ie its
meromorphic continuation to~$\De$) has a removable singularity at $\om$, because
$(L_1+L_2)(\om) = \int_0^\om \frac{\dd\xi}{\xi-\om_1} +
\int_0^\om \frac{\dd\xi}{\xi-\om_2} = 0$
in that case (by the change of variable $\xi \mapsto \om-\xi$ in one of the
integrals).
This is consistent with Lemma~\ref{lemstarshapedstb} (the set~$\De$ is clearly
star-shaped \wrt~$0$).
But it is easy to see that this does not happen for all the branches
of~$\hat\chi$: when considering all the paths~$\ga$ going from~$0$ to~$\om$ and
avoiding~$\om_1$ and~$\om_2$, we have
\[
\cont_\ga L_j(\om) = \int_\ga \frac{\dd\xi}{\xi-\om_j}, \qquad j=1,2,
\]
hence $\frac{1}{2\pi\I}\big(\cont_\ga L_1(\om) + \cont_\ga L_2(\om)\big)$ is the
sum of the winding numbers around~$\om_1$ and~$\om_2$ of the loop obtained by
concatenating~$\ga$ and the line segment $[\om,0]$;
elementary geometry shows that this sum of winding numbers can take any
integer value, but whenever this value is non-zero the corresponding
branch of~$\hat\chi$ does have a pole at~$\om$.

\smallskip

\begin{figure}
\begin{center}

\includegraphics[scale=1]{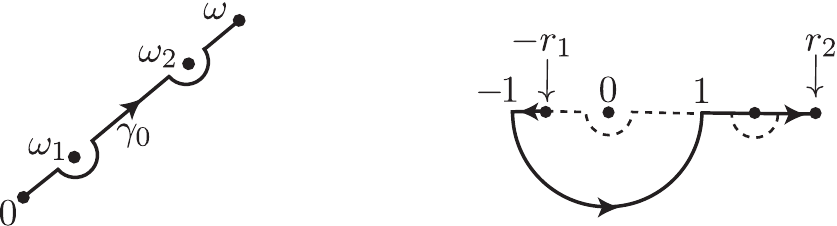}

\bigskip 


\caption{%
{Convolution of aligned poles.}
%
%
}
\label{figConvolPolesAl}

\end{center}
\end{figure}

-- The case $\om\notin\De$ is slightly different. Then we can write 
$\om_j = r_j\,\ee^{\I\th}$ with $r_1,r_2>0$ 
and consider the path~$\ga_0$ which follows the segment $[0,\om]$ except that it
circumvents $\om_1$ and~$\om_2$ by small half-circles travelled anti-clockwise
(notice that $\om_1$ and~$\om_2$ may coincide)---see the left part of Figure~\ref{figConvolPolesAl};
an easy computation yields
\[
\cont_{\ga_0} L_1(\om) = \int_{-r_1}^{-1} \frac{\dd\xi}{\xi}
+ \int_{\Ga_0} \frac{\dd\xi}{\xi}
+ \int_1^{r_2} \frac{\dd\xi}{\xi},
\]
where $\Ga_0$ is the half-circle from~$-1$ to~$1$ with radius~$1$ travelled
anti-clockwise (see the right part of Figure~\ref{figConvolPolesAl}),
hence $\cont_{\ga_0} L_1(\om) = \ln\frac{r_2}{r_1} + \I\pi$, 
similarly $\cont_{\ga_0} L_2(\om) = \ln\frac{r_1}{r_2} + \I\pi$, 
therefore $\cont_{\ga_0} L_1(\om) + \cont_{\ga_0} L_2(\om) = 2\pi\I$ is non-zero and this
again yields a branch of~$\hat\chi$ with a pole at~$\om$ (and infinitely many
others by using other paths than~$\ga_0$).
%
%
%

\medskip

In all cases, there are paths from~$0$ to $\om_1+\om_2$ which avoid~$\om_1$
and~$\om_2$ and which are not paths of analytic continuation for~$\hat\chi$.
This example thus shows that $\hat\gR_{\{\om_1,\om_2\}}$ is \emph{not} stable under
convolution: it contains $\hat\ph$ and~$\hat\psi$ but not $\hat\ph*\hat\psi$.

Now, whenever $\Om$ is not stable under addition, one can find
$\om_1,\om_2\in\Om$ such that $\om_1+\om_2\notin \Om$ and the previous example
then yields $\hat\ph,\hat\psi\in\hat\gR_\Om$ with $\hat\ph*\hat\psi \notin
\hat\gR_\Om$.
%


\section{Analytic continuation of a convolution product: the general case}
\label{sec_contconvolgen} 


\parag
The main result of this section is

\begin{thm}	\label{thmOmstbgROmstb}
Let $\Om$ be a non-empty closed discrete subset of~$\C$.
Then the space $\hat\gR_\Om$ is stable under convolution if and only if $\Om$ is stable under
addition. 
\end{thm}


The necessary and sufficient condition on~$\Om$ is satisfied by the typical
examples~$\Z$ or $2\pi\I\Z$, but also by~$\N^*$, $\Z+\I\Z$, $\N^*+\I\N$ or
$\{ m+n\sqrt{2} \mid m,n\in\N^* \}$ for instance.
An immediate consequence of Theorem~\ref{thmOmstbgROmstb} is
%

\begin{cor}   \label{corOmresursubalg}
Let $\Om$ be a non-empty closed discrete subset of~$\C$.
Then the space~$\ti\gR_\Om$ of $\Om$-resurgent formal series is a subalgebra of
$\C[[z\ii]]$ if and only if $\Om$ is stable under addition.
\end{cor}


The necessity of the condition on~$\Om$ was proved in
Section~\ref{sec_contconvexa}.
In the rest of this section we shall prove that the condition is
sufficient. However we shall restrict ourselves to the case where $0\in\Om$,
because this will allow us to give a simpler proof. The reader is referred to
\cite{stabiconv} for the proof in the general case.

\begin{figure}
\begin{center}

\includegraphics[scale=1]{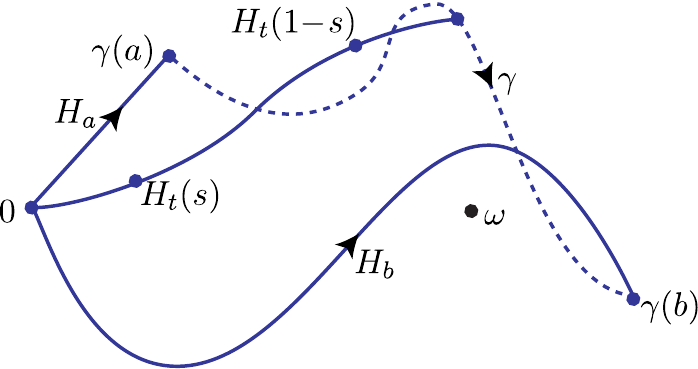} 

\bigskip 

\caption{%
\emph{Symmetric $\Om$-homotopy}
($H_a =$ initial path,
$H_b =$ final path,
$\ga =$ endpoint path~$\Ga_H$).}
\label{fig:SymHom}

\end{center}
\end{figure}

\parag
We thus fix~$\Om$ closed, discrete, containing~$0$ and stable under addition.
We begin with a new definition (see Figure~\ref{fig:SymHom}):

\begin{Def}	\label{defSymOmHom}
A continuous map
$H \col I \times J \to \C$, where $I=[0,1]$ and $J$ is a compact interval
of~$\R$, is called a \emph{symmetric $\Om$-homotopy} if, for each $t \in J$,
\[ s \in I \mapsto H_t(s) \defeq H(s,t) \]
defines a path which satisfies 
\begin{enumerate}
\item $H_t(0) = 0$,
\item $H_t\big( (0,1] \big) \subset \C\setminus\Om$,
\item $H_t(1) - H_t(s) = H_t(1-s)$ for every $s \in I$.
\end{enumerate}
We then call \emph{endpoint path} of~$H$ the path
\[
\Ga_H \col t \in J \mapsto H_t(1).
\]
Writing $J=[a,b]$, we call $H_a$ (\resp $H_b$) the \emph{initial path}
of~$H$ (\resp its \emph{final path}).
\end{Def}

The first two conditions imply that each path~$H_t$ is a path of analytic
continuation for any $\hat\ph\in\hat\gR_\Om$, in view of
Remark~\ref{reminitialpoint}.

We shall use the notation $H_{t|s}$ for the truncated paths $(H_t)_{|[0,s]}$,
$s\in I$, $t\in J$ (analogously to what we did when commenting
Lemma~\ref{lemeasyconvol}).
Here is a technical statement we shall use:
\begin{lemma}	\label{lemtechnic}
For a symmetric $\Om$-homotopy~$H$ defined on $I\times J$,
there exists $\de>0$ such that, for any
$\hat\ph\in\hat\gR_\Om$ and $(s,t) \in I\times J$, the radius of convergence of the
holomorphic germ $\cont_{H_{t|s}}\hat\ph$ at~$H_t(s)$ is at least~$\de$.
\end{lemma}


\begin{proof}
Let $\rho = {\rho(\Om)}$ (\cf Remark~\ref{reminitialpoint}).
Consider %
\[
U\defeq \big\{\, (s,t) \in I\times J \mid H\big( [0,s]\times\{t\} \big) 
\subset\Ddem \,\},
\qquad K \defeq I\times J \setminus U.
\]
Writing $K = 
\big\{\, (s,t) \in I\times J \mid \exists s'\in[0,s] \;\text{s.t.}\;
H(s',t) \in \C\setminus\Ddem \,\}$, 
we see that $K$ is a compact subset of $I\times J$ which is contained in
$(0,1]\times J$.
Thus $H(K)$ is a compact subset of $\C\setminus\Om$,
and $\de \defeq \min\big\{ \dist\big( H(K), \Om \big), \rho/2 \big\} > 0$.
Now, for any $s$ and~$t$, 
\begin{enumerate}[--]
\item either $(s,t)\in U$, then 
the truncated path~$H_{t|s}$ lies in~$\Ddem$, hence
$\cont_{H_{t|s}}\hat\ph$ is a holomorphic germ
at~$H_t(s)$ with radius of convergence $\ge \de$;
\item or $(s,t)\in K$, and then $\dist(H_t(s),\Om)\ge\de$, which yields the same
conclusion for the germ $\cont_{H_{t|s}}\hat\ph$.
\end{enumerate}
\end{proof}


The third condition in Definition~\ref{defSymOmHom} means that each path~$H_t$
is symmetric \wrt\ its midpoint~$\demi H_t(1)$.
Here is the motivation behind this requirement:
\begin{lemma}	\label{lemcontsymOm}
Let $\ga \col [0,1] \to \C \setminus \Om$ be a path such that $\ga(0) \in
\D_{\rho(\Om)}$ (\cf Remark~\ref{reminitialpoint}).
If there exists a symmetric $\Om$-homotopy whose endpoint path coincides
with~$\ga$ and whose initial path is contained in~$\D_{\rho(\Om)}$, 
then any convolution product $\hat\ph*\hat\psi$ with $\hat\ph,
\hat\psi\in\hat\gR_\Om$ can be analytically continued along~$\ga$.
\end{lemma}

\begin{proof}
We assume that $\ga = \Ga_H$ with a symmetric $\Om$-homotopy~$H$ defined on $I\times J$. 
Let $\hat\ph, \hat\psi\in\hat\gR_\Om$ and, for $t \in J$, consider the formula
\beglabel{eqpremform}
\hat\chi_t(\ze) = \int_{H_t} \hat\ph(\xi) \hat\psi(\ze-\xi)\,\dd\xi
+ \int_{\ga(t)}^\ze \hat\ph(\xi) \hat\psi(\ze-\xi)\,\dd\xi
\elabel
(recall that $\ga(t) = H_t(1)$). We shall check that $\hat\chi_t$ is a
well-defined holomorphic germ at~$\ga(t)$ and that it provides the analytic
continuation of $\hat\ph*\hat\psi$ along~$\ga$.

\medskip

\noindent \textbf{a)}
The idea is that when $\xi$ moves along~$H_t$, $\xi = H_t(s)$ with $s \in I$,
we can use for ``$\hat\ph(\xi)$'' the analytic continuation of~$\hat\ph$
along the truncated path~$H_{t|s}$;
correspondingly, if~$\ze$ is close to~$\ga(t)$, then $\ze-\xi$ is close to
$\ga(t)-\xi = H_t(1) - H_t(s) = H_t(1-s)$, thus for
``$\hat\psi(\ze-\xi)$'' we can use the analytic continuation of~$\hat\psi$
along~$H_{t|1-s}$.
In other words, setting $\ze=\ga(t)+\sig$, we wish to interpret~\eqref{eqpremform} as
\begin{multline}	\label{eqsecform}
\hat\chi_t(\ga(t)+\sig) \defeq \int_0^1 
(\cont_{H_{t|s}}\hat\ph)(H_t(s)) (\cont_{H_{t|1-s}}\hat\psi)(H_t(1-s)+\sig)
H_t'(s)\,\dd s \\
+ \int_0^1 (\cont_{H_t}\hat\ph)(\ga(t)+u\sig) \hat\psi((1-u)\sig)\sig\,\dd u
\end{multline}
(in the last integral, we have performed the change variable $\xi =
\ga(t)+u\sig$; it is the germ of~$\hat\psi$ at the origin that we use
there).

Lemma~\ref{lemtechnic} provides $\de>0$ such that, by regular dependence of the
integrals upon the parameter~$\sig$, the \rhs\ of~\eqref{eqsecform} is
holomorphic for $|\sig|<\de$.
We thus have a family of analytic elements $(\hat\chi_t,D_t)$, $t \in J$, with
$D_t \defeq \{\, \ze\in\C \mid |\ze-\ga(t)| < \de \,\}$.

\medskip

\noindent \textbf{b)}
For $t$ small enough, the path $H_t$ is contained in~$\D_{\rho(\Om)}$ which is open and simply
connected; then, for $|\ze|$ small enough, the line segment $[0,\ze]$ and
the concatenation of~$H_t$ and $[\ga(t),\ze]$ are homotopic in~$\D_{\rho(\Om)}$, hence the
Cauchy theorem implies $\hat\chi_t(\ze) = \hat\ph*\hat\psi(\ze)$.

\medskip

\noindent \textbf{c)}
By uniform continuity, there exists $\eps>0$ such that, for any $t_0,t \in J$,
\beglabel{inequnifcont}
|t-t_0| \le \eps 
\quad\Longrightarrow\quad
|H_t(s)-H_{t_0}(s)| < \de/2 \quad \text{for all $s \in I$}.
\elabel
To complete the proof, we check that, for any $t_0, t$ in~$J$ such that
$t_0\le t \le t_0+\eps$, we have $\hat\chi_{t_0} \equiv \hat\chi_{t}$ in
$D\big(\ga(t_0),\de/2)$ (which is contained in $D_{t_0}\cap D_{t}$).

Let $t_0,t\in J$ be such that $t_0\le t \le t_0+\eps$ and let $\ze \in D\big(\ga(t_0),\de/2)$.
By Lemma~\ref{lemtechnic} and~\eqref{inequnifcont}, we have for every $s \in I$
\begin{align*}
&\cont_{H_{t|s}}\hat\ph\big(H_t(s)\big) = \cont_{H_{t_0|s}}\hat\ph\big(H_t(s)\big), \\
&\cont_{H_{t|1-s}}\hat\psi\big(\ze-H_t(s)\big) =
\cont_{H_{t_0|1-s}}\hat\psi\big(\ze-H_t(s)\big)
\end{align*}
(for the latter identity, write
$\ze-H_t(s) = H_t(1-s) + \ze - \ga(t) = H_{t_0}(1-s) + \ze-\ga(t_0) +
H_{t_0}(s)-H_t(s)$, thus this point belongs to $D\big(H_t(1-s),\de)
\cap D\big(H_{t_0}(1-s),\de)$).
Moreover, $[\ga(t),\ze] \subset D\big(\ga(t_0),\de/2)$ by convexity,
hence $\cont_{H_{t}}\hat\ph \equiv \cont_{H_{t_0}}\hat\ph$ on this line segment,
and we can write
\begin{multline*}	
\hat\chi_t(\ze) = \int_0^1 
(\cont_{H_{t_0|s}}\hat\ph)(H_t(s)) (\cont_{H_{t_0|1-s}}\hat\psi)(\ze-H_t(s))
H_t'(s)\,\dd s \\
+ \int_{\ga(t)}^\ze (\cont_{H_{t_0}}\hat\ph)(\xi) \hat\psi(\ze-\xi)\,\dd\xi.
\end{multline*}
We then get $\hat\chi_{t_0}(\ze) = \hat\chi_{t}(\ze)$ from the Cauchy theorem by
means of the homotopy induced by~$H$ between 
the concatenation of~$H_{t_0}$ and $[\ga(t_0),\ze]$ 
and the concatenation of~$H_{t}$ and $[\ga(t),\ze]$.
\end{proof}


\begin{rem}
\label{remIntFinalPath}
With the notation of Definition~\ref{defSymOmHom}, when the initial path~$H_a$ is a
line segment contained in~$\D_{\rho(\Om)}$, the final path~$H_b$ is what
\'Ecalle calls a ``symmetrically contractible path'' in \cite{Eca81}.
The proof of Lemma~\ref{lemcontsymOm} shows that the analytic continuation of
$\hat\ph*\hat\psi$ until the endpoint $H_b(1)=\Ga_H(b)$ can be computed by the usual
integral taken over~$H_b$:
\beglabel{eqcontgaphpsiHb}
\cont_{\ga} (\hat\ph*\hat\psi)(\ze) =
\int_{H_b} \hat\ph(\xi) \hat\psi(\ze-\xi) \, \dd\xi,
\qquad
\ga = \Ga_H, \ens \ze = \ga(b)
\elabel
(with appropriate interpretation, as in~\eqref{eqsecform}).
However, it usually cannot be computed as the same integral over
$\ga=\Ga_H$ itself, even when the latter integral is well-defined).
\end{rem}


\parag
In view of Lemma~\ref{lemcontsymOm}, the proof of Theorem~\ref{thmOmstbgROmstb}
will be complete if we prove the following purely geometric result:


\begin{lemma}	\label{lemkey}
For any path $\ga \col I= [0,1]\to \C\setminus\Om$ such that $\ga(0)\in \D^*_ {\rho(\Om)}$,
%
there exists a symmetric $\Om$-homotopy~$H$ on $I\times I$ whose endpoint path is~$\ga$ and
whose initial path is a line segment, \ie $\Ga_H=\ga$ and $H_0(s)\equiv s\ga(0)$.
\end{lemma}


\begin{proof}
%
%
Assume that $\ga$ is given as in the hypothesis of Lemma~\ref{lemkey}.
We are looking for a symmetric $\Om$-homotopy whose initial path is imposed: it
must be
\[
s \in I \mapsto H_0(s) \defeq s\ga(0),
\]
which satisfies the three requirements of Definition~\ref{defSymOmHom} at $t=0$:
\begin{enumerate}[(i)]
\item
$H_0(0)=0$, 
\item
$H_0\big( (0,1] \big) \subset \C\setminus \Om$, 
\item
$H_0(1) - H_0(s) = H_0(1-s)$ for every $s \in I$.
\end{enumerate}
The idea is to define a family of maps $(\Psi_t)_{t\in[0,1]}$ so that
\beglabel{eqdefHsPsis}
H_t(s) \defeq \Psi_t\big( H_0(s) \big), \qquad s \in I,
\elabel
yield the desired homotopy. For that, it is sufficient that $(t,\ze) \in[0,1]
\times \C \mapsto \Psi_t(\ze)$ be continuously differentiable (for the structure
of real two-dimensional vector space of~$\C$), $\Psi_0=\id$ and, for each $t\in
[0,1]$,
\begin{enumerate}[(i')]
\item
$\Psi_t(0) = 0$,
\item
$\Psi_t(\C\setminus\Om) \subset \C\setminus\Om$,
\item
$\Psi_t\big( \ga(0) - \ze \big) = \Psi_t\big( \ga(0) \big) - \Psi_t(\ze)$ for
all $\ze\in\C$,
\item
$\Psi_t\big(\ga(0)\big) = \ga(t)$.
\end{enumerate}
In fact, the properties (i')--(iv') ensure that any initial path $H_0$
satisfying (i)--(iii) and ending at~$\ga(0)$ produces through~\eqref{eqdefHsPsis}
a symmetric $\Om$-homotopy whose endpoint path is~$\ga$. Consequently, we may
assume without loss of generality that~$\ga$ is~$C^1$ on $[0,1]$ 
(then, if $\ga$ is only piecewise~$C^1$, we just need to concatenate the
symmetric $\Om$-homotopies associated with the various pieces).

The maps~$\Psi_t$ will be generated by the flow of a non-autonomous vector
field $X(\ze,t)$ associated with~$\ga$ that we now define.
We view $(\C,|\,\cdot\,|)$ as a real $2$-dimensional Banach space and pick%
\footnote{
For instance pick a $C^1$ function $\ph_0 \col \R\to[0,1]$ such that 
$\{\, x\in\R \mid \ph_0(x) = 1 \,\} = \{0\}$ and $\ph_0(x)=0$ for $|x|\ge1$,
and a bijection $\om\col\N\to\Om$;
then set $\de_k \defeq \dist\big( \om(k), \Om\setminus\{\om(k)\} \big) >0$ and
$\sig(\ze) \defeq \sum_k \ph_0\big( \frac{4|\ze-\om(k)|^2}{\de_k^2} \big)$:
for each $\ze\in\C$ there is at most one non-zero term in this series 
(because $k\neq\ell$, $|\ze-\om(k)| < \de_k/2$ and $|\ze-\om(\ell)| <
\de_\ell/2$ would imply $|\om(k)-\om(\ell)|< (\de_k+\de_\ell)/2$, which would
contradict $|\om(k)-\om(\ell)| \ge \de_k$ and~$\de_\ell$),
thus $\sig$ is $C^1$, takes its values in $[0,1]$ and satisfies
$\{\, \ze\in\C \mid \sig(\ze) = 1 \,\} = \Om$, therefore $\eta\defeq 1-\sig$
will do.
%
%
}
a $C^1$ function $\eta \col \C \to [0,1]$ such that
\[
\{\,\ze\in\C \mid \eta(\ze) = 0 \,\} = \Om.
\]
Observe that $D(\ze,t) \defeq \eta(\ze) + \eta\big( \ga(t)-\ze \big)$ defines
a $C^1$ function of $(\ze,t)$
which satisfies
\[
D(\ze,t) > 0 \quad \text{for all $\ze\in\C$ and $t\in[0,1]$}
\]
\emph{because $\Om$ is stable under addition};
indeed, $D(\ze,t) = 0$ would imply $\ze\in\Om$ and $\ga(t)-\ze\in\Om$, hence
$\ga(t)\in\Om$, which would contradict our assumptions.
Therefore, the formula
\begin{equation}	\label{eqdefXnonaut}
X(\ze,t) \defeq \frac{\eta(\ze)}{\eta(\ze) + \eta\big( \ga(t)-\ze \big)} \ga'(t)
\end{equation}
defines a non-autonomous vector field, which is continuous in $(\ze,t)$ on
$\C\times[0,1]$, $C^1$ in~$\ze$ and has its partial derivatives continuous in
$(\ze,t)$.
The Cauchy-Lipschitz theorem on the existence and uniqueness of solutions to
differential equations applies to $\frac{\dd\ze}{\dd t} = X(\ze,t)$:
for every $\ze\in\C$ and $t_0\in[0,1]$ there is a unique solution $t\mapsto \Phi^{t_0,t}(\ze)$
such that $\Phi^{t_0,t_0}(\ze)=\ze$. 
The fact that the vector field~$X$ is bounded implies that $\Phi^{t_0,t}(\ze)$ is
defined for all $t\in[0,1]$
and the classical theory guarantees that $(t_0,t,\ze)\mapsto\Phi^{t_0,t}(\ze)$
is $C^1$ on $[0,1]\times[0,1]\times\C$.

Let us set $\Psi_t \defeq \Phi^{0,t}$ for $t\in[0,1]$ and check that this family of
maps satisfies (i')--(iv').
We have
\begin{gather}
\label{eqvanishX}
X(\om,t) = 0 \quad\text{for all $\om\in\Om$,}\\
\label{eqsymX}
X\big( \ga(t)-\ze, t \big) = \ga'(t) - X(\ze,t)
\quad\text{for all $\ze\in\C$}
\end{gather}
for all $t\in[0,1]$
(by the very definition of~$X$). Therefore
\begin{itemize}
\item
(i') and (ii') follow from~\eqref{eqvanishX} which yields
$\Phi^{t_0,t}(\om)=\om$ for every $t_0$ and~$t$, whence $\Psi_t(0)=0$ since
$0\in\Om$,
and from the non-autonomous flow property $\Phi^{t,0}\circ\Phi^{0,t}=\id$
(hence $\Psi_t(\ze)=\om$ implies $\ze=\Phi^{t,0}(\om)=\om$);
\item
(iv') follows from the fact that $X\big( \ga(t),t \big)=\ga'(t)$,
by~\eqref{eqvanishX} and~\eqref{eqsymX} with $\ze=0$, using again that
$0\in\Om$, hence $t\mapsto \ga(t)$ is a solution of~$X$;
\item
(iii') follows from~\eqref{eqsymX}: for any solution $t\mapsto\ze(t)$, the curve
$t\mapsto \xi(t)\defeq \ga(t)-\ze(t)$ satisfies $\xi(0) = \ga(0)-\ze(0)$ and
$\xi'(t) = \ga'(t) - X\big( \ze(t),t \big) = X\big( \xi(t),t \big)$,
hence it is a solution:
$\xi(t) = \Psi_t\big(\ga(0)-\ze(0)\big)$.
\end{itemize}

As explained above, formula~\eqref{eqdefHsPsis} thus produces the desired
symmetric $\Om$-homotopy.
\end{proof}


\parag
\noindent \emph{Note on this section:}
The presentation we adopted is influenced by \cite{CNP} (the example of
Section~\ref{sec_contconvexa} is taken from this book).
Lemma~\ref{lemkey}, which is the key to the proof of
Theorem~\ref{thmOmstbgROmstb}
and which essentially relies on the use of the flow of the non-autonomous vector
field~\eqref{eqdefXnonaut}, arose as an attempt to understand a related but more
complicated (somewhat obscure!) construction which can be found in an appendix
of \cite{CNP}.
See \cite{Eca81} and~\cite{Y_Ou} for other approaches to the stability under
convolution of the space of resurgent functions.

\begin{figure}
\begin{center}

\includegraphics[scale=1]{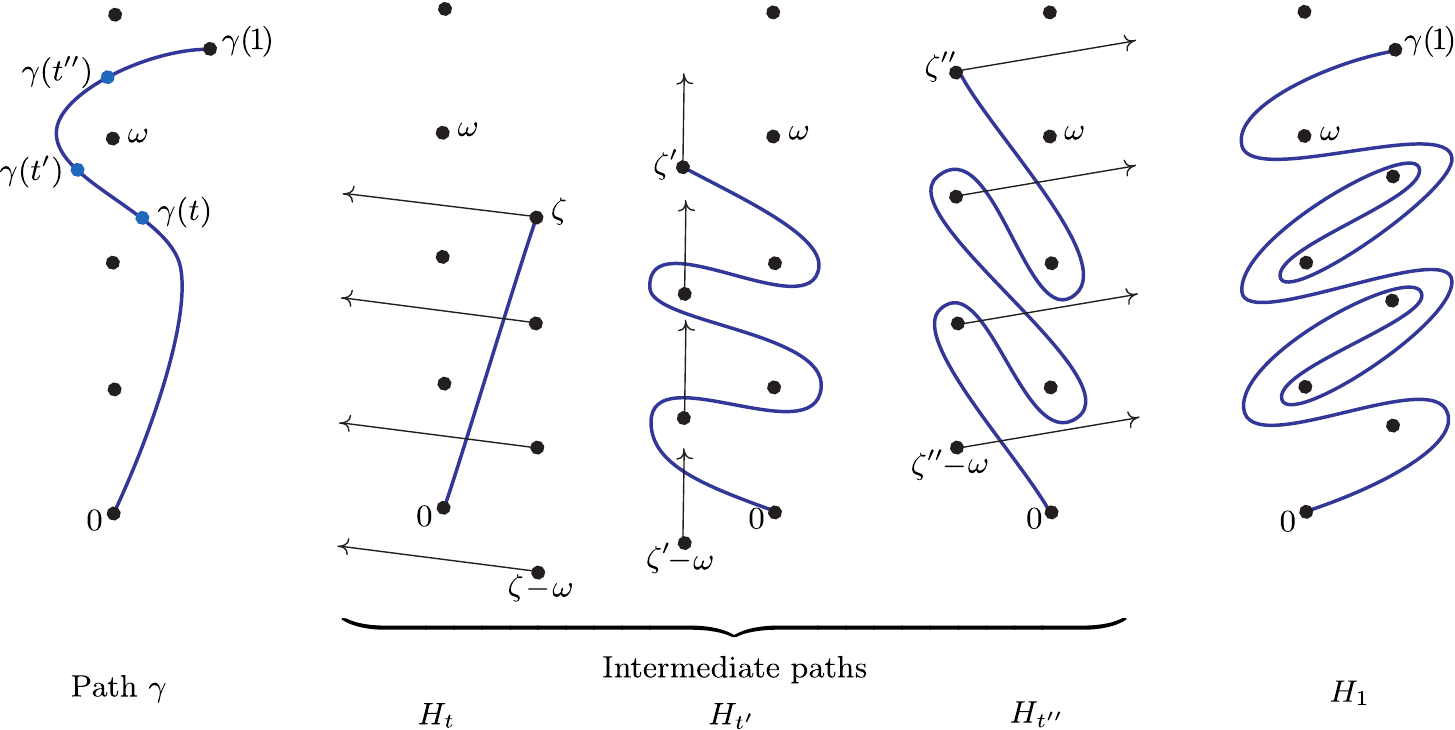}


\bigskip

\caption{From~$\ga$ to the integration path~$H_1$ used for 
$\cont_{\ga} (\hat\ph*\hat\psi)(\ga(1))$.}
\label{fignails}

\end{center}
\end{figure}

For the proof of Lemma~\ref{lemkey}, according to \cite{Eca81} and
\cite{CNP}, one can visualize the realization of a given path~$\ga$ as
the enpoint path~$\Ga_H$ of a symmetric $\Om$-homotopy as follows:
\label{secnails}
Let a point $\ze=\ga(t)$ move along~$\ga$ (as $t$ varies from~$0$ to~$1$) and
remain connected to~$0$ by an extensible thread, 
with moving nails pointing downwards at each point of~$\ze-\Om$,
while fixed nails point upwards at each point of~$\Om$
(imagine for instance that the first nails are fastened to a moving rule and the
last ones to a fixed rule).
As $t$ varies, the thread is progressively stretched but it has to meander
between the nails.
The path~$H_1$ used as integration path for 
$\cont_{\ga} (\hat\ph*\hat\psi)(\ga(1))$
in formula~\eqref{eqcontgaphpsiHb} is given by the
thread in its final form, when~$\ze$ has reached the extremity
of~$\ga$; the paths~$H_t$ correspond to the thread at intermediary
stages.
See Figure~\ref{fignails} (or Figure~5 of \cite{kokyu}).
The point is that none of the moving nails $\ze-\om' \in \ze-\Om$ will ever
collide with a fixed nail $\om''\in\Om$ because we assumed that $\ga$
avoids $\{ \om'+\om''  \} \subset \Om$.

%
%

%


\parag
\noindent \emph{Asymmetric version of the result.}
Theorem~\ref{thmOmstbgROmstb} admits a useful generalization,
concerning the convolution product of two resurgent germs which do not
belong to the same space of $\Om$-continuable germs:


\begin{thm}	\label{thmAsymConv}
Let $\Om_1$ and~$\Om_2$ be non-empty closed discrete subsets
of~$\C$. Let
\[
\Om \defeq \Om_1 \cup \Om_2 \cup (\Om_1+\Om_2),
\]
where $\Om_1+\Om_2 \defeq \{\, \om_1 + \om_2 \mid 
\om_1 \in \Om_1, \; \om_2 \in \Om_2 \,\}$.
If $\Om$ is closed and discrete, then
\[
\hat\ph \in \hat\gR_{\Om_1} \ens\text{and}\ens
\hat\psi \in \hat\gR_{\Om_2}
\quad \Longrightarrow \quad
\hat\ph * \hat\psi \in \hat\gR_\Om.
\]
\end{thm}


We shall content ourselves with giving hints about the proof when
both~$\Om_1$ and~$\Om_2$ are assumed to contain~$0$, 
in which case
\[
\Om = \Om_1 + \Om_2
\]
since both~$\Om_1$ and~$\Om_2$ are contained in $\Om_1+\Om_2$
(the general case is obtained by adapting the arguments of \cite{stabiconv}).
Assuming this, we generalize Definition~\ref{defSymOmHom} and
Lemma~\ref{lemcontsymOm}:


\begin{Def}	\label{defAsymHom}
A continuous map
$H \col I \times J \to \C$, where $I=[0,1]$ and $J$ is a compact interval
of~$\R$, is called an \emph{$(\Om_1,\Om_2)$-homotopy} if, for each $t \in J$,
the paths $s \in I \mapsto H_t(s) \defeq H(s,t)$ and 
$s \in I \mapsto H^*_t(s) \defeq H_t(1) - H_t(1-s)$
satisfy
\begin{enumerate}
\item $H_t(0) = 0$,
\item $H_t\big( (0,1] \big) \subset \C\setminus\Om_1$,
\item $H^*_t\big( (0,1] \big) \subset \C\setminus\Om_2$.
\end{enumerate}
We then call $t \in J \mapsto H_t(1)$ the \emph{endpoint path} of~$H$.
\end{Def}


\begin{lemma}	\label{lemcontASymOm}
Let $\ga \col [0,1] \to \C \setminus \Om$ be a path such that $\ga(0) \in
\D_{\rho(\Om_1)} \cap \D_{\rho(\Om_2)}$.
Suppose that there exists an $(\Om_1,\Om_2)$-homotopy whose endpoint
path coincides with~$\ga$ and such that $H_0(I) \subset \D_{\rho(\Om_1)}$
and $H^*_0(I) \subset \D_{\rho(\Om_2)}$.
Then any convolution product $\hat\ph*\hat\psi$ with $\hat\ph \in
\hat\gR_{\Om_1}$ and $\hat\psi\in\hat\gR_{\Om_2}$ can be analytically
continued along~$\ga$.
\end{lemma}

\begin{proof}[Idea of the proof of Lemma~\ref{lemcontASymOm}]
Mimick the proof of Lemma~\ref{lemcontsymOm}, replacing the \rhs\
of~\eqref{eqsecform} with
\begin{multline*}
\int_0^1 
(\cont_{H_{t|s}}\hat\ph)(H_t(s)) (\cont_{H^*_{t|1-s}}\hat\psi)(H^*_t(1-s)+\sig)
H_t'(s)\,\dd s \\
+ \int_0^1 (\cont_{H_t}\hat\ph)(\ga(t)+u\sig) \hat\psi((1-u)\sig)\sig\,\dd u
\end{multline*}
and showing that this expression is the value at $\ga(t)+\sig$ of a
holomorphic germ, which is
$\cont_{\ga|t} (\hat\ph*\hat\psi)$.
\end{proof}


To conclude the proof of Theorem~\ref{thmAsymConv}, it is thus
sufficient to show
\begin{lemma}	\label{lemASymKey}
For any path $\ga \col [0,1]\to \C$ such that 
$\ga(0)\in \D^*_ {\rho(\Om_1)} \cap \D^*_ {\rho(\Om_2)}$
and $\ga\big( (0,1] \big) \subset \C\setminus\Om$,
%
%
there exists an $(\Om_1,\Om_2)$-homotopy~$H$ on $I\times [0,1]$ whose
endpoint path is~$\ga$ 
and such that $H_0(s)\equiv s\ga(0)$.
\end{lemma}

Indeed, if this lemma holds true, then all such paths~$\ga$ will be,
by virtue of Lemma~\ref{lemcontASymOm},
paths of analytic continuation for our convolution products
$\hat\ph*\hat\psi$, which is the content of Theorem~\ref{thmAsymConv}.

\begin{proof}[Idea of the proof of Lemma~\ref{lemASymKey}]
It is sufficient to construct a family of maps $(\Psi_t)_{t\in[0,1]}$ such that
$(t,\ze) \in[0,1] \times \C \mapsto \Psi_t(\ze) \in \C$ be continuously
differentiable (for the structure of real two-dimensional vector space
of~$\C$), $\Psi_0=\id$ and, for each $t\in [0,1]$,
\begin{enumerate}[(i')]
\item
$\Psi_t(0) = 0$,
\item
$\Psi_t(\C\setminus\Om_1) \subset \C\setminus\Om_1$,
\item
the map $\ze \in \C \mapsto \Psi^*_t(\ze) \defeq \ga(t) - \Psi_t(\ze)$ 
satisfies
$\Psi^*_t(\C\setminus\Om_2) \subset \C\setminus\Om_2$,
\item
$\Psi_t\big(\ga(0)\big) = \ga(t)$.
\end{enumerate}
Indeed, the formula
$H_t(s) \defeq \Psi_t\big( s \ga(0) \big)$
then yields the desired homotopy, with
$H^*_t(s) = \Psi^*_t\big( (1-s)\ga(0) \big)$.

As in the proof of Lemma~\ref{lemkey}, the maps~$\Psi_t$ will be
generated by the flow of a non-autonomous vector field associated
with~$\ga$.
We view $(\C,|\,\cdot\,|)$ as a real $2$-dimensional Banach space and pick
$C^1$ functions $\eta_1,\eta_2 \col \C \to [0,1]$ such that
\[
\{\,\ze\in\C \mid \eta_j(\ze) = 0 \,\} = \Om_j,
\qquad j=1,2.
\]
Observe that $D(\ze,t) \defeq \eta_1(\ze) + \eta_2\big( \ga(t)-\ze \big)$ defines
a $C^1$ function of $(\ze,t)$
which satisfies
\[
D(\ze,t) > 0 \quad \text{for all $\ze\in\C$ and $t\in[0,1]$,}
\]
since $D(\ze,t) = 0$ would imply $\ze\in\Om_1$ and $\ga(t)-\ze\in\Om_2$, hence
$\ga(t)\in\Om_1+\Om_2$, which would contradict our assumptions.
Therefore, the formula
\begin{equation}	
X(\ze,t) \defeq \frac{\eta_1(\ze)}{\eta_1(\ze) + \eta_2\big( \ga(t)-\ze \big)} \ga'(t)
\end{equation}
defines a non-autonomous vector field and
the Cauchy-Lipschitz theorem applies to $\frac{\dd\ze}{\dd t} = X(\ze,t)$:
for every $\ze\in\C$ and $t_0\in[0,1]$ there is a unique solution
$t\in[0,1]\mapsto \Phi_X^{t_0,t}(\ze)$ such that $\Phi_X^{t_0,t_0}(\ze)=\ze$;
the flow map $(t_0,t,\ze)\mapsto\Phi_X^{t_0,t}(\ze)$
is $C^1$ on $[0,1]\times[0,1]\times\C$.

Setting $\Psi_t \defeq \Phi_X^{0,t}$ for $t\in[0,1]$, one can check that this family of
maps satisfies (i')--(iv') by mimicking the arguments in the proof of
Lemma~\ref{lemkey} and using the fact that the corresponding family of
maps~$(\Psi^*_t)$ in~(iii') can be obtained from the identity
\[
\ga(t) - \Phi_X^{0,t}(\ze) = 
\Phi_{X^*}^{0,t}\big( \ga(0)-\ze \big),
\]
where we denote by $(t_0,t,\ze)\mapsto\Phi_{X^*}^{t_0,t}(\ze)$ the
flow map of the non-autonomous vector field
\[
X^*(\ze,t) \defeq \ga'(t) - X\big( \ga(t) - \ze, t \big) 
= \frac{\eta_2(\ze)}{\eta_1\big( \ga(t)-\ze \big) + \eta_2(\ze)} \ga'(t).
\]
\end{proof}


\section{Non-linear operations with resurgent formal series}
\label{secNLopRES}


From now on, we give ourselves a non-empty closed discrete subset~$\Om$ of~$\C$
which is stable under addition.

We already mentioned the stability of~$\ti\gR_\Om$ under certain linear
difference/differential operators in Lemma~\ref{lemelemstab}. 
Now, with our assumption that $\Om$ is stable under addition, we can
obtain the stability of $\Om$-resurgent formal series under the non-linear
operations which were studied in Sections~\ref{secnonlinsumm}
and~\ref{secgrponesummdiffeos}.
However this requires quantitative estimates for iterated convolutions whose
proof is beyond the scope of the present text, we thus quote without proof the
following


\begin{lemma}	\label{lemestimiterconv}
Let $\ga$ be a path of $\C\setminus\Om$ starting from a point $\ze_0
\in\D_{\rho(\Om)}$ and ending at a point~$\ze_1$. Let $R>0$ be such that
$\ov{D(\ze_1,R)} \subset \C\setminus\Om$.
Then there exist a positive number~$L$ and a set~$\gC$ of paths parametrized by
$[0,1]$ and contained in $\D_L\setminus\Om$ such that,
for every $\hat\ph\in\hat\gR_\Om$, 
the number
\[
\norm{\hat\ph}_\gC \defeq \sup_{\ti\ga\in\gC} \abs*{\cont_{\ti\ga}\hat\ph\big(\ti\ga(1)\big)}
\]
is finite, and there exist $A,B>0$ such that,
for every $k\ge1$ and $\hat\ph,\hat\psi\in\hat\gR_\Om$, the iterated convolution products
\[
\hat\ph^{*k} \defeq \underbrace{\hat\ph * \cdots * \hat\ph}_\text{$k$ factors} 
\]
and $\hat\psi*\hat\ph^{*k}$
(which admit analytic continuation along~$\ga$, according to Theorem~\ref{thmOmstbgROmstb})
satisfy
\begin{align*}
\abs{\cont_\ga \hat\ph^{*k}(\ze) } &\le A \frac{B^k}{k!}
\big(\norm{\hat\ph}_\gC\big)^k, \\[1ex]
\abs{\cont_\ga (\hat\psi * \hat\ph^{*k})(\ze) } &\le A \frac{B^k}{k!} 
\norm{\hat\psi}_\gC \big(\norm{\hat\ph}_\gC\big)^k,
\end{align*}
for every $\ze \in \ov{D(\ze_1,R)}$.
\end{lemma}

The proof can be found in~\cite{NL_resur}.
Taking this result for granted, we can show


\label{secSubstResCV}

\begin{thm}	\label{thmresOmstbNL}
Suppose that $\ti\ph(z), \ti\psi(z), \ti\chi(z) \in \ti\gR_\Om$
and that $\ti\chi(z)$ has no constant term.
Let $H(t)\in\C\{t\}$.
Then 
\beglabel{eqdeflamu}
\ti\psi\circ(\id+\ti\ph) \in\ti\gR_\Om, 
\qquad
H\circ\ti\chi \in\ti\gR_\Om.
\elabel
\end{thm}


\begin{proof}
We can write 
$\ti\ph = a + \ti\ph_1$, $\ti\psi = b + \ti\psi_1$, where $a,b\in\C$
and $\ti\ph_1$ and~$\ti\psi_1$ have no constant term.
With notations similar to those of the proof of Theorem~\ref{thmcompatsumcirc},
we write the first formal series in~\eqref{eqdeflamu} as $b+\ti\la(z)$ and the second one as $c+\ti\mu(z)$,
where $c = H(0)$.
Since $\ti\la = (T_a\ti\psi_1)\circ(\id+\ti\ph_1)$, where
$T_a\ti\psi_1$ is $\Om$-resurgent (by Lemma~\ref{lemelemstab}) and has no constant term, 
we see that it is sufficient to deal with the case $a=b=0$; from now
on we thus suppose $\ti\ph=\ti\ph_1$ and $\ti\psi=\ti\psi_1$. Then
\[
\ti\la = \ti\psi\circ(\id+\ti\ph) = 
\sum_{k\ge0} \frac{1}{k!} (\pa^k\ti\psi) \ti\ph^k,
\qquad
\ti\mu = \sum_{k\ge1} h_k \ti\chi^k
\]
where $H(t) = c + \sum_{k\ge1} h_k t^k$ with $|h_k| \le C D^{k}$ for some
$C, D>0$ independent of~$k$, and the corresponding formal Borel transforms are
\[
\hat\la = 
\sum_{k\ge0} \frac{1}{k!} \big((-\ze)^k\hat\psi\big) * \hat\ph^{*k},
\qquad
\hat\mu = \sum_{k\ge1} h_k \hat\chi^{*k}.
\]
These can be viewed as formally convergent series of elements of $\C[[\ze]]$, in
which each term belongs to~$\hat\gR_\Om$ (by virtue of
Theorem~\ref{thmOmstbgROmstb}).
They define holomorphic germs in $\D_{\rho(\Om)}$ because they can also be seen as
normally convergent series of holomorphic functions in any compact disc
contained in $\D_{\rho(\Om)}$ (by virtue of inequalities~\eqref{ineqhatphkD}
and~\eqref{ineqhatchikD}). 

To conclude, it is sufficient to check that, given a path
$\ga\col[0,1]\to\C\setminus\Om$ starting in $\D_{\rho(\Om)}$,
for every $t\in[0,1]$ and $R_t>0$ such that $\ov{D\big(\ga(t),R_t\big)} \subset \C\setminus\Om$
the series of holomorphic functions
\[
\sum \frac{1}{k!} \cont_{\ga_{|[0,t]}} \Big(\big((-\ze)^k\hat\psi\big) * \hat\ph^{*k}\Big)
\quad \text{and} \quad
\sum h_k \cont_{\ga_{|[0,t]}} \big(\hat\chi^{*k}\big)
\]
are normally convergent on $\ov{D\big(\ga(t),R_t\big)}$ (indeed, this will
provide families of analytic elements which analytically continue~$\hat\la$ and~$\hat\mu$).
This follows from Lemma~\ref{lemestimiterconv}.
\end{proof}


\begin{exa}	\label{exaexpoStir}
In view of Lemma~\ref{lemEulPoinStrilRES}, since $2\pi\I\,\Z$ is stable under
addition, this implies that the exponential of the Stirling series $\ti\la =
\ee^{\ti\mu}$ is $2\pi\I\,\Z$-resurgent.
\end{exa}


Recall that $\ti\gG = \id + \C[[z\ii]]$ is the topological group of formal
tangent-to-identity diffeomorphisms at~$\infty$ studied in
Section~\ref{secFormalDiffeos}.

\begin{Def}
We call \emph{$\Om$-resurgent tangent-to-identity diffeomorphism} any $\ti f = \id + \ti\ph
\in \ti\gG$ where $\ti\ph$ is an $\Om$-resurgent formal series.
We use the notations
\[
\gGR(\Om) \defeq \{\, \ti f = \id + \ti\ph \mid \ti\ph\in\ti\gR_\Om \,\},
\qquad
\gGR_\sig(\Om) \defeq \gGR(\Om) \cap \ti\gG_\sig
\ens\text{for $\sig\in\C$.}
\]
\end{Def}

Observe that $\gGR(\Om)$ is not a closed subset of~$\ti\gG$ for the topology
which was introduced in Section~\ref{secFormalDiffeos}; 
in fact it is dense, since it contains the subset~$\gG$ of holomorphic
tangent-to-identity germs of diffeomorphisms at~$\infty$, which itself is dense
in~$\ti\gG$.


\begin{thm}	\label{thmresdiffeo}
The set~$\gGR(\Om)$ is a subgroup of~$\ti\gG$,
the set~$\gGR_0(\Om)$ is a subgroup of~$\ti\gG_0$.
\end{thm}


\begin{proof}
The stability under group composition stems from Theorem~\ref{thmresOmstbNL},
since 
$(\id+\ti\psi)\circ(\id+\ti\ph) = \id + \ti\ph + \ti\psi\circ(\id+\ti\ph)$.

For the stability under group inversion, we only need to prove
\[
\ti h = \id + \ti\chi \in \gGR(\Om) 
\quad\Longrightarrow\quad
\ti h\ic \in \gGR(\Om).
\]
It is sufficient to prove this when $\ti\chi$ has no constant term, 
\ie when $\ti h \in \gGR_0(\Om)$,
since we can always write
$\ti h = (\id+\ti\chi_1)\circ(\id+a)$ with a formal series 
$\ti\chi_1 = T_{-a}(-a+\ti\chi) \in\ti\gR_{\Om}$  
which has no constant term (taking $a =$ constant term of~$\ti\chi$ and using
Lemma~\ref{lemelemstab})
and then $\ti h\ic = (\id+\ti\chi_1)\ic - a$.

We thus assume that $\ti\chi = \ti\chi_1 \in\ti\gR_{\Om}$ has no constant term and apply the
Lagrange reversion formula~\eqref{eqnLagrangeRev} to $\ti h = \id+\ti\chi$.
We get $\ti h\ic = \id-\ti\ph$ with the Borel transform of~$\ti\ph$ given by
\[
\hat\ph =  \sum_{k\ge1} \frac{\ze^{k-1}}{k!} \hat\chi^{*k},
\]
formally convergent series in $\C[[\ze]]$, in which each term belongs
to~$\hat\gR_\Om$.
The holomorphy of~$\hat\ph$ in $\D_{\rho(\Om)}$ and its analytic continuation along the paths
of $\C\setminus\Om$ are obtained by invoking inequalities~\eqref{ineqhatchikDR}
and Lemma~\ref{lemestimiterconv}, similarly to what we did at the end of the
proof of Theorem~\ref{thmresOmstbNL}.
\end{proof}



\vspace{.9cm}

\centerline{\Large\sc Simple singularities}
\addcontentsline{toc}{part}{\sc Simple singularities}

\vspace{.6cm}

\section{Singular points}

When the analytic continuation of a holomorphic germ~$\htb\ph(\ze)$ has
singularities (\ie $\htb\ph$ does not extend to an entire function), its inverse
formal Borel transform $\ti\ph = \cB\ii\htb\ph$ is a divergent formal series,
and the location and the nature of the singularities in the $\ze$-plane
influence the growth of the coefficients of~$\ti\ph$.
By analysing carefully the singularities of~$\htb\ph$, one may hope to be able
to deduce subtler information on~$\ti\ph$ and, if Borel-Laplace
summation is possible, on its Borel sums.

Therefore, we shall now develop a theory which allows one to study and
manipulate singularities (in the case of isolated singular points).

First, recall the definition of a \emph{singular point} in complex analysis: 
given $f$ holomorphic in an open subset~$U$ of~$\C$, a boundary point~$\om$
of~$U$ is said to be a singular point of~$f$ if one cannot find an open
neighbourhood~$V$ of~$\om$, a function~$g$ holomorphic in~$V$, and an open
subset~$U'$ of~$U$ such that $\om\in\pa U'$ and $f_{|U'\cap V} = g_{|U'\cap V}$.

Thus this notion refers to the imposssibility of extending locally the function:
even when restricting to a smaller domain~$U'$ to which~$\om$ is adherent, we
cannot find an analytic continuation in a full neighbourhood of~$\om$.
Think of the example of the principal branch of logarithm: it can be defined as
the holomorphic function
\beglabel{eqdefLogppalbr}
\Log\ze \defeq \int_1^\ze \frac{\dd\xi}{\xi}
\quad \text{for $\ze \in U = \C\setminus\R^-$.}
\elabel
Then, for $\om<0$, one cannot find a holomorphic extension of $f = \Log$ from~$U$
to any larger open set containing~$\om$ (not even a continuous extension!),
however such a point~$\om$ is not singular: if we first restrict, say, to the
upper half-plane $U' \defeq
\{ \IM\ze>0 \}$, then we can easily find an analytic continuation of
$\Log_{|U'}$ to $U' \cup V$, where $V$ is the disc $D(\om,|\om|)$:
define~$g$ by 
\[
g(\ze) = \bigg(\int_\ga + \int_\om^\ze\bigg) \frac{\dd\xi}{\xi}
\] 
with any path $\ga\col [0,1] \to \C$ such that $\ga(0)=1$,
$\ga\big((0,1)\big)\subset U'$ and $\ga(1) = \om$.
In fact, for the function $f = \Log$, the only singular point is~$0$, there is
no other local obstacle to analytic continuation, even though there is no holomorphic
extension of this function to~$\C^*$.

If $\om$ is an isolated%
\footnote{
As a rule, all the singular points that we shall encounter in resurgence theory
will be isolated even when the same holomorphic function~$f$ is considered in
various domains~$U$
(\ie no ``natural boundary'' will show up).
This does not mean that our functions will extend in punctured dics centred on
the singular points, because there may be ``monodromy'': leaving the original domain of
definition~$U'$ on one side of~$\om$ or the other may lead to different analytic
continuations.
}
singular point for a holomorphic function~$f$, we can wonder what
kind of \emph{singularity} occurs at this point.
There are certainly many ways for a point to be singular: maybe the function
near~$\om$ looks like $\log(\ze-\om)$ (for an appropriate branch of the
logarithm), or like a pole $\frac{C}{(\ze-\om)^m}$, and the reader can imagine
many other singular behaviours (square-root branching $(\ze-\om)^{1/2}$, powers
of logarithm $\big(\log(\ze-\om)\big)^m$, iterated logarithms
$\log\big(\log(\ze-\om)\big)$, etc.).
The singularity of~$f$ at~$\om$ will be defined as an equivalence class modulo
regular functions in Section~\ref{secformalismsing}.
Of course, by translating the variable, we can always assume $\om=0$.
Observe that, in this text, we make a distinction between singular points and
singularities (the former being the locations of the latter).

As a preliminary, we need to introduce a few notations in relation with the
Riemann surface of the logarithm.

\section{The Riemann surface of the logarithm}	\label{secClog}

The Riemann surface of the logarithm~$\Clog$ can be defined topologically (without any
reference to the logarithm!) as the universal cover of~$\C^*$ with base point
at~$1$.
This means that we consider the set~$\gP$ of all paths\footnote{%
In this section, ``path'' means any continuous $\C$-valued map defined
on a real interval.
} 
$\ga \col [0,1] \to \C^*$ with $\ga(0)=1$,
we put on~$\gP$ the equivalence relation~$\sim$ of ``homotopy with fixed endpoints'', \ie
\begin{multline*}
\ga \sim \ga_0
\quad\Longleftrightarrow\quad
\exists H \col [0,1]\times[0,1] \to \C^* \;\text{continuous, such that
$H(0,\cdot\,) = \ga_0$, $H(1,\cdot\,) = \ga$,}\\[1ex]
\text{$H(s,0)=\ga_0(0)$ and $H(s,1)=\ga_0(1)$ for each $s\in[0,1]$,}
\end{multline*}
and we define~$\Clog$ as the set of all equivalence classes,
\[
\Clog \defeq \gP / \sim.
\]


Observe that, if $\ga\sim\ga_0$, then $\ga(1) = \ga_0(1)$: the endpoint $\ga(1)$
does not depend on~$\ga$ but only on its equivalence class~$[\ga]$.
We thus get a map
\[
\pi \col \Clog \to \C^*, \qquad
\pi(\ze) = \ga(1) \ens\text{for any $\ga\in\gP$ such that $[\ga]=\ze$}
\]
(recall that the other endpoint is the same for all paths $\ga\in\gP$:
$\ga(0)=1$).


Among all the representatives of an equivalence class $\ze\in\Clog$, there is a
canonical one: there exists a unique pair $(r,\th) \in (0,+\infty) \times\R$ such
that~$\ze$ is represented by the concatenation of the paths
$t\in[0,1] \mapsto \ee^{\I t \th}$
and $t\in[0,1] \mapsto \big(1 + t(r-1)\big) \ee^{\I \th}$.
In that situation, we use the notations
\[
\ze = r \, \eel^{\I\th},
\qquad r = \abs{\ze},
\qquad \th = \arg\ze,
\]
so that we can write $\pi(r \, \eel^{\I\th}) = r \, \ee^{\I\th}$.
Heuristically, one may think of $\th \mapsto \eel^{\I\th}$ as of a non-periodic exponential:
it keeps track of the number of turns around the origin, not only of the angle
$\th$ modulo $2\pi$.


There is a simple way of defining a Riemann surface structure on~$\Clog$. 
One first defines a Hausdorff topology on~$\Clog$ by taking as a basis $\big\{\,
\ti D(\ze,R) \mid \ze\in\Clog, \; 0 < R < \abs{\pi(\ze)} \,\big\}$, where
$\ti D(\ze,R)$ is the set of the equivalence classes of all paths~$\ga$ obtained as concatenation of a
representative of~$\ze$ and a line segment starting from~$\pi(\ze)$ and contained
in $D(\pi(\ze),R)$ (\cf notation~\eqref{eqdefDst}).
(Exercise: check that this is legitimate, \ie that $\big\{\, \ti D(\ze,R) \,
\big\}$ is a collection of subsets of~$\Clog$ which meets the necessary
conditions for being the basis of a topology, and check that the resulting
topology satisfies the Hausdorff separation axiom.)
It is easy to check that, for each basis element, the projection~$\pi$
induces a homeomorphism $\pi_{\ze,R} \col \ti D(\ze,R) \to D(\pi(\ze),R)$ and
that, for each pair of basis elements with non-empty intersection, the
transition map $\pi_{\ze',R'}\circ\pi_{\ze,R}\ii$ is the identity map on
$D(\pi(\ze),R) \cap D(\pi(\ze'),R') \subset\C$,
hence we get an atlas $\{ \pi_{\ze,R} \}$ which defines a Riemann surface
structure on~$\Clog$, \ie a $1$-dimensional complex manifold structure (because
the identity map is holomorphic!).


Now, why do we call~$\Clog$ the Riemann surface of the logarithm?
This is not so apparent in the presentation that was adopted here, but in fact the
above construction is related to a more general one, in which one starts
with an arbitrary open connected subset~$U$ of~$\C$ and a holomorphic
function~$f$ on~$U$, and one constructs (by quotienting a certain set of paths)
a Riemann surface in which $U$ is embedded and on which~$f$ has a holomorphic
extension.
We shall not give the details, but content ourselves with checking the last
property for $U = \C\setminus\R^-$ and $f=\Log$ defined
by~\eqref{eqdefLogppalbr},
defining a holomorphic function $\gL \col \Clog \to \C$ and explaining why it
deserves to be considered as a holomorphic extension of the logarithm. 

We first observe that $\ti U \defeq \pi\ii(U)$ is an open subset of~$\Clog$ with
infinitely many connected components,
\[
\ti U_m \defeq \{\, r\,\eel^{\I\th} \in \ti\C \mid r>0,\; 
2\pi m-\pi < \th < 2\pi m+\pi \,\},
\qquad m\in\Z.
\]
By restriction, the projection~$\pi$ induces a biholomorphism 
\[
\pi_0 \col \ti U_0 \xrightarrow{\sim} U 
\]
(it does so for any $m\in\Z$ but, quite arbitrarily, we choose $m=0$
here).
The \emph{principal sheet of the Riemann surface of the logarithm} is defined to be
the set $\ti U_0 \subset \ti\C$,
which is identified to the cut plane $U \subset \C$ by means of~$\pi_0$.

On the other hand, since the function $\xi\mapsto 1/\xi$ is holomorphic
on~$\C^*$, the Cauchy theorem guarantees that, for any $\ga\in\gP$, 
the integral $\int_\ga\frac{\dd\xi}{\xi}$ depends only on the equivalence
class~$[\ga]$, we thus get a function
\[
\gL \col \Clog \to \C, \qquad
\gL\big([\ga]\big) \defeq \int_\ga \frac{\dd\xi}{\xi}.
\]
This function is holomorphic on the whole of~$\Clog$, because its expression in any
chart domain $\ti D(\ze_0,R)$ is
\[
\gL(\ze) = \gL(\ze_0) + \int_{\pi(\ze_0)}^{\pi(\ze)} \frac{\dd\xi}{\xi},
\]
which is a holomorphic function of~$\pi(\ze)$.

Now, since any $\ze \in \ti U_0$ can be represented by a line segment starting
from~$1$, we have 
\[
\gL_{|\ti U_0} = \Log \circ \pi_0.
\]
In other words, if we identify $U$ and~$\ti U_0$ by means of~$\pi_0$, we can
view~$\gL$ as a holomorphic extension of~$\Log$ to the whole of~$\Clog$.

The function~$\gL$ is usually denoted $\log$.
Notice that $\gL(r\,\eel^{\I\th}) = \ln r + \I\th$ for all $r>0$ and $\th\in\R$,
and that $\gL = \log$ is a biholomorphism $\Clog \to \C$ 
(with our notations: $\gL\ii(x+\I y) = \ee^x \, \eel^{\I y}$).
Notice also that there is a natural multiplication 
$(r_1\,\eel^{\I\th_1}, r_2\,\eel^{\I\th_2}) \mapsto
r_1r_2\,\ee^{\I(\th_1+\th_2)}$
in~$\Clog$, inherited from the addition in~$\C$.

\section{The formalism of singularities}	\label{secformalismsing}

\begin{figure}
\begin{center}

\includegraphics[scale=.9]{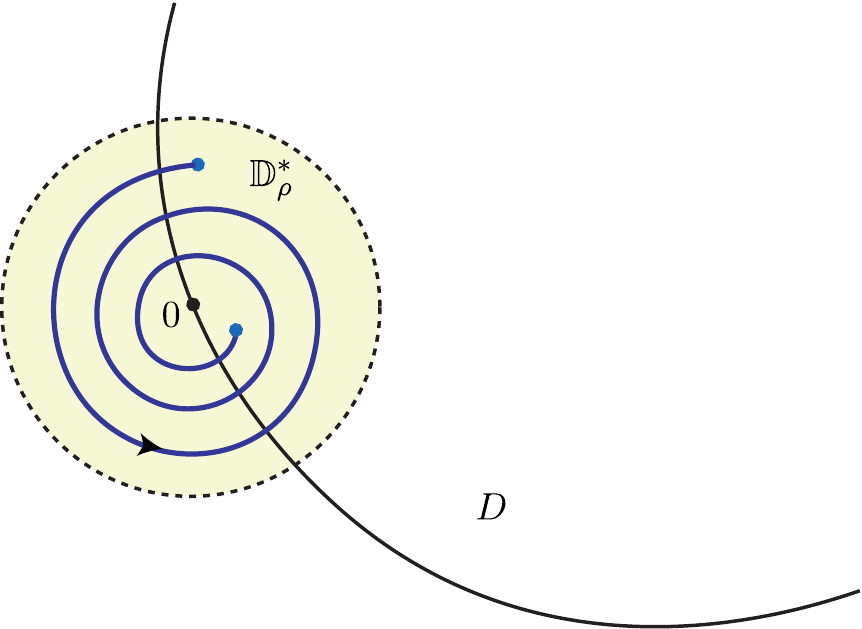} 

\bigskip 

\caption{%
The function~$f$ is holomorphic in~$D$ and has spiral continuation around~$0$.}
\label{fig:spiral_dom}
\end{center}
\end{figure}

We are interested in holomorphic functions~$f$ for which 
the origin is locally the only singular point in the following sense:


\begin{Def}
We say that a function~$f$ \emph{has spiral continuation around~$0$} if
it is holomorphic in an open disc~$D$ to which~$0$ is adherent and, 
for every $L>0$, there exists $\rho>0$ such that $f$ can be analytically
continued along any path of length $\le L$ starting from $D\cap
\D^*_\rho$ and staying in~$\D^*_\rho$
(recall the notation~\eqref{eqdefDDst}).
See Figure~\ref{fig:spiral_dom}.
\end{Def}

In the following we shall need to single out one of the connected components of
$\pi\ii(D)$ in~$\ti\C$, but there is no canonical choice in general.
(If one of the connected components is contained in the principal sheet of~$\ti\C$, we
may be tempted to choose this one, but this does not happen when the centre of~$D$ has
negative real part and we do not want to eliminate a priori this case.)
We thus choose $\ze_0\in\Clog$ such that $\pi(\ze_0)$ is the centre of~$D$, then
the connected component of $\pi\ii(D)$ which contains~$\ze_0$ is a domain~$\ti
D$ of the form $\ti D(\ze_0,R_0)$ (notation of the previous section) and this will be
the connected component that we single out.

Since $\pi$ induces a biholomorphism $\ti D \xrightarrow{\sim} D$,
we can identify~$f$ with $\ch f \defeq f \circ\pi$ viewed as a holomorphic function
on~$\ti D$.
Now, the spiral continuation property implies that~$\ch f$ extends analytically to a
domain of the form
\[
\cV(h) \defeq \{\, \ze = r \, \eel^{\I\th} \mid
0 < r < h(\th), \; \th \in \R \,\}
\subset \Clog,
\]
with a continuous function $h\col\R\to(0,+\infty)$, but in fact the precise
function~$h$ is of no interest to us.%
\footnote{
\label{footchgesheet}
Observe that there is a countable infinity of choices for~$\ze_0$ (all the possible
``lifts'' of the centre of~$D$ in~$\Clog$) thus, a priori, infinitely many
different functions~$\ch f$ associated with the same function~$f$;
they are all of the form $\ch f(\ze\,\eel^{2\pi\I m})$, $m\in\Z$, where $\ch
f(\ze)$ is one of them, so that if~$\ch f$ is holomorphic in a domain of the
form $\cV(h)$ then each of them is holomorphic in a domain of this form.
}
We are thus led to 


\begin{Def}
We define the space~$\ANA$ of all \emph{singular germs} as follows:
on the set of all pairs $(\ch f,h)$, where $h\col\R\to(0,+\infty)$ is continuous
and $\ch f \col \cV(h) \to \C$ is holomorphic, we put the equivalence relation 
\[
(\chfsub1,h_1) \sim (\chfsub2,h_2) 
\quad \xLeftrightarrow{\text{\ens def\ens}} \quad
\chfsub1 \equiv \chfsub2 \;\text{on $\cV(h_1) \cap \cV(h_2)$,}
\]
and we define~$\ANA$ as the quotient set. 
\end{Def}

Heuristically, one may think of a singular germ as of a ``germ of holomorphic
function at the origin of~$\Clog$'' (except that $\Clog$ has no origin!).
We shall usually make no notational difference between an element of~$\ANA$ and
any of its representatives.
As explained above, the formula $f = \ch f \circ \pi$ allows one to identify a
singular germ~$\ch f$ with a function~$f$ which has spiral continuation
around~$0$; however, one must be aware that this presupposes an identification,
by means of~$\pi$, between a simply connected domain~$D$ of~$\C^*$ (\eg an open
disc) and a subset~$\ti D$ of a domain of the form~$\cV(h)$ (and, given~$D$, there are
countably many choices for~$\ti D$).


\begin{exa}	\label{exaLaurSer}
Suppose that $f$ is holomorphic in the punctured disc $\D_\rho^*$, for some
$\rho>0$; in particular, it is holomorphic in $D =
D(\frac{\rho}{2},\frac{\rho}{2})$ and we can apply the above construction.
Then, for whatever choice of a connected component of $\pi\ii(D)$ in~$\Clog$, we
obtain the same $\ch f \defeq f\circ\pi$ holomorphic in $\cV(h)$ with a constant function
$h(\th) \equiv \rho$.
The corresponding element of~$\ANA$ identifies itself with the Laurent series
of~$f$ at~$0$, which is of the form 
\beglabel{eqLaurSer}
\sum_{n\in\Z} a_n \ze^n = S(1/\ze) + R(\ze), 
\elabel
with $R(\ze) \defeq \sum_{n\ge0} a_n\ze^n$ of radius of convergence $\ge \rho$ and
$S(\xi) \defeq \sum_{n>0} a_{-n}\xi^n$ of infinite radius of convergence.
Heuristically, the ``singularity of~$f$'' is encoded by the sole term $S(1/\ze)$;
Definition~\ref{DefSingul} will formalize the idea of discarding the regular term
$R(\ze)$.
\end{exa}


\begin{exa}	\label{exaLogSing}
Suppose that $f$ is of the form $f(\ze) = \htb\ph(\ze) \Log\ze$, where $\htb\ph$
is holomorphic in the disc~$\D_\rho$, for some $\rho>0$, and we are using the
principal branch of the logarithm.
Then we may define $\ch f(\ze) \defeq \htb\ph\big(\pi(\ze)\big) \log\ze$ for $\ze\in\cV(h)$
with a constant function $h(\th) \equiv \rho$;
this corresponds to the situation described above with 
$D = D(\frac{\rho}{2},\frac{\rho}{2})$
and $\ti D =$ the connected component of $\pi\ii(D)$ which is contained in the
principal sheet of~$\Clog$ 
(choosing some other connected component for~$\ti D$ would have resulted in
adding to the above~$\ch f$ an integer multiple of $2\pi\I\,\htb\ph\circ\pi$).
The corresponding element of~$\ANA$ identifies itself with 
\[
\bigg( \sum_{n\ge0} a_n\ze^n\bigg) \log\ze,
\] 
where $\sum_{n\ge0} a_n\ze^n$ is the Taylor series of~$\htb\ph$ at~$0$ (which
has radius of convergence $\ge\rho$).
\end{exa}


\begin{exa}	\label{exazeal}
For $\al\in\C^*$ we define ``the principal branch of~$\ze^\al$'' as
$\ee^{\al\Log\ze}$ for $\ze\in\C\setminus\R^-$.
If we choose~$D$ and~$\ti D$ as in Example~\ref{exaLogSing}, then the
corresponding singular germ is 
\[
\ze^\al \defeq \ee^{\al\log\ze},
\] 
which extends holomorphically to the whole of~$\Clog$.
One can easily check that $0$ is a singular point for~$\ze^\al$ if and only if
$\al\not\in\N$.
\end{exa}


\begin{exo}
Consider a power series $\sum_{n\ge0} a_n\xi^n$ with \emph{finite} radius of
convergence $R>0$ and denote by $\Phi(\xi)$ its sum for $\xi\in\D_R$.
Prove that there exists $\rho>0$ such that 
\[
f(\ze) \defeq \Phi(\ze\,\Log\ze)
\]
is holomorphic in the half-disc $\D_\rho \cap \{\RE\ze>0\}$ and that $0$ is a
singular point.
Prove that $f$ has spiral continuation around~$0$.
Consider any function~$\ch f$ associated with~$f$ as above; prove that one cannot
find a \emph{constant} function~$h$ such that $\ch f$ is holomorphic in $\cV(h)$.
\end{exo}


\begin{exo}	\label{exozeal}
Let $\al\in\C^*$ and
\[
f_\al(\ze) \defeq \frac{1}{\ze^{\al}-\ze^{-\al}}
\]
(notation of Example~\ref{exazeal}).
Prove that $f_\al$ has spiral continuation around~$0$ if and only if $\al \not\in \I\R$.
Suppose that $\al$ is not real nor pure imaginary and consider any 
function~$\chfsub\al$ associated with~$f_\al$ as above; prove that one cannot find
a constant function~$h$ such that $\chfsub\al$ is holomorphic in $\cV(h)$.
\end{exo}


The set~$\ANA$ is clearly a linear space which contains~$\C\{\ze\}$,
in the sense that there is a natural injective linear map $\C\{\ze\}
\hookrightarrow \ANA$
(particular case of Example~\ref{exaLaurSer} with $f$ holomorphic in a disc~$\D_\rho$).
We can thus form the quotient space:


\begin{Def}    \label{DefSingul}
We call \emph{singularities} the elements of the space $\SING \defeq \ANA /
\C\{\ze\}$.
The canonical projection is denoted by $\sing_0$ and we use the notation
\[
\sing_0 \col \left\{ \begin{aligned}
\ANA &\to \SING \\[1ex]
\ch f \hspace{.65em} &\mapsto \: \tr f = \sing_0\bigl(\ch f(\ze)\bigr).
\end{aligned} \right.
\]
Any representative~$\ch f$ of a singularity~$\tr f$ is called a \emph{major}
of~$\tr f$.
\end{Def}


The idea is that singular germs like $\log\ze$ and $\log\ze +
\frac{1}{1-\ze}$ have the same singular behaviour near~$0$: they are different
majors for the same singularity (at the origin).
Similarly, in Example~\ref{exaLaurSer}, the singularity 
$\sing_0\big(\ch f(\ze)\big)$ coincides with $\sing_0\big(S(1/\ze)\big)$.
The simplest case is that of a simple pole or a pole of higher order, for which we
introduce the notation
\beglabel{eqdefdek}
\de \defeq \sing_0\Big( \frac{1}{2\pi\I\ze} \Big), \qquad
\de\pp k \defeq \sing_0\Big( \frac{(-1)^k k!}{2\pi\I\ze^{k+1}} \Big)
\ens\text{for $k\ge0$}.
\elabel
The singularity of Example~\ref{exaLaurSer} can thus be written
$\dst 2\pi\I \sum_{k=0}^\infty \tfrac{(-1)^k}{k!} a_{-k-1} \de\pp k$.

\begin{rem}	\label{remambiglog}
In Example~\ref{exaLogSing}, a singular germ~$\ch f$ was defined from
$f(\ze) = \htb\ph(\ze) \Log\ze$, with $\htb\ph(\ze) \in \C\{\ze\}$, 
by identifying the cut plane $U = \C\setminus\R^-$ with the principal sheet~$\ti
U_0$ of~$\Clog$, and we can now regard~$\ch f$ as a major.
Choosing some other branch of the logarithm or identiying~$U$ with some other
sheet~$\ti U_m$ would yield another major \emph{for the same singularity}, 
because this modifies the major by an integer multiple of $2\pi\I\,\htb\ph(\ze)$
which is regular at~$0$.
The notation
\beglabel{eqdefbem}
\bem\htb\ph \defeq \sing_0\Big( \htb\ph(\ze) \frac{\log\ze}{2\pi\I} \Big)
\elabel
is sometimes used in this situation.
Things are different if we replace~$\htb\ph$ by the Laurent series of a function
which is holomorphic in a punctured disc $\D_\rho^*$ and not regular at~$0$;
for instance, if we denote by $\Llog\ze$ a branch of the logarithm in the
half-plane $V\defeq\{\RE\ze<0\}$, the function $\dfrac{1}{2\pi\I\ze}\Llog\ze$
defines a singular germ, hence a singularity, for any choice of a connected
component~$\ti V$ of $\pi\ii(V)$ in~$\Clog$,
but we change the singularity by an integer multiple of $2\pi\I\,\de$ if we change
the branch of the logarithm or the connected component~$\ti V$.
\end{rem}


\begin{exa}
Let us define
\beglabel{eqdefIsig}
\trn I_\sig \defeq \sing_0\bigl(\chn I_\sig\bigr), \qquad
\chn I_\sig(\ze) \defeq \frac{\ze^{\sig-1}}{(1-\ee^{-2\pi\I\sig})\Ga(\sig)}
\quad\text{for $\sig \in \C\setminus\Z$}
\end{equation}
(notation of Example~\ref{exazeal}).
For $k\in\N$, in view of the poles of the Gamma function (\cf
\eqref{eqmeromcontGa}), we have
$(1-\ee^{-2\pi\I\sig})\Ga(\sig) \xrightarrow[{\sig\to -k}]{}
2\pi\I (-1)^k/k!$,
which suggests to extend the definition by setting 
\[
\chn I_{-k}(\ze) \defeq \frac{(-1)^k k!}{2\pi\I \ze^{k+1}},
\qquad
\trn I_{-k} \defeq \de\pp k
\]
(we could have noticed as well that the reflection
formula~\eqref{eqreflectionGa} yields
$\chn I_\sig(\ze) = 
\tfrac{1}{2\pi\I} \ee^{\pi\I\sig} \Ga(1-\sig) \ze^{\sig-1}$, which yields the
same $\chn I_{-k}$ when $\sig = -k$).
If $n\in \N^*$, there is no limit for~$\chn I_\sig$ as $\sig\to n$, however
$\trn I_\sig$ can also be represented by the equivalent major
$\dst \frac{\ze^{\sig-1}-\ze^{n-1}}{(1-\ee^{-2\pi\I\sig})\Ga(\sig)}$
which tends to the limit
\[
\chn I_n(\ze) \defeq \frac{\ze^{n-1}}{(n-1)!} \frac{\log\ze}{2\pi\I},
\]
therefore we set $\dst \trn I_n \defeq 
\sing_0 \Big( \frac{\ze^{n-1}}{(n-1)!} \frac{\log\ze}{2\pi\I} \Big)$.
We thus get a family of singularities $\big( \trn I_\sig \big)_{\sig\in\C}$.

Observe that
\beglabel{eqobservIsig}
\sing_0(\ze^{\sig-1}) = (1-\ee^{-2\pi\I\sig})\Ga(\sig) \trn I_\sig,
\qquad \sig\in\C,
\elabel
with the convention
$(1-\ee^{-2\pi\I\sig})\Ga(\sig) = 2\pi\I (-1)^k/k!$ if $\sig = -k \in
-\N$
(and this singularity is~$0$ if and only if $\sig = n \in\N^*$).
\end{exa}


We shall not investigate deeply the structure of the space $\SING$, but let us
mention that there is a natural algebra structure on it:
one can define a commutative associative product~$\tr*$ on $\SING$,
for which~$\de$ is a unit,
and which is compatible with the convolution law of $\C\{\ze\}$ defined by
Lemma~\ref{lemconvoldisc} in the sense that
\beglabel{eqconvsingcompat}
\sing_0\Big(\htb\ph(\ze) \frac{\log\ze}{2\pi\I}\Big)
\, \tr* \, 
\sing_0\Big(\hta\psi(\ze) \frac{\log\ze}{2\pi\I}\Big)
= \sing_0\Big( \big(\htb\ph*\hta\psi\big)(\ze) \frac{\log\ze}{2\pi\I}\Big)
\elabel
for any $\htb\ph,\hta\psi\in\C\{\ze\}$.
See \cite{Eca81}, \cite[\S3.1--3.2]{kokyu} for the details.
The differentiation operator $\frac{\dd\,}{\dd\ze}$ passes to the quotient and
the notation~\eqref{eqdefdek} is motivated by the relation
$\de\pp k = \big(\frac{\dd\,}{\dd\ze}\big)^k \de$.
Let us also mention that $\de\pp k$ can be considered as the Borel transform
of~$z^k$ for $k\in\N$, and more generally $\trn I_\sig$ as the Borel transform
of~$z^{-\sig}$ for any $\sig\in\C$:
there is in fact a version of the formal Borel transform operator with values
in~$\SING$, which is defined on a class of formal objects much broader than
formal expansions involving only integer powers of~$z$.


There is a well-defined monodromy%
\footnote{%
The operator $\ch f(\ze) \in \ANA \mapsto \ch f(\ze\,\eel^{-2\pi\I}) \in \ANA$
reflects analytic continuation along a clockwise loop around the origin
for any function~$f$ holomorphic in a disc $D\subset\C^*$ and such that $\ti
f = f \circ \pi$ on one of the connected components of $\pi\ii(D)$.
}
operator $\ch f(\ze) \in \ANA \mapsto \ch f(\ze\,\eel^{-2\pi\I}) \in \ANA$
(recall that multiplication is well-defined in~$\Clog$), 
and the variation map 
$\ch f(\ze) \mapsto \ch f(\ze) - \ch f(\ze\,\eel^{-2\pi\I})$
obviously passes to the quotient:
\begin{Def}
The linear map induced by the variation map 
$\ch f(\ze) \mapsto \ch f(\ze) - \ch f(\ze\,\eel^{-2\pi\I})$
is denoted by
\[
\var \col \left\{ \begin{aligned}
\SING \hspace{.65em} &\to \hspace{.65em} \ANA \\[1ex]
\tr f=\sing_0\bigl(\ch f\bigr) &\mapsto 
\hta f(\ze) = \ch f(\ze) - \ch f(\ze\,\eel^{-2\pi\I}).
\end{aligned} \right.
\]
The germ $\hta f = \var\tr f$ is called the \emph{minor} of the
singularity~$\tr f$. 
\end{Def}


A simple but important example is
\beglabel{eqvarsinglog}
\var \bigg( \sing_0 \Big( \htb\ph(\ze) \frac{\log\ze}{2\pi\I} \Big) \bigg) = \htb\ph(\ze),
\elabel
for any $\htb\ph$ holomorphic in a punctured disc~$\D_\rho^*$.
Another example is provided by the singular germ of~$\ze^\al$ (notation of
Example~\ref{exazeal}): we get
$\var\big( \sing_0(\ze^\al) \big) = (1-\ee^{-2\pi\I\al}) \sing_0(\ze^\al)$,
hence
\[
\var\trn I_\sig = \frac{\ze^{\sig-1}}{\Ga(\sig)}
\quad \text{for all $\sig \in \C\setminus(-\N)$},
\]
but $\var\trn I_{-k} = \var\de\pp k = 0$ for $k\in\N$.
Clearly, the kernel of the linear map~$\var$ consists of the singularities
defined by the convergent Laurent series $\sum_{n\in\Z} a_n\ze^n$ of Example~\ref{exaLaurSer}.


\section{Simple singularities at the origin}


\parag
We retain from the previous section that, starting with a function~$f$ that
admits spiral continuation around~$0$, by identifying a part of the domain
of~$f$ with a subset of~$\Clog$, we get a function~$\ch f$ holomorphic in a
domain of~$\Clog$ of the form~$\cV(h)$ and then a singular germ, still denoted
by~$\ch f$ (by forgetting about the precise function~$h$);
we then capture the singularity of~$f$ at~$0$ by modding out by the regular
germs.

The space~$\SING$ of all singularities is huge.
In this text, we shall almost exclusively deal with singularities of a special kind:


\begin{Def}
We call \emph{simple singularity} any singularity of the form
\[
\trb\ph = a\,\de + \sing_0\Big( \htb\ph(\ze) \frac{\log\ze}{2\pi\I} \Big)
\]
with $a\in\C$ and $\htb\ph(\ze) \in \C\{\ze\}$.
The subspace of all simple singularities is denoted by $\SING\simp$.
We say that a function~$f$ \emph{has a simple singularity at~$0$} if it has spiral
continuation around~$0$ and, for any choice of a domain $\ti D\subset\Clog$
which projects injectively onto a part of the domain of~$f$, the formula 
$\ch f \defeq f \circ \pi_{|\ti D}$ defines the major of a simple singularity.
\end{Def}

In other words, $\SING\simp$ is the range of the $\C$-linear map
\beglabel{defisomSINGsimp}
a\,\de + \htb\ph(\ze) \in \C\,\de \oplus \C\{\ze\} 
\;\mapsto\;
 a\,\de + \sing_0\Big( \htb\ph(\ze) \frac{\log\ze}{2\pi\I} \Big) \in \SING,
\elabel
and a function~$f$ defined in an open disc~$D$ to which $0$ is adherent has a
simple singularity at~$0$ if and only if it can be written in the form
\beglabel{eqformsimplesing}
f(\ze) = \frac{a}{2\pi\I\ze} + \htb\ph(\ze) \frac{\Llog\ze}{2\pi\I} + R(\ze),
\qquad \ze\in D,
\elabel
where $a\in\C$, $\htb\ph(\ze)\in\C\{\ze\}$,
$\Llog\ze$ is any branch of the logarithm in~$D$,
and $R(\ze) \in \C\{\ze\}$.
Notice that we need not worry about the choice of the connected component~$\ti
D$ of $\pi\ii(D)$ in this case: the various singular germs defined from~$f$
differ from one another by an integer multiple of~$\htb\ph$ and thus define the
same singularity (\cf Remark~\ref{remambiglog}).


The map~\eqref{defisomSINGsimp} is injective
(exercise\footnote{Use~\eqref{eqvarsinglog}.}); 
it thus induces a $\C$-linear isomorphism 
\beglabel{eqisomSINGsimp}
\C\,\de\oplus\C\{\ze\} \xrightarrow{\sim} \SING\simp,
\elabel
which is also an algebra isomorphism if one takes into account the algebra
structure on the space of singularities which was alluded to earlier
(in view of~\eqref{eqconvsingcompat}).
This is why we shall identify $\sing_0\big( \htb\ph(\ze) \frac{\log\ze}{2\pi\I}
\big)$ with~$\htb\ph$ and use the notation
\beglabel{eqabusenot}
\sing_0\big(f(\ze)\big) = \trb\ph = a\,\de + \htb\ph(\ze)
\in \C\,\de \oplus \C\{\ze\} \simeq \SING\simp
\elabel
in the situation described by~\eqref{eqformsimplesing},
instead of the notation $a\,\de+\bem\htb\ph$ which is sometimes used in other texts.
(Observe that there is an abuse of notation in the \lhs\ of~\eqref{eqabusenot}:
we should have specified a major~$\ch f$ holomorphic in a subset of~$\Clog$ and
written $\sing_0\big(\ch f(\ze)\big)$, but there is no ambiguity here, as explained
above.)
The germ~$\htb\ph$ is the \emph{minor} of the singularity ($\htb\ph =
\var\trb\ph$) and the complex number~$a$ is called the \emph{constant term}
of~$\trb\ph$.

\parag   \label{paragLaplMaj}
The convolution algebra $\C\,\de \oplus \C\{\ze\}$ was studied in
Section~\ref{sec_convols} as the Borel image of the algebra $\C[[z\ii]]_1$ of
$1$-Gevrey formal series.
Then, in Section~\ref{secvarydir}, we defined its subalgebras 
$\C\,\de \oplus \cN(\eith\R^+)$ and $\C\,\de \oplus \cN(I)$,
Borel images of the subalgebras consisting of formal series $1$-summable in a
direction~$\th$ or in the directions of an open interval~$I$,
and studied the corresponding Laplace operators.

\begin{figure}
\begin{center}

\includegraphics[scale=1]{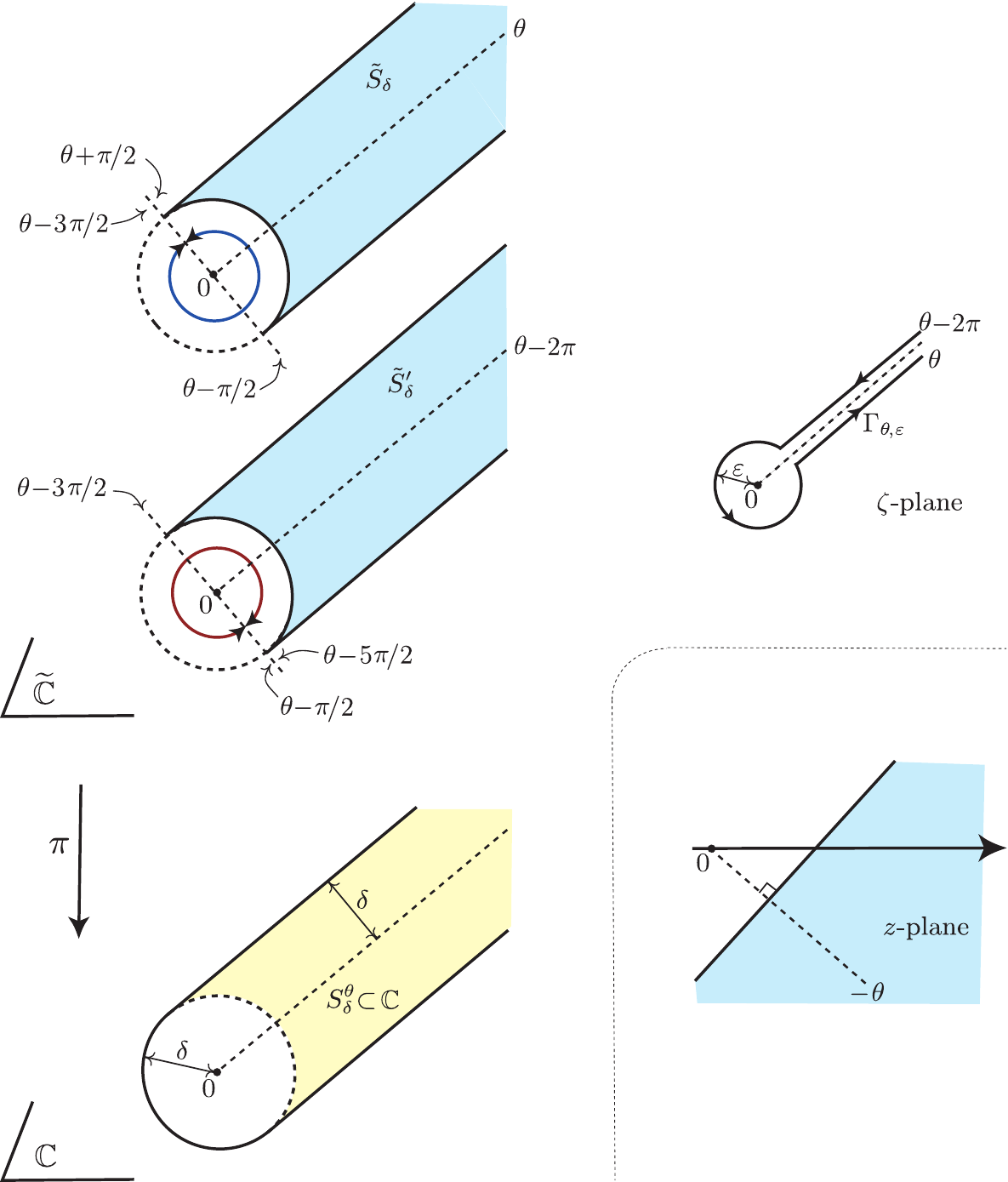} 

\bigskip 

\caption{%
\emph{Laplace transform of a major.}
Left: the domain of~$\protect\Clog$ where $\protect\chb\ph$ must be holomorphic and its
projection~$S_\de^\th$ in~$\C$.
Right: the contour~$\Ga_{\th,\eps}$ for~$\ze$ (above) and the domain where~$z$
belongs (below).}
\label{figLaplTrsfMaj}

\end{center}
\end{figure}

It is interesting to notice that the Laplace transform of a simple singularity
$\trb\ph = a\,\de + \htb\ph(\ze) \in \C\,\de \oplus \cN(\eith\R^+)$ 
can be defined in terms of a major of~$\trb\ph$:
we choose $\chb\ph(\ze) =$ the \rhs\ of~\eqref{eqformsimplesing} with
$R(\ze)=0$, 
or any major~$\chb\ph$ of~$\trb\ph$ for which there exist $\de,\ga>0$ such that this major
extends analytically to
\[
\big\{\, \ze\in\Clog \mid \th-\tfrac{5\pi}{2} < \arg\ze < \th + \tfrac{\pi}{2}
\;\text{and}\; \abs{\ze} < \de \,\big\}
\cup \ti S_\de \cup \ti S'_\de,
\]
where $\ti S_\de$ and~$\ti S'_\de$ are the connected components of 
$\pi\ii(S_\de^\th\setminus\D_\de) \subset \Clog$ which contain~$\eel^{\I\th}$
and~$\eel^{\I(\th-2\pi)}$ (see Figure~\ref{figLaplTrsfMaj}),
and satisfies 
\[
\abs{\chb\ph(\ze)} \le A\, \ee^{\ga\abs{\ze}},
\qquad \ze \in \ti S_\de \cup \ti S'_\de
\]
for some positive constant~$A$;
then, for $0<\eps<\de$ and 
\[
z \in \eel^{-\I\th}\big\{\, z_0\in\Clog \mid \RE z_0 >\ga \;\text{and}\;
\arg z_0 \in (-\tfrac{\pi}{2},\tfrac{\pi}{2}) \,\big\},
\]
we have
\beglabel{eqdefLaplMaj}
(\gS^\th \cB\ii\trb\ph)(z)=
a + (\cL^\th\htb\ph)(z) = \int_{\Ga_{\th,\eps}} 
\ee^{-z\ze} \chb\ph(\ze) \,\dd\ze,
\elabel
with an integration contour~$\Ga_{\th,\eps}$ which comes from infinity along
$\eel^{\I(\th-2\pi)}[\eps,+\infty)$, encircles the origin by following
counterclokwise the circle of radius~$\eps$, and go back to infinity along
$\eel^{\I\th}[\eps,+\infty)$
(a kind of ``Hankel contour''---see Figure~\ref{figLaplTrsfMaj}).
The proof is left as an exercise.\footnote{%
Use~\eqref{eqvarsinglog} for $\ze \in \eel^{\I\th}[\eps,+\infty)$ and then the
dominated convergence theorem for $\eps\to0$.
}


The \rhs\ of~\eqref{eqdefLaplMaj} is the ``Laplace transform of majors''. It
shows why the notation $\trb\ph = c\,\de + \htb\ph$ is consistent with the
notations used in the context of $1$-summability and suggests far-reaching
extensions of $1$-summability theory, which however we shall not pursue in this
text (the interested reader may consult \cite{Eca81} or \cite[\S3.2]{kokyu}).


\vspace{.9cm}

\centerline{\Large\sc Alien calculus for simple resurgent functions}
\addcontentsline{toc}{part}{\sc Alien calculus for simple resurgent functions}

\vspace{.6cm}

\section{Simple $\Om$-resurgent functions and alien operators}

We now leave aside the summability issues and come back to resurgent functions.
Let $\Om$ be a non-empty closed discrete subset of~$\C$ (for the moment we do
not require it to be stable under addition).
From now on, we shall always consider $\Om$-resurgent functions as simple
singularities (taking advantage of~\eqref{eqisomSINGsimp} and~\eqref{eqabusenot}):
\[
\C\,\de \oplus \hat\gR_\Om  \, \subset \,  \C\,\de \oplus \C\{\ze\} \simeq \SING\simp,
\]
where the germs of~$\hat\gR_\Om$ are characterized by $\Om$-continuability.


More generally, at least when $0\in\Om$, we define the space $\SING_\Om$ of
\emph{$\Om$-resurgent singularities} as the space of all $\trb\ph\in\SING$ whose
minors $\htb\ph = \var\trb\ph \in \ANA$ are $\Om$-continuable in the following sense:
denoting by $\cV(h) \subset \Clog$ a domain where~$\htb\ph$ defines a
holomorphic function, $\htb\ph$ admits analytic continuation along any path
$\ti\ga$ of~$\Clog$ starting in~$\cV(h)$ such that $\pi\circ\ti\ga$ is contained
in~$\C\setminus\Om$.
We then have the following diagram:
\[
\xymatrix{
& \rule[-1.25ex]{0ex}{0ex} \makebox[14em]{$\C\,\de\oplus\hat\gR_\Om = \SING\simp\cap\SING_\Om$} \ar@{^{(}->}[r]
& \rule[-1.25ex]{0ex}{0ex} \makebox[11em]{$\C\,\de\oplus\C\{\ze\} = \SING\simp$} 
\ar@{^{(}->}[r] & \makebox[4em]{$\SING$} \\
\makebox[4.5em]{$\C\{z\ii\}$} \ar@{^{(}->}[r]
& \rule[2.5ex]{0ex}{0ex} \makebox[4em]{$\ti\gR_\Om$} \ar@{^{(}->}[r] \ar[u]^-*[@]{\sim}_-*+{\cB} 
& \rule[2.5ex]{0ex}{0ex} \makebox[6em]{$\C[[z\ii]]_1$} \ar[u]^-*[@]{\sim}_-*+{\cB} 
&
}
\]


\begin{figure}
\begin{center}

\includegraphics[scale=1]{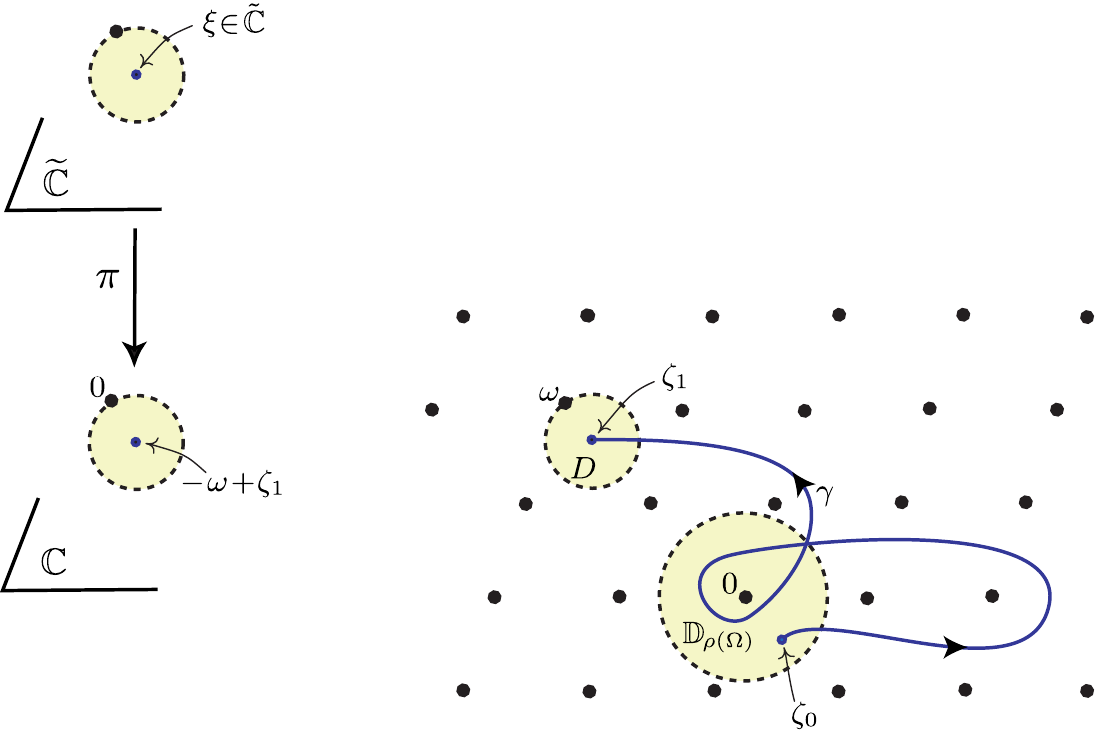} 

\bigskip 

\caption{%
The alien operator $\Al^{\ga,\xi}_\om$ measures the singularity at~$\om$ for the analytic
continuation along~$\ga$ of an $\Om$-continuable germ.}
\label{fig:Alien_op}

\end{center}
\end{figure}


\begin{Def}	\label{defalienomga}
Suppose that $\om\in\Om$, $\ga$ is a path of $\C\setminus\Om$ starting at a
point $\ze_0\in\D^*_{\rho(\Om)}$ and ending at a point~$\ze_1$ such that there
exists an open disc $D \subset \C\setminus\Om$ centred at~$\ze_1$ to which~$\om$
is adherent, and $\xi\in\Clog$ satisfies $\pi(\xi) = -\om+\ze_1$.
We then define a linear map, called the \emph{alien operator
  associated with $(\om,\ga,\xi)$},
\[
\Al^{\ga,\xi}_\om \col \C\,\de\oplus\hat\gR_\Om \to \SING
\]
by the formula
\beglabel{eqdefAlgaxiom}
\Al^{\ga,\xi}_\om (a\,\de + \htb\ph) \defeq \sing_0\big( \ch f(\ze)\big),
\qquad \ch f(\ze) = \cont_\ga\htb\ph\big(\om+\pi(\ze)\big)
\ens\text{for $\ze \in \ti D(\xi)$},
\elabel
where $\ti D(\xi) \subset\Clog$ is the connected component of~$\pi\ii(-\om+D)$ which
contains~$\xi$. See Figure~\ref{fig:Alien_op}.
\end{Def}


This means that we follow the analytic continuation of~$\htb\ph$
along~$\ga$ and get a function $\cont_\ga\htb\ph$ which is holomorphic
in the disc~$D$ centred at~$\ze_1$, of which $\om\in\pa D$ is possibly a
singular point;
we then translate this picture and get a function
\[
\ze \mapsto f(\ze) \defeq \cont_\ga\htb\ph(\om+\ze)
\] 
which is holomorphic in the disc $-\om+D$ centred at $-\om+\ze_1 =
\pi(\xi)$, of which $0 \in \pa(-\om+D)$ is possibly a singular point;
the function~$f$ has spiral continuation around~$0$ because $\htb\ph$
is $\Om$-continuable:
choosing $\eps>0$ small enough so that $D(\om,\eps)\cap\Om = \{\om\}$
(which is possible since $\Om$ is discrete), we see that $\cont_\ga\htb\ph$ can
be continued analytically along any path starting from~$\ze_1$ and staying in
$D(\om,\eps)\cup D$, hence $\ch f$ is holomorphic in $\cV(h)\cup\ti D$ with
$h(\th)\equiv\eps$ and formula~\eqref{eqdefAlgaxiom} makes sense.


\begin{rem}
\label{remunchangeAl}
It is clear that the operator $\Al^{\ga,\xi}_\om$ does not change
if~$\ga$ is replaced with a path which is homotopic (in
$\C\setminus\Om$, with fixed endpoints) to~$\ga$, nor if the endpoints
of~$\ga$ are modified in a continuous way (keeping satisfied the
assumptions of Definition~\ref{defalienomga}) provided that~$\xi$ is
modified accordingly.
On the other hand, modifying~$\xi$ while keeping~$\ga$ unchanged
results in an elementary modification of the result, in line with
footnote~\ref{footchgesheet} on p.~\pageref{footchgesheet}. 
\end{rem}


In a nutshell, the idea is to measure the singularity at~$\om$ for the
analytic continuation along~$\ga$ of the minor~$\htb\ph$.
Of course, if $\om$ is not a singular point for $\cont_\ga \htb\ph$, then
$\Al^{\ga,\xi}_\om\trb\ph = 0$.
In fact, the intersection of the kernels of all the operators
$\Al^{\ga,\xi}_\om$ is $\C\,\de \oplus \gO(\C)$, where $\gO(\C)$ is the set of
all entire functions.
In particular,
\[
\cB\ii\trb\ph \in \C\{z\ii\} 
\quad \Longrightarrow \quad
\Al^{\ga,\xi}_\om\trb\ph = 0.
\]


\begin{exa}
We had $\htb\ph_\al(\ze) \defeq (1+\ze)^{\al-1}$ with $\al\in\C$ in
Example~\ref{exaIncGa}, in connection with the incomplete Gamma function.
Here we can take any $\Om$ containing~$-1$ and we have
$\Al^{\ga,\xi}_\om\htb\ph_\al = 0$ whenever $\om \neq -1$, 
since $-1$ is the only possible singular point of a branch of the
analytic continuation of~$\htb\ph$.
For $\om = -1$, the value of $\Al^{\ga,\xi}_{-1}\htb\ph_\al$ depends on~$\ga$ and~$\xi$.
If $\ga$ is contained in the interval $(-1,0)$, then we find
$f(\ze) = \htb\ph_\al(-1+\ze) =$ the principal branch of $\ze^{\al-1}$ and,
if we choose~$\xi$ in the principal sheet of~$\Clog$, then
\[
\Al^{\ga,\xi}_{-1}\htb\ph_\al = (1-\ee^{-2\pi\I\al})\Ga(\al) \trn I_\al,
\]
which is~$0$ if and only if $\al\in\N^*$ (\cf \eqref{eqobservIsig}).
If~$\ga$ turns $N$ times around~$-1$, keeping the same endpoints for~$\ga$ and
the same~$\xi$, then this result is multiplied by $\ee^{2\pi\I N\al}$;
if we multiply~$\xi$ by~$\eel^{2\pi\I m}$, then the result is multiplied by
$\ee^{-2\pi\I m\al}$.
(In both cases the result is unchanged if $\al\in\Z$.)
\end{exa}


\begin{exa}	\label{exavarLogPole}
Let
$\htb\ph(\ze) = \frac{1}{\ze}\Log(1+\ze)$ (variant of Example~\ref{exaLogPole})
and $\Om = \{-1,0\}$; we shall describe the logarithmic singularity
which arises at~$-1$ and the simple pole at~$0$ for every branch of
the analytic continuation of~$\htb\ph$.
Consider first a path $\ga$ contained in the interval $(-1,0)$ and ending at $\ze_1 =
-\dem$.
For any~$\xi$ projecting onto $0+\ze_1 = -\dem$, we find
$\Al^{\ga,\xi}_{0}\htb\ph = 0$ (no singularity at the origin for the
principal branch),
while  for~$\xi$ projecting onto $1+\ze_1 = \dem$,
\[
\xi = \dem\, \eel^{2\pi\I m} \quad\Longrightarrow\quad
\Al^{\ga,\xi}_{-1}\htb\ph = \sing_0 \Big( \frac{1}{-1+\ze}(-2\pi\I m + \log\ze) \Big)
= -\frac{2\pi\I}{1-\ze}
\]
(using the notation~\eqref{eqabusenot}).
If~$\ga$ turns $N$ times around~$-1$, then the analytic continuation of~$\htb\ph$ is
augmented by $\frac{2\pi\I N}{\ze}$, which is regular at~$-1$ but singular at~$0$,
hence
$\Al^{\ga,\xi}_{-1}\htb\ph$ still coincides with $-\frac{2\pi\I}{1-\ze}$
(the logarithmic singularity at~$-1$ is the same for every branch)
but 
\[
\xi = -\dem\, \eel^{2\pi\I m} \quad\Longrightarrow\quad
\Al^{\ga,\xi}_{0}\htb\ph = \sing_0 \Big( \frac{2\pi\I N}{\ze} \Big) 
= (2\pi\I)^2 N \de.
\]
\end{exa}


\begin{exo}
Consider
$\htb\ph(\ze) = -\frac{1}{\ze}\Log(1-\ze)$ as in Example~\ref{exaLogPole},
with $\Om = \{0,1\}$,
and a path $\ga$ contained in $(0,1)$ and ending at $\ze_1 = \dem$.
Prove that
$\Al^{\ga,\xi}_{1}\htb\ph = -\frac{2\pi\I}{1+\ze}$ for any~$\xi$
projecting onto $-1+\ze_1 = -\dem$.
Compute $\Al^{\ga,\xi}_{0}\htb\ph$ for~$\ga$ turning $N$ times
around~$1$ and~$\xi$ projecting onto $0+\ze_1 = \dem$.
\end{exo}


Examples~\ref{exaLogPole} and~\ref{exavarLogPole} (but not
Example~\ref{exaIncGa} if $\al\notin\N$) are particular cases of

\begin{Def}
  We call \emph{simple $\Om$-resurgent function} any $\Om$-resurgent
  function~$\trb\ph$ such that, for all $(\om,\ga,\xi)$ as in
  Definition~\ref{defalienomga},
$\Al^{\ga,\xi}_\om\trb\ph$ is a simple singularity.
The set of all simple $\Om$-resurgent functions is denoted by
\[ \C\,\de \oplus \hat\gR_\Om\simp, \]
where
$\hat\gR_\Om\simp$ is the set of all simple $\Om$-resurgent functions without
constant term.
We call \emph{simple $\Om$-resurgent series} any element of
\[
\ti\gR_\Om\simp \defeq \cB\ii\big( \C\,\de \oplus \hat\gR_\Om\simp \big)
\subset \ti\gR_\Om.
\]
\end{Def}


\begin{lemma}	\label{lemAlgaxiAlga}
Let $\om,\ga,\xi$ be as in Definition~\ref{defalienomga}. Then
\begin{align*}
\trb\ph \in \C\,\de\oplus\hat\gR_\Om 
\quad &\Longrightarrow \quad
\Al^{\ga,\xi}_\om\trb\ph \in \SING_{-\om+\Om} \\[1ex]
\trb\ph \in \C\,\de \oplus \hat\gR_\Om\simp
\quad &\Longrightarrow \quad
\Al^{\ga,\xi}_\om\trb\ph \in \C\,\de \oplus \hat\gR_{-\om+\Om}\simp.
\end{align*}
Moreover, in the last case, $\Al^{\ga,\xi}_\om\trb\ph$ does not depend on the choice
of~$\xi$ in $\pi\ii(-\om+\ze_1)$; denoting it by
$\Al^\ga_\om\trb\ph$, we thus define an operator
$\Al^\ga_\om \col \C\,\de \oplus \hat\gR_\Om\simp 
\to \C\,\de \oplus \hat\gR_{-\om+\Om}\simp$.
\end{lemma}

\begin{proof}[Proof of Lemma~\ref{lemAlgaxiAlga}]
Let $\trb\ph \in \C\,\de\oplus \hat\gR_\Om$ and $\tr\psi \defeq
\Al^{\ga,\xi}_\om\trb\ph \in \SING$, $\hta\psi \defeq \var\tr\psi \in \ANA$. 
With the notations of Definition~\ref{defalienomga} and~$\eps$ as in the
paragraph which follows it,
we consider the path~$\ga'$ obtained by concatenating~$\ga$ and a loop of
$D(\om,\eps)\cup D$ that starts and ends at~$\ze_1$ and encircles~$\om$
clockwise. We then have $\hta\psi = \ch f - \ch g$, with
\[
\ch f(\ze) \defeq \cont_\ga\htb\ph\big(\om+\pi(\ze)\big)
\quad\text{and}\quad
\chg g(\ze) \defeq \cont_{\ga'}\htb\ph\big(\om+\pi(\ze)\big)
\quad \text{for $\ze\in\ti D$},
\]
where $\ti D$ is the connected component of $\pi\ii(D)$ which contains~$\xi$.

For any path~$\ti\la$ of~$\Clog$ which starts at~$\xi$ and whose projection $\la
\defeq \pi\circ\ti\la$ is contained in $\C\setminus(-\om+\Om)$, the analytic
continuation of~$\ch f$ and~$\chg g$ along~$\ti\la$ exists and is given by
\[
\cont_{\ti\la}\ch f(\ze) = \cont_\Ga\htb\ph\big(\om+\pi(\ze)\big), \qquad
\cont_{\ti\la}\chg g(\ze) = \cont_{\Ga'}\htb\ph\big(\om+\pi(\ze)\big),
\]
where $\Ga$ is obtained by concatenating~$\ga$ and~$\om+\la$, and~$\Ga'$ by
concatenating~$\ga'$ and~$\om+\la$. 
Hence the analytic continuation of~$\hta\psi$ along any such path~$\ti\la$
exists, and this is sufficient to ensure that $\tr\psi\in\SING_\Om$, which was
the first statement to be proved.

If we suppose $\trb\ph \in \C\,\de \oplus \hat\gR_\Om\simp$, then $\trb\psi \in
\SING\simp$, and the second statement follows from Example~\ref{exaLaurSer} and
Remark~\ref{remambiglog}: changing~$\xi$ amounts to adding to~$\ch f$ an integer
multiple of~$\hta\psi$ which is now assumed to be regular at the origin, and
hence does not modify $\sing_0\big(\ch f(\ze)\big)$.
Putting these facts together, we obtain  
$\Al^{\ga,\xi}_\om\trb\ph
\in \SING_{-\om+\Om} \cap \SING\simp 
= \C\,\de \oplus\hat\gR_{-\om+\Om}\simp$
independent of~$\xi$.
\end{proof}


In other words, an $\Om$-resurgent function~$\trb\ph$ is simple if and only if
all the branches of the analytic continuation of the minor $\htb\ph =
\var\trb\ph$ have only simple singularities;
the relation $\Al^\ga_\om \trb\ph = a\,\de + \hta\psi(\ze)$ then means 
\beglabel{eqmeaningAlgaom}
\cont_\ga\htb\ph(\om+\ze) =
\frac{a}{2\pi\I\ze} + \hta\psi(\ze) \frac{\Llog\ze}{2\pi\I} + R(\ze)
\elabel
for $\ze$ close enough to~$0$,
where $\Llog\ze$ is any branch of the logarithm
and $R(\ze) \in \C\{\ze\}$.

\begin{nota}
We just defined an operator $\Al^\ga_\om \col
\C\,\de \oplus \hat\gR_\Om\simp \to \C\,\de \oplus \hat\gR_{-\om+\Om}\simp$.
We shall denote by the same symbol the counterpart of this operator in spaces of
formal series:
\[
\xymatrix{
{\rule[-.5ex]{0ex}{0ex} \makebox[5.5em]{$\C\,\de \oplus \hat\gR_\Om\simp$}} \ar[r]^*+{\Al^\ga_\om} 
& {\rule[-.5ex]{0ex}{0ex} \makebox[6em]{$\C\,\de \oplus \hat\gR_{-\om+\Om}\simp$}} \\
{\makebox[5em]{$\ti\gR_\Om\simp$}} \ar[r]^*+{\Al^\ga_\om} \ar[u]^-*[@]{\sim}_-*+{\cB} 
& {\makebox[5em]{$\ti\gR_{-\om+\Om}\simp$}} \ar[u]^-*[@]{\sim}_-*+{\cB} 
}
\]
\end{nota}


\begin{Def}
Let $\om\in\Om$.
We call \emph{alien operator at~$\om$} any linear combination of composite
operators of the form
\[
\Al^{\ga_r}_{\om-\om_{r-1}} \circ \cdots 
\circ \Al^{\ga_2}_{\om_2-\om_1}
\circ \Al^{\ga_1}_{\om_1}
\]
(viewed as operators
$\C\,\de \oplus \hat\gR_\Om\simp \to 
\C\,\de \oplus \hat\gR_{-\om+\Om}\simp$
or, equivalently,
$\ti\gR_\Om\simp \to \ti\gR_{-\om+\Om}\simp$)
with any $r\ge1$, $\om_1,\ldots,\om_{r-1} \in\Om$,
$\ga_j$ being any path of $\C\setminus(-\om_{j-1}+\Om)$ starting in
$\D^*_{\rho(-\om_{j-1}+\Om)}$ and ending in a disc 
$D_j \subset \D\setminus(-\om_{j-1}+\Om)$ 
to which $\om_j-\om_{j-1}$ is adherent, 
with the conventions $\om_0=0$ and $\om_r = \om$,
so that $\Al^{\ga_j}_{\om_j-\om_{j-1}} \col 
\ti\gR_{-\om_{j-1}+\Om}\simp \to \ti\gR_{-\om_j+\Om}\simp$
is well defined.
\end{Def}


Clearly $\C\,\de \oplus \gO(\C) \subset \C\,\de \oplus \hat\gR_\Om\simp$
(since an entire function has no singularity at all!), hence 
\[
\C\{z\ii\} \subset \ti\gR_\Om\simp,
\]
and of course all alien operators act trivially on such resurgent
functions.
Another easy example of simple $\Om$-resurgent function is provided by any
meromorphic function~$\htb\ph$ of~$\ze$ which is regular at~$0$ and whose poles are all
simple and located in~$\Om$.
In this case $\Al^{\ga}_\om\htb\ph$ does not depend on~$\ga$:
its value is $2\pi\I c_\om \de$, where $c_\om$ is the residuum of~$\htb\ph$ at~$\om$.

\begin{exa}          \label{exasimpleEPS}
By looking at the proof of Lemma~\ref{lemEulPoinStrilRES}, we see that we have
meromorphic Borel transforms for the formal series associated with the names of
Euler, Poincar\'e and Stirling, hence
\[
\ti\ph\eul \in \ti\gR_{\{-1\}}\simp,  \qquad
\ti\ph\poin \in \ti\gR_{s + 2\pi\I\,\Z}\simp, \qquad 
\ti\mu \in \ti\gR_{2\pi\I\,\Z^*}\simp,
\]
and we can compute
\[
\Al^{\ga}_{-1} \htb\ph\Eul = 2\pi\I\de, \qquad
\Al^{\ga}_{s+2\pi\I k} \htb\ph\Poin = 2\pi\I\de, \qquad
\Al^{\ga}_{2\pi\I m} \htb\mu = \frac{1}{m}\de,
\]
for $k\in\Z$, $m\in\Z^*$ with any~$\ga$
(and correspondingly $\Al^{\ga}_{-1} \ti\ph\eul =
\Al^{\ga}_{s+2\pi\I k} \ti\ph\poin = 2\pi\I$,
$\Al^{\ga}_{2\pi\I m} \ti\mu = \frac{1}{m}$).
A less elementary example is $\ti\la = \ee^{\ti\mu}$;
we saw that $\ti\la \in \ti\gR_{2\pi\I\,\Z}$ in
Example~\ref{exaexpoStir}, 
we shall see in Section~\ref{secAlienOpNL} that it belongs to
$\ti\gR_{2\pi\I\,\Z}\simp$ and that any alien operator maps~$\ti\la$
to a multiple of~$\ti\la$.
\end{exa}


Here is a variant of Lemma~\ref{lemelemstab} adapted to the case of simple resurgent functions:
\begin{lemma}	\label{lemSimplestab}
Let $\Om$ be any non-empty closed discrete subset of~$\C$.
\begin{enumerate}[--]
\item
If $\htn B$ is an entire function, then multiplication by~$\htn B$ leaves
$\hat\gR_\Om\simp$ invariant, with
\beglabel{eqalienopmult}
\Al^\ga_\om \htb\ph = a\,\de + \hta\psi(\ze)
\ens\Longrightarrow\ens
\Al^\ga_\om (\htn B \htb\ph) = \htn B(\om)a\, \de + \htn B(\om+\ze)\hta\psi(\ze).
\elabel
\item
  As a consequence, for any $c\in \C$, the operators~$\hat\pa$
  and~$\hat T_c$ (defined by~\eqref{eqdefhatpa}
  and~\eqref{eqdefcounterTc}) leave $\C\,\de \oplus \hat\gR_\Om\simp$
  invariant or, equivalently,
$\ti\gR_\Om\simp$ is stable by $\pa = \frac{\dd\,}{\dd z}$ and~$T_c$; one has
\beglabel{eqalienoppartcases}
\ti\ph_0 \in \ti\gR_\Om\simp
\ens\Longrightarrow\ens
\cA_\om^\ga(\pa\ti\ph_0) = (-\om+\pa)\cA_\om^\ga\ti\ph_0
\ens\text{and}\ens 
\cA_\om^\ga(T_c\ti\ph_0) = \ee^{-c\om} T_c (\cA_\om^\ga\ti\ph_0).
\elabel
\item
If $\ti\psi\in z^{-2}\C\{z\ii\}$, then the solution in
$z\ii\C[[z\ii]]$ of the difference equation
\[
\ti\ph(z+1) - \ti\ph(z) = \ti\psi(z)
\]
belongs to $\ti\gR_{2\pi\I\,\Z^*}\simp$, with
$\Al^\ga_\om \htb\ph = -2\pi\I \hta\psi(\om) \,\de$
for all $(\om,\ga)$ with $\om\in 2\pi\I\,\Z^*$.
\end{enumerate}
\end{lemma}

\begin{proof}
Suppose that $\Al^\ga_\om \htb\ph = a\,\de + \hta\psi(\ze)$.
Since multiplication by~$\htn B$ commutes with analytic continuation, the
relation~\eqref{eqmeaningAlgaom} implies
\[
\cont_\ga\big(\htn B\htb\ph\big)(\om+\ze) = 
\htn B(\om+\ze) \cont_\ga\htb\ph(\om+\ze) =
\frac{\htn B(\om)a}{2\pi\I\ze} + 
\htn B(\om+\ze)\hta\psi(\ze) \frac{\Llog\ze}{2\pi\I} + R^*(\ze)
\]
with $R^*(\ze) = R(\ze) + a \dfrac{\htn B(\om+\ze)-\htn B(\om)}{2\pi\I\ze} \in \C\{\ze\}$,
hence $\Al^\ga_\om (\htn B \htb\ph) = \htn B(\om)a\, \de + \htn B(\om+\ze)\hta\psi(\ze)$.

Suppose now that $\ti\ph_0 \in \ti\gR_\Om\simp$ has Borel transform 
$\trb\ph_0 = \al\,\de + \htb\ph$ with $\al\in\C$ and $\htb\ph$ as above.
According to~\eqref{eqdefhatpa} and~\eqref{eqdefcounterTc}, we have
$\hat\pa\trb\ph_0 = -\ze\htb\ph(\ze)$ and
$\hat T_c \trb\ph_0 = \al\,\de + \ee^{-c\ze}\htb\ph(\ze)$;
we see that both of them belong to $\C\,\de \oplus \hat\gR_\Om\simp$ by applying the first
statement with $\htn B(\ze) = -\ze$ or $\ee^{-c\ze}$, and
\begin{align*}
\cA_\om^\ga(\hat\pa\trb\ph_0) &= 
-\om a \de + (-\om-\ze) \hta\psi(\ze)
= (-\om+\hat\pa)\cA_\om^\ga\trb\ph_0 \\[1ex]
\cA_\om^\ga(\hat T_c\trb\ph_0) &= 
\ee^{-c\om} a \de + \ee^{-c(\om+\ze)} \hta\psi(\ze)
= \ee^{-c\om} \hat T_c (\cA_\om^\ga\trb\ph_0),
\end{align*}
which is equivalent to~\eqref{eqalienoppartcases}.
 
For the last statement, we use Corollary~\ref{coreqdifflin}, according to which $\htb\ph =
\htn B \hta\psi$ with $\htn B(\ze) = \frac{1}{\ee^{-\ze}-1}$ and $\hta\psi(\ze) \in
\ze\gO(\C)$:
the function~$\htb\ph$ is meromorphic on~$\C$ and all its poles are simple and located in
$\Om = 2\pi\I\,\Z^*$, therefore it is a simple $\Om$-resurgent function and we get the
values of $\cA_\om^\ga\trb\ph$ by computing the residues of~$\htb\ph$
(\cf the paragraph just before Example~\ref{exasimpleEPS}).
\end{proof}


\begin{exo}	\label{exoHurwitz}
Given $s\in\C$ with $\RE s >1$, the Hurwitz zeta function\footnote{
Notice that $\ze(s,1)$ is the Riemann zeta value~$\ze(s)$.} 
is defined as
\[
\ze(s,z) = \sum_{k=0}^\infty \frac{1}{(z+k)^s}, 
\qquad z\in\C\setminus\R^-
\]
(using the principal branch of $(z+k)^s$ for each~$k$).
Show that, for $s\in\N$ with $s\ge2$, 
\[
\ti\ph\hurw_s(z) \defeq
\frac{1}{(s-1) z^{s-1}} + \frac{1}{2 z^s} +
\sum_{k=1}^\infty \binom{s+2k-1}{s-1} \frac{B_{2k}}{(s+2k-1) z^{s+2k-1}} 
\]
(where the Bernoulli numbers~$B_{2k}$ are defined in Exercise~\ref{exoBernoulli})
is a simple $2\pi\I\,\Z^*$-resurgent formal series which is $1$-summable in the
directions of $I = (-\frac{\pi}{2},\frac{\pi}{2})$, with
\[
\ze(s,z) = (\gS^I\ti\ph\hurw_s)(z) \sim_1 \ti\ph\hurw_s(z).
\]
Hint: Use Lemma~\ref{lemSimplestab} and prove that $\ze(s,z)$ coincides with the
Laplace transform of
\beglabel{eqBorelHurwitz}
\htb\ph\Hurw_s(\ze) = \frac{\ze^{s-1}}{\Ga(s)(1-\ee^{-\ze})}.
\elabel
\end{exo}

\begin{rem}
If $s\in\C\setminus\N$ has $\RE s > 1$, then \eqref{eqBorelHurwitz} is not
regular at the origin but still provides an example of $2\pi\I\,\Z$-continuable minor (in the
sense of the definition given in the paragraph just before Definition~\ref{defalienomga}).
In fact, there is an extension of $1$-summability theory in which the Laplace
transform of~$\htb\ph\Hurw_s$ in the directions of $(-\frac{\pi}{2},\frac{\pi}{2})$ is
still defined and coincides with~$\ze(s,z)$ (see \cite{Eca81}, \cite[\S3.2]{kokyu}).
\end{rem}


We end this section with a look at the action of alien operators on convolution
products in the ``easy case'' considered in Section~\ref{sec_contconvoleasy}.

\begin{thm}   \label{thmalieneasycase}
Suppose that $\trn B_0 \in \C\,\de\oplus\hat\gR_\Om\simp$ with $\htn B \defeq
\var\trn B_0$ entire. 
Then, for any $\om \in \Om$, all the alien operators 
$\C\,\de\oplus\hat\gR_\Om\simp \to \C\,\de\oplus\hat\gR_{-\om+\Om}\simp$ 
commute with the operator of convolution with~$\trn B_0$.
\end{thm}

\begin{proof}
It suffices to show that, for any $\ga \subset \C\setminus\Om$ 
starting at a point $\ze_0 \in \D_{\rho(\Om)}^*$ 
and ending at the centre~$\ze_1$ of a disc $D\subset\C\setminus\Om$ to
which~$\om$ is adherent,
and for any 
$\trb\ph_0 \in \C\,\de\oplus\hat\gR_\Om\simp$, 
\[
\Al^\ga_\om(\trn B_0*\trb\ph_0) = \trn B_0*(\Al^\ga_\om \trb\ph_0).
\]
We can write 
\[
\trn B_0 = b\de + \htn B, \quad
\trb\ph_0 = c\de + \htb\ph, \quad
\Al^\ga_\om\trb\ph_0 = a\de + \hta\psi.
\]
Then $\trn B_0*\trb\ph_0 = b \trb\ph_0 + c \htn B + \htn B*\htb\ph$
and
$\Al^\ga_\om(\trn B_0*\trb\ph_0) = b \Al^\ga_\om \trb\ph_0 +
\Al^\ga_\om(\htn B*\htb\ph)$, hence we just need to prove that
$\Al^\ga_\om(\htn B*\htb\ph) = \htn B*\Al^\ga_\om\htb\ph$, \ie that
\[
\Al^\ga_\om(\htn B*\htb\ph) = a \htn B + \htn B*\hta\psi.
\]

According to Lemma~\ref{lemeasyconvol}, we have
\[
\cont_\ga(\htn B*\htb\ph)(\om+\ze) = 
\int_0^{\ze_0} \htn B(\om+\ze-\xi) \htb\ph(\xi) \,\dd\xi 
+ \int_\ga \htn B(\om+\ze-\xi) \htb\ph(\xi) \,\dd\xi 
+ \int_{\ze_1}^{\om+\ze} \htn B(\om+\ze-\xi) \htb\ph(\xi) \,\dd\xi
\]
for $\ze \in -\om+D$, 
where it is understood that $\htb\ph(\xi)$ represents the value
at~$\xi$ of the appropriate branch of the analytic continuation
of~$\htb\ph$ (which is $\cont_\ga\htb\ph$ for the third integral).
The standard theorem about an integral depending holomorphically on a
parameter ensures that the sum $R_1(\ze)$ of the first two integrals
extends to an entire function of~$\ze$.
Let $\De \defeq -\om+D$ (a disc to which $0$ is adherent).
Performing the change of variable $\xi \to \om+\xi$ in the third
integral, we get
\[
\cont_\ga(\htn B*\htb\ph)(\om+\ze) = R_1(\ze) +
\int_{{-\om+\ze_1}}^{\ze} \htn B(\ze-\xi) \cont_\ga\htb\ph(\om+\xi) \,\dd\xi,
\qquad \ze \in \De.
\]

Now, according to~\eqref{eqmeaningAlgaom}, we can write
\[
\cont_\ga\htb\ph(\om+\xi) = S(\xi) + R_2(\xi), 
\qquad \xi \in \De \cap \D_\rho^*,
\]
where 
$\dst
S(\xi) = \frac{a}{2\pi\I\xi} + \hta\psi(\xi) \frac{\Llog\xi}{2\pi\I}
$,
$\Llog\xi$ being a branch of the logarithm holomorphic in~$\De$,
$R_2(\xi) \in \C\{\xi\}$,
and $\rho>0$ is smaller than the radii of convergence of~$\hta\psi$ and~$R_2$.
Let us pick $\sig \in \De \cap \D_\rho^*$ and set
\[
R(\ze) \defeq R_1(\ze) + 
\int_{{-\om+\ze_1}}^{\sig} \htn B(\ze-\xi) \cont_\ga\htb\ph(\om+\xi) \,\dd\xi,
\]
so that
\[
\cont_\ga(\htn B*\htb\ph)(\om+\ze) = R(\ze) +
\int_{\sig}^{\ze} \htn B(\ze-\xi) \cont_\ga\htb\ph(\om+\xi) \,\dd\xi,
\qquad \ze \in \De.
\]
We see that $R(\ze)$ extends to an entire function of~$\ze$ and, for
$\ze\in\De\cap\D_\rho^*$, the last integral can be written
\[
\int_{\sig}^{\ze} \htn B(\ze-\xi) \cont_\ga\htb\ph(\om+\xi) \,\dd\xi
= f(\ze) + R_3(\ze),
\qquad f(\ze) \defeq
\int_{\sig}^{\ze} \htn B(\ze-\xi) S(\xi) \,\dd\xi,
\]
with $R_3(\ze)$ defined by an integral involving $R_2(\xi)$ and thus
extending holomorphically for $\ze \in \D_\rho$.
The only possibly singular term in $\cont_\ga(\htn
B*\htb\ph)(\om+\ze)$ is thus $f(\ze)$, which is seen to 
admit analytic continuation along every path~$\Ga$
starting from~$\sig$ and contained in~$\D_\rho^*$; indeed,
\beglabel{eqconfGaf}
\cont_\Ga f(\ze) = \int_\Ga \htn B(\ze-\xi) S(\xi) \,\dd\xi.
\elabel
In particular, $f$ has spiral continuation around~$0$. We now show
that it defines a simple singularity, which is none other than 
$a \htn B + \htn B*\hta\psi$.

Let us first compute the difference $g \defeq f^+-f$, where we denote by~$f^+$
the branch of the analytic continuation of~$f$ obtained by starting
from $\De\cap\D_\rho^*$, turning anticlockwise around~$0$ and coming
back to $\De\cap\D_\rho^*$.
We have 
\[
g(\ze) = \int_{C_\ze} \htn B(\ze-\xi) S(\xi) \, \dd\xi,
\qquad \ze \in \De\cap\D_\rho^*,
\]
where $C_\ze$ is the circular path $t\in[0,2\pi] \mapsto \ze\,\ee^{\I t}$.
For any $\eps \in (0,1)$, by the Cauchy theorem,
\[
g(\ze) = 
\int_{\eps\ze}^\ze \htn B(\ze-\xi) \hta\psi(\xi) \, \dd\xi +
\int_{C_{\eps\ze}} \htn B(\ze-\xi) S(\xi) \, \dd\xi,
\]
because $S^+-S = \hta\psi$.
Keeping~$\ze$ fixed, we let $\eps$ tend to~$0$:
the first integral clearly tends to $\htn B*\hta\psi(\ze)$ and the
second one can be written
\begin{multline*}
a \int_{C_{\eps\ze}} \frac{\htn B(\ze-\xi)}{2\pi\I\xi} \, \dd\xi
+ \int_{C_{\eps\ze}} \htn B(\ze-\xi)\hta\psi(\xi)\frac{\Llog\xi}{2\pi\I} \, \dd\xi
= \\[1ex] a \htn B(\ze) +
\int_0^{2\pi} \htn B(\ze-\eps\ze\,\ee^{\I t}) \hta\psi(\eps\ze\,\ee^{\I t})
\frac{\ln\eps + \Llog\ze + \I t}{2\pi\I} \I \eps \ze \, \ee^{\I t} \, \dd t
\end{multline*}
(because the analytic continuation of~$\Llog$ is explicitly known),
which tends to $a \htn B(\ze)$ since the last integral is bounded in modulus by
$C \eps \big(C' + \abs{\ln\eps} \big)$ with appropriate constants $C,C'$.
We thus obtain
\[
g(\ze) = a \htn B(\ze) + \htn B*\hta\psi(\ze). 
\]
Since this function is regular at the origin and holomorphic
in~$\D_\rho$, we can reformulate this result on $f^+-f$ by saying that
the function
\[
\ze \in \De\cap\D_\rho^* \mapsto h(\ze) \defeq f(\ze) - g(\ze) \frac{\Llog\ze}{2\pi\I}
\]
extends analytically to a (single-valued) function holomorphic in
$\D_\rho^*$, \ie it can be represented by a Laurent series~\eqref{eqLaurSer}.

We conclude by showing that the above function~$h$ is in fact regular at
the origin. For that, it is sufficient to check that, in $\D_{\abs{\sig}}^*$,
it is bounded by $C \big(C' + \ln\frac{1}{\abs\ze} \big)$ with appropriate constants
$C,C'$
(indeed, this will imply $\ze h(\ze) \xrightarrow[\ze\to0]{} 0$, 
thus the origin will be a removable singularity for~$h$).
Observe that every point of~$\D_{\abs{\sig}}^*$ can be written in the form
$\ze = \sig u \, \ee^{\I v}$ with 
$0 < u \defeq { \abs\ze } / { \abs{\sig} } < 1$ 
and $0 \le v < 2\pi$, hence it can be reached by starting from~$\sig$
and following the concatenation~$\Ga_\ze$ of the circular path
$ 
t \in [0,v] \mapsto \sig \ee^{\I t}$
%
%
and the line segment 
$ 
t \in [0,1] \mapsto \sig \ee^{\I v} x(t)$
with $x(t) \defeq 1 - t (1-u) >0$,
hence
\begin{multline*}
(\cont_{\Ga_\ze} h)(\ze) =
(\cont_{\Ga_\ze} f)(\ze) - 
\frac{1}{2\pi\I} g(\ze)(\cont_{\Ga_\ze} \Llog)(\ze) \\[1ex]
= \int_{\Ga_\ze} \htn B(\ze-\xi) S(\xi) \,\dd\xi -
\frac{1}{2\pi\I} g(\ze)(\Llog\sig + \ln u + \I v)
\end{multline*}
(using~\eqref{eqconfGaf} for the first term).
The result follows from the existence of a constant $M>0$ such
that
$\abs{\htn B} \le M$ on $\D_{2\rho}$,
$\abs g \le M$ on $\D_{\rho}$
and $\abs{S(\xi)} \le M / \abs\xi$ for $\xi \in \D_{\rho}$,
because the first term in the above representation of~$h(\ze)$ has
modulus $\le$
\[
\abs*{ \int_0^v \htn B(\ze-\sig\ee^{\I t}) S(\sig\ee^{\I t})
  \sig\I \ee^{\I t} \, \dd t
+ \int_0^1 \htn B\big( \ze - \sig \ee^{\I v} x(t) \big) S\big( \sig \ee^{\I v} x(t) \big)
\sig \ee^{\I v} x'(t)\,\dd t }
\le M^2 v + M^2 \ln \frac{1}{u}.
\]
\end{proof}

\section{The alien operators $\De^+_\om$ and $\De_\om$}
\label{secalopDepomDeom}

We still denote by~$\Om$ a non-empty closed discrete subset of~$\C$.
We now define various families of alien operators acting on simple
$\Om$-resurgent functions, among which the
most important will be $(\De^+_\om)_{\om\in\Om\setminus\{0\}}$
and $(\De_\om)_{\om\in\Om\setminus\{0\}}$.

\subsection{Definition of $\Al_{\om,\eps}^\Om$, $\De^+_\om$ and $\De_\om$}

\begin{figure}
\begin{center}
\includegraphics[scale=1]{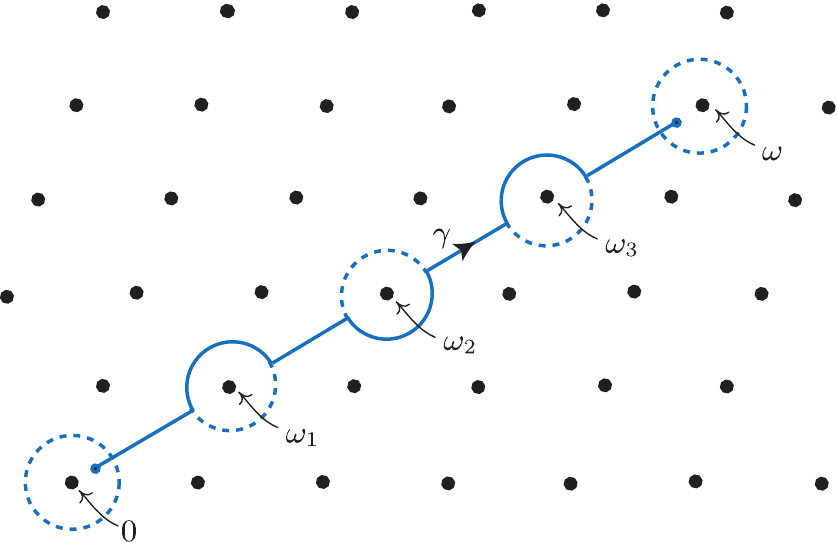}
\bigskip 

\caption{%
  An example of path~$\ga$ used in the definition
  of~$\Al_{\om,\eps}^\Om$, here with $\eps = (-,+,-)$.}
\label{figpathgaeps}

\end{center}
\end{figure}

\begin{Def}
Let $\om \in \Om\setminus\{0\}$.
We denote by $\prec$ the total order on $[0,\om]$ induced by
$t\in[0,1] \mapsto t\,\om \in [0,\om]$
and write
\beglabel{eqdefrcapOm}
[0,\om] \cap \Om = \{\om_0, \om_1, \ldots, \om_{r-1}, \om_r\},
\qquad 0 = \om_0 \prec \om_1 \prec \cdots \prec \om_{r-1} \prec \om_r = \om
\elabel
(with $r \in \N^*$ depending on~$\om$ and~$\Om$).
With any $\eps = (\eps_1,\ldots,\eps_{r-1}) \in \{+,-\}^{r-1}$ we
associate an alien operator at~$\om$
\beglabel{eqdefAlomeps}
\Al_{\om,\eps}^\Om \col \ti\gR_\Om\simp \to \ti\gR_{-\om+\Om}\simp
\elabel
defined as $\Al_{\om,\eps}^\Om = \Al_\om^{\ga}$
for any path~$\ga$ chosen as follows:
we pick $\de>0$ small enough so that
%
%
the closed discs
$D_j \defeq \ov D(\om_j,\de)$, $j=0,1,\ldots r$,
are pairwise disjoint and satisfy
$D_j \cap \Om = \{\om_j\}$,
and we take a path~$\ga$ connecting ${]0,\om[} \cap D_0$
and ${]0,\om[} \cap D_r$ by following the line segment ${]0,\om[}$ except
that, for $1\le j\le r-1$,
the subsegment ${]0,\om[} \cap D_j$ is replaced by one of the two
half-circles which are the connected components of 
${]0,\om[} \cap \pa D_j$:
the path~$\ga$ must circumvent~$\om_j$ 
to the right if $\eps_j = +$, to the left if $\eps_j = -$.
See Figure~\ref{figpathgaeps}.
\end{Def}


Observe that the notation~\eqref{eqdefAlomeps} is justified by the
fact that, in view of Remark~\ref{remunchangeAl}, the operator
$\Al_{\om,\eps}^\Om$ does not depend on~$\de$ nor on the endpoints of~$\ga$.


\begin{Def}
\label{DefDeompDeom}
For any $\om \in \Om\setminus\{0\}$, we define two particular alien
operators at~$\om$
\[
\De^+_\om,\; \De_\om \col 
\ti\gR_\Om\simp \to \ti\gR_{-\om+\Om}\simp
\]
by the formulas
\beglabel{eqdefDepomDeom}
\De^+_\om \defeq \Al_{\om,(+,\ldots,+)}^\Om, \qquad
\De_\om \defeq \sum_{\eps \in \{+,-\}^{r-1}} 
\frac{p(\eps)!q(\eps)!}{r!} \Al_{\om,\eps}^\Om,
\elabel
where $r = r(\om,\Om)$ is defined by~\eqref{eqdefrcapOm} and
$p(\eps)$ and $q(\eps)$ represent the number of symbols~`$+$'
and~`$-$' in the tuple~$\eps$
(so that $p(\eps)+q(\eps) = r-1$).
\end{Def}

We thus have (still with notation~\eqref{eqdefrcapOm}), for $r=1,2,3$:
\begin{align*}
\De_{\om_1}^+ &= \Al_{\om_1,()}^\Om, 
&
\De_{\om_1} &= \Al_{\om_1,()}^\Om, \\
\De_{\om_2}^+ &= \Al_{\om_2,(+)}^\Om, 
&
\De_{\om_2} &= \demi\Al_{\om_2,(+)}^\Om + \demi\Al_{\om_2,(-)}^\Om, \\
\De_{\om_3}^+ &= \Al_{\om_3,(+,+)}^\Om, 
&
\De_{\om_3} &= \tiers\Al_{\om_3,(+,+)}^\Om +
\frac{1}{6}\Al_{\om_3,(+,-)}^\Om + \frac{1}{6}\Al_{\om_3,(-,+)}^\Om
+ \tiers\Al_{\om_3,(-,-)}^\Om.
\end{align*}
Of course, the operators $\De^+_\om, \De_\om, \Al_{\om,\eps}^\Om$ can all
be considered as operators
$\C\,\de\oplus\hat\gR_\Om\simp \to \C\,\de\oplus\hat\gR_{-\om+\Om}\simp$
as well.


\begin{rem}   \label{remAlienDeriv}
Later on, in Sections \ref{sec:ResurIter}--\ref{sec:Bridge}, we shall assume that~$\Om$ is an
additive subgroup of~$\C$, so $-\om+\Om = \Om$ and 
$\De^+_\om, \De_\om, \Al_{\om,\eps}^\Om$ are operators from~$\ti\gR_\Om\simp$
to itself;
we shall see in Section~\ref{secsubalg} that, in that case, $\ti\gR_\Om\simp$ is a subalgebra
of~$\ti\gR_\Om$ (which is itself a subalgebra of $\C[[z\ii]]$ by
Corollary~\ref{corOmresursubalg}) 
of which each~$\De_\om$ is a derivation (\ie it satisfies the Leibniz rule).
For that reason the operators~$\De_\om$ are called ``alien derivations''.
\end{rem}


Observe that, given $r\ge1$, 
there are $r$~possibilities for the value of $p=p(\eps)$ and,
for each~$p$, there are $\binom{r-1}{p}$ tuples~$\eps$ such that
$p(\eps)=p$;
since in the definition of $\De_\om$ the coefficient in front of
$\Al_{\om,\eps}^\Om$ is the inverse of $\dst{r}\binom{r-1}{p(\eps)}$,
it follows that the sum of all these $2^{r-1}$ coefficients is~$1$.
The resurgent function $\De_\om(c\,\de+\htb\ph)$ can thus be viewed as
an average of the singularities at~$\om$ of the branches of the
minor~$\htb\ph$ obtained by following the $2^{r-1}$ ``most direct''
paths from~$0$ to~$\om$.
The reason for this precise choice of coefficients will appear later (Theorem~\ref{thmrelDeomDeomp}). 

As a consequence, when the
minor~$\htb\ph$ is meromorphic, both $\De_\om(c\,\de+\htb\ph)$ and
$\De^+_\om(c\,\de+\htb\ph)$ coincide with $2\pi\I c_\om \de$, where
$c_\om$ is the residuum of~$\htb\ph$ at~$\om$ (\cf the remark just
before Example~\ref{exasimpleEPS}).
For instance, for the resurgent series associated with the names of
Euler, Poincar\'e, Stirling and Hurwitz,
\beglabel{eqDeomnames}
\De_{-1} \ti\ph\eul = 2\pi\I,  \quad
\De_{s+2\pi\I k} \ti\ph\poin = 2\pi\I, \quad 
\De_{2\pi\I m} \ti\mu = \frac{1}{m}, \quad
\De_{2\pi\I m} \ti\ph\hurw = 2\pi\I\frac{(2\pi\I m)^{s-1}}{\Ga(s)}, \quad
\elabel
for $k\in\Z$ and $\RE s < 0$ in the case of Poincar\'e, 
$m\in\Z^*$ for Stirling, and $s\in\N$ with $s\ge2$ for Hurwitz,
in view of Example~\ref{exasimpleEPS} and Exercise~\ref{exoHurwitz}.


We note for later use an immediate consequence of
formula~\eqref{eqalienoppartcases} of Lemma~\ref{lemSimplestab}:

\begin{lemma}   \label{lemcommutDeompaTc}
Let $\ti\ph \in \ti\gR_\Om\simp$ and $c\in \C$. Then
\begin{align}
\label{eqcommutDeompa}
\De^+_\om\pa\ti\ph &= (-\om+\pa)\De^+_\om\ti\ph, & 
\De_\om\pa\ti\ph &= (-\om+\pa)\De_\om\ti\ph, \\
\label{eqcommutDeomTc}
\De^+_\om T_c\ti\ph &= \ee^{-c\om} T_c\De^+_\om\ti\ph, &
\De_\om T_c\ti\ph &= \ee^{-c\om} T_c\De_\om\ti\ph, 
\end{align}
where $\pa = \frac{\dd\,}{\dd z}$ and~$T_c$ is defined by~\eqref{eqdefshift}.
\end{lemma}


\subsection{Dependence upon~$\Om$}


\begin{lemma}   \label{lemnodepOm}
Suppose that we are given $\om \in \C^*$ 
and~$\Om_1, \Om_2$ closed discrete such that $\om \in
\Om_1\cap\Om_2$.
Then there are two operators ``$\De_\om^+$'' defined
by~\eqref{eqdefDepomDeom}, 
an operator $\ti\gR_{\Om_1}\simp \to \ti\gR_{-\om+\Om_1}\simp$
and an operator $\ti\gR_{\Om_2}\simp \to \ti\gR_{-\om+\Om_2}\simp$,
but they act the same way on $\ti\gR_{\Om_1}\simp \cap \ti\gR_{\Om_2}\simp$.
The same is true of~``$\De_\om$".
\end{lemma}


The point is that the sets ${]0,\om[} \cap \Om_1$ and ${]0,\om[} \cap \Om_2$ may
differ, but their difference is constituted of points which
are artificial singularities for the minor of any
$\ti\ph \in \ti\gR_{\Om_1}\simp \cap \ti\gR_{\Om_2}\simp$, in the sense that no branch of
its analytic continuation is singular at any of these points.
So Lemma~\ref{lemnodepOm} claims that, in the above situation, we get the same resurgent series for
  $\De_\om^+\ti\ph$ whether computing it in $\ti\gR_{-\om+\Om_1}\simp$
  or in $\ti\gR_{-\om+\Om_2}\simp$.


\begin{proof}
Let $\Om \defeq \Om_1$, $\Om^* \defeq \Om_1 \cup \Om_2$ 
and $\ti\ph \in \ti\gR_{\Om_1}\simp \cap \ti\gR_{\Om_2}\simp$.
As in~\eqref{eqdefrcapOm}, we write
\[
{]0,\om[} \cap \Om =
\{ \om_1 \prec \cdots \prec \om_{r-1} \},
\qquad
{]0,\om[} \cap \Om^* =
\{ \om^*_1 \prec \cdots \prec \om^*_{s-1} \},
\]
with $1 \le r \le s$. Given $\eps^* \in \{+,-\}^{s-1}$, we have 
$\Al_{\om,\eps^*}^{\Om^*}\ti\ph = \Al_{\om,\eps}^{\Om}\ti\ph$
with $\eps \defeq \eps^*_{|\Om}$, \ie 
the tuple $\eps \in \{+,-\}^{r-1}$ is obtained from~$\eps^*$ by deleting the
symbols~$\eps^*_j$ corresponding to the fictitious singular points
$\om^*_j \in \Om^*\setminus\Om$.

In view of formula~\eqrefa{eqdefDepomDeom}, when $\eps^* = (+,\ldots,+)$
we get the same resurgent series for $\De_\om^+\ti\ph$ whether
computing it in $\ti\gR_{-\om+\Om^*}\simp$ or in
$\ti\gR_{-\om+\Om}\simp$,
which yields the desired conclusion by exchanging the roles of~$\Om_1$
and~$\Om_2$. 

We now compute $\De_\om\ti\ph$ in $\ti\gR_{-\om+\Om^*}\simp$ by
applying formula~\eqrefb{eqdefDepomDeom} and grouping together the
tuples~$\eps^*$ that have the same restriction~$\eps$: with the
notation $c \defeq s-r$, we get
\[
\sum_{\eps \in \{+,-\}^{r-1}} 
\sum_{\substack{ \eps^* \in \{+,-\}^{r+c-1} \\ \text{with}\; \eps^*_{|\Om}=\eps }} 
\frac{p(\eps^*)!q(\eps^*)!}{(r+c)!}  
\Al_{\om,\eps}^\Om\ti\ph
= \sum_{\eps \in \{+,-\}^{r-1}} \sum_{c=a+b} 
\binom{c}{a}
\frac{\big(p(\eps)+a\big)! \big(q(\eps)+b \big)!}{(r+c)!} 
\Al_{\om,\eps}^\Om\ti\ph,
\]
which yields the desired result because
\[
\sum_{c=a+b}  \binom{c}{a}
\frac{(p+a)! (q+b)!}{(r+c)!} 
= \frac{p!q!}{r!}
\]
for any non-negative integers $p,q,r$ with $r=p+q+1$, as is easily
checked by rewriting this identity as
\[
\sum_{c=a+b} \frac{(p+a)!}{a!p!} \frac{(q+b)!}{b!q!}
= \frac{(r+c)!}{c!r!}
\]
and observing that the generating series 
$\sum_{a\in\N} \frac{(p+a)!}{a!p!} X^a = (1-X)^{-p-1}$
satisfies $(1-X)^{-p-1} (1-X)^{-q-1} = (1-X)^{-r-1}$.
\end{proof}


\begin{rem}   \label{remnodeponOm}
  Given $\om\in\C^*$, we thus can compute $\De_\om^+\ti\ph$ or
  $\De_\om\ti\ph$ as soon as there exists~$\Om$ so that $\om\in\Om$ and
  $\ti\ph \in \ti\gR_\Om\simp$, and the result does not depend on~$\Om$.
We thus have in fact a family of operators
$\De^+_\om,\; \De_\om \col 
\ti\gR_\Om\simp \to \ti\gR_{-\om+\Om}\simp$,
indexed by the closed discrete sets~$\Om$ which contain~$\om$,
and there is no need that the notation for these operators depend explicitly
on~$\Om$.
\end{rem}

\subsection{The operators $\De_\om^+$ as a system of generators}

\begin{thm}
\label{thmDepomgen}
Let $\Om$ be a non-empty closed discrete subset of~$\C$ and let $\om\in\Om$.
Any alien operator at~$\om$ can be expressed as a linear combination
of composite operators of the form
\beglabel{eqdefcomposDep}
\bDe^+_{\eta_1,\ldots,\eta_s} \defeq
\De^+_{\eta_s-\eta_{s-1}} \circ \cdots 
\circ \De^+_{\eta_2-\eta_1} \circ \De^+_{\eta_1}
\elabel
with $s\ge1$, $\eta_1,\ldots,\eta_{s-1} \in\Om$, $\eta_s=\om$,
$\eta_1\neq0$ and $\eta_j\neq\eta_{j+1}$ for $1\le j<s$,
with the convention $\bDe^+_\om \defeq \De^+_\om$ for $s=1$
(viewing~$\bDe^+_{\eta_1,\ldots,\eta_s}$ as an operator
$\C\,\de \oplus \hat\gR_\Om\simp \to 
\C\,\de \oplus \hat\gR_{-\om+\Om}\simp$
or, equivalently,
$\ti\gR_\Om\simp \to \ti\gR_{-\om+\Om}\simp$).
\end{thm}

Observe that the composition~\eqref{eqdefcomposDep} is well defined
because, with the convention $\eta_0=0$, the operator
$\De^+_{\eta_j-\eta_{j-1}}$ maps
$\ti\gR\simp_{-\eta_{j-1}+\Om}$
into $\ti\gR\simp_{-\eta_j+\Om}$.
We shall not give the proof of this theorem, but let us indicate a few
examples: with the notation~\eqref{eqdefrcapOm},
\begin{gather*}
\Al_{\om_2,(+)}^\Om = \De^+_{\om_2}, \qquad
\Al_{\om_2,(-)}^\Om = 
\De^+_{\om_2} - \De^+_{\om_2-\om_1}\circ\De^+_{\om_1} \\[1ex]
\Al_{\om_3,(+,+)}^\Om = \De^+_{\om_3}, \qquad
\Al_{\om_3,(-,+)}^\Om = 
\De^+_{\om_3} - \De^+_{\om_3-\om_1}\circ\De^+_{\om_1}, \qquad
\Al_{\om_3,(+,-)}^\Om = 
\De^+_{\om_3} - \De^+_{\om_3-\om_2}\circ\De^+_{\om_2}, \\
\Al_{\om_3,(-,-)}^\Om =
\De^+_{\om_3} - \De^+_{\om_3-\om_1}\circ\De^+_{\om_1}
- \De^+_{\om_3-\om_2}\circ\De^+_{\om_2} +
\De^+_{\om_3-\om_2}\circ\De^+_{\om_2-\om_1}\circ\De^+_{\om_1}.
\end{gather*}


\begin{rem}   \label{remDeomgen}
One can omit the `$+$' in Theorem~\ref{thmDepomgen}, \ie the family
$\{\De_\eta\}$ is a system of generators as well.
This will follow from the relation~\eqref{eqforminverse} of next section.
\end{rem}


\begin{exo}
Suppose that $1 \le s \le r-1$ and $\eps,\eps^* \in \{+,-\}^{r-1}$
assume the form 
\[
\eps = a (-) b, \qquad \eps^* = a (+) b,
\qquad \text{with $a\in \{+,-\}^{s-1}$,}
\]
\ie $(\eps_1, \ldots, \eps_{s-1}) = (\eps^*_1, \ldots, \eps^*_{s-1}) =
a$,
$\eps_s = -$, $\eps^*_s = +$,
$(\eps_{s+1}, \ldots, \eps_{r-1}) = (\eps^*_{s+1}, \ldots,
\eps^*_{r-1}) = b$.
%
%
Prove that
\[
\Al_{\om_r,\eps}^\Om = \Al_{\om_r,\eps^*}^\Om -
\Al_{\om_r-\om_s,b}^\Om \circ \Al_{\om_s,a}^\Om
\]
with the notation~\eqref{eqdefrcapOm}.
Deduce the formulas given in example just above. 
\end{exo}


\begin{rem}
There is also a strong ``freeness'' statement for the operators~$\De^+_\eta$:
consider an arbitrary finite set~$F$ of finite sequences~$\etab$ of elements of~$\Om$,
so that each $\etab \in F$ is of the form $(\eta_1,\ldots,\eta_s)$ for
some $s\in \N$, with $\eta_1\neq0$ and $\eta_j\neq\eta_{j+1}$ for $1\le j<s$,
with the convention $\etab=\est$ and $\bDe^+_\est=\ID$
for $s=0$;
then, for any non-trivial family $\big(\ti\psi^\etab\big)_{\etab\in F}$ of
simple $\Om$-resurgent series,
\[
\ti\ph \in \gR_\Om\simp \mapsto \sum_{\etab\in F}
\ti\psi^\etab \cdot \bDe^+_\etab\ti\ph
\]
is a non-trivial linear map: one can construct a simple
$\Om$-resurgent series which is not annihilated by this operator.
There is a similar statement for the family $\{\De_\eta\}$.
See \cite[Vol.~3]{Eca81} or adapt \cite[\S 12]{mouldSN}.
\end{rem}

\section{The symbolic Stokes automorphism for a direction~$d$}

\subsection{Exponential-logarithm correspondence between $\{\De_\om^+\}$ and $\{\De_\om\}$}
\label{secexplogcorresp}

For any $\om\in\C^*$, we denote by $\prec$ the total order on
$[0,\om]$ induced by $t\in[0,1] \mapsto t\,\om \in [0,\om]$.
\begin{thm}     \label{thmrelDeomDeomp}
Let $\Om$ be a non-empty closed discrete subset of~$\C$.
Then, for any $\om\in\Om\setminus\{0\}$,
\begin{align}
\label{eqDeomlogDeomp}
\De_\om & = \sum_{s\in\N^*} \tfrac{(-1)^{s-1}}{s}
%
%
\sum_{ (\eta_1,\ldots,\eta_{s-1}) \in \Sig(s,\om,\Om) }
\De^+_{\om-\eta_{s-1}} \circ \cdots
\circ \De^+_{\eta_2-\eta_1} \circ \De^+_{\eta_1} \\[1ex]
\label{eqforminverse}
\De^+_\om & = \sum_{s\in\N^*} \tfrac{1}{s!}
%
%
\sum_{ (\eta_1,\ldots,\eta_{s-1}) \in \Sig(s,\om,\Om) }
\De_{\om-\eta_{s-1}} \circ \cdots
\circ \De_{\eta_2-\eta_1} \circ \De_{\eta_1}
\end{align}
where $\Sig(s,\om,\Om)$ is the set of all increasing sequences $(\eta_1,\ldots,\eta_{s-1})$
of ${]0,\om[}\cap\Om$,
\[
0 \prec \eta_1 \prec \cdots \prec \eta_{s-1} \prec \om,
\]
%
%
%
with the convention that the composite operator $\De^+_{\om-\eta_{s-1}} \circ \cdots
\circ \De^+_{\eta_2-\eta_1} \circ \De^+_{\eta_1}$ is reduced to~$\De^+_\om$
when $s=1$
(in which case $\Sig(1,\om,\Om)$ is reduced to the empty sequence)
and similarly for the composite operator appearing in~\eqref{eqforminverse}.
\end{thm}


With the notation~\eqref{eqdefrcapOm}, this means
\begin{gather*}
\De_{\om_1} = \De^+_{\om_1} \\[1.5ex]
\De_{\om_2} = \De^+_{\om_2} - 
\tfrac{1}{2}\De^+_{\om_2-\om_1} \circ \De^+_{\om_1} \\[1.5ex]
\De_{\om_3} = \De^+_{\om_3} - \tfrac{1}{2}\left( 
\De^+_{\om_3-\om_1} \circ \De^+_{\om_1}
+ \De^+_{\om_3-\om_2} \circ \De^+_{\om_2} \right) 
+ \tfrac{1}{3}\De^+_{\om_3-\om_2} \circ \De^+_{\om_2-\om_1} \circ \De^+_{\om_1} \\[.5ex]
\vdots\\[1ex]
%
%
\De^+_{\om_1} = \De_{\om_1} \\[1.5ex]
\De^+_{\om_2} = \De_{\om_2} +
\tfrac{1}{2!}\De_{\om_2-\om_1} \circ \De_{\om_1} \\[1.5ex]
\De^+_{\om_3} = \De_{\om_3} + \tfrac{1}{2!}\left( 
\De_{\om_3-\om_1} \circ \De_{\om_1}
+ \De_{\om_3-\om_2} \circ \De_{\om_2} \right) 
+ \tfrac{1}{3!}\De_{\om_3-\om_2} \circ \De_{\om_2-\om_1} \circ \De_{\om_1} \\[.5ex]
\vdots
\end{gather*}
%
%
We shall obtain Theorem~\ref{thmrelDeomDeomp} in next section as a consequence of
Theorem~\ref{thmDelogDep}, which is in fact an equivalent formulation
in term of series of homogeneous operators in a graded vector space.

\subsection{The symbolic Stokes automorphism and the symbolic Stokes
  infinitesimal generator}    \label{secSymbolStokes}

From now on, we fix $\Om$ and a ray
$d = \{ t\, \eith \mid t \ge 0 \}$, with some $\th\in\R$,
and denote by~$\prec$ the total order on~$d$ induced by $t\mapsto t\,\eith$.
We shall be interested in the operators $\De_\om^+$ and $\De_\om$ with
$\om\in d$.
Without loss of generality we can suppose that the set $\Om \cap d$ is
infinite and contains~$0$; indeed, if it is not the case, then we
can enrich~$\Om$ and replace it say with 
$\Om^* \defeq \Om \cup \{ N\,\eith \mid N\in\N \}$,
and avail ourselves of Remark~\ref{remnodeponOm}, observing that 
$\hat\gR_\Om\simp \hookrightarrow \hat\gR_{\Om^*}\simp$ and that any
relation proved for the alien operators in the larger space induces a
relation in the smaller, with $\De_{\om^*}^+$ and $\De_{\om^*}$
annihilating the smaller space when $\om^* \in \Om^*\setminus\Om$.

We can thus write $\Om \cap d$ as an increasing sequence
\beglabel{eqdefseqomm}
\Om \cap d = \{ \om_m \}_{m \in \N},
\qquad
\om_0 = 0 \prec \om_1 \prec \om_2 \prec \cdots
\elabel
For each $\om = \om_m \in \Om \cap d$, we define
\begin{itemize}
\item
  $\htn E_\om(\Om)$ as the space of all functions~$\hat\phi$ which are
  holomorphic at~$\om$, which can be analytically continued along
  any path of $\C\setminus\Om$ starting close enough to~$\om$, and
  whose analytic continuation has at worse simple singularities;
\item
$\trn E_\om(\Om)$ as the vector space $\C\,\de_\om \oplus \htn
E_\om(\Om)$, where each~$\de_\om$ is a distinct symbol\footnote{%
to be understood as a ``the translate of~$\de$ from~$0$ to~$\om$'', or
``the simple singularity at~$\om$ represented by $\frac{1}{2\pi\I(\ze-\om)}$''
} analogous to the convolution unit~$\de$;
\item
  $\chn E_\om(\Om,d)$ as the space of all functions~$\ch f$ holomorphic on the
  line segment $ ]\om_m, \om_{m+1}[ $ which can be analytically continued along any path of
  $\C\setminus\Om$ starting from this segment and whose analytic
  continuation has at worse simple singularities.
\end{itemize}
We shall often use abridged notations $\trn E_\om$ or~$\chn E_\om$.
Observe that there is a linear isomorphism
\beglabel{eqdefisomtaum}
\tau_\om \col \left|
\begin{array}{cccc}
\C\,\de \oplus \hat\gR_{-\om+\Om}\simp 
& \xrightarrow{\sim} &
\trn E_\om & \\[1ex]
a\,\de + \hat\ph 
& \mapsto &
a\,\de_\om + \hat\ph^\om, &
\qquad 
\hat\ph^\om(\ze) \defeq \hat\ph(\ze-\om),
\end{array} \right.
\elabel
and a linear map
\[ 
\ds \col \left|
\begin{array}{cccc}
\chn E_{\om_m}
& \to &
\trn E_{\om_{m+1}} & \\[1ex]
\ch f
& \mapsto &
\tau_{\om_{m+1}} \trb\ph, &
\qquad 
\trb\ph \defeq \sing_0 \big(\ch f(\om_{m+1}+\ze) \big).
\end{array} \right.
\]
The idea is that an element of $\trn E_\om(\Om)$ is nothing but a simple $\Om$-resurgent 
singularity ``based at~$\om$'' and that any element of $\chn E_\om(\Om,d)$
has a well-defined simple singularity ``at~$\om_{m+1}$'',
\ie we could have written
$\ds\ch f = \sing_{\om_{m+1}} \big( \ch f(\ze) \big)$
with an obvious extension of Definition~\ref{DefSingul}.

We also define a ``minor'' operator~$\mu$ and two ``lateral
continuation'' operators~$\dl+$ and~$\dl-$ by the formulas
\[
\mu \col \left|
\begin{array}{ccc}
\trn E_\om
& \to &
\chn E_\om \\[1ex]
a\,\de_\om + \hat\phi
& \mapsto &
\hat\phi_{| \, {]\om_m,\om_{m+1}[} }
\end{array} \right.\qquad\qquad
\dl\pm \col \left|
\begin{array}{ccc}
\chn E_\om
& \to &
\chn E_{\om_{m+1}} \\[1ex]
\ch f
& \mapsto &
\cont_{\ga_\pm}\ch f 
\end{array} \right.
\]
where~$\ga_+$, \resp $\ga_-$, is any path which connects 
$]\om_m,\om_{m+1}[$ and $]\om_{m+1},\om_{m+2}[$
staying in a neighbourhood of $]\om_m,\om_{m+2}[$ whose intersection
with~$\Om$ is reduced to $\{\om_{m+1}\}$
and circumventing~$\om_{m+1}$ to the right, \resp to the left.

Having done so for every $\om\in\Om\cap d$, we now ``gather'' the vector
spaces~$\trn E_\om$ or~$\chn E_\om$ and consider the completed graded
vector spaces
\[
\trn E(\Om,d) \defeq \bigop \trn E_\om(\Om), \qquad
\chn E(\Om,d) \defeq \bigop \chn E_\om(\Om,d)
\]
(we shall often use the abridged notations $\trn E$ or $\chn E$).
This means that, for instance, $\chn E$ is the cartesian product of
all spaces~$\chn E_\om$, but with additive notation for its
elements:
they are infinite series
\beglabel{eqdefchbphhomogcomp}
\chb\ph = \sum_{\om\in\Om\cap d}\chbph \om \in \chn E,
\qquad \chbph \om \in \chn E_\om \ens\text{for each $\om\in\Om\cap d$}.
\elabel
This way $\chn E_{\om_m} \hookrightarrow \chn E$ can be considered as the
subspace of homogeneous elements of degree~$m$ for each~$m$.
Beware that $\chb\ph\in\chn E$ may have infinitely many non-zero
homogeneous components~$\chbph \om$---this is the difference with the
direct sum\footnote{   \label{footdist}
One can define translation-invariant distances which
make~$\chn E$ and~$\trn E$ complete metric spaces as follows:
let $\ord \col \chn E \to \N\cup\{\infty\}$ be the ``order function''
associated with the decomposition in homogeneous components, \ie
$\ord\chb\ph \defeq \min\{\, m \in \N \mid \chbph{\om_m}\neq0 \,\}$ if
$\chb\ph\neq0$ and $\ord0=\infty$,
and let
$\dist(\chb\ph_1,\chb\ph_2) \defeq 2^{-\ord(\chb\ph_2-\chb\ph_1)}$, 
and similarly for~$\trn E$.
This allows one to consider a series of homogeneous components as the
limit of the sequence of its partial sums;
we thus can say that a series like~\eqref{eqdefchbphhomogcomp} is
convergent for the topology of the formal convergence (or
``formally convergent'').
Compare with Section~\ref{parKrull}.}
$\dst\bigoplus_{\om\in\Om\cap d}\chn E_\om$.

We get homogeneous maps
\[
\mu \col \trn E \to \chn E,
\qquad
\ds \col \chn E \to \trn E,
\qquad
\dl\pm \col \chn E \to \chn E
\]
by setting, for instance, 
\[
\dl+ \bigg( \sum_{\om\in\Om\cap d}\chbph \om \bigg) \defeq 
\sum_{\om\in\Om\cap d} \dl+ \chbph \om.
\]
The maps $\dl+, \dl-$ and~$\ds$ are $1$-homogeneous, in the sense
that for each~$m$ they map homogeneous elements of degree~$m$ to
homogeneous elements of degree $m+1$, while $\mu$ is $0$-homogeneous.
Notice that
\beglabel{eqmudsdlpdlm}
\mu \circ \ds = \dl+ - \dl-.
\elabel
For each $r\in\N^*$, let us define two $r$-homogeneous operators
$\dDep r, \, \dDe r \col \trn E \to \trn E$ 
by the formulas
\beglabel{eqdefdDeprdDer}
\dDep r \defeq \ds \circ \dll{+}{r-1} \circ \mu, \qquad
\dDe r \defeq  \sum_{\eps \in \{+,-\}^{r-1}} 
\frac{p(\eps)!q(\eps)!}{r!} \:
\ds \circ \dl{\eps_{r-1}} \circ \cdots \circ \dl{\eps_1} \circ \mu,
\elabel
with notations similar to those of~\eqref{eqdefDepomDeom}.


\begin{thm}   \label{thmDelogDep}
(i)
For each $m\in\N$ and $r\in\N^*$, the diagrams 
\[
\xymatrix @!0 @C=12em @R=10ex {
\makebox[6.5em]{$\C\,\de \oplus \hat\gR_{-\om_m+\Om}\simp$}
\ar[r]^*+{\De^+_{\om_{m+r}-\om_m}} \ar[d]_*+{\tau_{\om_m}}
& \makebox[7.5em]{$\C\,\de \oplus \hat\gR_{-\om_{m+r}+\Om}\simp$}
\ar[d]^*+{\tau_{\om_{m+r}}} \\
\makebox[5em]{$\trn E_{\om_m}(\Om,d)$} \ar[r]^*+{\dDep r}
& \makebox[6em]{$\trn E_{\om_{m+r}}(\Om,d)$}
}  \qquad
\xymatrix @!0 @C=12em @R=10ex {
\makebox[6.5em]{$\C\,\de \oplus \hat\gR_{-\om_m+\Om}\simp$} 
\ar[r]^*+{\De_{\om_{m+r}-\om_m}} \ar[d]_*+{\tau_{\om_m}} 
& \makebox[7.5em]{$\C\,\de \oplus \hat\gR_{-\om_{m+r}+\Om}\simp$} 
\ar[d]^*+{\tau_{\om_{m+r}}} \\
\makebox[5em]{$\trn E_{\om_m}(\Om,d)$} \ar[r]^*+{\dDe r}
& \makebox[6em]{$\trn E_{\om_{m+r}}(\Om,d)$}
}
\]
commute.

\medskip

(ii)
The formulas
$\DD+ d \defeq \ID + \sum_{r\in\N^*} \dDep r$
and $\DD{}d \defeq \sum_{r\in\N^*} \dDe r$
define two operators 
\[ \DD+ d,\,\DD{}d \col \trn E(\Om,d) \to \trn E(\Om,d), \]
the first of which has a well-defined logarithm which coincides with
the second,
\ie
\beglabel{eqDelogDep}
\sum_{r\in\N^*} \dDe r = 
\sum_{s\in\N^*} \frac{(-1)^{s-1}}{s}
\bigg(
\sum_{r\in\N^*} \dDep r
\bigg)^{\!\!s}.
\elabel

\medskip

(iii)
The operator~$\DD{}d$
has a well-defined exponential which coincides with~$\DD+ d$,
\ie
\beglabel{eqDepexpDe}
\ID + \sum_{r\in\N^*} \dDep r = 
\sum_{s\in\N} \frac{1}{s!}
\bigg(
\sum_{r\in\N^*} \dDe r
\bigg)^{\!\!s}.
\elabel
\end{thm}


\begin{proof}[Proof of Theorem~\ref{thmDelogDep}]
\begin{enumerate}[(i)]
\item
Put together~\eqref{eqdefDepomDeom}, \eqref{eqdefisomtaum}
and~\eqref{eqdefdDeprdDer}.
\item
The fact that $\DD+ d \col \trn E \to \trn E$ and its logarithm are
well-defined series of operators stems from the $r$-homogeneity
of~$\dDep r$ for every $r\ge 1$, which ensures formal convergence.

The \rhs\ of~\eqref{eqDelogDep} can be written
\[
\sum_{\substack{r_1,\ldots,r_s\ge1 \\ s\ge1}} \tfrac{(-1)^{s-1}}{s} 
\dDep{r_1} \cdots \dDep{r_s}
= 
\sum_{\substack{m_1,\ldots,m_s\ge0 \\ s\ge1}} \tfrac{(-1)^{s-1}}{s} 
\ds\, \dll{+}{m_1} \mu \ds\, \dll{+}{m_2} 
\cdots \mu \ds\, \dll{+}{m_s} \mu
= \sum_{r\ge1} \ds B_r \mu,
\]
where we have omitted the composition symbol ``$\circ$'' to lighten
notations, made use of~\eqref{eqdefdDeprdDer}, and availed ourselves
of~\eqref{eqmudsdlpdlm} to introduce the $(r-1)$-homogeneous operators
\[
B_r \defeq
\sum_{\substack{m_1+\cdots+m_s+s=r \\ m_1,\ldots,m_s\ge0, \ s\ge1}} 
\tfrac{(-1)^{s-1}}{s} \,
\dll{+}{m_1} (\dl+-\dl-) \, \dll{+}{m_2} 
\cdots (\dl+-\dl-) \, \dll{+}{m_s},
\]
with the convention $B_1=\ID$.
It is an exercise in non-commutative algebra to check that
\[
B_r = \sum_{\eps \in \{+,-\}^{r-1}} 
\frac{p(\eps)!q(\eps)!}{r!} \:
\dl{\eps_{r-1}} \cdots \dl{\eps_1}
\]
(viewed as an identity between polynomials in two non-commutative
variables $\dl+$ and~$\dl-$),
hence \eqref{eqdefdDeprdDer} shows that $\ds B_r \mu = \dDe r$ and we
are done.

\item
Clearly equivalent to~(ii).
\end{enumerate}
\end{proof}


\begin{Def}
\begin{enumerate}[--]
\item
The elements of $\trn E(\Om,d)$ are called \emph{$\Om$-resurgent symbols
with support in~$d$}.
\item
The operator~$\DD+ d$ is called the \emph{symbolic Stokes automorphism for
the direction~$d$}.
\item
The operator~$\DD{}d$ is called the \emph{symbolic Stokes
  infinitesimal generator for the direction~$d$}.
\end{enumerate}
\end{Def}

The connection between~$\DD+ d$ and the Stokes phenomenon will be
explained in next section.
This operator is clearly a linear invertible map, but there is a
further reason why it deserves the name ``automorphism'': we shall see
in Section~\ref{secActionStokesPr} that, when $\Om$ is stable under
addition, there is a natural algebra structure for which~$\DD+d$ is an
algebra automorphism.


\begin{proof}[Theorem~\ref{thmDelogDep} implies Theorem~\ref{thmrelDeomDeomp}]
Given~$\Om$ and a ray~$d$, Theorem~\ref{thmDelogDep}(i) says that
\beglabel{eqdDermDeom}
\DDe{r}{\om_m} = \tau_{\om_{m+r}} \circ \De_{\om_{m+r}-\om_m} \circ \tau_{\om_m}\ii,
\qquad
\DDep{r}{\om_m} = \tau_{\om_{m+r}} \circ \De^+_{\om_{m+r}-\om_m} \circ \tau_{\om_m}\ii
\elabel
for every~$m$ and~$r$.
By restricting the
identity~\eqref{eqDelogDep} to~$\trn E_0$ and extracting homogeneous
components we get the identity
\[
\DDe{r}{0} = \sum_{s\in\N^*} \tfrac{(-1)^{s-1}}{s} 
\sum_{r_1+\cdots+r_s= r}
\DDep{r_s}{\om_{r_1+\cdots+r_{s-1}}} \circ \cdots \circ \DDep{r_2}{\om_{r_1}} \circ \DDep{r_1}{0}
\]
for each $r\in\N^*$, which is equivalent, by~\eqref{eqdDermDeom}, to
\beglabel{eqequivDeomlogDeomp}
\De_{\om_r} = \sum_{s\in\N^*} \tfrac{(-1)^{s-1}}{s} 
\sum_{r_1+\cdots+r_s= r}
\De^+_{ \om_r - \om_{r_1+\cdots+r_{s-1}} } \circ \cdots 
\circ \De^+_{ \om_{r_1+r_2} - \om_{r_1} } \circ \De^+_{ \om_{r_1} }.
\elabel
Given $\om \in \Om\setminus\{0\}$, we can apply this with the ray
$\{ t \om \mid t\ge0 \}$:
the notations~\eqref{eqdefrcapOm} and~\eqref{eqdefseqomm} agree for
$1\le m < r$, with $r\in\N^*$ defined by $\om = \om_r$;
the identity~\eqref{eqequivDeomlogDeomp} is then seen to be equivalent
to~\eqref{eqDeomlogDeomp} by the change of indices
\[
\eta_1 = \om_{r_1}, \ens \eta_2 = \om_{r_1+r_2}, 
\ens \ldots \ens, \ens
\eta_{s-1} = \om_{r_1+\cdots+r_{s-1}}.
\]
The identity~\eqref{eqforminverse} is obtained the same way from~\eqref{eqDepexpDe}.
\end{proof}


\begin{exo}   \label{exoinvDDp}
Show that, for each $r\in\N^*$, the $r$-homogeneous component of 
\[ \DD-d \defeq \exp(-\DD{}d) = \big(\DD+d\big)\ii \] 
is
$\dDem r \defeq - \ds \circ \dll{-}{r-1} \circ \mu$,
giving rise to the family of operators 
$\De^-_\om \defeq - \Al_{\om,(-,\ldots,-)}^\Om$, $\om\in\Om\setminus\{0\}$.
\end{exo}

\subsection{Relation with the Laplace transform and the Stokes
  phenomenon}
\label{secrelLaplStokesPhen}

We keep the notations of the previous section, in particular
$d = \{ t\, \eith \mid t \ge 0 \}$ with $\th\in\R$ fixed.
With a view to use Borel-Laplace summation, we suppose that~$I$ is an
open interval of length~$<\pi$ which contains~$\th$, such that the
sector $\{\, \xi\,\ee^{\I\th'} \mid \xi>0,\; \th' \in I \,\}$
intersects~$\Om$ only along~$d$:
\[
\Om \cap \{\, \xi\,\ee^{\I\th'} \mid \xi\ge0,\; \th' \in I \,\} = \{ \om_m \}_{m \in \N}
\subset d,
\qquad
\om_0 = 0 \prec \om_1 \prec \om_2 \prec \cdots
\]
We then set
\[
I^+ \defeq \{\, \th^+ \in I \mid \th^+ < \th \,\}, \qquad
I^- \defeq \{\, \th^- \in I \mid \th^- > \th \,\}
\]
(mark the somewhat odd convention: the idea is that the directions
of~$I^+$ are to the right of~$d$, and those of~$I^-$ to the left)
and
\[
I_\epsilon^+ \defeq \{\, \th^+ \in I \mid \th^+ < \th-\epsilon \,\}, \qquad
I_\epsilon^- \defeq \{\, \th^- \in I \mid \th^- > \th+\epsilon \,\}
\]
{for $0 < \epsilon < \min\big\{ \tfrac{\pi}{2}, \dist(\th,\pa I) \big\}$.}

Let us give ourselves a locally bounded function $\ga \col I^+\cup I^-
\to \R$.
Recall that in Section~\ref{paragcNIga} we have defined the spaces
$\cN(I^\pm,\ga)$, consisting of holomorphic germs at~$0$ which extend
analytically to the sector $\{\, \xi\,\ee^{\I\th^\pm} \mid \xi>0,\; \th^\pm \in I \,\}$,
with at most exponential growth along each ray $\R^+\ee^{\I\th^\pm}$
as prescribed by~$\ga(\th^\pm)$,
and that according to Section~\ref{paraggDIga},
the Laplace transform gives rise to two operators $\cL^{I^+}$ and~$\cL^{I^-}$
defined on $\C\,\de\oplus\cN(I^+,\ga)$ and $\C\,\de\oplus\cN(I^-,\ga)$,
producing functions holomorphic in the domains $\gD(I^+,\ga)$ or
$\gD(I^-,\ga)$.

The domains $\gD(I^+,\ga)$ and $\gD(I^-,\ga)$ are sectorial
neighbourhoods of~$\infty$ which overlap:
their intersection is a sectorial neighbourhood of~$\infty$ centred on
the ray $\arg z = -\th$, with aperture~$\pi$.
For a formal series~$\ti\ph$ such that $\cB\ti\ph \in 
\C\,\de\oplus \big( \cN(I^+,\ga) \cap \cN(I^-,\ga) \big)$,
the Borel sums $\gS^{I^+}\ti\ph = \cL^{I^+}\cB\ti\ph$ and 
$\gS^{I^-}\ti\ph = \cL^{I^-}\cB\ti\ph$ may
differ, but their difference is exponentially small on
$\gD(I^+,\ga) \cap \gD(I^-,\ga)$.
We shall investigate more precisely this difference when $\cB\ti\ph$
satisfies further assumptions.

To get uniform estimates, we shall restrict to a domain of the form
$\gD(I_\epsilon^+,\ga+\epsilon) \cap \gD(I_\epsilon^-,\ga+\epsilon)$,
which is a sectorial neighbourhood of~$\infty$ of aperture
$\pi-2\epsilon$ centred on the ray $\arg z = -\th$.


\begin{nota}    \label{notatrigaE}
For each $m\in\N$, we set
$\trn E(\Om,d,m) \defeq \dst\bigoplus_{j=0}^m \trn E_{\om_j}(\Om)$
and denote by $[\,\cdot\,]_m$ the canonical
projection
\[
\Phi = \sum_{\om\in\Om\cap d} \Phi^\om \in \trn E(\Om,d) 
\mapsto 
[\Phi]_m \defeq \sum_{j=0}^m \Phi^{\om_j} \in \trn E(\Om,d,m).
\]
For each $\om \in \Om\cap d$, we set
\[
\trigaE{\om}(\Om) \defeq
\tau_\om \Big( \C\,\de \oplus \big( \hat\gR_{-\om+\Om}\simp
\cap \cN(I^+,\ga) \cap \cN(I^-,\ga) \big) \Big)
\subset \trn E_\om(\Om).
\]
We also define
$\trigaE{}(\Om,d,m) \defeq \dst\bigoplus_{j=0}^m \trigaE{\om_j}(\Om)
\subset \trn E(\Om,d,m)$,
on which we define the ``Laplace operators''~$\cL^+$ and~$\cL^-$ by
\begin{multline*}
\Phi = \sum_{j=0}^m \Phi^{\om_j} 
\mapsto \text{$\cL^\pm\Phi$ holomorphic in $\gD(I\pm,\ga)$,} \quad
\cL^\pm\Phi(z) \defeq \sum_{j=0}^m 
\ee^{-\om_j z} \cL^{I^\pm}( \tau_{\om_j}\ii \Phi^{\om_j} )(z).
\end{multline*}
\end{nota}


\begin{thm}   \label{thmSymbStokesLapl}
Consider $m\in\N$ and real numbers $\rho$ and~$\epsilon$ such that 
$\abs{\om_m} < \rho < \abs{\om_{m+1}}$
and $0 < \epsilon < \min\big\{ \tfrac{\pi}{2}, \dist(\th,\pa I) \big\}$.
Then, for any $\Phi \in \trigaE{}(\Om,d,m)$ such that
$[\DD+d\Phi]_m \in \trigaE{}(\Om,d,m)$,
one has
\beglabel{eqSymbolStokes}
\cL^+\Phi(z) = \cL^- [\DD+d\Phi]_m(z) + O(\ee^{-\rho\RE(\eith z)})
\elabel
uniformly for $z \in \gD(I_\epsilon^+,\ga+\epsilon) \cap \gD(I_\epsilon^-,\ga+\epsilon)$.
\end{thm}


\begin{figure}
\begin{center}

\includegraphics[scale=1]{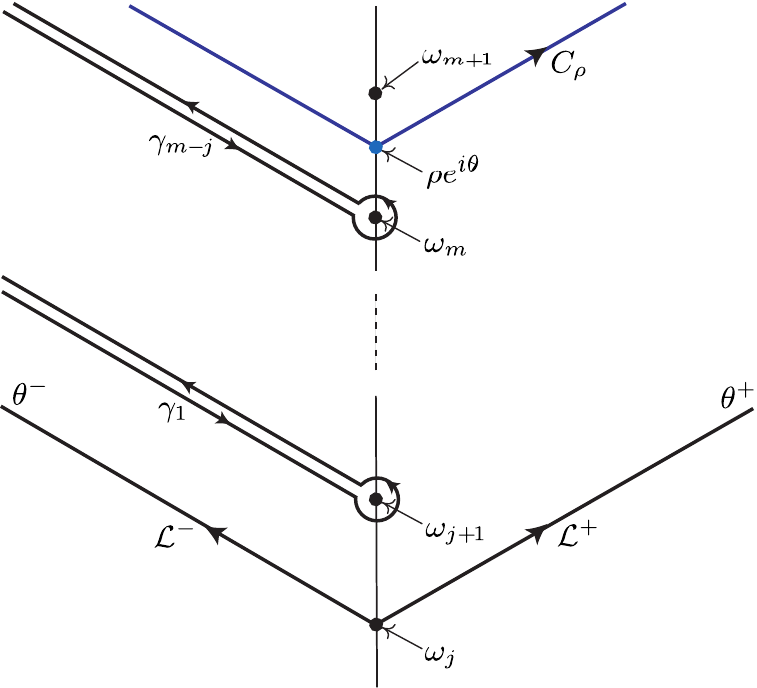}

\bigskip 

\caption{%
  From the Stokes phenomenon to the symbolic Stokes automorphism.}
\label{figStokes}

\end{center}
\end{figure}


\begin{proof}
It is sufficient to prove it for each homogeneous component of~$\Phi$, so
we can assume $\Phi = a\,\de_{\om_j}+\hat\ph \in \trigaE{\om_j}(\Om)$, with $0\le j \le m$.
Given $z \in \gD(I_\epsilon^+,\ga+\epsilon) \cap \gD(I_\epsilon^-,\ga+\epsilon)$, we choose $\th^+\in I^+$
and $\th^-\in I^-$ so that $\ze \mapsto \ee^{-z\ze}$ is exponentially
decreasing on the rays $\R^+\ee^{\I\th^\pm}$.
Then $\cL^\pm\Phi(z)$ can be written
$a\,\ee^{-\om_j z} + \int_{\om_j}^{\ee^{\I\th^\pm}\!\!\infty} 
\ee^{-z\ze} \hat\ph(\ze)\,\dd\ze$
(by the very definition of $\tau_{\om_j}$).
Decomposing the integration path as indicated on
Figure~\ref{figStokes}, we get
\begin{multline*}
\cL^+\Phi(z) = a\,\ee^{-\om_j z} + 
\Big( \int_{\om_j}^{\ee^{\I\th^-}\!\!\infty} + \int_{\ga_1} + 
\cdots + \int_{\ga_{m-j}} + \int_{C_\rho} \Big)
\ee^{-z\ze} \hat\ph(\ze)\,\dd\ze \\
= \cL^-\Phi(z) +
\sum_{r=1}^{m-j} \int_{\ga_r} \ee^{-z\ze} \,\dll{+}{r-1}\mu\Phi(\ze)
\,\dd\ze +
\int_{C_\rho} \ee^{-z\ze} \,\dll{+}{m-j-1}\mu\Phi(\ze)\,\dd\ze,
\end{multline*}
where the contour~$C_\rho$ consists of the negatively oriented half-line
$[\rho\,\eith,\ee^{\I\th^-}\!\!\infty)$ followed by the positively oriented half-line
$[\rho\,\eith,\ee^{\I\th^+}\!\!\infty)$.
We recognize in the $m-j$ terms of the sum in the \rhs\ the Laplace
integral of majors (\cf Section~\ref{paragLaplMaj}) applied to the
homogeneous components of $[\DD+d\Phi]_m$; all these integrals are
convergent by virtue of our hypothesis that $[\DD+d\Phi]_m \in
\trigaE{}(\Om,d,m)$, and also the last term in the \rhs\ is seen to be
a convergent integral which yields an $O(\ee^{-\rho\RE(\eith z)})$
error term.
\end{proof}


Observe that the meaning of~\eqref{eqSymbolStokes} for $\Phi =
\hat\ph \in \trigaE{0}(\Om)$, \ie
$\hat\ph \in \hat\gR_\Om\simp \cap \cN(I^+,\ga)\cap \cN(I^-,\ga)$,
is
\[
\cL^{\th^+}\hat\ph = \cL^{\th^-}\hat\ph +
\ee^{-\om_1 z} \cL^{\th^-}\De^+_{\om_1}\hat\ph(z) +
\cdots +
\ee^{-\om_m z} \cL^{\th^-}\De^+_{\om_m}\hat\ph(z) +
O(\ee^{-\rho\RE(\eith z)}).
\]
The idea is that the action of the symbolic Stokes automorphism yields
the exponentially small corrections needed to pass from the Borel sum
$\cL^{\th^+}\hat\ph$ to the Borel sum $\cL^{\th^-}\hat\ph$.
It is sometimes possible to pass to the limit $m\to\infty$ and to get
rid of any error term, in which case one could be tempted to write
\beglabel{eqtemptedheuristic}
\text{``}
\cL^+ = \cL^- \circ \DD+d
\text{''}.
\elabel


\begin{exa}
The simplest example of all is again provided by the Euler series, for
which there is only one singular ray, $d = \R^-$.
Taking any $\Om \subset \R^-$ containing~$-1$, we have
\beglabel{eqDDphatEul}
\DD+{\R^-}\hat\ph\eul = \hat\ph\eul + 2\pi\I\,\de_{-1}
\elabel
(in view of~\eqref{eqDeomnames}).
If we set $I^+ \defeq (\frac{\pi}{2},\pi)$ and $I^- \defeq
(\pi,\frac{3\pi}{2})$, then the functions $\ph^+ = \cL^+\hat\ph\eul$ and
$\ph^- = \cL^-\hat\ph\eul$ coincide with those of Section~\ref{secEulStokes}.
Recall that one can take
$\ga=0$ in this case, so $\ph^\pm$ is holomorphic in $\gD(I^\pm,0)$
(at least) and
the intersection $\gD(I^+,0) \cap \gD(I^-,0)$ is the half-plane $\{
\RE z < 0 \}$, on which Theorem~\ref{thmSymbStokesLapl} implies
\[
\ph^+ = \ph^- + 2 \pi\I \,\ee^z,
\]
which is consistent with formula~\eqref{eqmonodEul}.
\end{exa}


\begin{exa}
Similarly, for Poincar\'e's example with parameter $s\in\C$ of
negative real part, according to Section~\ref{sec:ReturnPoinc}, the
singular rays are $d_k \defeq \R^+\ee^{\I\th_k}$, $k\in\Z$, with 
$\om_k = s+2\pi\I k$ and 
$\th_k \defeq \arg\om_k \in (\frac{\pi}{2},\frac{3\pi}{2})$.
We take any $\Om$ contained the union of these rays and containing $s+2\pi\I\Z$.
For fixed~$k$, we can set $I^+ \defeq J_{k-1} = (\arg\om_{k-1},\arg\om_k)$,
$I^- \defeq J_k = (\arg\om_k,\arg\om_{k+1})$, and $\ga(\th)\equiv\cos\th$.
Then, according to Theorem~\ref{thmSommaPoin}, the Borel sums
$\cL^+\hat\ph\poin = \gS^{J_{k-1}}\ti\ph\poin$ and 
$\cL^-\hat\ph\poin = \gS^{J_k}\ti\ph\poin$ are well defined.
In view of~\eqref{eqDeomnames}, we have
\[
\DD+{d_k}\hat\ph\poin = \hat\ph\poin + 2\pi\I\,\de_{\om_k},
\quad \text{hence}\quad
\cL^+\hat\ph\poin = \cL^-\hat\ph\poin + 2\pi\I\,\ee^{-\om_k z}
\]
by Theorem~\ref{thmSymbStokesLapl},
which is consistent with~\eqref{eqdiffsumpoin}.
\end{exa}


\begin{exa}
The asymptotic expansion $\ti\ph\hurw_s(z)$ of the Hurwitz zeta
function was studied in Exercise~\ref{exoHurwitz}.
For $s\ge2$ integer, with $I = (-\frac{\pi}{2},\frac{\pi}{2})$, we
have \[ \ph^+(z) \defeq \sum_{k\in\N}(z+k)^{-s} = \gS^I\ti\ph\hurw_s(z)\]
for $z \in \gD(I,0) = \C\setminus\R^-$.
With the help of the difference equation $\ph(z)-\ph(z+1) = z^{-s}$,
it is an exercise to check that 
\[
\ph^-(z) \defeq -\sum_{k\in\N^*}(z-k)^{-s}
\]
coincides with the Borel sum $\gS^J\ti\ph\hurw_s$
defined on $\gD(J,0) = \C\setminus\R^+$, with 
$J = (\frac{\pi}{2},\frac{3\pi}{2})$
or $(-\frac{3\pi}{2},-\frac{\pi}{2})$.
In this case, we can take $\Om = 2\pi\I\Z^*$ and we have two singular
rays, $\I\R^+$ and~$\I\R^-$, for each of which the symbolic Stokes
automorphism yields infinitely many non-trivial homogeneous
components:
indeed, according to~\eqref{eqDeomnames},
\[
\DD+{\I\R^+}\hat\ph\hurw_s = \hat\ph\hurw_s +
\frac{2\pi\I}{\Ga(s)} \sum_{m=1}^\infty (2\pi\I m)^{s-1} \de_{2\pi\I  m},
\qquad
\DD+{\I\R^-}\hat\ph\hurw_s = \hat\ph\hurw_s +
\frac{2\pi\I}{\Ga(s)} \sum_{m=1}^\infty (-2\pi\I m)^{s-1} \de_{-2\pi\I  m}.
\]
Applying Theorem~\ref{thmSymbStokesLapl} with $I^+ =
(0,\frac{\pi}{2})$ and $I^- = (\frac{\pi}{2},\pi)$, or with $I^+ =
(-\pi,-\frac{\pi}{2})$ and $I^- = (-\frac{\pi}{2},0)$, for each $m\in \N$ we get
\begin{gather}
\label{eqdiffHurwbas}
\IM z < 0 \ens\Longrightarrow\ens
\ph^+(z) = \ph^-(z) +
\frac{2\pi\I}{\Ga(s)} \sum_{j=1}^m (2\pi\I j)^{s-1} \ee^{-2\pi\I jz}
+ O(\ee^{-2\pi(m+\demi)\abs{\IM z}}), \\[1ex]
\label{eqdiffHurwhaut}
\IM z > 0 \ens\Longrightarrow\ens
\ph^-(z) = \ph^+(z) +
\frac{2\pi\I}{\Ga(s)} \sum_{j=1}^m (-2\pi\I j)^{s-1} \ee^{2\pi\I jz}
+ O(\ee^{-2\pi(m+\demi)\abs{\IM z}})
\end{gather}
(and the constants implied in these estimates are uniform provided one
restricts oneself to $\abs{z} >\epsilon$ and 
$\abs{ \arg z \pm \I\frac{\pi}{2} } < \frac{\pi}{2} - \epsilon$).
In this case we see that we can pass to the limit $m\to\infty$ because
the finite sums involved in
\eqref{eqdiffHurwbas}--\eqref{eqdiffHurwhaut} are the partial sums of
convergent series.
In fact this could be guessed in advance: since~$\ph^+$ and~$\ph^-$
satisfy the same difference equation $\ph(z)-\ph(z+1) = z^{-s}$, their
difference yields $1$-periodic functions holomorphic in the
half-planes $\{ \IM z < 0 \}$ and $\{ \IM z > 0 \}$,
which thus have convergent Fourier series of the form\footnote{
See Section~\ref{secFourierHalfplane}.}
\[
(\ph^+ - \ph^-)_{|\{ \IM z < 0 \}} = \sum_{m\ge0} A_m \ee^{-2\pi\I mz},
\qquad
(\ph^+ - \ph^-)_{|\{ \IM z > 0 \}} = \sum_{m\ge0} B_m \ee^{2\pi\I mz},
\]
but the finite sums in \eqref{eqdiffHurwbas}--\eqref{eqdiffHurwhaut}
are nothing but the partial sums of these series (up to sign for the second).
So, in this case, the symbolic Stokes automorphism delivers the
Fourier coefficients of the diffence between the two Borel sums:
\[
\sum_{k\in\Z}(z+k)^{-s} = \begin{cases}
\dst \tfrac{2\pi\I}{\Ga(s)} \sum_{m=1}^\infty (2\pi\I m)^{s-1} \ee^{-2\pi\I mz}
& \text{for $\IM z < 0$,} \\[2ex]
\dst \tfrac{2\pi\I}{\Ga(s)} \sum_{m=1}^\infty (-1)^s(2\pi\I m)^{s-1} \ee^{2\pi\I mz}
& \text{for $\IM z > 0$.}
\end{cases}
\]
\end{exa}


\begin{exa}   \label{exaSymbStokesStir}
The case of the Stirling series~$\ti\mu$ studied in
Section~\ref{sec:Stirling} is somewhat similar,
with~\eqref{eqDeomnames} yielding
\beglabel{eqDDphatStirl}
\DD+{\I\R^+}\hat\mu = \hat\mu + \sum_{m\in\N^*} \frac{1}{m} \, \de_{2\pi\I m},
\qquad
\DD+{\I\R^-}\hat\mu = \hat\mu -\sum_{m\in\N^*} \frac{1}{m} \, \de_{-2\pi\I m}.
\elabel
Here we get
\begin{align}
\label{eqdiffStirlbas}
\IM z < 0 \ens&\Longrightarrow\ens
\mu^+(z) = \mu^-(z) +
\sum_{m=1}^\infty \tfrac{1}{m} \ee^{-2\pi\I mz}
= \mu^-(z) - \log(1-\ee^{-2\pi\I z}),\\[1ex]
\label{eqdiffStirlhaut}
\IM z > 0 \ens&\Longrightarrow\ens
\mu^-(z) = \mu^+(z)
- \sum_{m=1}^\infty \tfrac{1}{m} \ee^{2\pi\I mz}
= \mu^+(z) + \log(1-\ee^{2\pi\I z})
\end{align}
(compare with Exercise~\ref{exoStirlprepStokes}). 
\end{exa}

\subsection{Extension of the inverse Borel transform to
  $\Om$-resurgent symbols}
\label{secExtBorelSymb}

In the previous section, we have defined the Laplace
operators~$\cL^+$ and~$\cL^-$ on
\[
\trigaE{}(\Om,d,m) \subset \trn E(\Om,d,m) \subset \trn E(\Om,d),
\]
\ie the $\Om$-resurgent symbols to which they can be applied are
subjected to two constraints:
finitely many non-trivial homogeneous components, with at most
exponential growth at infinity for their minors.
There is a natural way to define on the whole of $\trn E(\Om,d)$ a
formal Laplace operator, which is an extension of the inverse Borel
transform~$\cB\ii$ on $\C\,\de \oplus \hat\gR_\Om\simp$.
Indeed, replacing the function $\ee^{-\om z} = \cL^\pm \de_\om$ by a
symbol $\ee^{-\om z}$, we define
\beglabel{eqdeftiEomOm}
\ti E_\om(\Om) \defeq \ee^{-\om z} \ti\gR_{-\om+\Om}\simp
\quad\text{for}\ens
\om\in\Om\cap d,
\qquad
\ti E(\Om,d) \defeq \bigop \ti E_\om(\Om),
\elabel
\ie we take the completed graded vector space obtained as cartesian
product of the spaces $\ti\gR_{-\om+\Om}\simp$, representing its elements
by formal expressions of the form
$\ti\Phi = \sum_{\om\in\Om\cap d} \ee^{-\om z}\ti\Phi_\om(z)$,
where each $\ti\Phi_\om(z)$ is a formal series and $\ee^{-\om z}$ is
just a symbol meant to distinguish the various homogeneous components.
We thus have for each $\om\in\Om\cap d$ a linear isomorphism 
\[
\ti\tau_\om \col \ti\ph(z)\in\ti\gR_{-\om+\Om}\simp 
\mapsto \ee^{-\om z}\ti\ph(z) \in \ti E_\om(\Om),
\]
which allow us to define 
\[
\cB_\om \defeq \tau_\om \circ \cB \circ \ti\tau_\om\ii
\col \ti E_\om(\Om) \xrightarrow{\sim} \trn E_\om(\Om).
\]
The map~$\cB_0$ can be identified with the Borel transform~$\cB$
acting on simple $\Om$-resurgent series;
putting together the maps~$\cB_\om$, $\om\in\Om\cap d$, we get a
linear isomorphism
\[
\cB \col \ti E(\Om,d) \xrightarrow{\sim} \trn E(\Om,d),
\]
which we can consider as the Borel transform acting on ``$\Om$-resurgent
transseries in the direction~$d$'',
and whose inverse can be considered as the formal Laplace transform
acting on $\Om$-resurgent symbols in the direction~$d$.


Observe that, if $\ee^{-\om z}\ti\ph(z) \in \ti E_\om(\Om)$ is such that
$\ti\ph(z)$ is $1$-summable in the directions of $I^+\cup I^-$,
then $\cB(\ee^{-\om z}\ti\ph) \in \trigaE{\om}(\Om)$ and
\[
\cL^\pm \cB(\ee^{-\om z}\ti\ph) = \ee^{-\om z} \gS^{I^\pm}\ti\ph.
\]
Beware that in the above identity, $\ee^{-\om z}$ is a \emph{symbol}
in the \lhs, whereas it is a \emph{function} in the \rhs.


Via~$\cB$, the operators~$\DD+d$ and~$\DD{}d$ give rise to operators
which we denote with the same symbols:
\[
\DD+d, \DD{}d \col \ti E(\Om,d) \to \ti E(\Om,d),
\]
so that we can \eg rephrase~\eqref{eqDDphatEul} as
\beglabel{eqDDptiEul}
\DD+{\R^-}\ti\ph\eul = \ti\ph\eul + 2\pi\I\,\ee^z
\elabel
or~\eqref{eqDDphatStirl} as
\beglabel{eqDDptiStirl}
\DD+{\I\R^+}\ti\mu = \ti\mu + \sum_{m\in\N^*} \frac{1}{m} \, \ee^{-2\pi\I mz},
\qquad
\DD+{\I\R^-}\ti\mu = \ti\mu -\sum_{m\in\N^*} \frac{1}{m} \, \ee^{2\pi\I mz}.
\elabel
In Section~\ref{secActionStokesPr},
we shall see that, if $\Om$ is stable under addition, then $\trn
E(\Om,d)$ and thus also $\ti E(\Om,d)$ have algebra structures, for
which it is legitimate to write
\beglabel{eqlogsymbolic}
-\log(1-\ee^{-2\pi\I z}) = \sum_{m\in\N^*} \frac{1}{m} \, \ee^{-2\pi\I
  mz}, \quad
\log(1-\ee^{2\pi\I z}) = -\sum_{m\in\N^*} \frac{1}{m} \, \ee^{2\pi\I
  mz}. 
\elabel


\begin{rem}   \label{rempadotted}
One can always extend the definition of $\pa = \frac{\dd\,}{\dd z}$ to
$\ti E(\Om,d)$ by setting
\[
\ti\ph \in \ti\gR_{-\om+\Om}\simp
\ens\Longrightarrow\ens
\pa(\ee^{-\om z}\ti\ph) \defeq \ee^{-\om z} (-\om+\pa)\ti\ph.
\]
(When~$\Om$ is stable under addition $\pa$ will be a derivation of the
algebra $\ti E(\Om,d)$, which will thus be a differential algebra.)

On the other hand, writing as usual $\Om\cap d = \{
0 = \om_0 \prec \om_1 \prec \om_2 \prec \cdots \}$,
we see that the homogeneous components of~$\DD+d$ and~$\DD{}d$ acting
on $\ti E(\Om,d)$
(Borel counterparts of the operators~$\dDep r, \, \dDe r \col \trn E
\to \trn E$ defined by~\eqref{eqdefdDeprdDer}) 
act as follows on $\ti E_\om(\Om)$ for each $\om = \om_m \in \Om \cap d$:
\beglabel{eqdottedalienop}
\ti\ph \in \ti\gR_{-\om_m+\Om}\simp
\ens\Longrightarrow\ens
\left\{ \begin{aligned} 
\dDep r(\ee^{-\om_m z}\ti\ph) &= \ee^{-\om_{m+r} z} \De^+_{\om_{m+r}-\om_m}\ti\ph,\\[1ex]
\ens
\dDe r(\ee^{-\om_m z}\ti\ph) &= \ee^{-\om_{m+r} z} \De_{\om_{m+r}-\om_m}\ti\ph.
\end{aligned} \right. 
\elabel
Formula~\eqref{eqcommutDeompa} then says
\[
\dDep r\pa\phi = \pa\dDep r\phi, \quad
\dDe r\pa\phi = \pa\dDe r\phi
\]
for every $\phi \in \ti E(\Om,d)$,
whence $\DD+d \circ \pa = \pa \circ \DD+d$ and $\DD{}d \circ \pa = \pa
\circ \DD{}d$.
\end{rem}

\section{The operators $\De_\om$ are derivations}
\label{secopDeomDeriv}

We now investigate the way the operators~$\De_\om$
and~$\De^+_\om$ act on a product of two terms (convolution product or
Cauchy product, according as one works with formal series or their
Borel transforms).

Let $\Om'$ and~$\Om''$ be non-empty closed discrete subsets of~$\C$
such that
\beglabel{eqassumpOm}
\Om \defeq \Om' \cup \Om'' \cup (\Om'+\Om'')
\elabel
is also closed and discrete.
Recall that, according to Theorem~\ref{thmAsymConv}, 
\[
\ti\ph \in \ti\gR_{\Om'} \ens\text{and}\ens
\ti\psi \in \ti\gR_{\Om''}
\quad \Longrightarrow \quad
\ti\ph \ti\psi \in \ti\gR_\Om.
\]

\subsection{Generalized Leibniz rule for the operators $\De^+_\om$}

We begin with the operators~$\De^+_\om$.


\begin{thm}   \label{ThmDepomAutom}
Let $\ti\ph \in \ti\gR_{\Om'}\simp$ and $\ti\psi \in \ti\gR_{\Om''}\simp$. Then
$\ti\ph \ti\psi \in \ti\gR_\Om\simp$ and, for every $\om \in \Om\setminus\{0\}$,
\beglabel{eqDeompLeibnm}
\De^+_\om(\ti\ph \ti\psi) = (\De^+_\om\ti\ph) \ti\psi + 
\sum_{ \substack{ \om = \om'+\om'' \\ \om'\in\Om'\cap{]0,\om[}, \, \om''\in\Om''\cap{]0,\om[} } }
(\De^+_{\om'}\ti\ph) (\De^+_{\om''}\ti\psi)
\quad + \; \ti\ph (\De^+_\om \ti\psi).
\elabel
%
%
%
\end{thm}


\begin{proof}
%
%
\textbf{a)}
The fact that $\ti\ph \ti\psi \in \ti\gR_\Om\simp$ follows from the
proof of formula~\eqref{eqDeompLeibnm} and Theorem~\ref{thmDepomgen},
we omit the details.

\medskip

\noindent \textbf{b)}
To prove formula~\eqref{eqDeompLeibnm}, we define
\[
\Sig_\om \defeq \{\, \eta\in {]0,\om[} \mid \eta \in \Om'\cup\Om''
\;\text{or}\; \om-\eta \in \Om'\cup\Om'' \,\}
\]
and write
$\cB\ti\ph = a\,\de + \hat\ph$, $\cB\ti\psi = b\,\de + \hat\psi$,
with $a, b \in \C$, $\hat\ph \in \hat\gR_{\Om'}\simp$,
$\hat\psi \in \hat\gR_{\Om''}\simp$,
\begin{align}
\label{eqDepetaph}
\cB \De^+_\eta \ti\ph &= a_\eta\,\de + \hat\ph_\eta,
& & a_\eta\in\C, \ens \hat\ph_\eta \in \hat\gR_{-\eta+\Om'}\simp,
& & \eta \in \Sig_\om \cup \{\om\} \\[1ex]
\label{eqDepommetapsi}
\cB \De^+_{\om-\eta} \ti\psi &= b_{\om-\eta}\,\de + \hat\psi_{\om-\eta},
& & b_{\om-\eta}\in\C, \ens \hat\psi_{\om-\eta} \in \hat\gR_{-(\om-\eta)+\Om''}\simp,
& & \eta \in \{0\} \cup \Sig_\om.
\end{align}
Since 
$\cB\De^+_\om(\ti\ph \ti\psi) = b \De^+_\om\hat\ph + a \De^+_\om\hat\psi
+ \De^+_\om(\hat\ph * \hat\psi)$,
formula~\eqref{eqDeompLeibnm} is equivalent to
\beglabel{eqlefttoprv}
\De^+_\om(\hat\ph * \hat\psi) =
\sum_{\eta \in \SSig}
(a_\eta \, \de + \hat\ph_\eta) * (b_{\om-\eta} \, \de + \hat\psi_{\om-\eta}),
\elabel
with the convention $a_0=0$, $\hat\ph_0=\hat\ph$ 
and $b_0=0$, $\hat\psi_0=\hat\psi$.

Consider a neighbourhood of $[0,\om]$
of the form $U_\de = \{\, \ze \in \C \mid \dist\big(\ze,[0,\om]\big) <\de \,\}$ 
with $\de>0$ small enough so that $U_\de\setminus[0,\om]$ does not
meet~$\Om$.
Let $u \defeq \om \,\ee^{-\I\al}$ with $0<\al<\frac{\pi}{2}$,
$\al$~small enough so that $u \in U_\de$ and the line segment $\ell \defeq
[0,u]$ can be considered as a path issuing from~$0$ circumventing to
the right all the points of
$ {]0,\om[} \cup \Om$.
We must show that $\cont_\ell(\hat\ph*\hat\psi)(\om+\ze)$ has a simple
singularity at~$0$ and compute this singularity.

\medskip

\noindent \textbf{c)}
We shall show that, when all the numbers~$a_\eta$
and~$b_{\om-\eta}$ vanish, 
\beglabel{eqcasevanish}
f(\ze) \defeq
\cont_\ell(\hat\ph*\hat\psi)(\om+\ze) =
\Bigg(
\sum_{\eta \in \SSig}
 \hat\ph_\eta * \hat\psi_{\om-\eta}
\Bigg)
\frac{\Llog\ze}{2\pi\I}
+ R(\ze),
\elabel
where $\Llog\ze$ is a branch of the logarithm and $R(\ze) \in \C\{\ze\}$.
This is sufficient to conclude, because in the general case we can
write
\[
\hat\ph * \hat\psi = \big( \tfrac{\dd\,}{\dd\ze} \big)^2 (\hat\ph^* *
\hat\psi^*),
\qquad \hat\ph^* \defeq 1*\hat\ph, \ens \hat\psi^* \defeq 1*\hat\psi,
\]
%
%
and, by Theorem~\ref{thmalieneasycase},
the anti-derivatives~$\hat\ph^*$ and~$\hat\psi^*$ satisfy
\[
\De^+_\eta \hat\ph^* = a_\eta + 1*\hat\ph_\eta, \qquad
\De^+_{\om-\eta} \hat\psi^* = b_{\om-\eta} + 1*\hat\psi_{\om-\eta}
\]
instead of \eqref{eqDepetaph}--\eqref{eqDepommetapsi},
thus we can apply~\eqref{eqcasevanish} to them and get
\begin{multline*}
\cont_\ell(\hat\ph^* * \hat\psi^*)(\om+\ze) =
\Bigg(
\sum_{\eta \in \SSig}
a_\eta b_{\om-\eta} \, \ze + a_\eta \, \ze * \hat\psi_{\om-\eta}
+ b_{\om-\eta} \, \ze * \hat\ph_\eta + \ze * \hat\ph_\eta * \hat\psi_{\om-\eta}
\Bigg)
\frac{\Llog\ze}{2\pi\I} \\[1ex]
+ R(\ze),
\end{multline*}
whence, by differentiating twice, a formula whose interpretation is precisely~\eqref{eqlefttoprv}
(because $(\frac{\dd\,}{\dd\ze} (\ze*A))/\ze$ and $(\ze*A)/{\ze^2}$
are regular at~$0$ for whatever regular germ~$A$).

\medskip

\noindent \textbf{d)}
From ow on, we thus suppose that all the numbers~$a_\eta$
and~$b_{\om-\eta}$ vanish.
Our aim is to prove~\eqref{eqcasevanish}.
We observe that
$D^+ \defeq \{\, \ze \in D(\om,\de) \mid
\IM(\ze/\om) < 0 \,\}$
is a half-disc such that, for all $\ze\in D^+$, the line segment
$[0,\ze]$ does not meet $\Om\setminus\{0\}$, hence
$\cont_\ell(\hat\ph*\hat\psi)(\ze) = \int_0^\ze
\hat\ph(\xi)\hat\psi(\ze-\xi) \, \dd\xi$
for such points.
We know by Section~\ref{sec_contconvolgen} that~$f$ has spiral
continuation around~$0$. 
Following the ideas of Section~\ref{secformalismsing},
we choose a determination of~$\arg\om$ and lift the half-disc
$-\om+D^+$ to the Riemann surface of logarithm by setting
$\ti D^+ \defeq \{\, \ze = r\,\eel^{\I\th} \in \ti\C \mid r < \de,\ 
\arg\om-\pi < \th <\arg\om \,\}$.
This way we can write $f = \ch f \circ \pi$, where~$\ch f$ is a
representative of a singular germ, explicitly defined on~$\ti D^+$ by
\beglabel{eqchfDp}
\ze \in \ti D^+ 
\ens\Longrightarrow\ens
\ch f(\ze) = \int_{\ell_{\pi(\ze)}}
\hat\ph(\xi)\hat\psi(\om+\pi(\ze)-\xi) \, \dd\xi
\quad\text{with $\ell_{\pi(\ze)} \defeq [0, \om+\pi(\ze)]$.}
\elabel
The analytic continuation of~$\ch f$ in
\[
\ti D^- \defeq \{\, \ze = r\,\eel^{\I\th} \in \ti\C \mid r < \de,\ 
\arg\om-3\pi < \th \le \arg\om-\pi \,\}
\]
is given by
\beglabel{eqchfDm}
\ze \in \ti D^-
\ens\Longrightarrow\ens
\ch f(\ze) = \int_{L_{\pi(\ze)}}
\hat\ph(\xi)\hat\psi(\om+\pi(\ze)-\xi) \, \dd\xi,
\elabel
where the symmetrically contractible path~$L_{\pi(\ze)}$ is obtained by
following the principles expounded in Section~\ref{sec_contconvolgen}
(\cf particularly~\eqref{eqcontgaphpsiHb});
this is illustrated in Figure~\ref{figPathsColDep}.


\begin{figure} 
\begin{center}

\begin{subfigure}{0.4\textwidth}
\centering
\includegraphics[scale=1]{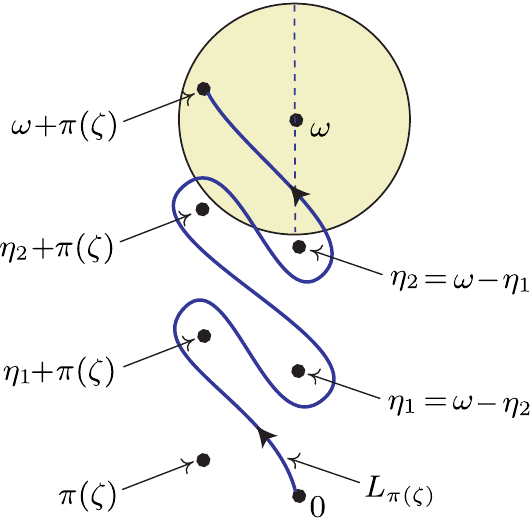}
\end{subfigure}
%
\begin{subfigure}{0.4\textwidth}
\centering
\includegraphics[scale=1]{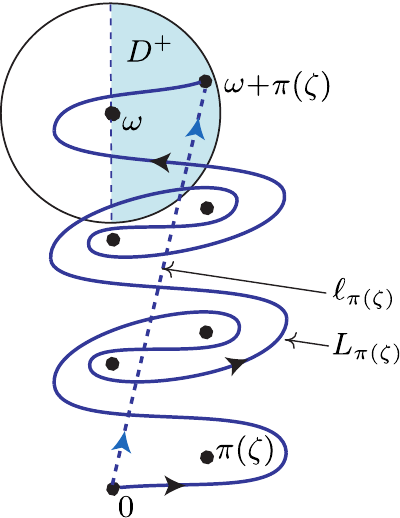}
\end{subfigure}

\bigskip

\caption{%
\emph{Integration paths for $\hat\ph*\hat\psi$.}
  Left: $L_{\pi(\ze)}$ for $\arg\om-2\pi < \arg\ze \le \arg\om-\pi$.\\
  Right: $\ell_{\pi(\ze)}$ and $L_{\pi(\ze)}$ for $\arg\om-3\pi < \arg\ze \le \arg\om-2\pi$.
  (Case when $\Sig_\om$ has two elements.)}
\label{figPathsColDep}

\end{center} 
\end{figure}


We first show that
\beglabel{eqvarchf}
\ze \in \ti D^+ 
\ens\Longrightarrow\ens
\ch f(\ze) - \ch f(\ze\,\eel^{-2\pi\I}) = 
\sum_{\eta \in \SSig}
 \hat\ph_\eta * \hat\psi_{\om-\eta}.
\elabel
The point is that $\Sig_\om$ is symmetric \wrt\ its midpoint
$\frac{\om}{2}$, thus of the form $\{ \eta_1 \prec \cdots \eta_{r-1}
\}$ with $\eta_{r-j} = \om - \eta_j$ for each~$j$,
and when~$\ze$ travels along a small circle around~$\om$, the ``moving
nail'' $\ze-\eta_j$ turns around the ``fixed nail''~$\eta_{r-j}$, 
to use the language of Section~\ref{secnails}.
Thus, for $\ze \in \ti D^+$, we can decompose the difference of paths
$\ell_{\pi(\ze)} - L_{\pi(\ze)}$ as on Figure~\ref{figComputColDep}
and get
\[
\ch f(\ze) - \ch f(\ze\,\eel^{-2\pi\I}) = 
\bigg( 
\int_{\pi(\ze)-\ga} + \int_{\om+\ga} +
\sum_{\eta \in \Sig_\om} \int_{\eta+\Ga}
\bigg) \hat\ph(\xi)\hat\psi(\om+\pi(\ze)-\xi) \, \dd\xi,
\]
where $\ga$ goes from $\ze\,\eel^{-2\pi\I}$ to~$\ze$ by turning
anticlockwise around~$0$,
whereas $\Ga$ goes from $\ze\,\eel^{-2\pi\I}$ to~$\ze$ the same way but
then comes back to $\ze\,\eel^{-2\pi\I}$ (see Figure~\ref{figComputColDep}).
%
%
\begin{figure} 
\begin{center}  

\includegraphics[scale=1]{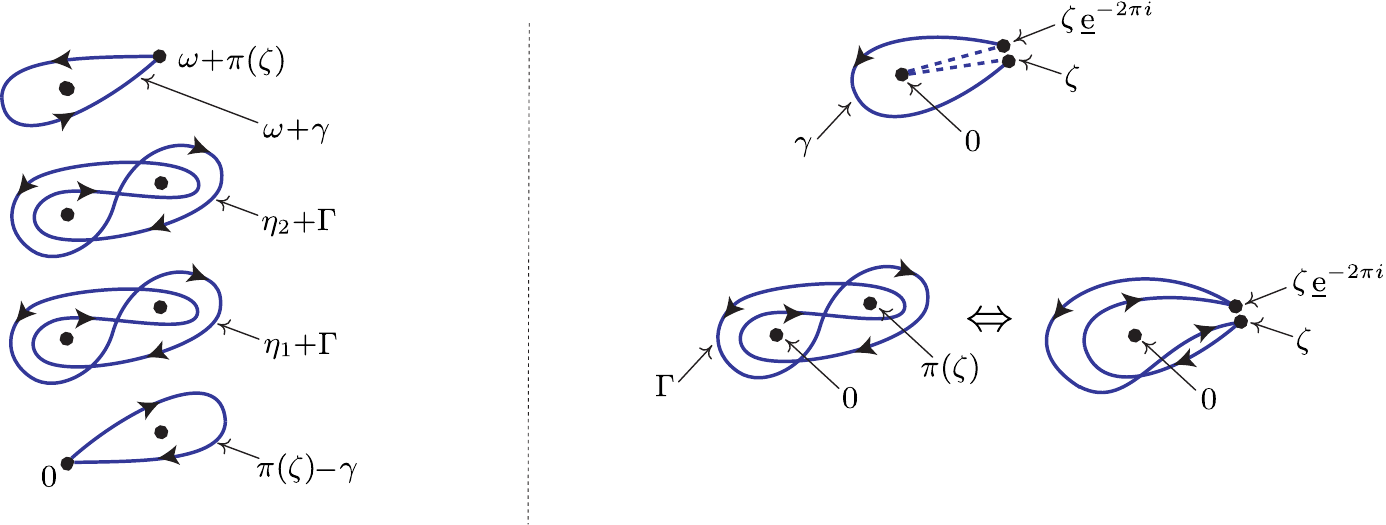}




\bigskip

\caption{%
Computation of the variation of the singularity at~$\om$ of $\hat\ph*\hat\psi$.}
\label{figComputColDep}

\end{center} 
\end{figure}
%
%
With an appropriate change of variable in each of these integrals,
this can be rewritten as
\begin{multline*}
\ch f(\ze) - \ch f(\ze\,\eel^{-2\pi\I}) = 
\int_\ga \hat\ph(\pi(\ze)-\xi) \hat\psi(\om+\xi)\,\dd\xi +
\int_\ga \hat\ph(\om+\xi) \hat\psi(\pi(\ze)-\xi)\,\dd\xi \\[1ex] +
\sum_{\eta \in \Sig_\om} \int_\Ga 
\hat\ph(\eta+\xi)\hat\psi(\om-\eta+\pi(\ze)-\xi) \,\dd\xi.
\end{multline*}
In the first two integrals, since
\[
\hat\psi(\om+\xi) = 
\frac{1}{2\pi\I} \hat\psi_\om(\xi) \Llog\xi
+ R'(\xi), \qquad
\hat\ph(\om+\xi) = 
\frac{1}{2\pi\I} \hat\ph_\om(\xi) \Llog\xi
+ R''(\xi),
\]
with $R'$ and~$R''$ regular at~$0$,
and we can diminish~$\de$ so that $\xi$ and $\pi(\ze)-\xi$ stay in a
neighbourhood of~$0$ where $\hat\ph$, $\hat\psi_\om$, $R'$,
$\hat\ph_\om$, $R''$ and~$\hat\psi$ are holomorphic, the Cauchy
theorem cancels the contribution of~$R'$ and~$R''$, while the
contribution of the logarithms can be computed by collapsing~$\ga$
onto the line segment $[\ze\,\eel^{-2\pi\I},0]$ followed by $[0,\ze]$,
hence the sum of the first two integrals is $\hat\ph*\hat\psi_\om + \hat\ph_\om*\hat\psi$.
Similarly, by collapsing~$\Ga$ as indicated on Figure~\ref{figComputColDep},
\[
\int_\Ga \hat\ph(\eta+\xi)
\hat\psi(\om-\eta+\pi(\ze)-\xi) \,\dd\xi
= \frac{1}{2\pi\I} \int_\Ga \hat\ph(\eta+\xi)
\big( \hat\psi_{\om-\eta}(\pi(\ze)-\xi) \Llog\xi + R_{\om-\eta}(\xi) \big)
\,\dd\xi
\]
(with some regular germ~$R_{\om-\eta}$)
is seen to coincide with
$\int_\ga \hat\ph(\eta+\xi) \hat\psi_{\om-\eta}(\pi(\ze)-\xi)
\,\dd\xi$,
which is itself seen to coincide with $\hat\ph_\eta*\hat\psi_{\om-\eta}(\ze)$
by arguing as above.
So~\eqref{eqvarchf} is proved.

\medskip

\noindent \textbf{e)}
We now observe that, since $g(\ze) \defeq \ch f(\ze) - \ch
f(\ze\,\eel^{-2\pi\I})$ is a regular germ at~$0$,
\[
R(\ze) \defeq f(\ze) - g(\ze) \frac{\Llog\ze}{2\pi\I}
\]
extends analytically to a (single-valued) function holomorphic in a
punctured disc, \ie it can be represented by a Laurent series~\eqref{eqLaurSer}.
But~$R(\ze)$ can be bounded by $C \big(C' + \ln\frac{1}{\abs\ze} \big)$ with appropriate constants
$C,C'$
(using~\eqref{eqchfDp}--\eqref{eqchfDm} to bound the analytic
continuation of~$f$),
thus the origin is a removable singularity for~$R$, which is thus regular at~$0$.
The proof of~\eqref{eqcasevanish} is now complete.
\end{proof}

\subsection{Action of the symbolic Stokes automorphism on a product}
\label{secActionStokesPr}

Theorem~\ref{ThmDepomAutom} can be rephrased in terms of the symbolic
Stokes automorphism~$\DD+d$ of Section~\ref{secSymbolStokes}.
Let us fix a ray $d = \{ t\, \eith \mid t \ge 0 \}$, with total
order~$\prec$ defined as previously.
Without loss of generality we can assume that both $\Om'\cap d$ and
$\Om''\cap d$ are infinite and contain~$0$.
With the convention $\De^+_0\defeq\ID$, formula~\eqref{eqDeompLeibnm} can be rewritten
\beglabel{eqrewrDepLeibn}
\De^+_\sig(\ti\ph \ti\psi) = 
\sum_{ \substack{ \sig = \sig'+\sig'' \\ \sigp, \,
    \sigpp } }
(\De^+_{\sig'}\ti\ph) (\De^+_{\sig''}\ti\psi),
\qquad \sig \in \Om\cap d.
\elabel
For every $\omp$ and $\ompp$ we have
commutative diagrams
\[
\xymatrix @!0 @C=10.5em @R=10ex {
\makebox[6.5em]{$\C\,\de \oplus \hat\gR_{-\om'+\Om'}\simp$}
\ar@{^{(}->}[r] \ar[d]_*+{\tau_{\om'}}
& \makebox[7.5em]{$\C\,\de \oplus \hat\gR_{-\om'+\Om}\simp$}
\ar[d]^*+{\tau_{\om'}} \\
\makebox[5.5em]{$\trn E_{\om'}(\Om',d)$} \ar@{^{(}->}[r]
& \makebox[6em]{$\trn E_{\om'}(\Om,d)$}
}  \qquad
\xymatrix @!0 @C=10.5em @R=10ex {
\makebox[6.5em]{$\C\,\de \oplus \hat\gR_{-\om''+\Om''}\simp$}
\ar@{^{(}->}[r] \ar[d]_*+{\tau_{\om''}}
& \makebox[7.5em]{$\C\,\de \oplus \hat\gR_{-\om''+\Om}\simp$}
\ar[d]^*+{\tau_{\om''}} \\
\makebox[5.5em]{$\trn E_{\om''}(\Om'',d)$} \ar@{^{(}->}[r]
& \makebox[6em]{$\trn E_{\om''}(\Om,d)$}
}%
\]
hence $\trn E(\Om',d) = \bigopp \trn E_{\om'}(\Om')$ 
and $\trn E(\Om'',d) = \bigoppp \trn E_{\om''}(\Om'')$
can be viewed as subspaces of $\trn E(\Om,d) \defeq \bigop \trn E_\om(\Om)$.
We shall often abbreviate the notations, writing for instance
\[
\Ep \hookrightarrow \trn E, \qquad \Epp \hookrightarrow \trn E.
\]
The convolution law 
$\big( \C\,\de \oplus \hat\gR_{-\om'+\Om'}\simp \big) \times
\big( \C\,\de \oplus \hat\gR_{-\om''+\Om''}\simp \big) \to
\C\,\de \oplus \hat\gR_{-(\om'+\om'')+\Om}\simp$
induces a bilinear map~$*$ defined by
\beglabel{eqdefconvolPhiPsi}
(\Phi,\Psi) =
\Big( \sum_{\omp} \ph^{\om'}, \sum_{\ompp} \psi^{\om''} \Big)   
\in \Ep \times \Epp \mapsto 
\sum_{\omp,\ \ompp} \ph^{\om'} * \psi^{\om''}
\in \trn E,
\elabel
where
\beglabel{eqdefconvolompompp}
(\ph,\psi) \in \Ep_{\om'}\times \Epp_{\om''} \ens\Longrightarrow\ens
\ph*\psi \defeq \tau_{\om'+\om''}\big( \tau_{\om'}\ii \ph * \tau_{\om''}\ii \psi \big)
\in \trn E_{\om'+\om''}.
\elabel


\begin{thm}    \label{thmStokesAlgMorphism}
With the above notations and definitions,
\beglabel{eqDepAutom}
(\Phi,\Psi) \in \Ep(\Om',d) \times \Epp(\Om'',d)
\quad\Longrightarrow\quad 
\DD+ d(\Phi*\Psi) = (\DD+ d\Phi) * (\DD+ d\Psi).
\elabel
\end{thm}

%
\begin{proof}
It is sufficient to prove~\eqref{eqDepAutom} for
$(\Phi,\Psi) = (\ph,\psi) \in \Ep_{\om'}\times \Epp_{\om''}$,
with $(\om',\om'') \in \Om' \times \Om''$. Recall that
\beglabel{eqrecallDep}
\DD+ d\ph = \sum_{ \eta' \succeq \om',\, \eta'\in\Om'\cap d }
\tau_{\eta'} \De^+_{ \eta'-\om' } \tau_{\om'}\ii \ph, \qquad
\DD+ d\psi = \sum_{ \eta'' \succeq \om'',\, \eta''\in\Om''\cap d }
\tau_{\eta''} \De^+_{ \eta''-\om'' } \tau_{\om''}\ii \psi.
\elabel
Let $\om \defeq \om'+\om''$, so that $\ph*\psi \in \trn E_\om$.
%
%
We have
\[
\DD+ d(\ph*\psi) = \sum_{ \eta\succeq\om,\, \eta\in\Om\cap d}
\tau_{ \eta } \De^+_{ \eta-\om } \tau_\om\ii(\ph*\psi)
= \sum_{ \eta\succeq\om,\, \eta\in\Om\cap d}
\tau_{ \eta } \De^+_{ \eta-\om } 
\big( (\tau_{\om'}\ii\ph)*(\tau_{\om''}\ii\psi) \big)
\]
by definition of~$\DD+ d$ and~$*$.
For each~$\eta$, applying~\eqref{eqrewrDepLeibn} to $\sig = \eta-\om$,
$\tau_{\om'}\ii\ph \in \C\,\de \oplus \hat\gR_{-\om'+\Om'}\simp$,
$\tau_{\om''}\ii\psi \in \C\,\de \oplus \hat\gR_{-\om''+\Om''}\simp$,
we get
\[
\De^+_{ \eta-\om } \big( (\tau_{\om'}\ii\ph)*(\tau_{\om''}\ii\psi) \big)
= \sum_{ \substack{
\eta-\om = \sig' + \sig'' \\
\sig' \in (-\om'+\Om')\cap d, \, \sig'' \in (-\om''+\Om'')\cap d
} }
\big( \De^+_{\sig'} \tau_{\om'}\ii\ph \big) *
\big( \De^+_{\sig''} \tau_{\om''}\ii\psi \big).
\]
With the change of indices $(\sig',\sig'') \mapsto (\eta',\eta'') =
(\om'+\sig',\om''+\sig'')$, this yields
\[
\DD+ d(\ph*\psi) = \sum_{ \substack{ \eta\in\Om\cap d \\ \eta\succeq\om} }
\, \sum_{ \substack{
\eta = \eta' + \eta'' \\
\eta' \in \Om'\cap d,\, \eta'' \in \Om''\cap d \\
\om'\preceq\eta',\, \om''\preceq\eta''
} }
\tau_{\eta} \big( (\De^+_{\eta'-\om'} \tau_{\om'}\ii\ph) *
(\De^+_{\eta''-\om''} \tau_{\om''}\ii\psi) \big).
\]
By Fubini, this is
\[
\DD+ d(\ph*\psi) = \sum_{ \substack{
\eta' \in \Om'\cap d,\, \eta'' \in \Om''\cap d \\
\om'\preceq\eta',\, \om''\preceq\eta''
} }
\tau_{\eta'+\eta''} \big( (\De^+_{\eta'-\om'} \tau_{\om'}\ii\ph) *
(\De^+_{\eta''-\om''} \tau_{\om''}\ii\psi) \big) =
(\DD+ d\ph) * (\DD+ d\psi)
\]
by definition of~$*$ and~\eqref{eqrecallDep}.
Hence~\eqref{eqDepAutom} is proved.
\end{proof}


\begin{rem}   \label{remStokesAlgMorphism}
When $\Om$ is stable under addition, one can take $\Om' = \Om'' =
\Om$. In that case, the operation~$*$ makes $\trn E(\Om,d)$ an algebra
and Theorem~\ref{thmStokesAlgMorphism} implies that~$\DD+d$ is an algebra automorphism.
At a heuristical level, this could be guessed from~\eqref{eqtemptedheuristic},
since both~$\cL^+$ and~$\cL^-$ take convolution products to pointwise products.
\end{rem}


\begin{rem}   \label{remBilinConc}
Via the linear isomorphism 
$\cB \col \ti E(\Om,d) \xrightarrow{\sim} \trn E(\Om,d)$ of
Section~\ref{secExtBorelSymb}, the bilinear map~$*$ gives rise to a
bilinear map
$\conc\, \col \ti E(\Om',d) \times \ti E(\Om'',d) \to \ti E(\Om,d)$
which, for homogeneous components, is simply
$\ee^{-\om' z}\ti\ph(z) \conc\, \ee^{-\om'' z}\ti\psi(z) =
\ee^{-(\om'+\om'')z} \ti\ph(z) \ti\psi(z)$.
This justifies~\eqref{eqlogsymbolic}.
\end{rem}

\subsection{Leibniz rule for the symbolic Stokes infinitesimal
  generator and the operators $\De_\om$}

From~$\DD+d$ we now wish to move on to its logarithm $\DD{}d$, which will give
us access to the way the operators~$\De_\om$ act on products.
We begin with a purely algebraic result,
according to which, roughly speaking, ``the logarithm of an automorphism is a derivation''.


\begin{lemma}   \label{lemPureAlg}
Suppose that $E$ is a vector space over~$\Q$, on which we have a
translation-invariant distance~$d$ which makes it a complete metric space,
and that $T \col E \to E$ is a $\Q$-linear contraction, so that
$D \defeq \log(\ID+T) = \sum_{s\ge1} \frac{(-1)^{s-1}}{s} T^s$ is well
defined.

Suppose that $E'$ and~$E''$ are $T$-invariant closed subspaces and
that $* \col E'\times E'' \to E$ is $\Q$-bilinear, with
$d(\Phi*\Psi,0) \le C d(\Phi,0) d(\Psi,0)$
for some $C>0$, and
\beglabel{eqIDpTautomGr}
(\Phi,\Psi) \in E'\times E''
\quad\Longrightarrow\quad 
(\ID+T)(\Phi*\Psi) = \big((\ID+T)\Phi \big) * \big((\ID+T)\Psi \big).
\elabel
Then 
\beglabel{eqDeDerivGr}
(\Phi,\Psi) \in E'\times E''
\quad\Longrightarrow\quad 
D(\Phi*\Psi) = (D\Phi) * \Psi + \Phi * (D\Psi).
\elabel
\end{lemma}


\begin{proof}
By~\eqref{eqIDpTautomGr}, $T(\Phi*\Psi) = (T\Phi)*\Psi + \Phi*(T\Psi)
+ (T\Phi) * (T\Psi)$.
Denoting by $N(s',s'',s)$ the coefficient of $X^{s'}Y^{s''}$ in the polynomial
$(X+Y+XY)^s \in \Z[X,Y]$
for any $s', s'', s \in \N$, we obtain by induction
\[
T^s(\Phi*\Psi) = \sum_{s',s''\in\N} N(s',s'',s) (T^{s'}\Phi) *(T^{s''}\Psi)
\]
for every $s\in\N$, whence
\[
D(\Phi*\Psi) = \sum_{s',s''\in\N} \sum_{s\in\N} \tfrac{(-1)^{s-1}}{s} N(s',s'',s) (T^{s'}\Phi) *(T^{s''}\Psi).
\]
The result follows from the fact that, for every $s',s''\in\N$, the number
$\sum \frac{(-1)^{s-1}}{s} N(s',s'',s)$
is the coefficient of $X^{s'}Y^{s''}$ in the formal series
$\sum \frac{(-1)^{s-1}}{s}(X+Y+XY)^s = \log(1+X+Y+XY) = \log(1+X) + \log(1+Y)
\in \Q[[X,Y]]$.
\end{proof}


The main result of this section follows easily:

\begin{thm}   \label{ThmDeomDeriv}
Under the assumption~\eqref{eqassumpOm}, one has for every direction~$d$
%
%
\beglabel{eqDeDirectDeriv}
(\Phi,\Psi) \in \Ep(\Om',d) \times \Epp(\Om'',d)
\quad\Longrightarrow\quad 
\DD{}d(\Phi*\Psi) = (\DD{}d\Phi) * \Psi + \Phi * (\DD{}d\Psi)
\elabel
and, for every $\om \in \Om\setminus\{0\}$,
\beglabel{eqDeomderiv}
(\ti\ph, \ti\psi) \in \ti\gR_{\Om'}\simp \times \ti\gR_{\Om''}\simp
\quad\Longrightarrow\quad 
\De_\om(\ti\ph \ti\psi) = (\De_\om\ti\ph) \ti\psi + \ti\ph (\De_\om \ti\psi).
\elabel
\end{thm}


\begin{proof}
The requirements of Lemma~\ref{lemPureAlg} are satisfied by $T
\defeq \DD+ d-\ID$ and the distance on~$\trn E$ indicated in
footnote~\ref{footdist};
since $\log\DD+d = \DD{}d$, this yields~\eqref{eqDeDirectDeriv}.

One gets~\eqref{eqDeomderiv} by evaluating~\eqref{eqDeDerivGr} with
$\Phi = \tau_0 \cB\ti\ph \in \Ep_0$ and
$\Psi = \tau_0 \cB\ti\psi \in \Epp_0$,
and extracting the homogeneous component 
$\tau_\om \De_\om (\cB\ti\ph * \cB\ti\psi) \in\trn E_\om$.
\end{proof}

\subsection{The subalgebra of simple $\Om$-resurgent functions}
\label{secsubalg}

We now suppose that $\Om$ is stable under addition, so that, by Corollary~\ref{corOmresursubalg},
$\ti\gR_\Om$ is a subalgebra of~$\C[[z\ii]]_1$ and
$\C\,\de\oplus\hat\gR_\Om$ is a subalgebra of the convolution algebra $\C\,\de\oplus\C\{\ze\}$.
Taking $\Om' = \Om'' = \Om$ in Theorem~\ref{ThmDepomAutom}, we get
\begin{cor}
If $\Om$ is stable under addition, then $\ti\gR_\Om\simp$ is a subalgebra of~$\ti\gR_\Om$ and
$\C\,\de\oplus\hat\gR_\Om\simp$ is a subalgebra of $\C\,\de\oplus\hat\gR_\Om$.
\end{cor}


As anticipated in Remark~\ref{remStokesAlgMorphism},
there is also for each ray~$d$ an algebra structure on $\trn E(\Om,d)$
given by the operation~$*$ defined in~\eqref{eqdefconvolPhiPsi},
for which the symbolic Stokes automorphism~$\DD+d$ is an algebra
automorphism;
the symbolic Stokes infinitesimal generator~$\DD{}d$ now appears
as a derivation, in view of formula~\eqref{eqDeDirectDeriv} of
Theorem~\ref{ThmDeomDeriv}
(for that reason~$\DD{}d$ is sometimes called ``directional alien derivation'').


\begin{rem}   \label{remcommutDDexpomz}
In particular, for each $\om\in\Om$ and $\Phi \in \ti E(\Om,d)$, we
have $\ee^{-\om z}\Phi \in \ti E(\Om,d)$,
\beglabel{eqcommutDDexpomz}
\DD+d (\ee^{-\om z}\Phi) = \ee^{-\om z} \DD+d \Phi, \qquad
\DD{}d (\ee^{-\om z}\Phi) = \ee^{-\om z} \DD{}d \Phi 
\elabel
(because $\ee^{-\om z}$ is fixed by~$\DD+d$ and annihilated by~$\DD{}d$).
\end{rem}


As indicated in formula~\eqref{eqDeomderiv} of Theorem~\ref{ThmDeomDeriv},
the homogeneous components~$\De_\om$ of~$\DD{}d$ inherit the Leibniz
rule, however it is only if $-\om+\Om \subset \Om$ that $\De_\om \col \ti\gR_\Om\simp
\to \ti\gR_\Om\simp$ is a derivation of the algebra~$\ti\gR_\Om\simp$,
and this is the case for all $\om\in\Om\setminus\{0\}$ when $\Om$ is an additive subgroup of~$\C$.
As anticipated in Remark~\ref{remAlienDeriv}, the operators~$\De_\om$
are called ``alien derivations'' for that reason.

Let us investigate farther the rules of ``alien calculus'' for
non-linear operations.


\begin{thm}   \label{thmComposSimpRes}
Suppose that $\Om$ is stable under addition.
Suppose that $\ti\ph(z), \ti\psi(z), \ti\chi(z) \in \ti\gR_\Om\simp$
and that $\ti\chi(z)$ has no constant term.
Let $H(t)\in\C\{t\}$.
Then
\[
\ti\psi\circ(\id+\ti\ph) \in \ti\gR_\Om\simp,
\quad H\circ\ti\chi \in \ti\gR_\Om\simp
\]
and, for any $\om \in \Om\setminus\{0\}$,
$(\De_\om\ti\psi) \circ (\id+\ti\ph) \in \ti\gR_{-\om+\Om}\simp$ and
\begin{align}
\label{eqDeomCompos}
\De_\om\big( \ti\psi\circ(\id+\ti\ph) \big) &=
(\pa\ti\psi)\circ(\id+\ti\ph) \cdot \De_\om\ti\ph
+ \ee^{-\om\ti\ph} \cdot (\De_\om\ti\psi) \circ (\id+\ti\ph), \\[1ex]
\label{eqDeomSubst}
\De_\om( H\circ\ti\chi ) &= (\tfrac{\dd H}{\dd t} \circ \ti\chi) \cdot \De_\om\ti\chi.
\end{align}
\end{thm}


The proof requires the following technical statement.
\begin{lemma}   \label{lemConvSeriesSimpRes}
Let $U \defeq \{\, r\,\eel^{\I\th} \in \ti\C \mid 0 < r < R,\ \th\in I \,\}$, 
where $I$ is an open interval of~$\R$ of length $>4\pi$ and $R>0$.
Suppose that, for each $k\in \N$, we are given a function~$\chb\ph_k$
which is holomorphic in~$U$ and is the major of a simple singularity
$a_k\,\de + \htb\ph_k$, 
and that the series $\sum\chb\ph_k$ converges normally on every compact
subset of~$U$.

Then the numerical series $\sum a_k$ is absolutely convergent, 
the series of functions $\sum \htb\ph_k$ converges normally on every
compact subset of $\D_R$,
and the function $\chb\ph \defeq \sum_{k\in\N} \chb\ph_k$, which is
holomorphic in~$U$, is the major of the simple singularity
$\big(\sum_{k\in\N} a_k \big) \de + \sum_{k\in\N} \htb\ph_k$.
\end{lemma}


\begin{proof}[Proof of Lemma~\ref{lemConvSeriesSimpRes}]
Pick~$\th_0$ such that $[\th_0,\th_0+4\pi] \subset I$ and let $J
\defeq [\th_0+2\pi,\th_0+4\pi]$.
For any $R'<R$,
writing $\htb\ph_k\big( \pi(\ze) \big) = \chb\ph_k(\ze)-\chb\ph_k(\ze\,\eel^{-2\pi\I})$
for $\ze \in U$ with $\arg\ze \in J$ and $\abs\ze \le R'$, we get the normal convergence of
$\sum\htb\ph_k$ on~$\ov\D_{R'}$.

Now, for each~$k$, $\chn L_k(\ze) \defeq \chb\ph_k(\ze) -
\htb\ph_k\big(\pi(\ze)\big)\frac{\log\ze}{2\pi\I}$ is a major of
$a_k\,\de$ and is holomorphic in~$U$; its monodromy is trivial, thus
$\chn L_k = L_k \circ \pi$ with~$L_k$ holomorphic in~$\D^*_R$.
For any circle~$C$ centred at~$0$, contained in~$\D_R$ and positively
oriented, we have $a_k = \int_C L_k(\ze) \, \dd\ze$.
The normal convergence of $\sum\chb\ph_k$ and $\sum\htb\ph_k$ implies
that of $\sum L_k$, hence the absolute convergence of $\sum a_k$.
Moreover, for every $n\in\N^*$, $\int_C L_k(\ze) \ze^{-n}\, \dd\ze = 0$, hence 
$L \defeq \sum_{k\in\N} L_k$ satisfies $\int_C L(\ze)  \ze^{-n}\, \dd\ze = 0$,
whence $\sing_0\big(L(\ze)\big) = \big(\sum_{k\in\N} a_k \big) \de$.

We conclude by observing that $\chb\ph(\ze) = L\big(\pi(\ze)\big) + 
\big( \sum_{k\in\N} \htb\ph_k(\pi(\ze)) \big) \frac{\log\ze}{2\pi\I}$.
\end{proof}


\begin{proof}[Proof of Theorem~\ref{thmComposSimpRes}]
We proceed as in the proof of Theorem~\ref{thmresOmstbNL}, writing 
$\ti\ph = a + \ti\ph_1$, $\ti\psi = b + \ti\psi_1$, where $a,b\in\C$
and $\ti\ph_1$ and~$\ti\psi_1$ have no constant term, 
and $H(t) = \sum_{k\ge0} h_k t^k$. Thus
\beglabel{eqdefSIMPlamu}
\ti\psi\circ(\id+\ti\ph) = b + \ti\la 
\quad\text{with}\ens
\ti\la \defeq T_a\ti\psi_1 \circ (\id+\ti\ph_1),
\quad
H\circ\ti\chi = h_0 + \ti\mu
\quad\text{with}\ens
\ti\mu \defeq \sum_{k\ge1} h_k \ti\chi^k.
\elabel
Both~$\ti\la$ and~$\ti\mu$ are naturally defined as formally
convergent series of formal series without constant term:
\[
\ti\la = \sum_{k\ge0} \ti\la_k
\quad\text{with}\ens
\ti\la_k \defeq \frac{1}{k!} (\pa^k T_a \ti\psi_1) \ti\ph_1^k,
\qquad
\ti\mu = \sum_{k\ge1} \ti\mu_k
\quad\text{with}\ens
\ti\mu_k \defeq h_k \ti\chi^k.
\]
By Lemma~\ref{lemcommutDeompaTc} and Theorem~\ref{ThmDepomAutom}, each
Borel transform
\[
\hat\la_k = \frac{1}{k!} \big((-\ze)^k\ee^{-a\ze}\hat\psi_1\big) * \hat\ph_1^{*k},
\qquad
\hat\mu_k = h_k \hat\chi^{*k}
\]
belongs to~$\hat\gR_\Om\simp$,
and we have checked in the proof of Theorem~\ref{thmresOmstbNL} that
their sums~$\hat\la$ and~$\hat\mu$ belong to~$\hat\gR_\Om$,
with their analytic continuations along the paths of $\C\setminus\Om$
given by the sums of the analytic continuations of the
functions~$\hat\la_k$ or~$\hat\mu_k$.
The argument was based on Lemma~\ref{lemestimiterconv}; we use it
again to control the behaviour of $\cont_\ga\hat\la$ and
$\cont_\ga\hat\mu$ near an arbitrary $\om \in \Om$, for a path
$\ga \col [0,1] \to \C\setminus\Om$ starting close to~$0$ and ending
close to~$\om$.
Choosing a lift~$\xi$ of $\ga(1)-\om$ in~$\ti\C$, we shall then apply
Lemma~\ref{lemConvSeriesSimpRes} to the functions
$\chb\ph_k(\ze)$
defined by $\cont_\ga\hat\la_k\big( \om + \pi(\ze) \big)$
or $\cont_\ga\hat\mu_k\big( \om + \pi(\ze) \big)$
for $\ze\in\ti\C$ close to~$\xi$.

Without loss of generality, we can suppose that $\abs{\ga(1)-\om} =
R/2$ with $R>0$ small enough so that $D(\om,R) \cap \Om = \{0\}$.
Let us extend~$\ga$ by a circle travelled twice, setting 
$\ga(t) \defeq \om + \big(\ga(1)-\om\big) \ee^{2\pi\I(t-1)}$ for $t\in[1,3]$.
For every $t\in[1,3]$ and $R_t < R/2$, we can apply
Lemma~\ref{lemestimiterconv} and get the normal convergence of
$\sum \cont_{ \ga_{|[0,t]} }\hat\la_k$ and $\sum \cont_{ \ga_{|[0,t]} }\hat\mu_k$
on $\ov{D\big( \ga(t),R_t \big)}$.
Now Lemma~\ref{lemConvSeriesSimpRes} shows that $\cont_\ga\hat\la$ and
$\cont_\ga\hat\mu$ have simple singularities at~$\om$.
Hence $\hat\la, \hat\mu \in \hat\gR_\Om\simp$.

A similar argument shows that
$(\De_\om\ti\psi) \circ (\id+\ti\ph) \in \ti\gR_{-\om+\Om}\simp$.

Lemma~\ref{lemConvSeriesSimpRes} also shows that
$\De_\om\hat\la = \sum_{k\ge0} \De_\om\hat\la_k$
and $\De_\om\hat\mu = \sum_{k\ge1} \De_\om\hat\mu_k$.
By means of~\eqref{eqDeomderiv}, we compute easily 
$\De_\om\ti\mu_k = k h_k \ti\chi^{k-1} \De_\om\ti\chi$,
whence $\De_\om\ti\mu = (\frac{\dd H}{\dd t}\circ\ti\chi) \cdot \De_\om\ti\chi$,
which yields~\eqref{eqDeomSubst} since~\eqref{eqdefSIMPlamu} shows that $\De_\om\ti\mu = \De_\om( H\circ\ti\chi )$.
By means of~\eqref{eqcommutDeompa}--\eqref{eqcommutDeomTc} and~\eqref{eqDeomderiv}, we compute 
\[
\De_\om\ti\la_k = A_k + B_k, \qquad
A_k \defeq \frac{k}{k!} (\pa^k T_a \ti\psi_1) \ti\ph_1^{k-1}
\De_\om\ti\ph_1,
\quad
B_k \defeq \frac{\ee^{-a\om}}{k!} \big( (-\om+\pa)^k T_a
\De_\om\ti\psi_1 \big) \ti\ph_1^k,
\]
$\sum_{k\ge0} A_k =  (\pa T_a \ti\psi_1) \circ (\id+\ti\ph_1)
 \cdot \De_\om\ti\ph_1
=  (\pa \ti\psi_1) \circ (\id+\ti\ph)
 \cdot \De_\om\ti\ph_1
=  (\pa \ti\psi) \circ (\id+\ti\ph)
 \cdot \De_\om\ti\ph$,
and
\begin{multline*}
\sum_{k\ge0} B_k = \ee^{-a\om} \sum_{k',k''\ge0} 
\frac{ (-\om)^{k'} }{ k'! k''! } (\pa^{k''} T_a \De_\om \ti\psi_1)
\ti\ph_1^{k'+k''} \\[1ex]
= \ee^{-a\om} \sum_{k'\ge0} \frac{ (-\om)^{k'} }{ k'! } \ti\ph_1^{k'}
\sum_{k''\ge0} \frac{1}{k''!} (\pa^{k''} T_a \De_\om \ti\psi_1) \ti\ph_1^{k''}
= \exp(-a\om -\om\ti\ph_1) \cdot 
( T_a \De_\om \ti\psi_1 ) \circ ( \id + \ti\ph_1 ) \\[1ex]
= \ee^{-\om\ti\ph} \cdot ( \De_\om \ti\psi_1 ) \circ ( \id + \ti\ph )
= \ee^{-\om\ti\ph} \cdot ( \De_\om \ti\psi ) \circ ( \id + \ti\ph ),
\end{multline*}
which yields~\eqref{eqDeomCompos} since~\eqref{eqdefSIMPlamu} shows that 
$\De_\om\ti\la = \De_\om \big( \ti\psi\circ(\id+\ti\ph) \big)$.
\end{proof}


\begin{exa}
\label{secAlienOpNL}
As promised in Example~\ref{exasimpleEPS}, we can now study the
exponential of the Stirling series $\ti\mu \in \ti\gR_{2\pi\Z^*}$. 
Since $2\pi\I\Z^*$ is not stable under addition, we need to take at
least $\Om = 2\pi\I\Z$ to ensure $\ti\la = \exp \ti\mu \in
\ti\gR_\Om\simp$.
Formulas \eqref{eqDeomnames} and~\eqref{eqDeomSubst} yield 
\beglabel{eqDeomtila}
\De_{2\pi\I m}\ti\la = \frac{1}{m} \ti\la,
\qquad m\in\Z^*.
\elabel
In view of Remark~\ref{remDeomgen}, this implies that 
\emph{any alien operator maps~$\ti\la$ to a multiple of~$\ti\la$}.
This clearly shows that the analytic continuation of the Borel
transform $\cB(\ti\la-1)$ is multiple-valued, since \eg
\eqref{eqDeomtila} with $m=\pm1$ says that the singularity
at~$\pm2\pi\I$ of the principal branch has a non-trivial minor.
Let us show that
\beglabel{eqDeptila}
\De^+_{2\pi\I m}\ti\la = \begin{cases}
\ti\la & \text{for $m = -1, +1, +2, +3, \ldots$} \\
0  & \text{for $m = -2, -3, \ldots$}
\end{cases}
\elabel
(notice that the last formula implies that the analytic continuation
of $\cB(\ti\la-1)$ from the line segment $(-2\pi\I,2\pi\I)$ to $(-2\pi\I,-4\pi\I)$ obtained by
circumventing $-2\pi\I$ to the right is free of singularity in the
rest of~$\I\R^-$: it extends analytically to
$\C\setminus[-2\pi\I,+\I\infty)$, but that this is not the case of the
analytic continuation to the left!)

Formula~\eqref{eqDeptila} could probably be obtained from the relation
$\De^+_{2\pi\I m}\ti\mu = \frac{1}{m}$ by repeated use
of~\eqref{eqDeompLeibnm}, but it is simpler to
use~\eqref{eqforminverse} and~\eqref{eqDeomtila}, 
and even better to perform the computation at the level of the symbolic
Stokes automorphism and its infinitesimal generator.
This time, we manipulate the multiplicative counterpart of $\DD+{\I\R^\pm}$
and~$\DD{}{\I\R^\pm}$ obtained through~$\cB$ as indicated in
Section~\ref{secExtBorelSymb} and Remark~\ref{remBilinConc}, writing for
instance
\[
\DD{}{\I\R^+}\ti\la = \sum_{m\in\N^*} \frac{1}{m} \, \ee^{-2\pi\I mz} \ti\la
= - \log(1-\ee^{-2\pi\I z}) \ti\la,
\quad
\DD{}{\I\R^-}\ti\la = -\sum_{m\in\N^*} \frac{1}{m} \, \ee^{2\pi\I mz} \ti\la
= \log(1-\ee^{2\pi\I z}) \ti\la.
\]
By exponentiating in $\ti E(\Om,\I\R^+)$ or $\ti E(\Om,\I\R^-)$, with
the help of~\eqref{eqcommutDDexpomz}, we get
\beglabel{eqSymbStokesExpStirl}
\DD+{\I\R^+}\ti\la = (1-\ee^{-2\pi\I z})\ii \ti\la =
\sum_{m\in\N} \ee^{-2\pi\I mz} \ti\la,
\qquad
\DD+{\I\R^-}\ti\la = (1-\ee^{2\pi\I z}) \ti\la 
= \ti\la - \ee^{2\pi\I z} \ti\la.
\elabel
One gets~\eqref{eqDeptila} by extracting the homogeneous components of
these identities.


The Stokes phenomenon for the two Borel sums~$\ti\la(z)$ can be described as follows:
with $I \defeq (-\frac{\pi}{2},\frac{\pi}{2})$, we have $\la^+ \defeq
\la = \gS^I\ti\la$ holomorphic in $\C\setminus\R^-$,
and with $J \defeq (\frac{\pi}{2},\frac{3\pi}{2})$, we have $\la^- \defeq
\gS^J\ti\la$ holomorphic in $\C\setminus\R^+$;
by adapting the chain of reasoning of Example~\ref{exaSymbStokesStir},
one can deduce from~\eqref{eqSymbStokesExpStirl} that
\[
\IM z < 0 \ens\Longrightarrow\ens
\la^+(z) = (1-\ee^{-2\pi\I z})\ii \la^-(z),
\qquad
\IM z > 0 \ens\Longrightarrow\ens
\la^-(z) = (1-\ee^{2\pi\I z}) \la^+(z)
\]
(one can also content oneself with exponentiating \eqref{eqdiffStirlbas}--\eqref{eqdiffStirlhaut}),
getting thus access to the exponentially small discrepancies between
both Borel sums.

Observe that it follows that~$\la^\pm$ admits a multiple-valued meromorphic continuation
which gives rise to a function meromorphic in the whole of~$\ti\C$:
for instance, since $\la^+_{| \{\IM z>0 \} }$ coincides with
$(1-\ee^{2\pi\I z})\ii\la^-$, it can be meromorphically continued to
$\C\setminus\R^-$ and its anticlockwise continuation to $\{ \IM z < 0
\}$ is given by
$(1-\ee^{2\pi\I z})\ii\la^-_{| \{\IM z<0 \} }$, which coincides with 
$(1-\ee^{2\pi\I z})\ii (1-\ee^{-2\pi\I z}) \la^+_{| \{\IM z<0 \} }$,
and can thus be anticlockwise continued to $\{ \IM z > 0 \}$:
we find
\[
\la^+(\eel^{2\pi\I} z) = 
(1-\ee^{2\pi\I z})\ii (1-\ee^{-2\pi\I z}) \la^+(z) = 
-\ee^{-2\pi\I z}\la^+(z)
\]
(compare with Remark~\ref{remStokeslaStirl}).
Since $z^{-\demi+z} = \ee^{(-\demi+z)\log z}$ gets multiplied by $-\ee^{2\pi\I z}$
after one anticlockwise turn around~$0$,
we can deduce that the product $\sqrt{2\pi} \, \ee^{-z}
z^{-\demi+z}\la^+(z)$ is single-valued,
not a surprise in view of~\eqref{eqleftStirlingf}:
this product function is none other than Euler's gamma function, which
is known to be meromorphic in the whole complex plane!
\end{exa}

\section{A glance at a class of non-linear differential equations}

We give here a brief account of some results based on \'Ecalle's works
that one can find in \cite{mouldSN}. 
Our purpose is to illustrate alien calculus on the example of the
simple $\Z$-resurgent series which appear when dealing with a
non-linear generalization of the Euler equation.
In this section, we omit most of the proofs but try to acquaint the reader with concrete computations with
alien operators.

%
\parag
Let us give ourselves 
\[
B(z,y) = \sum_{n\in\N} b_n(z) y^n \in \C\{z\ii,y\},
\ens\text{with $b_1(z) = 1+O(z^{-2})$ and $b_n(z) = O(z\ii)$ if $n\neq1$,}
\]
and consider the differential equation
\beglabel{eqEulerNL}
\frac{\dd\ti\phi}{\dd z} = B(z,\ti\phi) =
b_0(z) + b_1(z) \ti\phi + b_2(z) \ti\phi^2 + \cdots
\elabel
(one recovers the Euler equation for $B(z,y) = -z\ii + y$).
Observe that if $\ti\phi(z) \in \zcz$ then $B(z,\ti\phi(z))$ is given
by a formally convergent series, so~\eqref{eqEulerNL} makes sense.


\begin{thm}   \label{thmSimpResPhiZero}
Equation~\eqref{eqEulerNL} admits a unique formal solution $\ti\phi_0
\in \zcz$.
This formal series is $1$-summable in the directions of $(-\pi,\pi)$
and 
\[
\ti\phi_0(z) \in \ti\gR_{\Z^*_-}\simp,
\quad \text{where $\Z^*_- \defeq \{ -1, -2, -3, \ldots \}$.}
\]
Its Borel sum $\gS^{(-\pi,\pi)}\ti\phi_0$ is a particular solution of
Equation~\eqref{eqEulerNL}, defined and holomorphic in a domain of the
form
$\ti\gD\big( (-\pi,\pi), \ga \big) \subset \ti\C$.
\end{thm}


We omit the proof, which can be found in \cite{mouldSN}.
Let us only give a hint on why one must take $\Om=\Z^*_-$.
Writing $B(z,y) - y = \sum a_n(z) y^n$, we have $a_n(z) \in
z\ii\C\{z\ii\}$ for all $n\in\N$, thus~\eqref{eqEulerNL} can be
rewritten
$\frac{\dd\ti\phi}{\dd z} -\ti\phi = \sum a_n \ti\phi^n$, which
via~$\cB$ is equivalent to
\[
\hat\phi_0(\ze) = \frac{-1}{1+\ze} (\hat a_0 +
\hat a_1 * \hat\phi + \hat a_2 * \hat\phi^{*2} + \cdots).
\]
The Borel transforms~$\hat a_n$ are entire functions, thus it is only
the division by $1+\ze$ which is responsible for the appearance of
singularities in the Borel plane: a pole at~$-1$ in the first place,
but also, because of repeated convolutions, a simple singularity
at~$-1$ rather than only a simple pole and other simple singularities
at all points of the additive semigroup generated by~$-1$.

%
\parag   \label{paragCalculaEnrich}
The next question is: what about the Stokes phenomenon for~$\ti\phi_0$
and the action of the alien operators?
Let us first show how, taking for granted that $\ti\phi_0 \in
\ti\gR_{\Z^*_-}\simp$, one can by elementary alien calculus see that
$\De_\om\ti\phi_0 = 0$ for $\om\neq-1$ and compute $\De_{-1}\ti\phi_0$ up
to a multiplicative factor.
We just need to enrich our ``alien toolbox'' with two lemmas.


\begin{nota}
Since $\pa = \frac{\dd\,}{\dd z}$ increases the standard valuation by
at least one unit (\cf \eqref{ineqvald}), the operator $\mu+\pa \col
\C[[z\ii]] \to \C[[z\ii]]$ is invertible for any $\mu\in\C^*$ and its
inverse $(\mu+\pa)\ii$ is given by the formally convergent series of
operators $\sum_{p\ge0} \mu^{-p-1} (-\pa)^p$
(and its Borel counterpart is just division by $\mu-\ze$).
For $\mu=0$, we \emph{define}~$\pa\ii$ as the unique operator
$\pa\ii \col z^{-2}\C[[z\ii]] \to \zcz$ 
such that $\pa\circ\pa\ii$ on $z^{-2}\C[[z\ii]]$
(its Borel counterpart is division by $-\ze$).
\end{nota}


\begin{lemma}   \label{lemAlDerivInvLinDiff}
Let $\Om$ be any non-empty closed discrete subset of~$\C$.
Let $\ti\ph \in \ti\gR_\Om\simp$ and $\mu\in\Om$.
If $\mu=0$ we assume $\ti\ph\in z^{-2}\C[[z\ii]]$;
if $\mu\neq0$ we assume $\De_\mu\ti\ph\in z^{-2}\C[[z\ii]]$.
Then $(\mu+\pa)\ii\ti\ph \in \ti\gR_\Om\simp$ and
\[
\om \in \Om\setminus\{0,\mu\}
\ens\Longrightarrow\ens
\De_\om (\mu+\pa)\ii\ti\ph = (\mu-\om+\pa)\ii \De_\om\ti\ph,
\]
while, if $\mu\neq0$, there exists $C\in\C$ such that
\[
\De_\mu (\mu+\pa)\ii\ti\ph = C + \pa\ii \De_\mu\ti\ph.
\]
\end{lemma}


\begin{lemma}   \label{lemNonAutSubst}
Let $B(z,y) \in \C\{z\ii,y\}$.
Suppose that $\Om$ is stable under addition and
$\ti\ph(z) \in \ti\gR_\Om\simp$ has no constant term.
Then $B\big( z, \ti\ph(z) \big) \in \ti\gR_\Om\simp$ and, for every $\om\in\Om\setminus\{0\}$,
\[
\De_\om B\big( z, \ti\ph(z) \big) = \pa_y B\big( z, \ti\ph(z) \big)
\cdot \De_\om\ti\ph.
\]
\end{lemma}


\noindent
The proofs of Lemmas~\ref{lemAlDerivInvLinDiff}
and~\ref{lemNonAutSubst} are left to the reader.


Let us come back to the solution~$\ti\phi_0$ of~\eqref{eqEulerNL}.
For $\om\in\Z^*_-$, we derive a differential equation for $\ti\psi
= \De_\om\ti\phi_0$ by writing on the one hand 
$\De_\om \pa_z\ti\phi_0 = \pa_z\ti\psi - \om\ti\psi$ (by~\eqref{eqcommutDeompa})
and, on the other hand,
$\De_\om\big( B(z,\ti\phi_0) \big) = \pa_y B(z,\ti\phi_0) \cdot \ti\psi$
by Lemma~\eqref{lemNonAutSubst}, thus 
alien differentiating Equation~\eqref{eqEulerNL} yields
\beglabel{eqfirstderivNL}
\frac{\dd\ti\psi}{\dd z} = \big(\om + \pa_y B(z,\ti\phi_0)\big) \cdot \ti\psi.
\elabel
Since $\om + \pa_y B(z,\ti\phi_0) = \om+1 +O(z^{-2})$, it is immediate that
the only solution of this equation in $\zcz$ is~$0$ when
$\om\neq-1$. This proves 
\[
\om\neq-1
\ens\Longrightarrow\ens
\De_\om\ti\phi_0=0.
\]
For $\om=-1$, Equation~\eqref{eqfirstderivNL} reads
\beglabel{eqLinHomgPhiun}
\frac{\dd\ti\psi}{\dd z} = \ti\be_1 \ti\psi
\elabel
with $\ti\be_1(z) \defeq -1 + \pa_y B(z,\ti\phi_0(z)) \in \ti\gR_{\Z^*_-}\simp$ 
(still by Lemma~\ref{lemNonAutSubst}).
Since $\ti\be_1(z) =O(z^{-2})$, Lemma~\ref{lemAlDerivInvLinDiff} implies
$\ti\al \defeq \pa\ii\ti\be_1 \in \ti\gR_{\Z_-}\simp$
(beware that we must replace $\Z^*_-$ with $\Z_- = \{0\} \cup \Z^*_-$
because a priori only the principal branch of
$\hat\al \defeq - \frac{1}{\ze} \hat\be_1(\ze)$ is regular at~$0$).
Then
\[
\ti\phi_1 \defeq \ee^{\pa\ii\ti\be_1} = 1 + O(z\ii) \in \ti\gR_{\Z_-}\simp
\]
is a non-trivial solution of~\eqref{eqLinHomgPhiun}.
This implies that
\[
\De_{-1}\ti\phi_0 = C \ti\phi_1,
\]
with a certain $C\in\C$.

%
\parag   \label{paragCalculaGoOn}
We go on with the computation of the alien derivatives of~$\ti\phi_1$.
Let \[ \ti\be_2(z) \defeq \pa_y^2 B\big(z,\ti\phi_0(z)\big) \in
\ti\gR_{\Z^*_-}\simp,\]
so that $\De_{-1}\ti\be_1 = C \ti\be_2 \ti\phi_1(z)$ and $\De_\om\ti\be_1 = 0$ for
$\om\neq-1$ (by Lemma~\ref{lemNonAutSubst}).
Computing $\De_\om(\pa\ii\ti\be_1)$ by Lemma~\ref{lemAlDerivInvLinDiff}
and then $\De_\om\ti\phi_1$ by~\eqref{eqDeomSubst}, we get
\beglabel{eqdefphideux}
\De_{-1}\ti\phi_1 = 2 C\ti\phi_2, \qquad
\ti\phi_2 \defeq \demi \ti\phi_1 \cdot (1+\pa)\ii(\ti\be_2\ti\phi_1)
\in \ti\gR_{\Z_-\cup\{1\}}\simp
\elabel
and $\De_\om\ti\phi_1=0$ for $\om\neq-1$.

By the same kind of computation, we get at the next step
$\De_\om\ti\phi_2=0$ for $\om\notin\{-1,1\}$,
\[
\De_{-1}\ti\phi_2 = 3 C\ti\phi_3, \qquad
\ti\phi_3 \defeq \frac{1}{3} \ti\phi_2 \cdot (1+\pa)\ii(\ti\be_2\ti\phi_1)
+ \frac{1}{6}\ti\phi_1 \cdot (2+\pa)\ii(\ti\be_3\ti\phi_1^2 + 2 \ti\be_2\ti\phi_2)
\in \ti\gR_{\Z_-\cup\{1,2\}}\simp
\]
with $\ti\be_3 \defeq \pa_y^3 B\big( z, \ti\phi_0(z) \big)$.
A new undetermined constant appears for $\om=1$:
Lemma~\ref{lemAlDerivInvLinDiff} yields a $C'\in\C$ such that
$\De_1 (1+\pa)\ii (\ti\be_2\ti\phi_1) = C' + \pa\ii
\De_1(\ti\be_2\ti\phi_1) = C'$, hence~\eqref{eqdefphideux} implies
\[
\De_1\ti\phi_2 = C'\ti\phi_3.
\]


We see that Equation~\eqref{eqEulerNL} generates not only the formal
solution~$\ti\phi_0$ but also a sequence of resurgent series
$(\ti\phi_n)_{n\ge1}$,
in which~$\ti\phi_1$ was constructed as the unique solution of the
linear homogeneous differential equation~\eqref{eqLinHomgPhiun} whose
constant term is~$1$;
the other series in the sequence can be characterized by linear
non-homogeneous equations:
alien differentiating~\eqref{eqLinHomgPhiun}, we get
$(1+\pa)\De_{-1}\ti\psi = \De_{-1}\pa\ti\psi = \De_{-1}(\ti\be_1\ti\psi)
= \ti\be_1 \De_{-1}\ti\psi + C\ti\be_2\ti\phi_1\ti\psi$,
thus $\pa(\De_{-1}\ti\phi_1) = (-1+\ti\be_1)\De_{-1}\ti\phi_1 + C\ti\be_2\ti\phi_1^2$,
and it is not a surprise that~$\ti\phi_2$ is the unique formal
solution of
\beglabel{eqLinNonHomgPhideux} 
\pa\ti\phi_2 = (-1+\ti\be_1)\ti\phi_2 + \demi\ti\be_2\ti\phi_1^2.
\elabel
Similarly, $\ti\phi_3$ is the unique formal solution of
\beglabel{eqLinNonHomgPhitrois} 
\pa\ti\phi_3 = (-2+\ti\be_1)\ti\phi_3 + \ti\be_2\ti\phi_1\ti\phi_2
+\frac{1}{6}\ti\be_3\ti\phi_1^3.
\elabel
%


\parag
The previous calculations can be put into perspective with the notion
of \emph{formal integral}, \ie a formal object which solves
Equation~\eqref{eqEulerNL} and is more general than a formal series
like~$\ti\phi_0$.
Indeed, both sides of~\eqref{eqEulerNL} can be evaluated on an
expression of the form
\beglabel{eqdefFormIntForm}
\ti\phi(z,u) = \sum_{n\in\N} u^n \,\ee^{n z} \ti\phi_n(z)
= \ti\phi_0(z) + u\,\ee^z \ti\phi_1(z)  + u^2\,\ee^{2z} \ti\phi_1(z) + \ldots
\elabel
if $(\ti\phi_n)_{n\in\N}$ is any sequence of formal series such
that~$\ti\phi_0$ has no constant term,
simply by treating~$\ti\phi(z,u)$ as a formal series in~$u$ whose
coefficients are transseries of a particular form and writing
\begin{align*}
\frac{\pa\ti\phi}{\pa z}(z,u) &= \sum_{n\in\N} u^n\,\ee^{nz} (n+\pa)\ti\phi_n \\[1ex]
B\big( z, \ti\phi(z,u) \big) &= B\big( z, \ti\phi_0(z) \big) +
\sum_{r\ge1} \frac{1}{r!}\pa_y^r B\big( z, \ti\phi_0(z) \big) 
\sum_{n_1,\ldots,n_r\ge1} u^{n_1+\cdots+n_r} \ee^{ (n_1+\cdots+n_r) z}
\ti\phi_{n_1} \cdots \ti\phi_{n_r}.
\end{align*}
This is equivalent to setting 
$\ti Y(z,y) = \sum_{n\in\N} y^n \ti\phi_n(z)$, 
so that $\ti\phi(z,u) = \ti Y(z,u\,\ee^z)$,
and to considering the equation
\beglabel{eqEulNLY}
\pa_z \ti Y + y \pa_y \ti Y = B\big( z, \ti Y(z,y) \big).
\elabel
for an unknown double series $\ti Y\in\C[[z\ii,y]]$ without constant term.

For an expression~\eqref{eqdefFormIntForm}, Equation~\eqref{eqEulerNL}
is thus equivalent to the sequence of equations
\begin{align}
\tag{$E_0$} & \qquad\qquad \pa\ti\phi_0 = B(z,\ti\phi_0) \\[1ex]
\tag{$E_1$} & (1+\pa)\ti\phi_1 - \pa_y B(z,\ti\phi_0)\cdot\ti\phi_1 = 0 \\[1ex]
\tag{$E_n$} & (n+\pa)\ti\phi_n - \pa_y B(z,\ti\phi_0)\cdot\ti\phi_n = 
\sum_{r\ge2} \frac{1}{r!}\pa_y^r B(z,\ti\phi_0) 
\sum_{ \substack{n_1,\ldots,n_r\ge1 \\ n_1+\cdots+n_r = n} }
\ti\phi_{n_1} \cdots \ti\phi_{n_r} \quad \text{for $n\ge2$.}
\end{align}
Of course~($E_0$) is identical to Equation~\eqref{eqEulerNL} for a
formal series without constant term. 
The reader may check that Equation~($E_1$) coincides
with~\eqref{eqLinHomgPhiun}, 
($E_2$) with~\eqref{eqLinNonHomgPhideux} 
and~($E_3$) with~\eqref{eqLinNonHomgPhitrois}.


\begin{thm}   \label{thmFIsumma}
Equation~\eqref{eqEulerNL} admits a unique solution of the
form~\eqref{eqdefFormIntForm} for which the constant term
of~$\ti\phi_0$ is~$0$ and the constant term
of~$\ti\phi_1$ is~$1$, called ``Formal Integral''.
The coefficients~$\ti\phi_n$ of the formal integral are $1$-summable
in the directions of $(-\pi,0)$ and $(0,\pi)$, and
\beglabel{eqPhinResurZn}
\ti\phi_n(z) \in \ti\gR_{\Z^*_- \cup \{0,1,\ldots,n-1\}}\simp,
\qquad n\in\N.
\elabel

The dependence on~$n$ in the exponential bounds for the Borel
transforms~$\hat\phi_n$ is controlled well enough to ensure 
the existence of locally bounded functions
$\ga$ and~$R >0$ on $(-\pi,0) \cup (0,\pi)$
such that, for $I = (-\pi,0)$ or $(0,\pi)$,
\[ Y^I(z,y) \defeq \sum_{n\in\N} y^n \gS^I\ti\phi_n(z) \]
is holomorphic in 
\[
\gD(I,\ga,R) \defeq \big\{\,
(z,y) \in \C\times\C \mid
\exists \th \in I \;\text{such that}\; 
\RE(z\,\eith) > \ga(\th)
\;\text{and}\;
\abs{y} < R(\th)
\,\}.
\]
Correspondingly, the function
\[ \phi^I(z,u) \defeq \sum_{n\in\N} (u\,\ee^z)^n \gS^I\ti\phi_n(z) \]
is holomorphic in $\{\, (z,u) \in \gD(I,\ga) \times \C \mid
(z, u\,\ee^z) \in \gD(I,\ga,R) \,\}$.

The Borel sums~${\phi^{(-\pi,0)}}_{|u=0}$ and~${\phi^{(0,\pi)}}_{|u=0}$
both coincide with the particular solution of Equation~\eqref{eqEulerNL}
mentioned in Theorem~\ref{thmSimpResPhiZero}.
For $I = (-\pi,0)$ or $(0,\pi)$ and for each $u\in\C^*$, the function $\phi^I(\,.\,,u)$ is a
solution of~\eqref{eqEulerNL} holomorphic in $\{\, z\in\gD(I,\ga) \mid
\RE z < \ln\frac{R}{\abs{u}} \,\}$.
\end{thm}


The reader is once more referred to \cite{mouldSN} for the proof. 

Observe that when we see the formal integral~$\ti\phi(z,u)$ as a
solution of~\eqref{eqEulerNL}, we must think of~$u$ as of an
\emph{indeterminate}, the same way as~$z$ (or rather~$z\ii$) is an
\emph{indeterminate} when we manipulate ordinary formal series;
after Borel-Laplace summation of each~$\ti\phi_n$, we get holomorphic
functions of the \emph{variable} $z\in\gD(I,\ga)$, coefficients of a formal expression
$\sum u^n \ee^{n z}\gS^I\ti\phi_n(z)$;
Theorem~\ref{thmFIsumma} says that, for each $z\in\gD(I,\ga)$, this
expression is a convergent formal series, Taylor expansion of the
function obtained by substituting the indeterminate~$u$ with a
\emph{variable} $u \in \D_{R'\,\ee^{-\RE z}}$ (with $R'>0$ small
enough depending on~$z$).

If we think of~$z$ as of the main variable, the interpretation of the
indeterminate/variable~$u$ is that of a free parameter in the solution
of a first-order differential equation:
$\ti\phi(z,u)$ appears as a formal $1$-parameter family of formal
solutions,
$\phi^{(-\pi,0)}$ and~$\phi^{(0,\pi)}$ as two $1$-parameter families
of analytic solutions.


As for the Borel sum $Y^I(z,y)$, it is an analytic solution of
Equation~\eqref{eqEulNLY} in its domain $\gD(I,\ga,R)$;
this means that the vector field\footnote{
If we change the variable~$z$ into $x\defeq -z\ii$, the vector
field~$X_B$ becomes $x^2\frac{\pa\,}{\pa x} + B(z,Y)
\frac{\pa\,}{\pa Y}$, which has a saddle-node singularity at $(0,0)$.
} $X_B\defeq \frac{\pa\,}{\pa z} + B(z,Y)
\frac{\pa\,}{\pa Y}$ is the direct image of $N\defeq \frac{\pa\,}{\pa z} + y
\frac{\pa\,}{\pa y}$ by the diffeormophism 
\[ \Th^I \col (z,y) \mapsto (z,Y) = \big(z,Y^I(z,y)\big). \]
We may consider~$N$ as a normal form for~$X_B$
and~$\Th^{(-\pi,0)}$ and~$\Th^{(0,\pi)}$ as two sectorial
normalizations. 


The results of the alien calculations of Sections
\ref{paragCalculaEnrich}--\ref{paragCalculaGoOn} are contained in
following statement (extracted from Section~10 of \cite{mouldSN}):


\begin{thm}   \label{thmBridgeSN} 
There are uniquely determined complex numbers $C_{-1},C_1, C_2,\ldots$ such that, for each $n\in\N$,
\begin{align}
\label{eqresurmoinsdeux}
& \De_{-m} \ti\phi_n = 0 \quad \text{for $m\ge 2$,} \\[1ex]
\label{eqresurmoinsun}
& \De_{-1} \ti\phi_n = (n+1) C_{-1} \ti\phi_{n+1}, \\[1ex]
& \De_m \ti\phi_n = (n-m) C_m \ti\phi_{n-m} 
\quad \text{for $1\le m\le n-1$.}
\label{eqresurmun}
\end{align}
Equivalently, letting act the alien derivation~$\De_\om$ on an
expression like~$\ti\phi(z,u)$ or~$\ti Y(z,y)$ by declaring that it
commutes with multiplication by~$u$, $\ee^z$ or~$y$, on has
\beglabel{eqBridgeEulerNL}
\De_m\ti\phi = C_m u^{m+1} \,\ee^{m z} \frac{\pa\ti\phi}{\pa u}
\quad\text{or}\quad
\De_m\ti Y = C_m y^{m+1} \frac{\pa\ti Y}{\pa y},
\qquad \text{for $m=-1$ or $m\ge2$.}
\elabel
\end{thm}


Equation~\eqref{eqBridgeEulerNL} (either for~$\ti\phi$ or for~$\ti Y$)
was baptized ``Bridge Equation'' by \'Ecalle, in view of the bridge it
establishes between ordinary differential calculus (involving $\pa_u$
or~$\pa_y$) and alien calculus (when dealing with the solution of an
analytic equation like~$\ti\phi$ or~$\ti Y$).


\begin{proof}[Proof of Theorem~\ref{thmBridgeSN}]
Differentiating~\eqref{eqEulNLY} \wrt~$y$, we get
\[
(\pa_z + y\pa_y) \pa_y\ti Y = \big( -1 + \pa_y B(z,\ti Y) \big) \pa_y\ti Y.
\]
Alien differentiating~\eqref{eqEulNLY}, we get (in view of~\eqref {eqcommutDeompa})
\[
(\pa_z + y\pa_y) \De_m\ti Y = \big( m + \pa_y B(z,\ti Y) \big)  \De_m\ti Y.
\]
Now $\pa_y\ti Y = 1 +O(z\ii,y)$ is invertible and
we can consider 
$\ti\chi \defeq (\pa_y\ti Y)\ii \De_m\ti Y \in \C[[z\ii,y]]$,
for which we get $(\pa_z + y\pa_y)\ti\chi = (m+1)\ti\chi$,
and this implies the existence of a unique $C_m\in\C$ such that $\ti\chi = C_m
y^{m+1}$.
This yields the second part of~\eqref{eqBridgeEulerNL}, from which the
first part follows, and also
\eqref{eqresurmoinsdeux}--\eqref{eqresurmun} by expanding the formula.
\end{proof}


\parag
The Stokes phenomenon for $\ti\phi(z,u)$ takes the form of two
connection formulas, one for $\RE z<0$, the other for $\RE z>0$,
between the two families of solutions $\phi^{(-\pi,0)}$
and~$\phi^{(0,\pi)}$.
For $\RE z < 0$, it is obtained by analyzing the action of~$\DD+{\R^-}$, the
symbolic Stokes automorphism for the direction~$\R^-$.

Let $\Om \defeq \Z^*_-$.
Since $\ti\phi_n \in \ti\gR_{n+\Om}\simp$ (by~\eqref{eqPhinResurZn}),
the formal integral~$\ti\phi$ can be considered as an $\Om$-resurgent
symbol with support in~$\R^-$ at the price of a slight extension of
the definition: we must allow our resurgent symbols to depend on the
indeterminate~$u$, so we replace~\eqref{eqdeftiEomOm} with
\[
\ti E(\Om,d) \defeq \bigg\{\,
\sum_{\om\in (\Om\cup\{0\}) \cap d} \ee^{-\om z} \ti\ph_\om(z,u) \mid
\ti\ph_\om(z,u) \in \ti\gR_{-\om+\Om}[u]
\,\bigg\}
\]
(thus restricting ourselves to a polynomial dependence on~$u$ for each
homogeneous component).
Then $\ti\phi(z,u) = \sum_{n\in\N} u^n \,\ee^{n z} \ti\phi_n(z) \in \ti E(\Om,\R^-)$.
According to \eqref{eqresurmoinsdeux}, only
one homogeneous component of~$\DD{}{\R^-}$ needs to be taken into
account, and~\eqref{eqdottedalienop} yields
$\DD{}{\R^-} \ti\phi(z,u) = \ee^z\De_{-1}\ti\phi(z,u)$,
whence, by~\eqref{eqresurmoinsun}, 
\[
\DD{}{\R^-} \ti\phi(z,u) = 
\sum_{n\ge0} (n+1) C_{-1} u^n \ee^{(n+1)z} \ti\phi_{n+1}(z)
= C_{-1} \frac{\pa\ti\phi}{\pa u}(z,u).
\]
It follows that
\[
\DD+{\R^-} \ti\phi(z,u) = 
\ti\phi(z,u+C_{-1}) =
\sum_{n\ge0} (u+C_{-1})^n \ee^{n z} \ti\phi_n(z)
\]
and one ends up with
\begin{thm} 
For $z \in \gD\big( (-\pi,0), \ga \big) \cap \gD\big(
(0,\pi), \ga \big)$ with $\RE z < 0$,
\[
\phi^{(0,\pi)}(z,u) \equiv \phi^{(-\pi,0)}(z,u+C_{-1}),
\qquad
Y^{(0,\pi)}(z,y) \equiv Y^{(-\pi,0)}(z,y+C_{-1}\,\ee^z).
\]
\end{thm}


\parag
For $\RE z >0$, we need to inquire about the action of~$\DD+{\R^+}$,
however the action of this operator is not defined on the space of
resurgent symbols with support in~$\R^-$.
Luckily, we can view~$\ti\phi(z,u)$ as a member of the space 
$\ti F(\Z,\R^-) = \ti F_0 \supset \ti F_1 \supset \ti F_2 \supset \cdots$, where
\begin{multline*}
\ti F_p \defeq \bigg\{\,
\sum_{n\in \N} u^{n+p} \ee^{n z} \ti\ph_n(z,u) \mid
\ti\ph_n(z,u) \in \ti\gR_\Z[[u]] \;\text{and}\\[-1ex]
\De_{m_r}\cdots\De_{m_1}\ti\ph_n = 0  \;\text{for $m_1,\ldots,m_r\ge1$
with $m_1+\cdots+m_r > n$}
\,\bigg\}
\end{multline*}
for each $p \in \N$.
One can check that the operator $\DD{}{\R^+} = \sum_{m\ge1} \ee^{-m z}\De_m$ is well
defined on $\ti F(\Z,\R^+)$ and maps~$\ti F_p$ in~$\ti F_{p+1}$,
with 
\[
\DD{}{\R^+} \Big(\sum_{n\ge0} u^{n+p} \ee^{n z} \ti\ph_n(z,u) \Big)
= \sum_{n\ge0} u^{n+p+1} \ee^{n z} \ti\psi_n(z,u),
\qquad
\ti\psi_n(z,u) \defeq \sum_{m\ge1} u^{m-1}\De_m\ti\ph_{m+n}(z,u),
\]
therefore its exponential is well defined and coincides with~$\DD+{\R^+}$.

In the case of the formal integral~$\ti\phi(z,u)$, thanks to~\eqref{eqresurmun}, we find
\[
\DD{}{\R^+} \ti\phi(z,u) = 
\sum_{n\ge0,\, m\ge1} n C_m u^{n+m} \ee^{n z} \ti\phi_n =
\gC \ti\phi(z,u)
\]
with a new operator
$\gC \defeq \dst\sum_{m\ge1} C_m u^{m+1}\dfrac{\pa\,}{\pa u}$.

One can check that $\ti F(\Z,\R^-)$ is an algebra and its
multiplication maps $\ti F_p \times \ti F_q$ to~$\ti F_{p+q}$.
Since~$\gC$ is a derivation which maps~$\ti F_p$ to~$\ti F_{p+1}$, its
exponential~$\exp\gC$ is well defined and is an automorphism (same argument as
for Lemma~\ref{lemPureAlg}).
Reasoning as in Exercise~\ref{exoAllAlgEndo}, one can see that there
exists $\xi(u) \in u\C[[u]]$ such that $\exp\gC$ coincides with the
composition operator associated with
$(z,u) \mapsto \big( z, \xi(u) \big)$:
\[
\ti\ph(z,u) \in \ti F(\Z,\R^-) 
\ens\Longrightarrow\ens
(\exp\gC)\ti\ph(z,u) = \ti\ph\big( z, \xi(u) \big).
\]
In fact, there is an explicit formula
\[
\xi(u) = u + \sum_{m\ge1} \bigg(
\sum_{r\ge1} \sum_{ \substack{m_1,\ldots,m_r \ge 1 \\ m_1+\cdots+m_r = m} }
\frac{1}{r!} \be_{m_1,\ldots,m_r} C_{m_1} \cdots C_{m_r}
\bigg) u^{m+1}
\]
with the notations $\be_{m_1}=1$ and
$\be_{m_1,\ldots,m_r} = (m_1+1)(m_1+m_2+1)\cdots(m_1+\cdots+m_{r-1}+1)$.
We thus obtain
\beglabel{eqDDpRptiphi}
\DD+{\R^+}\ti\phi(z,u) = \ti\phi\big( z, \xi(u) \big).
\elabel

\begin{thm}
The series $\xi(u)$ has positive radius of convergence and,
for $z \in \gD\big( (-\pi,0), \ga \big) \cap \gD\big(
(0,\pi), \ga \big)$ with $\RE z > 0$,
\[
\phi^{(-\pi,0)}(z,u) \equiv \phi^{(0,\pi)}\big(z,\xi(u)\big),
\qquad
Y^{(-\pi,0)}(z,y) \equiv Y^{(0,\pi)}\big(z,\xi(y\,\ee^{-z})\,\ee^z\big).
\]
\end{thm}

\begin{proof}[Sketch of proof]
Let $I \defeq [\eps,\pi-\eps]$, $J \defeq [-\pi+\eps,-\eps]$, and
consider the diffeomorphism
$
\th \defeq \big[ \Th^{(0,\pi)} \big]\ii \circ \Th^{(-\pi,0)}
$ in $\{\, z \in \gD(I, \ga) \cap \gD(J, \ga) \mid \RE z > 0 \,\}
\times \{ \abs{y} < R'\}$ with $R'>0$ small enough.
It is of the form $\th(z,y)= \big( z, \chi^+(z,y) \big)$ with $\chi^+(z,0)\equiv0$.
The direct image of $N = \frac{\pa\,}{\pa z} + y \frac{\pa\,}{\pa y}$
by~$\th$ is~$N$,
this implies that 
$\chi^+ = N\chi^+$,
whence $\frac{1}{u\,\ee^z}\chi^+(z,u\,\ee^z)$ is independent of~$z$
and can be written $\frac{\xi^+(u)}{u}$ with $\xi^+(u) \in \C\{u\}$.
Thus $\chi^+(z,y) = \xi^+(y\,\ee^{-z})\, \ee^z$, \ie
\[
Y^{J}(z,y) \equiv
Y^{I}\big(z,\xi^+(y\,\ee^{-z})\,\ee^z\big).
\]
To conclude, it is thus sufficient to prove that the Taylor series
of~$\xi^+(u)$ is~$\xi(u)$.
This can be done using~\eqref{eqDDpRptiphi}, by arguing as in the proof of Theorem~\ref{thmSymbStokesLapl}.
\end{proof}


\begin{exo}
[\textbf{Analytic invariants}]
Assume we are given two equations of the form~\eqref{eqEulerNL}
and, correspondingly, two vector fields
$X_{B_1} =  \frac{\pa\,}{\pa z} + B_1(z,Y)\frac{\pa\,}{\pa Y}$
and $X_{B_2} =  \frac{\pa\,}{\pa z} + B_2(z,Y)\frac{\pa\,}{\pa Y}$
with the same assumptions as previously on $B_1,B_2 \in \C\{z\ii,y\}$.
Prove that there exists a formal series $\ti\chi(z,y) \in
\C[[z\ii,y]]$ such that the formula
$\th(z,y) \defeq \big( z, \ti\chi(z,y) \big)$ defines a formal
diffeomorphism which conjugates~$X_{B_1}$ and~$X_{B_2}$.
Prove that~$X_{B_1}$ and~$X_{B_2}$ are analytically conjugate, \ie
$\ti\chi(z,y) \in \C\{z\ii,y\}$, if and only both equations
give rise to the same sequence $(C_{-1},C_1,C_2,\ldots)$,
or, equivalently, to the same pair $\big( C_{-1},\xi(u) \big)$
(the latter pair is called the ``Martinet-Ramis modulus'').
\end{exo}


\begin{exo}
Study the particular case where~$B$ is of the form
$B(z,y) = b_0(z) + \big( 1+b_1(z) \big) y$,
with $b_0\in z\ii\C\{z\ii\}$, $b_1\in z^{-2}\C\{z\ii\}$.
Prove in particular that the Borel transform of $b_0\,\ee^{-\pa\ii
  b_1}$ is an entire function whose value at~$-1$ is
$-\frac{1}{2\pi\I}C_{-1}$ and that $C_m=0$ for $m\neq-1$ in that case.
\end{exo}


\begin{rem}
The numbers $C_m$, $m\in\{-1\}\cup\N^*$, which encode such a subtle
analytic information, are usually impossible to compute in closed
form.
An exception is the case of the ``canonical Riccati equations'', for
which
$B(z,y) = y - \frac{1}{2\pi\I}(B_-+B_+y^2) z\ii$,
with $B_-,B_+ \in \C$.
One finds $C_m =0$ for $m \notin \{-1,1\}$ and
\[
C_{-1} = B_- \sig(B_-B^+), \quad C_{1} = -B_+ \sig(B_-B^+)
\]
with $\sig(b) \defeq \frac{2}{b^{1/2}} sin \frac{b^{1/2}}{2}$.
See \cite{mouldSN} for the references.
\end{rem}


\newpage


\centerline{\Large\sc The resurgent viewpoint on  holomorphic}
\bigskip

\centerline{\Large\sc tangent-to-identity germs}
\addcontentsline{toc}{part}{\sc The resurgent viewpoint on holomorphic
tangent-to-identity germs} 


\bigskip

\bigskip

The last part of this text is concerned with germs of holomorphic
tangent-to-identity diffeomorphisms.
The main topics are the description of the local dynamics (describing the
local structure of the orbits of the discrete dynamical system induced
by a given germ)
and the description of the conjugacy classes (attaching to a given
germ quantities which characterize its analytic conjugacy class).
We shall give a fairly complete account of the results in the simplest
case, limiting ourselves to germs at~$\infty$ of the form
\beglabel{eqSimplestFormalClass}
f(z) = z + 1 + O(z^{-2})
\elabel
(corresponding to germs at~$0$ of the form $F(t) = t - t^2 + t^3 +
O(t^4)$ by \eqref{eqdefFtsigtau}--\eqref{eqdeffzFt}).
The reader is referred to 
\cite[Vol.~2]{Eca81}, \cite{Milnor}, \cite{Loray}, \cite{kokyu}, \cite{DSun}, \cite{DSdeux}
for more general studies.

It turns out that formal tangent-to-identity diffeomorphisms play a
prominent role, particularly those which are $1$-summable and
$2\pi\I\Z$-resurgent. So the ground was prepared in Sections
\ref{secGermsholdiffeos}--\ref{secgrponesummdiffeos}
and in Theorem~\ref{thmresdiffeo}.
In fact, because of the restriction~\eqref{eqSimplestFormalClass}, all the
resurgent functions which will appear will be simple;
we thus begin with a preliminary section.

\section{Simple $\Om$-resurgent tangent-to-identity diffeomorphisms}
\label{sec:GermsDeb}

Let us give ourselves a non-empty closed discrete subset~$\Om$ of~$\C$
which is stable under addition.
Recall that, according to Section~\ref{secNLopRES}, 
$\Om$-resurgent tangent-to-identity diffeomorphisms form a group
$\gGR(\Om)$ for composition (subgroup of the group $\ti\gG = \id +
\C[[z\ii]]$ of all formal tangent-to-identity diffeomorphisms
at~$\infty$).


\begin{Def}   \label{DefSimpOmresDiffeo}
We call \emph{simple $\Om$-resurgent tangent-to-identity diffeomorphism} any $\ti f = \id + \ti\ph
\in \gGR$ where $\ti\ph$ is a simple $\Om$-resurgent series.
We use the notations
\[
\gGS(\Om) \defeq \{\, \ti f = \id + \ti\ph \mid \ti\ph\in\ti\gR_\Om\simp \,\},
\qquad
\gGS_\sig(\Om) \defeq \gGS(\Om) \cap \ti\gG_\sig
\ens\text{for $\sig\in\C$.}
\]
  We define $\De_\om \col \gGS(\Om) \to \ti\gR_{-\om+\Om}\simp$ for any
  $\om\in\Om$ by setting \[\De_\om(\id+\ti\ph) \defeq \De_\om\ti\ph.\]
\end{Def}


Recall that, in Section~\ref{secFormalDiffeos}, $\pa\ti f$ was defined
as the invertible formal series $1+\pa\ti\ph$ for any $\ti f = \id +
\ti\ph \in \ti\gG$.
Clearly $\ti f \in \gGS(\Om) \imp \pa\ti f \in \gGS(\Om)$.


\begin{thm}	\label{thmSimpResdiffeo}
The set~$\gGS(\Om)$ is a subgroup of~$\gGR(\Om)$,
the set~$\gGS_0(\Om)$ is a subgroup of~$\gGR_0(\Om)$.
For any $\ti f, \ti g \in \gGS(\Om)$ and $\om\in\Om$, we have
\begin{gather}
\label{eqDeomComposDiffeo}
\De_\om(\ti g \circ \ti f) = (\pa\ti g)\circ \ti f \cdot \De_\om\ti f
+ \ee^{ -\om(\ti f-\id) }\cdot (\De_\om g)\circ\ti f, \\[1ex]
\label{eqDeomDiffeoInvers}
\ti h = \ti f\ic \Imp
\De_\om \ti h = - \ee^{ -\om(\ti h - \id) } \cdot 
(\De_\om\ti f) \circ \ti h\cdot \pa\ti h.
\end{gather}
\end{thm}

\begin{proof}
The stability under group composition stems from Theorem~\ref{thmComposSimpRes},
since 
$(\id+\ti\psi)\circ(\id+\ti\ph) = \id + \ti\ph + \ti\psi\circ(\id+\ti\ph)$.
The stability under group inversion is proved from Lagrange reversion
formula as in the proof of Theorem~\ref{thmresdiffeo}, adapting the
arguments of the proof of Theorem~\ref{thmComposSimpRes}.

Formula~\eqref{eqDeomComposDiffeo} results from~\eqref{eqDeomCompos},
and formula~\eqref{eqDeomDiffeoInvers} follows by choosing $g=f\ic$.
\end{proof}


\section{Simple parabolic germs with vanishing resiter} \label{sec:Abel}


We now come to the heart of the matter, giving ourselves a germ $F(t) \in\C\{t\}$ of
holomorphic tangent-to-identity diffeomorphism at~$0$ and the
corresponding germ $f(z) \defeq 1/F(1/z) \in \gG$ at~$\infty$.

The germ~$F$ gives rise to a discrete dynamical system $F \col U \to
\C$, where~$U$ is an open neighbourhood of~$0$ on which a
representative of~$F$ is holomorphic.
This means that for any $t_0 \in U$ we can define a finite or infinite \emph{forward
  orbit} $\{\, t_n = F^{\circ n}(t_0) \mid 0 \le n< N \,\}$, where $N\in
\N^*\cup\{\infty\}$ is characterized by $t_1 = F(t_0)\in U$, \ldots,
$t_{N-1} = F(t_{N-2}) \in U$ and
$t_N= F(t_{N-1})\notin U$ (so that apriori $t_{N+1}$ cannot be defined), 
and similarly a finite or infinite \emph{backward orbit} $\{\, t_{-n}
= F^{\circ(-n)}(t_0) \mid 0 \le n< M \,\}$ with $M\in\N^*\cup\{\infty\}$.

We are interested in the local structure of the orbits starting close
to~$0$, so the domain~$U$ does not matter.
Moreover, the qualitative study of a such a dynamical system is
insensitive to analytic changes of coordinate:
we say that~$G$ is \emph{analytically conjugate} to~$F$ if there exists an
invertible $H\in t\C\{t\}$ such that $G = H\ic \circ F \circ H$;
the germ~$G$ is then itself tangent-to-identity and it should be
considered as equivalent to~$F$ from the dynamical point of view
(because~$H$ maps the orbits of~$F$ to those of~$G$).
The description of the analytic conjugacy classes is thus dynamically relevant.


\emph{We suppose that~$F$ is non-degenerate in the sense that $F''(0)\neq0$.}
Observe that $G = H\ic \circ F \circ H \imp G''(0) = H'(0) F''(0)$,
thus we can rescale the variable~$w$ so as to make the second
derivative equal to~$-2$, \ie we assume from now on
$F(t) = t - t^2 + (\rho+1) t^3 + O(t^4)$ with a certain $\rho\in \C$,
and correspondingly
\beglabel{eqformfrho}
f(z) = z + 1 - \rho z\ii + O(z^{-2}) \in \gG_1.
\elabel
Such a germ~$F$ or~$f$ is called a \emph{simple parabolic germ}.

Once we have done that, we should only consider tangent-to-identity
changes of coordinate~$G$, so as to maintain the condition
$F''(0)=-2$.
In the variable~$z$, this means that we shall study the
$\gG$-conjugacy class 
$\{\, h\ic \circ f \circ h \mid h\in\gG \,\} \subset \gG_1$.


As already alluded to, the $\ti\gG$-conjugacy class of~$f$ in~$\ti\gG_1$
plays a role in the problem, \ie we must also consider the formal conjugacy equivalence
relation.
The point is that it may happen that two \emph{holomorphic} germs~$f$ and~$g$
are \emph{formally} conjugate (there exists $\ti h \in \ti\gG$ such
that $f \circ \ti h = \ti h \circ g$) without being \emph{analytically}
conjugate (there exists no $h\in \gG$ with the same property):
the $\gG$-conjugacy classes we are interested in form a finer
partition of~$\gG_1$ than the $\ti\gG$-conjugacy classes.

It turns out that the number~$\rho$ in~\eqref{eqformfrho} is invariant
by formal conjugacy and that two germs with the same~$\rho$ are always
formally conjugate (we omit the proof).
This number is called ``resiter''.


\emph{We suppose further that the resiter~$\rho$ is~$0$}, \ie we limit
ourselves to the most elementary formal conjugacy class.
This implies that our~$f$ is of the form~\eqref{eqSimplestFormalClass} and formally conjugate to 
$f_0(z)\defeq z+1$, the most elementary simple parabolic germ with
vanishing resiter,
which may be considered as a formal \emph{normal form} for all simple
parabolic germs with vanishing resiter.
The corresponding normal form at~$0$ is $F_0(t) \defeq \frac{t}{1+t}$.
The orbits of the normal form are easily computed:
we have $f_0^{\circ n} = \id + n$ and $F_0^{\circ n}(t) = \frac{t}{1+n t}$
for all $n\in\Z$,
thus the backward and forward orbits of a point $t_0\neq0$ are infinite
and contained either in~$\R$ (if $t_0\in\R$) or in a circle passing
through~$0$ centred at a point of $\I\R^*$.

In particular, all the forward orbits of~$F_0$ converge to~$0$ and
all its backward orbits converge in negative time to~$0$.
If the formal conjugacy between~$F$ and~$F_0$ happens to be
convergent, then such qualitative properties of the dynamics automatically
hold for the orbits of~$F$ itself (at least for those which start close
enough to~$0$).
We shall see that in general the picture is more complex\ldots


\section{Resurgence and summability of the iterators} 
\label{sec:ResurIter}
\label{sec:SummaIter}


\begin{nota}
Given $\ti g \in \ti\gG$, the operator of composition with~$\ti g$ is
denoted by
\[
C_{\ti g} \col \ti\ph \in \C[[z\ii]] \mapsto \ti\ph \circ \ti g \in \C[[z\ii]].
\]
The operator $C_{\id-1}-\ID$ induces an invertible map $\zcz
\to z^{-2}\C[[z\ii]]$ with Borel counterpart
$\hat\ph(\ze) \in \C[[\ze]] \mapsto (\ee^\ze-1)\hat\ph(\ze) \in \ze\C[[\ze]]$;
we denote by
\[
E \col z^{-2}\C[[z\ii]] \to \zcz, \qquad
\hat E \col \ze\C[[\ze]] \to \C[[\ze]]
\]
its inverse and the Borel counterpart of its inverse,
hence $(\hat E \hat\ph)(\ze) = \frac{1}{\ee^\ze-1} \hat\ph(\ze)$
(\cf Corollary~\ref{coreqdifflin}).
We also set
\beglabel{eqdefFormalMFfzero}
f_0 \defeq \id + 1 \in \gG_1.
\elabel
\end{nota}


The operator~$E$ will allow us to give a very explicit proof of the
existence of a formal conjugacy between a diffeomorphism with
vanishing resiter and the normal form~\eqref{eqdefFormalMFfzero}.


\begin{lemma}   \label{lemtivstsumtiphk}
Given a simple parabolic germ with vanishing resiter $f\in \gG_1$,
there is a unique $\ti v_* \in \ti\gG_0$ such that 
\beglabel{eqconjugtiv} 
\ti v_* \circ f = f_0 \circ \ti v_*.
\elabel
It can be written as a formally convergent series
\beglabel{eqdefvsttiphk}
\ti v_* = \id + 
%
%
\sum_{k\in\N} \ti\ph_k, \qquad
\ti\ph_k \defeq (EB)^k Eb \in z^{-2k-1}\C[[z\ii]]
\ens\text{for each $k\in\N$,}
\elabel
with a holomorphic germ $b \defeq f \circ f_0\ic - \id \in
z^{-2}\C\{z\ii\}$ and an operator $B \defeq C_{\id+b}-\ID$.

The solutions in~$\ti\gG$ of the conjugacy equation $\ti v \circ f =
f_0 \circ \ti v$ are the formal diffeomorphisms $\ti v = \ti v_*+c$
with arbitrary $c\in\C$.
\end{lemma}


\begin{proof}
The conjugacy equation can be written $\ti v\circ f = \ti v+1$ or,
equivalently (composing with $f_0\ic = \id-1$), 
$\ti v\circ(\id+b) = \ti v\circ(\id-1)+1$.
Searching for a formal solution in the form
$\ti v = \id + \ti\ph$ with $\ti\ph \in \C[[z\ii]]$,
we get $b + \ti\ph\circ(\id+b) = \ti\ph\circ(\id-1)$, \ie
\beglabel{eqconjugtiph}
(C_{\id-1}-\ID)\ti\ph = B\ti\ph + b.
\elabel
We have
$\val\big( (C_{\id-1}-\ID)\ti\ph \big) \ge \val(\ti \ph) + 1$
for the standard valuation~\eqref{eqstdval},
and $\val(B\ti \ph) \ge \val(\ti \ph) + 3$
(because $B$ can be written as the formally convergent series if
operators
$\sum_{r\ge1}\frac{1}{r!}b^r\pa^r$
with $\val(b)\ge2$ and $\val(\pa\ti \ph) \ge \val(\ti \ph) + 1$),
thus the difference between any two formal solutions
of~\eqref{eqconjugtiph} is a constant.
If we specify $\ti\ph\in\zcz$, then~\eqref{eqconjugtiph} is equivalent
to
\[
\ti\ph = EB\ti\ph + Eb,
\]
where $\val(EB\ti \ph) \ge \val(\ti \ph) + 2$,
thus the formal series~$\ti\ph_k$ of~\eqref{eqdefvsttiphk} have
valuation at least $2k+1$ and yield the unique formal solution without
constant term in the form $\ti\ph = \sum_{k\in\N} \ti\ph_k$.
\end{proof}


\begin{Def}
  The unique formal diffeomorphism $\ti v_* \in \ti\gG_0$ such that 
  $\ti v_* \circ f = f_0 \circ \ti v_*$
is called the ``iterator'' of~$f$.
Its inverse $\ti u_* \defeq \ti v_*\ic \in \ti\gG_0$ is called the
``inverse iterator'' of~$f$.
\end{Def}

We illustrate this in the following commutative diagram, including the
parabolic germ at~$0$ defined by $F(t) \defeq 1/f(1/t)$:
\[
\xymatrix{
{\rule[-.6ex]{0ex}{0ex} \makebox[1em]{$z$}} \ar@<-.5ex>[d]_*+{\ti u_*} \ar[rr] & 
& {\rule[-.6ex]{0ex}{0ex} \makebox[3em]{$z+1$}} \ar@<-.5ex>[d]_*+{\ti u_*}  \\
{\rule[-1ex]{0ex}{3ex} \makebox[1em]{$z$}} \ar@<-.5ex>[u]_*+{\ti v_*} 
%
\ar[dr];[]^*+{ \raisebox{-.5ex}[0ex][0ex]{$z=1/t$} } \ar[rr] & 
& {\rule[-1ex]{0ex}{3ex} \makebox[2em]{$f(z)$}} \ar[dr];[] \ar@<-.5ex>[u]_*+{\ti v_*} \\
& {t} \ar[rr] & & {F(t)}
}
\]
Observe that
\beglabel{eqDifftiust}
f \circ \ti u_* = \ti u_* \circ f_0,
\elabel
which can be viewed as a difference equation: 
$\ti u_*(z+1) = f\big( \ti u_*(z) \big)$.

\begin{figure}
\begin{center}

\includegraphics[scale=1]{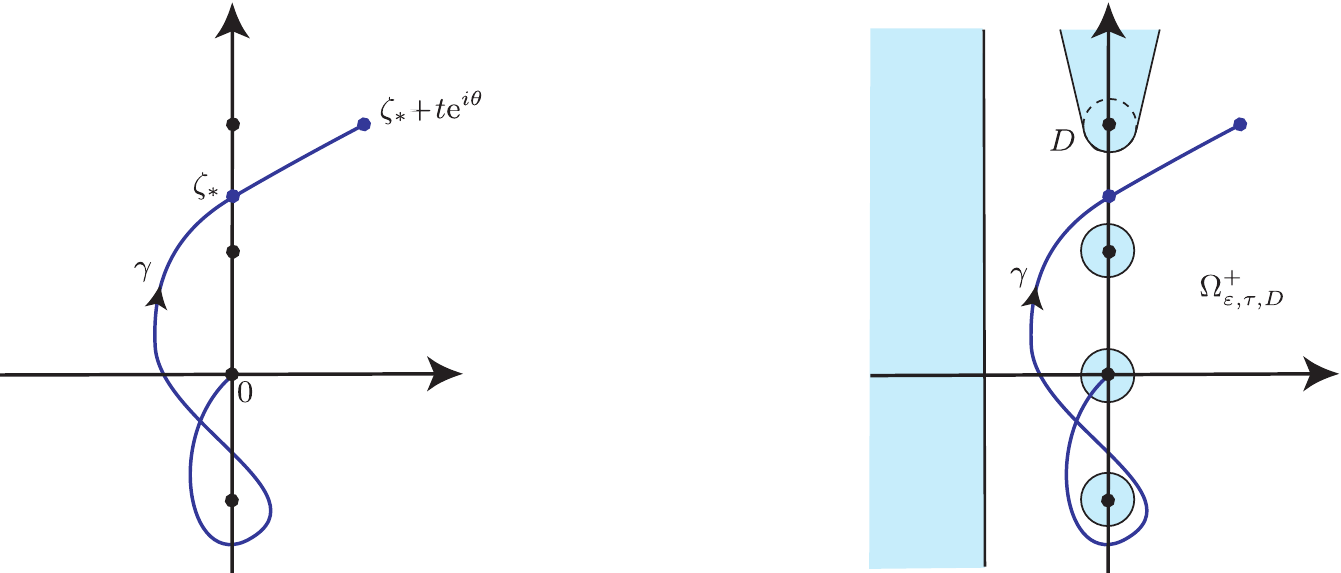}

 \bigskip 

\caption{%
\emph{Resurgence of the iterator (Theorem~\ref{thmIterSimpResurSumma}).}
  Left: A path of analytic continuation for~$\hat\ph$. Right: The
  domain $\Om^+_{\eps,\tau,D}$ of Lemma~\ref{lemcontrolhatphk}.}
\label{figxxaabb}

\end{center}
\end{figure}

\begin{thm}    \label{thmIterSimpResurSumma}
Suppose that $f\in \gG_1$ has vanishing resiter. Then its
iterator~$\ti v_*$ and its inverse interator~$\ti u_*$ belong to 
$\gGS_0(2\pi\I\Z) \cap \ti\gG_0(I^+) \cap \ti\gG_0(I^-)$ 
with $I^+ \defeq (-\frac{\pi}{2}, \frac{\pi}{2})$ 
and $I^- \defeq (\frac{\pi}{2}, \frac{3\pi}{2})$
(notations of Definitions~\ref{DefSummaDiffeo} and~\ref{DefSimpOmresDiffeo}).

Moreover, the iterator can be written $\ti v_* = \id+\ti\ph$ with a simple
$2\pi\I\Z$-resurgent series~$\ti\ph$ whose Borel transform satifies
the following:
for any path~$\ga$ issuing from~$0$ and then avoiding $2\pi\I\Z$ and
ending at a point $\ze_*\in\I\R$,
or for $\ga = \{0\}$ and $\ze_*=0$,
there exist locally bounded functions $\al,\be \col I^+ \cup I^- \to
\R^+$ such that
\beglabel{ineqboundcontgahatph}
\abs*{ \cont_\ga\hat\ph\big( \ze_* + t\,\ee^{\I\th} \big) }
\leq
\al(\th)\,\ee^{\be(\th) t}
\quad \text{for all $t\geq0$ and $\th \in I^+\cup I^-$}
\elabel
(see Figure~\ref{figxxaabb}a).
\end{thm}


Since $\ti\ph \defeq \ti v_* - \id$ is given by
Lemma~\ref{lemtivstsumtiphk} in the form of the formally convergent
series $\sum_{k\ge0} \ti\ph_k$, the statement can be proved by
controlling the formal Borel transforms~$\hat\ph_k$.


\begin{lemma}   \label{lemSimpReshatphk}
For each $k\in\N$ we have $\hat\ph_k \defeq \cB(\ti\ph_k) \in \hat\gR_{2\pi\I\Z}\simp$.
\end{lemma}


\begin{lemma}   \label{lemcontrolhatphk}
Suppose that $0 < \eps < \pi < \tau$, $0 < \ka \leq1$ and $D$ is a closed
disc of radius~$\eps$ centred at $2\pi\I m$ with $m \in \Z^*$,
and let
\beglabel{condtauka}
\Om^+_{\eps,\tau,D} \defeq \{\, \ze\in\C \mid
\RE \ze > -\tau, \;
\dist\big( \ze, 2\pi\I\Z^* \big) > \eps \,\}
\setminus
\{\, u\ze \in\C \mid u\in[1,+\infty), \;
\pm\ze \in D \,\}
\elabel
(see Figure~\ref{figxxaabb}b).
Then there exist $A,M,R>0$ such that, for any naturally parametrised path
$\ga \col [0,\ell] \to \Om^+_{\eps,\tau,D}$ with
%
%
\begin{gather}
\label{condeps}
s\in[0,\eps] \imp \abs{\ga(s)} = s,
\qquad
s>\eps \imp \abs{\ga(s)} > \eps,
\qquad
s\in[0,\ell] \imp
\abs{\ga(s)} > \ka s,\\[1ex]
\shortintertext{one has}
\label{ineqestimphk}
\abs*{ \cont_\ga\hat\ph_k\big( \ga(\ell) \big) }
\leq A\frac{ (M\ell)^k }{k!} \ee^{R\ell}
\quad \text{for every $k\geq0$}.
\end{gather}
%
%
\end{lemma}


\begin{proof}[Lemmas~\ref{lemSimpReshatphk} and~\ref{lemcontrolhatphk} imply Theorem~\ref{thmIterSimpResurSumma}]
According to notation~\eqref{eqdefDDst}, we denote by~$\D_{2\pi}$
or~$\D_1$ the open disc centred at~$0$ of radius $2\pi$ or~$1$.
Lemma~\ref{lemcontrolhatphk} implies that the series of holomorphic functions
$\sum \hat\ph_k$ converges normally in any compact subset of~$\D_{2\pi}$
(using paths~$\ga$ of the form $[0,\ze]$)
and that its sum, which is~$\hat\ph$, extends analytically along any
naturally parametrised path~$\ga$ which starts as the line segment
$[0,1]$ and then stays in $\C\setminus 2\pi\I\Z$:
indeed, taking $\eps,\ka$ small enough and $\tau, m$ large enough,
we see that Lemma~\ref{lemcontrolhatphk} applies to~$\ga$ and the neighbouring paths,
so that \eqref{ineqestimphk} yields the normal convergence of
$\sum_{k\geq0} \cont_\ga\hat\ph_k(\ga(t)+\ze) = \cont_\ga\hat\ph(\ga(t)+\ze)$
for all~$t$ and~$\ze$ with $\abs{\ze}$ small enough.
Therefore $\hat\ph$ is $2\pi\I\Z$-resurgent and,
combining Lemma~\ref{lemSimpReshatphk} with the estimates~\eqref{ineqestimphk},
we also get $\hat\ph \in \hat\gR_{2\pi\I\Z}\simp$ by Lemma~\ref{lemConvSeriesSimpRes}.


This establishes $\ti v_* \in \gGS_0(2\pi\I\Z)$, whence $\ti u_* \in
\gGS_0(2\pi\I\Z)$ by Theorem~\ref{thmSimpResdiffeo}.


For the part of~\eqref{ineqboundcontgahatph} relative to~$I^+$, we
give ourselves an arbitrary $n>1$ and set $\de_n\defeq \frac{\pi}{2n}$,
$I^+_n \defeq [-\frac{\pi}{2}+\de_n, \frac{\pi}{2}-\de_n]$.
Given~$\ga$ with endpoint $\ze_*\in\I\R$, we first
replace an initial portion of~$\ga$ with a line segment of length~$1$
(unless~$\ga$ stays in~$\D_1$, in which case the modification of the
arguments which follow is trivial)
and switch to its natural parametrisation $\ga \col [0,\ell] \to \C$.
We then choose~$\eps_n$ and~$\ka_n$ small enough:
\begin{align*}
\eps_n &< \min\Big\{ 1, \, \min_{[1,\ell]}\abs{\ga}, \,
\dist\big( \ga\big([0,\ell]\big), 2\pi\I\Z^* \big), \,
\dist\big(\ze_*,2\pi\I\Z\big)\cos\de_n \Big\}, \\[1ex]
\ka_n &< \min\Big\{
\min_{[0,\ell]} \tfrac{\abs{\ga(s)}}{s}, \,
\min_{t\geq0} \tfrac{\abs{ \ze_* + t\,\ee^{\pm\I\de_n} }}{\ell+t}
\Big\},
\end{align*}
and $\tau$ and~$m_n$ large enough:
\[
\tau > -\min\RE\ga, \qquad
m_n > \frac{1}{2\pi} \big(\eps_n + \max\abs{\IM\ga} \big),
\]
so that Lemma~\ref{lemcontrolhatphk} applies to the concatenation of paths 
$\Ga \defeq \ga + [\ze_*,\ze_*+t\,\ee^{\I\th}]$
for each $t\geq0$ and $\th\in I^+_n$;
since $\Ga$ has length $\ell+t$,
\eqref{ineqestimphk} yields
\[
\text{$t\ge0$ and $\th\in I^+_n$} \Imp
\abs*{ \cont_\ga\hat\ph\big( \ze_* + t\,\ee^{\I\th} \big) }
= \abs*{ \cont_\Ga\hat\ph\big( \Ga(\ell+t) \big) } \leq
A_n\,\ee^{(M_n+R_n)(\ell+t)},
\]
where $A_n$, $M_n$ and~$R_n$ depend on~$n$ and~$\ga$ but not on~$t$
or~$\th$.
We thus take
\[
\al^+(\th) \defeq \ee^{\ell\be^+(\th)} \max\big\{\, A_n \mid
\text{$n \ge 1$ s.t.\ $\th\in I^+_n$} \,\big\},
\quad
\be^+(\th) \defeq \max\big\{\, M_n+R_n \mid
\text{$n \ge 1$ s.t.\ $\th\in I^+_n$} \,\big\}
\]
for any $\th\in I^+$, and get
\[
\text{$t\ge0$ and $\th\in I^+$} \Imp
\abs*{ \cont_\ga\hat\ph\big( \ze_* + t\,\ee^{\I\th} \big) }
\le \al^+(\th)\,\ee^{\be^+(\th) t}.
\]
%
%
The part of~\eqref{ineqboundcontgahatph} relative to~$I^-$ follows
from the fact that $\hat\ph^-(\ze) \defeq \hat\ph(-\ze)$ satisfies all
the properties we just obtained for $\hat\ph(\ze)$, since
it is the formal Borel transform of
$\ti\ph^-(z) \defeq -\ti\ph(-z)$
which solves the equation
$C_{\id-1} \ti\ph^- = C_{\id+b^-} \ti\ph^- + b^-_*$
associated with the simple parabolic germ $f^-(z) \defeq -f\ii(-z) = z+1+b^-(z+1)$.


This establishes~\eqref{ineqboundcontgahatph}, which yields (in the
particular case $\ga=\{0\}$)
$\ti v_* \in \ti\gG_0(I^+) \cap \ti\gG_0(I^-)$, whence $\ti u_* \in
\ti\gG_0(I^+) \cap \ti\gG_0(I^-)$ by Theorem~\ref{thmSummaDiffeo}.
\end{proof}


\begin{proof}[Proof of Lemma~\ref{lemSimpReshatphk}]
Since $b(z) \in z^{-2}\C\{\ze\}$, its formal Borel transform is an
entire function~$\hat b(\ze)$ vanishing at~$0$, hence
\[
\hat\ph_0(\ze) = \frac{\hat b(\ze)}{\ee^\ze - 1} 
\in \hat\gR_{2\pi\I\Z}\simp
\]
(\cf Lemma~\ref{lemSimplestab}).

We proceed by induction on~$k$ and assume 
$k\ge1$ and $\ti\ph_{k-1} \in \ti\gR_{2\pi\I\Z}\simp$.
By Theorem~\ref{thmComposSimpRes} we get
$C_{\id+b}\ti\ph_{k-1} \in \ti\gR_{2\pi\I\Z}\simp$, 
thus $B\ti\ph_{k-1} \in \ti\gR_{2\pi\I\Z}\simp$, 
thus (since $\cB(B\ti\ph_{k-1})(\ze) \in \ze\C\{\ze\}$)
\[ 
\hat\ph_k(\ze) = \frac{1}{\ee^\ze-1} \cB(B\ti\ph_{k-1})(\ze) 
\in \hat\gR_{2\pi\I\Z},
\] 
but is it true that all the singularities of all the branches of the analytic continuation
of~$\hat\ph_k$ are simple?

By repeated use of~\eqref{eqDeomCompos}, we get
\[
\De_{\om_s} \cdots \De_{\om_1} C_{\id+b}\ti\ph_{k-1} =
\ee^{ -(\om_1+\cdots+\om_s)b } C_{\id+b}
\De_{\om_s} \cdots \De_{\om_1} \ti\ph_{k-1}
\]
for every $s\ge1$ and $\om_1,\ldots,\om_s \in 2\pi\I\Z^*$, hence
\[
\De_{\om_s} \cdots \De_{\om_1} B \ti\ph_{k-1} =
\bB_{\om_1,\ldots,\om_s}
\De_{\om_s} \cdots \De_{\om_1} \ti\ph_{k-1}
\quad\text{with}\ens
\bB_{\om_1,\ldots,\om_s} \defeq \ee^{ -(\om_1+\cdots+\om_s)b } C_{\id+b} - \ID.
\]
Now, for any $\ti\psi \in \C[[z\ii]]$, we have
$\bB_{\om_1,\ldots,\om_s}\ti\psi = 
\ee^{ -(\om_1+\cdots+\om_s)b } B \ti\psi 
+ (\ee^{ -(\om_1+\cdots+\om_s)b } -1)\ti\psi
\in z^{-2}\C[[z\ii]]$,
thus each of the simple $2\pi\I\Z$-resurgent series 
$\De_{\om_s} \cdots \De_{\om_1} B \ti\ph_{k-1}$ has valuation $\ge2$.
By Remark~\ref{remDeomgen}, the same is true of $\Al^\ga_\om B \ti\ph_{k-1}$
for every $\om\in2\pi\I\Z$ and every~$\ga$ starting close to~$0$ and ending
close to~$\om$: we have
\[
\cont_\ga \cB(B\ti\ph_{k-1})(\om+\ze) = \hat\psi(\ze)
\frac{\Log\ze}{2\pi} + R(\ze)
\]
with $\hat\psi \in \ze\C\{\ze\}$ and $R \in \C\{\ze\}$ depending on
$k$, $\om$, $\ga$, hence 
$\hat\chi(\ze) \defeq \frac{\hat\psi(\ze)}{\ee^\ze-1} \in \C\{\ze\}$
and (since $\ee^{\om+\ze} \equiv\ee^\ze$)
\[
\cont_\ga \hat\ph_k(\om+\ze) = 
\frac{c}{2\pi\I\ze} + \hat\chi(\ze)\frac{\Log\ze}{2\pi} + R^*(\ze),
\quad \text{with $c\defeq 2\pi\I R(0)$ and $R^*(\ze)\in\C\{\ze\}$.}
\]
Therefore $\hat\ph_k$ has only simple singularities.
\end{proof}


\begin{proof}[Proof of Lemma~\ref{lemcontrolhatphk}]
The set $\Om^+_{\eps,\tau,D}$ is such that we can find $M_0,L>0$ so
that
\beglabel{ineqFcnMeromOmP}
\ze \in \Om^+_{\eps,\tau,D} \Imp
\abs*{ \frac{\ze}{\ee^\ze-1} } \leq M_0 \, \ee^{-L \abs{\ze}}.
\elabel
On the other hand, we can find $C>L$ and $R_0>0$ such that the entire function~$\hat b$ satisfies
\[
\abs{\hat b(\ze)} \leq C \abs{\ze} \,\ee^{R_0\abs{\ze}}
\ens\text{for all $\ze\in\C$, hence}\ens
\abs{\hat b^{*k}(\ze)} \leq C^k \frac{\abs{\ze}^{2k-1}}{(2k-1)!} \ee^{R_0\abs{\ze}}
\ens\text{for all $\ze\in\C$ and $k\in\N^*$}
\]
by Lemma~\ref{lemestimconv}.


Let us give ourselves a naturally parametrised path
$\ga \col [0,\ell] \to \Om^+_{\eps,\tau,D}$
satisfying~\eqref{condeps}.
For any $2\pi\I\Z$-resurgent series~$\ti\psi$ with formal Borel
transform~$\hat\psi$, we have $B\ti\psi \in \ti\gR_{2\pi\I\Z}$ by
Theorem~\ref{thmresOmstbNL}, the proof of which shows that
$\hat B\hat\psi \defeq \cB(B\ti\psi)$
can be expressed as an integral transform
$\hat B\hat\psi(\ze) = \int_0^\ze K(\xi,\ze) \hat\psi(\xi) \,\dd\xi$
for $\ze$ close to~$0$,
with kernel function
\[
K(\xi,\ze) = \sum_{k\geq1} \frac{(-\xi)^k}{k!} \hat b^{*k}(\ze-\xi).
\]
The estimates available for~$\hat b^{*k}$ show that~$K$ is holomorphic
in $\C\times\C$, we can thus adapt the arguments of the ``easy''
Lemma~\ref{lemeasyconvol} and get
\[
\cont_\ga \hat B\hat\psi\big( \ga(s) \big) =
\int_0^s K\big( \ga(\sig),\ga(s)\big) \cont_\ga \hat\psi\big(\ga(\sig) \big)
\ga'(\sig)\,\dd\sig
\quad \text{for all $s\in[0,\ell]$.}
\]
The crude estimate
\[
\abs{K(\xi,\ze)} \leq C \abs{\xi} \, \ee^{\frac{C}{\mu}\abs{\xi} + (R_0+\mu)\abs{\ze-\xi}}
\quad\text{for all $(\xi,\ze)\in\C\times\C$,}
\]
with arbitrary $\mu>1$, will allow us to bound inductively 
$\cont_\ga\hat\ph_k = \cont_\ga\hat E\hat B\hat\ph_{k-1}$.


Indeed, the meromorphic function $\hat\ph_0 = \frac{\hat b}{\ee^\ze-1}$
satisfies~\eqref{ineqestimphk} with 
$A \defeq M_0 C$ and any $R\ge R_0$.
Suppose now that a $2\pi\I\Z$-resurgent function~$\hat\psi$ satisfies
\[
\abs*{\cont_\ga \hat\psi\big( \ga(s) \big)} \le \ee^{R s} \Psi(s)
\ens\text{for all $s\in[0,\ell]$,} \qquad
\text{with}\ens
R \defeq R_0 + \mu, \ens
\mu \defeq \frac{C}{\ka L},
\]
and a certain positive continuous function~$\Psi$.
Since $\abs{\ga(\sig)} \le \sig$ and $\abs{\ga(s)-\ga(\sig)} \le s-\sig$, we obtain
\[
\abs*{\cont_\ga \hat B\hat\psi\big( \ga(s) \big)} \le 
C s \, \ee^{ (\frac{C}{\mu} + R) s} \int_0^s\Psi(\sig)\,\dd\sig
\ens\text{for all $s\in[0,\ell]$,}
\]
whence
$\abs*{\cont_\ga \hat E\hat B\hat\psi\big( \ga(s) \big)} \le 
M \, \ee^{ R s} \int_0^s\Psi(\sig)\,\dd\sig$
with $M\defeq \frac{C M_0}{\ka}$
by~\eqref{ineqFcnMeromOmP},
using $\abs{\ga(s)} \ge \ka s$.
We thus get 
$\abs*{ \cont_\ga\hat\ph_k\big( \ga(s) \big) }
\leq A \, \ee^{Rs} \frac{ (Ms)^k }{k!}$
by induction on~$k$.
\end{proof}

\section{Fatou coordinates of a simple parabolic germ}


\parag
For every $R>0$ and $\de \in (0,\pi/2)$, we define
\[
\Sig^+_{R,\de} \defeq \{\, r\,\eith \in \C \mid
r>R,\; \abs{\th} < \pi-\de \,\},
\quad
\Sig^-_{R,\de} \defeq \{\, r\,\eith \in \C \mid
r>R,\; \abs{\th-\pi} < \pi-\de \,\}.
\]

\begin{Def}
A \emph{pair of Fatou coordinates at~$\infty$} is a pair $(v^+,v^-)$ of injective
holomorphic maps
\[
v^+ \col \Sig^+_{R,\de} \to \C, \qquad 
v^- \col \Sig^-_{R,\de} \to \C, 
\]
with some $R>0$ and $\de \in (0,\pi/2)$, 
such that
\[
v^+ \circ f = f_0 \circ v^+, \qquad
v^- \circ f = f_0 \circ v^-.
\]
\end{Def}

We still assume that $f \in \gG_1$ has vanishing resiter, with
iterator~$\ti v_*$ and inverse iterator~$\ti u_*$.
We still use the notations $I^+ = (-\frac{\pi}{2}, \frac{\pi}{2})$ 
and $I^- = (\frac{\pi}{2}, \frac{3\pi}{2})$.

\begin{thm}    \label{thmFatouCoord}
There exists locally bounded functions $\be, \be_1 \col I^+\cup I^- \to
(0,+\infty)$ such that $\be<\be_1$ and
\begin{itemize}
\item
$\ti v_* \in \ti\gG_0(I^+,\be) \cap \ti\gG_0(I^-,\be)$ and
$v_*^\pm \defeq \gS^{I^\pm}\ti v_*$ is injective on $\gD(I^\pm,\be)$
(notation of Definition~\ref{defglueLapl});
\item
$\ti u_* \in \ti\gG_0(I^+,\be_1) \cap \ti\gG_0(I^-,\be_1)$ and 
$u_*^\pm \defeq \gS^{I^\pm}\ti u_*$ is injective on $\gD(I^\pm,\be_1)$,
with 
\[
u_*^\pm\big( \gD(I^\pm,\be_1) \big) \subset \gD(I^\pm,\be)
\quad \text{and} \quad
v_*^\pm \circ u_*^\pm = \id
\ens\text{on $\gD(I^\pm,\be_1)$}.
\]
\end{itemize}
Moreover, the pairs of Fatou coordinates at~$\infty$ are the pairs $(v_*^++c^+,v_*^-+c^-)$
with arbitrary $c^+,c^-\in\C$.
\end{thm}

\begin{rem}
We may consider $(v_*^+,v_*^-)$ as a normalized pair of Fatou
coordinates.
Being obtained as Borel sums of a $1$-summable formal diffeomorphism,
they admit a $1$-Gevrey asymptotic expansion,
and the same is true of the inverse Fatou coordinates~$u_*^+$ and~$u_*^-$.
The first use of Borel-Laplace summation for obtaining Fatou coordinates is in \cite{Eca81}.
The asymptotic property without the Gevrey qualification can be found
in earlier works by G.~Birkhoff, G.~Szekeres, T.~Kimura and
J.~\'Ecalle---see \cite{Loray} and \cite{ML} for the references;
see \cite{LY} for a recent independent proof and an application to
numerical computations.
\end{rem}

\begin{proof}[Proof of Theorem~\ref{thmFatouCoord}]
The case $\ga=\{0\}$ of Theorem~\ref{thmIterSimpResurSumma} yields
locally bounded functions $\al,\be \col I^+ \cup I^- \to \R^+$ such
that $\ti v_* \in \ti\gG_0(I^\pm,\be,\al)$ (notation of
Definition~\ref{DefSummaDiffeo}).
In view of Theorem~\ref{thmSummaDiffeo}, we can replace~$\be$ by a
larger function so that~$v_*^\pm$ is injective on $\gD(I^\pm,\be)$.
We apply again Theorem~\ref{thmSummaDiffeo}: setting
\[
\be \;<\; 
\be^* \defeq \be + 2\sqrt{\al} \;<\;
\be_1 \defeq \be + (1+\sqrt2)\sqrt\al,
\]
we get $\ti u_* \in \ti\gG_0(I^\pm,\be^*)$, hence $\ti u_* \in
\ti\gG_0(I^\pm,\be_1)$, and all the desired properties follow.

By Lemma~\ref{lemSUMCVcase}, we have $f = \gS^{I^\pm}f$;
replacing the above function~$\be$ by a larger one if necessary so as
to take into account the domain of definition of~$f$,
Theorem~\ref{thmcompossummdiffeo} shows that 
$\gS^{I^\pm}(\ti v_*\circ f)=v_*^\pm\circ f$ and 
$\gS^{I^\pm}(f\circ\ti u_*)=f\circ u_*^\pm$.
In view of~\eqref{eqconjugtiv} and \eqref{eqDifftiust},
this yields
\beglabel{eqconjuupmvpm}
v_*^\pm \circ f = f_0 \circ v_*^\pm, \qquad
f \circ u_*^\pm = u_*^\pm \circ f_0.
\elabel
We see that for any $\de \in (0,\pi/2)$ there exists $R>0$ such that
$\Sig^\pm_{R,\de} \subset \gD(I^\pm,\be)$, therefore $(v_*^+,v_*^-)$
is a pair of Fatou coordinates.

Suppose now that $v^\pm$ is holomorphic and injective on
$\Sig^\pm_{R,\de}$.
Replacing the above function~$\be$ by a larger one if necessary, we
may suppose $\be\ge R$, then
$\gD(J^\pm,\be) \subset \Sig^\pm_{R,\de}$ with 
$J^+ \defeq (-\frac{\pi}{2}+\de, \frac{\pi}{2}-\de)$,
$J^- \defeq (\frac{\pi}{2}+\de, \frac{3\pi}{2}-\de)$.
By Theorem~\ref{thmSummaDiffeo}, we have
$u_*^\pm\big( \gD(J^\pm,\be_1) \big) \subset \gD(J^\pm,\be)$,
thus $\Phi^\pm \defeq v^\pm \circ u_*^\pm$ is holomorphic and injective on $\gD(J^\pm,\be_1)$.
In view of~\eqref{eqconjuupmvpm}, 
the equation $v^\pm\circ f = f_0\circ v^\pm$ is equivalent to
$f_0\circ\Phi^\pm = \Phi^\pm\circ f_0$, \ie $\Phi^\pm = \id+\Psi^\pm$
with~$\Psi^\pm$ $1$-periodic.
If $\Psi^\pm$ is a constant~$c^\pm$, then we find $v^\pm = v_*^\pm + c^\pm$.
In general, the periodicity of~$\Psi^\pm$ allows one to extend analytically~$\Phi^\pm$ to the whole
of~$\C$ and we get an injective entire function;
the Casorati-Weierstrass theorem shows that such a function must be of the form
$a z+c$, hence $\Psi^\pm$ is constant.
\end{proof}

\begin{figure} 
\begin{center}

\includegraphics[scale=1]{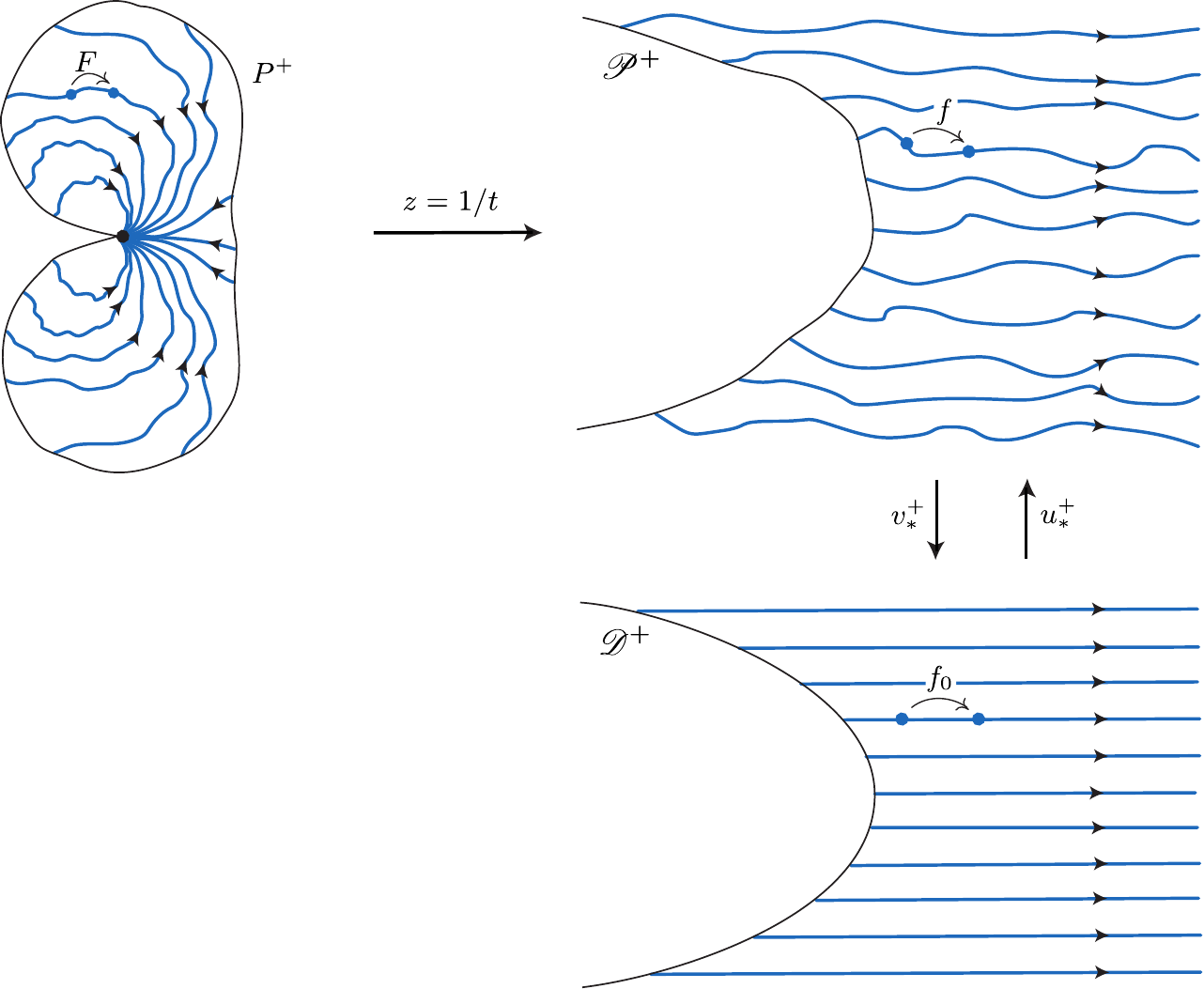}
%

\bigskip

\caption{%
\emph{The dynamics in the attracting petal viewed in three coordinates.}
}
\label{figzzABC}

\end{center} 
\end{figure}

\parag
Here are a few dynamical consequences of Theorem~\ref{thmFatouCoord}.
The domain $\gD^+ \defeq \gD(I^+,\be_1)$ is invariant by the normal
form $f_0 = \id+1$, while $\gD^- \defeq \gD(I^-,\be_1)$ is invariant
by the backward dynamics $f_0\ic = \id-1$,
hence 
\beglabel{eqdefgPpgPm}
\gP^+ \defeq u_*^+ (\gD^+) \;
\text{ is invariant by~$f$,}
\quad
\gP^- \defeq u_*^- (\gD^-) \;
\text{ is invariant by~$f\ic$,}
\elabel
and the conjugacy relations $f = u_*^+\circ f_0 \circ v_*^+$,
$f\ic = u_*^- \circ f_0\ic \circ v_*^-$ yield
\[
z \in \gP^+ \imp f^{\circ n}(z) = u_*^+\big( v_*^+(z)+n \big),
 \qquad
z \in \gP^- \imp f^{\circ(-n)}(z) = u_*^-\big( v_*^-(z)-n \big)
\]
for every $n\in\N$.
We thus see that all the forward orbits of~$f$ which start in~$\gP^+$ and all the
backward orbits of~$f$ which start in~$\gP^-$ are infinite and converge to the fixed
point at~$\infty$ (we could even describe the asymptotics \wrt\ the
discrete time~$n$)---see Figure~\ref{figzzABC}.


All this can be transferred to the variable $t=1/z$ and we get for the
dynamics of~$F$ a version of what is usually called the ``Leau-Fatou
flower theorem'':
we define the attracting and repelling ``petals'' by
\[
P^+ \defeq \{\, t \in \C^* \mid 1/t \in \gP^+ \,\}, \qquad
P^- \defeq \{\, t \in \C^* \mid 1/t \in \gP^- \,\},
\]
whose union is a punctured neighbourhood of~$0$,
and we see that all the forward orbits of~$F$ which start in~$P^+$ and all the
backward orbits of~$F$ which start in~$P^-$ are infinite and converge
to~$0$ (see Figure~\ref{figzzABC}).
Notice that~$P^+$ and~$P^-$ overlap, giving rise to two families of
bi-infinite orbits which are positively and negatively asymptotic to
the fixed point.

We can also define Fatou coordinates and inverse Fatou coordinates
at~$0$ as well as their formal counterparts by
\[
V_*^\pm(t) \defeq v_*^\pm(1/t), \quad
U_*^\pm(z) \defeq 1/u_*^\pm(z),
\qquad
\ti V_*(t) \defeq \ti v_*(1/t), \quad
\ti U_*(z) \defeq 1/\ti u_*(z),
\]
so that 
\begin{align}
\label{eqVpUp}
t\in P^+ &\Rightarrow V_*^+\big( F(t) \big) = V_*^+(t)+1, &
\quad z \in \gD^+ &\Rightarrow F\big( U_*^+(z) \big) = U_*^+(z+1), \\[1ex]
\label{eqVmUm}
t\in P^- &\Rightarrow V_*^-\big( F\ic(t) \big) = V_*^-(t)-1, &
\quad z \in \gD^- &\Rightarrow F\ic\big( U_*^-(z) \big) = U_*^-(z-1).
\end{align}
Observe that, with the notation $\ti v_*(z) = z + \sum_{k\ge1} a_k
z^{-k}$, we have
\[
V_*^\pm(t) \sim \ti V_*(t) = \frac{1}{t} + \sum_{k\ge1} a_k t^k,
\]
whereas $\ti u_*(z) = z + \ti\psi(z)$ with 
$\ti\psi(z) = \sum_{k\ge1} b_k z^{-k} \in \zcz$
implies
\[
\ti U_*(z) = z\ii \big( 1 + z\ii\ti\psi(z) \big)\ii \in \zcz.
\]
By Theorems~\ref{thmcompatsumcirc} and~\ref{thmComposSimpRes}, we see
that~$\ti U_*$ is a simple $2\pi\I\Z$-resurgent series, which is
$1$-summable in the directions of~$I^+$ and~$I^-$, with 
\[
U_*^\pm = \gS^{I^\pm}\ti U_*.
\]


\parag
Of course it may happen that one of the formal series $\ti v_*$, $\ti
u_*$, $\ti V_*$, $\ti U_*$ and thus all of them be convergent.
But this is the exception rather than the rule.


There is a case in which one easily proves that all of them are divergent.

\begin{lemma}
If~$F(t)$ or~$F\ic(t)$ extends to an entire function, then the formal series $\ti v_*$, $\ti
u_*$, $\ti V_*$, $\ti U_*$ are divergent.
\end{lemma}


\begin{proof}
Suppose that~$F$ is entire.
The function~$U_*^-(z)$, intially defined and holomorphic in~$\gD^-$,
which contains a left half-plane $\{\RE z < -c\}$, can be analytically
continued by repeated use of~\eqref{eqVmUm}:
for any $n\in\N^*$, the formula
\[
U_*^-(z) = F\big( U_*^-(z-1) \big) = \cdots = F^{\circ n}\big( U_*^-(z-n) \big)
\]
yields its analytic continuation in $\{ \RE(z) < -c+n \}$,
hence~$U_*^-$ extends to an entire function.
If~$\ti U_*$ had positive radius of convergence, then we would get
$U_*^- \sim_1 \ti U_*$ in a full neighbourhood of~$\infty$ by
Lemma~\ref{lemSUMCVcase}, in particular $U_*^-(z)$ would tend to~$0$ as
$\abs{z}\to\infty$ and thus be uniformly bounded; then the entire
function~$U_*^-$ would be constant by Liouville's theorem, which is
impossible because $\ti U_*(z) = z\ii + O(z^{-2})$.

If it is~$F\ic$ that extends to an entire function, then~$U_*^+$ extends to an entire function by virtue
of~\eqref{eqVpUp} and one can argue similarly to prove that~$\ti U_*$
is divergent.
\end{proof}


\section{The horn maps and the analytic
classification} \label{sec:HMAC}


In~\eqref{eqdefgPpgPm} we have defined~$\gP^+$ and~$\gP^-$ so
that~$v_*^+$ induces a biholomorphism $\gP^+ \xrightarrow{\sim} \gD^+$
and~$u_*^-$ induces a biholomorphism $\gD^- \xrightarrow{\sim} \gP^-$.
We can thus define a holomorphic function 
\beglabel{eqdefhornmaph}
h \defeq v_*^+ \circ u_*^- \col
\gD^- \cap (u_*^-)\ii(\gP^+) \to
\gD^+ \cap v_*^+(\gP^-),
\quad\text{such that}\ens
h \circ f_0 = f_0 \circ h
\elabel
(the fact that~$h$ conjugates~$f_0$ with itself stems
from~\eqref{eqconjuupmvpm}).

Let us define, for any $R>0$ and $\de \in (0,\pi/2)$,
\[
\gV\upp_{R,\de} \defeq \{\, r\,\eith \mid
r>R,\; \de < \th < \pi-\de \,\},
\quad
\gV\low_{R,\de} \defeq \{\, r\,\eith \mid
r>R,\; \pi+\de < \th < 2\pi-\de \,\}.
\]
Since~$v_*^+$ and~$u_*^-$ are close to identity near~$\infty$, there
exists $R>0$ such that the domain of definition of~$h$ has a connected
component which contains~$\gV\upp_{R,\pi/4}$ and a connected component
which contains~$\gV\low_{R,\pi/4}$,
so that in fact formula~\eqref{eqdefhornmaph} defines a
function~$h\upp$ and a function~$h\low$.


\begin{lemma}   \label{lemFourierHorn}
There exists $\sig>0$ such that the function~$h\upp$ extends analytically
to the upper half-plane $\{ \IM z > \sig \}$
and the function~$h\low$ extends analytically
to the lower half-plane $\{ \IM z < -\sig \}$.
The functions $h\upp-\id$ and $h\low-\id$ are $1$-periodic and admit
convergent Fourier expansions
\beglabel{eqFourierHorn}
h_*\upp(z) - z = \sum_{m=1}^{+\infty} A_{-m} \ee^{2\pi\I m z},
\qquad
h_*\low(z) - z = \sum_{m=1}^{+\infty} A_m \ee^{-2\pi\I m z},
\elabel
with $A_m = O( \ee^{\la\abs m} )$ for every $\la>2\pi \sig$.
\end{lemma}


\begin{proof}
The conjugacy relation $h\uplow\circ f_0 = f_0\circ h\uplow$ implies that $h\uplow$ is of the form
$\id+P\uplow$ with a $1$-periodic holomorphic function $P\uplow \col
\gV\uplow_{R,\pi/4} \to \C$.
By $1$-periodicity, $P\uplow$ extends analytically to an upper/lower
half-plane and can be written as $\chi(\ee^{\pm2\pi\I z})$, with~$\chi$
holomorphic in the punctured disc $\D^*_{2\pi\sig}$. 
The asymptotic behaviour of~$v_*^+$ and~$u_*^-$ at~$\infty$ in $\gD\big(
(-\frac{\pi}{4},\frac{\pi}{4}), \be_1 \big)$ shows that $h\uplow(z)=z
+o(1)$, hence $\chi(Z) \xrightarrow[Z\to0]{} 0$. 
Thus~$\chi$ is holomorphic in~$\D_{2\pi\sig}$ and vanishes at~$0$; 
\label{secFourierHalfplane}
its Taylor expansions yields the Fourier series of~$P\uplow$.
\end{proof}


\begin{Def}
We call $(h\upp,h\low)$ the pair of lifted horn maps of~$f$.
We call the coefficients of the sequence $(A_m)_{m\in\Z^*}$ the
\'Ecalle-Voronin invariants of~$f$.
\end{Def}


\begin{thm}   \label{thmAnalyticInv}
Two simple parabolic germs at~$\infty$ with vanishing resiter, $f$
and~$g$, are analytically conjugate if and only if
there exists $c\in\C$ such that their pairs of lifted horn maps
$(h_f\upp,h_f\low)$ and $(h_g\upp,h_g\low)$
are related by
\begin{gather}
\label{eqrelhghf}
h_g\upp(z) \equiv h_f\upp(z+c)-c, \qquad
h_g\low(z) \equiv h_f\low(z+c)-c,\\
\shortintertext{or, equivalently,}
A_m(g) = \ee^{-2\pi\I mc} A_m(f)
\quad \text{for every $m\in\Z^*$.}
\end{gather}
\end{thm}


\begin{proof}
  We denote by~$\ti v_f$, $v_f^\pm$, $\ti u_f$, $u_f^\pm$ the iterator
  of~$f$, its Borel sums and their inverses, and similarly~$\ti v_g$,
  $v_g^\pm$, $\ti u_g$, $u_g^\pm$ for~$g$.

Suppose that $f$ and~$g$ are analytically conjugate, so there exists
$h\in\gG$ (convergent!) such that $g\circ h = h \circ f$.
It follows that $\ti v_f\circ h\ic \circ g = f_0 \circ \ti v_f\circ
h\ic$, hence there exists $c\in\C$ such that $\ti v_f\circ h\ic = \ti
v_g + c$ by Lemma~\ref{lemtivstsumtiphk}.
Let $\tau \defeq \id + c$. We have $\ti v_g = \tau\ii \circ \ti v_f
\circ h\ic$ and $\ti u_g = h \circ \ti u_f \circ \tau$, whence
$v_g^+ = \tau\ii\circ v_f^+ \circ h\ic$ and 
$u_g^- = h \circ u_f^- \circ \tau$
by Theorem~\ref{thmcompossummdiffeo} and Lemma~\ref{lemSUMCVcase}.
This implies $v_g^+\circ u_g^- = \tau\ii\circ v_f^+ \circ u_f^- \circ \tau$,
\ie $h_g\uplow = \tau\ii \circ h_f\uplow \circ \tau$,
as desired.


Suppose now that there exists $c\in\C$
satisfying~\eqref{eqrelhghf}. We rewrite this relation as
\[
h_g\upp = \tau\ii \circ h_f\upp \circ \tau, \qquad
h_g\low = \tau\ii \circ h_f\low \circ \tau,
\]
with $\tau = \id+c$.
This implies
\[
\tau \circ v_g^+ \circ u_g^- = v_f^+ \circ u_f^- \circ \tau
\quad \text{on $\gV_{R,\de}\upp \cup \gV_{R,\de}\low$}
\]
with, say, $\de = 3\pi/4$ and $R$ large enough.
Therefore
\[
u_f^+ \circ \tau \circ v_g^+ = u_f^- \circ \tau \circ v_g^-
\quad \text{on $\gV_{R',\pi/4}\upp \cup \gV_{R',\pi/4}\low$}.
\]
This indicates that the functions $u_f^+ \circ \tau \circ v_g^+$ and
$u_f^- \circ \tau \circ v_g^-$ can be glued to form a function~$h$
holomorphic in punctured neighbourhood of~$\infty$;
the asymptotic behaviour then shows that~$h$ is holomorphic
at~$\infty$, with Taylor series $\ti u_f\circ \tau \circ \ti v_g$.
The conjugacy relations
$\ti u_g = g\circ\ti u_g\circ f_0\ic$ and
$\tau\circ\ti v_f\circ f = f_0\circ \tau\circ\ti v_f$
imply
$\ti u_g \circ \tau\circ\ti v_f\circ f = g\circ\ti u_g\circ
\tau\circ\ti v_f$,
hence $f$ and~$g$ are analytically conjugate by~$h$.
\end{proof}


Theorem~\ref{thmAnalyticInv} is just one part of \'Ecalle-Voronin's
classification result in the case of simple parabolic germs with
vanishing resiter.
The other part of the result (more difficult) says that
any pair of Fourier series of the form
$\big( \sum_{m\geq1} A_{-m} \ee^{2\pi\I m z}, \sum_{m\geq1} A_m
\ee^{-2\pi\I m z} \big)$,
where the first (\resp second) one is holomorphic in an upper (\resp lower) half-plane,
can be obtained as
$( h_*\upp-\id, h_*\low-\id )$
for a simple parabolic germ~$f$ with vanishing resiter.


\section{The Bridge Equation and the action of the symbolic Stokes
  automorphism} \label{sec:Bridge}


\parag
Let us give ourselves a simple parabolic germ at~$\infty$ with
vanishing resiter,~$f$.
So far, we have only exploited the summability statement contained in
Theorem~\ref{thmFatouCoord} and we have see that a deep information on
the analytic conjugacy class of~$f$ is encoded by the discrepancy
between the Borel sums~$v_*^+$ and~$v_*^-$, \ie by the lifted horn
maps.
Let us now see how the analysis of this discrepancy lends itself to
alien calculus, \ie to the study of the singularities in the Borel
plane.

We first use the operators~$\De_\om$ of Sections
\ref{secalopDepomDeom}--\ref{secopDeomDeriv} with $\om\in 2\pi\I\Z^*$.
They are derivations of the algebra $\ti\gR_{2\pi\I\Z}\simp$, and
they induce operators $\De_\om \col \ti\gG_{2\pi\I\Z}\simp \to \ti\gR_{2\pi\I\Z}\simp$ defined by
$\De_\om(\id+\ti\ph) \equiv \De_\om\ti\ph$.


\begin{thm}
There exists a sequence of complex numbers
$(C_\om)_{\om\in2\pi\I\Z^*}$ such that
\beglabel{eqBEDeomvu}
\De_\om\ti u_* = C_\om \pa\ti u_*,
\qquad 
\De_\om \ti v_* = - C_\om \, \ee^{-\om(\ti v_*-\id)}
\elabel
for each $\om\in2\pi\I\Z^*$.
\end{thm}


\begin{proof}
Let us apply $\De_\om$ to both sides of the conjugacy
equation~\eqref{eqDifftiust}:
by Theorem~\ref{thmSimpResdiffeo}, since $\De_\om f$ and $\De_\om f_0$
vanish, we get
\[
(\pa f)\circ \ti u_* \cdot \De_\om\ti u_* = (\De_\om\ti u_*)\circ f_0
\]
(we also used the fact that $\ee^{-\om(f_0-\id)} = 1$, since $\om\in 2\pi\I\Z^*$).
By applying~$\pa$ to~\eqref{eqDifftiust}, we also get
\[
(\pa f)\circ \ti u_* \cdot \pa\ti u_* = (\pa\ti u_*)\circ f_0.
\]
Since $\pa \ti u_* = 1+O(z^{-2})$, this implies that the formal series
$\ti C\defeq \frac{\De_\om\ti u_*}{\pa \ti u_*} \in \C[[z\ii]]$ satisfies
$\ti C = \ti C\circ f_0$.
Writing $\ti C\circ f_0-\ti C = \pa\ti C + \frac{1}{2!}\pa^2\ti C +
\cdots$ and reasoning on the valuation of~$\pa\ti C$, we see that $\ti
C$ must be constant.

We have $(\De_\om\ti u_*)\circ \ti v_* \cdot \pa\ti v_* = \ti C
(\pa\ti u_*)\circ \ti v_* \cdot \pa\ti v_* = \ti C \,\pa(\ti u_*\circ
\ti v_*) = \ti C$, hence
Formula~\eqref{eqDeomDiffeoInvers} yields 
$\De_\om\ti v = - \ti C\,\ee^{-\om(\ti v_*-\id)}$.
\end{proof}


The first equation in~\eqref{eqBEDeomvu} is called ``the Bridge
Equation for simple parabolic germs'': like
Equation~\eqref{eqBridgeEulerNL}, it yields a bridge between ordinary
differential calculus (here involving $\pa$) and alien calculus (when
dealing with the solution~$\ti u$ of the conjugacy equation~\eqref{eqDifftiust}).


\parag
From the operators~$\De_\om$ we can go the operators~$\De_\om^+$ by
means of formula~\eqref{eqDeomlogDeomp} of
Theorem~\ref{thmrelDeomDeomp}, according to which,
if one sets $\Om \defeq 2\pi\I\N^*$ or $\Om \defeq - 2\pi\I\N^*$, then
\beglabel{eqrecallDepom}
\De^+_\om = \sum_{s\ge1} \tfrac{1}{s!} \sum_{ \substack{ \om_1,\ldots,\om_s \in
    \Om \\ \om_1 + \cdots + \om_s = \om} }
\De_{\om_s} \circ \cdots \circ \De_{\om_2} \circ \De_{\om_1}
\quad\text{for $\om \in \Om$.}
\elabel
We also define
\beglabel{eqDefDemom}
\De^-_\om \defeq \sum_{s\ge1} \tfrac{(-1)^s}{s!} \sum_{ \substack{ \om_1,\ldots,\om_s \in
    \Om \\ \om_1 + \cdots + \om_s = \om} }
\De_{\om_s} \circ \cdots \circ \De_{\om_2} \circ \De_{\om_1}
\quad\text{for $\om \in \Om$.}
\elabel
The latter family of operators is related to Exercise~\ref{exoinvDDp}: they
correspond to the homogeneous components of $\exp(-\DD{}{\I\R^\pm})$
the same way the operators~$\De^+_\om$ correspond to the homogeneous
components of $\exp(\DD{}{\I\R^\pm})$---see formulas
\eqref{eqdottedalienop}--\eqref{eqHomCompStInvAutmo}. 


\begin{cor}   \label{cordefAmp}
Let $\Om \defeq 2\pi\I\N^*$ or $\Om \defeq - 2\pi\I\N^*$.
For each $\om \in\Om$, define
\[
S^+_\om \defeq 
- \sum_{s\ge 1} \tfrac{1}{s!}  \sum_{ \substack{ \om_1,\ldots,\om_s \in
    \Om \\ \om_1 + \cdots + \om_s = \om} }
\Ga_{\om_1,\ldots,\om_s}
C_{\om_1} \cdots C_{\om_s},
\quad
S^-_\om \defeq 
\sum_{s\ge 1} \tfrac{(-1)^{s-1}}{s!} \sum_{ \substack{ \om_1,\ldots,\om_s \in
    \Om \\ \om_1 + \cdots + \om_s = \om} }
\Ga_{\om_1,\ldots,\om_s}
C_{\om_1} \cdots C_{\om_s} 
\]
with $\Ga_{\om_1} \defeq 1$ and
$\Ga_{\om_1,\ldots,\om_s} \defeq \om_1 (\om_1+\om_2) \cdots (\om_1+\cdots+\om_{s-1})$.
Then 
\beglabel{eqBEDepomv}
\De^+_\om \ti v_* = S^+_\om\, \ee^{-\om(\ti v_*-\id)},
\qquad
\De^-_\om \ti v_* = S^-_\om\, \ee^{-\om(\ti v_*-\id)}.
\elabel
\end{cor}


\begin{proof}
Let $\ti\ph \defeq \ti v_*-\id$, so that the second equation
in~\eqref{eqBEDeomvu} reads
$\De_\om\ti\ph = - C_\om \, \ee^{-\om\ti\ph}$.
By repeated use of formula~\eqref{eqDeomSubst} of
Theorem~\ref{thmComposSimpRes}, we get
$\De_{\om_2}\De_{\om_1}\ti\ph = 
\om_1 C_{\om_1} \,\ee^{-\om_1\ti\ph} \De_{\om_2}\ti\ph 
= - \om_1 C_{\om_1} C_{\om_2} \,
\ee^{-(\om_1+\om_2)\ti\ph}$, 
$\De_{\om_3}\De_{\om_2}\De_{\om_1}\ti\ph = \ldots\,$, etc. The general
formula is
\[
\De_{\om_s} \cdots \De_{\om_1}\ti\ph = 
- \Ga_{\om_1,\ldots,\om_s}
C_{\om_1} \cdots C_{\om_s}\,
\ee^{-(\om_1+\cdots+\om_s)\ti\ph},
\]
whence the conclusion follows with the help of \eqref{eqrecallDepom}--\eqref{eqDefDemom}.
\end{proof}


In fact, in view of Remark~\ref{remDeomgen}, the above proof shows that, for every
$\om\in2\pi\I\Z$ and for every path~$\ga$ which starts close to~$0$
and ends close to~$\om$, there exists $S_\om^\ga \in \C$ such that
$\Al_\om^\ga\ti v_* = S_\om^\ga \,\ee^{-\om(\ti v_*-\id)}$.


\parag
We now wish to compute the action of the symbolic Stokes automorphism $\DD+{\I\R^\pm}$
on~$\ti v_*$ and
to describe the Stokes phenomenon in the spirit of
Section~\ref{secrelLaplStokesPhen},
so as to recover the horn maps of Section~\ref{sec:HMAC}.
We shall make use of the spaces
\[
\ti E^\pm \defeq \ti E(2\pi\I\Z,\I\R^\pm) = \bigpmN \ee^{-\om z} \ti\gR_{2\pi\I\Z}\simp
\]
introduced in Section~\ref{secExtBorelSymb};
since $2\pi\I\Z$ is an additive subgroup of~$\C$, these spaces are
differential algebras,
\[
\ti E^+ = \ti\gR_{2\pi\I\Z}\simp[[ \ee^{-2\pi\I z} ]],
\qquad
\ti E^- = \ti\gR_{2\pi\I\Z}\simp[[ \ee^{2\pi\I z} ]],
\qquad
\pa = \frac{\dd\,}{\dd z},
\]
on which are defined the directional alien
derivation~$\DD{}{\I\R^\pm}$ and the symbolic Stokes
automorphism $\DD+{\I\R^\pm} = \exp(\DD{}{\I\R^\pm})$.
According to Remark~\ref{rempadotted}, both operators commute with the
differential~$\pa$.
So does the ``inverse symbolic Stokes automorphism''
$\DD-{\I\R^\pm} \defeq \exp(-\DD{}{\I\R^\pm})$.


We find it convenient to modify slightly the notation for their
homogeneous components: from now on, we set
\beglabel{eqMODIFdottedalienop}
\om \in 2\pi\I\Z, \ens m\in\Z,\ens
\ti\ph \in \ti\gR_{2\pi\I\Z}\simp
\Imp
\left\{ \begin{aligned} 
\dDe\om(\ee^{-2\pi\I m z}\ti\ph) & \defeq \ee^{-(2\pi\I m+\om) z}
\De_\om\ti\ph,\\[1ex]
\dDepm\om(\ee^{-2\pi\I m z}\ti\ph) & \defeq \ee^{-(2\pi\I m+\om) z} \De^\pm_\om\ti\ph,
\end{aligned} \right. 
\elabel
so that
\begin{align}
\notag
\DD{}{\I\R^+} &= \sum_{\om\in 2\pi\I\N^*} \dDe\om
\quad\text{on $\ti E^+$,}
& & & 
\DD{}{\I\R^-} &= \sum_{\om\in -2\pi\I\N^*} \dDe\om
\quad\text{on $\ti E^-$,}\\[1ex]
\label{eqHomCompStInvAutmo}
\DD\pm{\I\R^+} &= \exp(\pm\DD{}{\I\R^+}) = \ID + \sum_{\om\in 2\pi\I\N^*} \dDepm\om,
& & & 
\DD\pm{\I\R^-} &= \exp(\pm\DD{}{\I\R^-}) = \ID + \sum_{\om\in -2\pi\I\N^*} \dDepm\om.
\end{align}


We may consider~$\ti v_*$ as an element of
$ \id + \ti\gR_{2\pi\I\Z}\simp \subset \id + \ti E^\pm$.
We thus set $\dDe\om\id \defeq 0$ and $\dDepm\om\id \defeq 0$ so that
the previous operators induce
\[
\DD{}{\I\R^+}, \DD+{\I\R^+}, \DD-{\I\R^+} \col 
\id + \ti E^+ \to \ti E^+,
\qquad
\DD{}{\I\R^-}, \DD+{\I\R^-}, \DD-{\I\R^-} \col 
\id + \ti E^- \to \ti E^-.
\]
This way \eqref{eqBEDeomvu} yields 
$\dDe\om\ti v_* = - C_\om \,\ee^{-\om\ti v_*}$
and~\eqref{eqBEDepomv} yields
$\dDepm\om\ti v_* = S^\pm_\om \,\ee^{-\om\ti v_*}$,
and we can write
\[
\DD\pm{\I\R^+} \ti v_* = \ti v_* +
\sum_{\om\in 2\pi\I\N^*} S^\pm_\om \,\ee^{-\om\ti v_*},
\qquad
\DD\pm{\I\R^-} \ti v_* = \ti v_* +
\sum_{\om\in -2\pi\I\N^*} S^\pm_\om \,\ee^{-\om\ti v_*}.
\]


\begin{thm}    \label{thmFromBEtoStokes}
We have
\begin{align}
\label{eqSpomhlow}
z + \sum_{\om\in 2\pi\I\N^*} S^+_\om \,\ee^{-\om z} &\equiv h_*\low(z),
&
\DD+{\I\R^+}\ti v_* &= h_*\low \circ \ti v_*, \\[1ex]
\label{eqSmomhlow}
z + \sum_{\om\in 2\pi\I\N^*} S^-_\om \,\ee^{-\om z} &\equiv (h_*\low)\ic(z),
&
\DD-{\I\R^+}\ti v_* &= (h_*\low)\ic \circ \ti v_*, \\[1ex]
\label{eqSpomhupp}
z + \sum_{\om\in -2\pi\I\N^*} S^+_\om \,\ee^{-\om z} &\equiv (h_*\upp)\ic(z),
&
\DD+{\I\R^-}\ti v_* &= (h_*\upp)\ic \circ \ti v_*, \\[1ex]
\label{eqSmomhupp}
z + \sum_{\om\in -2\pi\I\N^*} S^-_\om \,\ee^{-\om z} &\equiv h_*\upp(z),
&
\DD-{\I\R^-}\ti v_* &= h_*\upp \circ \ti v_*.
\end{align}
In particular the \'Ecalle-Voronin invariants $(A_m)_{m\in\Z^*}$ of
Lemma~\ref{lemFourierHorn} are given by
\beglabel{eqFromStoA}
A_{-m} =  S^-_{-2\pi\I m}, \quad A_m =  S^+_{2\pi\I m}, 
\qquad m\in\N^*.
\elabel
\end{thm}


\begin{rem}
  The ``exponential-like'' formulas which define $(S^\pm_\om)_{\om\in
    2\pi\I\N^*}$ from $(C_\om)_{\om\in 2\pi\I\N^*}$ in
  Corollary~\ref{cordefAmp} are clearly invertible, and similarly
$(C_\om)_{\om\in -2\pi\I\N^*} \mapsto (S^\pm_\om)_{\om\in
    2\pi\I\N^*}$ is invertible.
It follows that the coefficients~$C_\om$ of the Bridge
Equation~\eqref{eqBEDeomvu} are analytic conjugacy invariants too.
However there is an important difference between the $C$'s and the $S$'s:
Theorem~\ref{thmFromBEtoStokes} implies that there exists $\la>0$ such
that $S^\pm_{2\pi\I m} = O(\ee^{\la\abs m})$, but there are in general no
estimates of the same kind for the coefficients~$C_{2\pi\I m}$ of the
Bridge Equation.

\end{rem}


\begin{proof}
Let $I\defeq (0,\pi)$ and $\th\defeq\frac{\pi}{2}$, so that
$I^+ = (0,\frac{\pi}{2})$ and $I^- =(\frac{\pi}{2},\pi)$ with the notations of Section~\ref{secrelLaplStokesPhen}.
Let us pick $R>0$ large enough so that $h\low$ is defined by $v_*^+
\circ (v_*^-)\ic$ in $\gV\low_{R,\pi/4}$ (\cf \eqref{eqdefhornmaph}).

For any $m \in \N$, we deduce from the relation 
$\DD+{\I\R^+} \ti v_* = \ti v_* +
\sum_{\om\in 2\pi\I\N^*} S^+_\om \,\ee^{-\om\ti v_*}$
that
\[
[\DD+{\I\R^+} \ti v_*]_m = \ti v_* +
\sum_{j=0}^m S^+_{2\pi\I j} \,\ee^{-2\pi\I j\ti v_*}
\]
with notation~\ref{notatrigaE}.
Each term $\ee^{-2\pi\I j\ti v_*}$ is $2\pi\I\Z$-resurgent and
$1$-summable in the directions of~$I^\pm$, with
Borel sums
$\gS^{I^\pm} (\ee^{-2\pi\I j\ti v_*}) = \ee^{-2\pi\I j v^\pm_*}$,
hence Theorem~\ref{thmSymbStokesLapl} implies that
\[
z \in \gV\low_{R,\pi/4} \Imp
v_*^+(z) = v_*^-(z) + 
\sum_{j=0}^m S^+_{2\pi\I j} \,\ee^{-2\pi\I j v_*^-(z)} 
+ O(\ee^{-\rho\abs{\IM z}})
\]
for any $\rho \in (2\pi m, 2\pi(m+1))$.
It follows that
\[
z \in \gV\low_{R,\pi/4} \Imp
h_*\low(z) = z +
\sum_{j=0}^m S^+_{2\pi\I j} \,\ee^{-2\pi\I j z} 
+ O(\ee^{-\rho\abs{\IM z}})
\]
for any $\rho \in (2\pi m, 2\pi(m+1))$, whence~\eqref{eqSpomhlow} follows.

Formula~\eqref{eqSmomhlow} is obtained by the same chain of reasoning,
using a variant of Theorem~\ref{thmSymbStokesLapl} relating $\gS^-\ti
v_*$ and $\gS^+[\DD+{\I\R^+}\ti v_*]_m$.

Formulas~\eqref{eqSpomhupp} and~\eqref{eqSmomhupp} are obtained the
same way, using $I^+ \defeq (-\pi,-\frac{\pi}{2})$ and
$I^- \defeq (-\frac{\pi}{2},0)$,
but this time 
$\gS^{I^+} \ti v_* = v_*^-$ and $\gS^{I^-} \ti v_* = v_*^+$.
\end{proof}


\parag
We conclude by computing the action of the symbolic Stokes
automorphism $\DD+{\I\R^\pm}$ on~$\ti u_*$.

\begin{Def}
The derivation of~$\ti E^\pm$
\[
D_{\I\R^\pm} \defeq C_{\I\R^\pm}(z) \pa,
\quad\text{where}\ens
C_{\I\R^\pm}(z) = \sum_{\om\in\pm2\pi\I\N^*} C_\om \ee^{-\om z},
\]
is called the ``formal Stokes vector field'' of~$f$.
\end{Def}


Such a derivation $D_{\I\R^\pm}$ has a well-defined exponential, for the same reason
by which $\DD{}d$ had one according to Theorem~\ref{thmDelogDep}(iii):
it increases homogeneity by at least one unit.

\begin{lemma}
For any $\ti\phi\in\ti\gR_{2\pi\I\Z}\simp$,
\begin{align*}
\exp\big(C_{\I\R^\pm}(z) \pa\big) \ti\phi &=
\ti\phi \circ P_{\I\R^\pm}
&
\hspace{-4em}\text{with}\ens P_{\I\R^\pm}(z) &\defeq 
z + \sum_{\om\in\pm2\pi\I\N^*} S^-_{\om} \ee^{-\om z} \\[1ex]
\exp\big(-C_{\I\R^\pm}(z) \pa\big) \ti\phi &=
\ti\phi \circ Q_{\I\R^\pm}
&
\hspace{-4em}\text{with}\ens Q_{\I\R^\pm}(z) &\defeq 
z + \sum_{\om\in\pm2\pi\I\N^*} S^+_{\om} \ee^{-\om z}.
\end{align*}
\end{lemma}

\begin{proof}
Let $\Om = 2\pi\I\N^*$ or $\Om = -2\pi\I\N^*$ and,
accordingly, $C = C_{\I\R^+}$ or $C = C_{\I\R^-}$,
$D = D_{\I\R^+}$ or $D = D_{\I\R^-}$.
We have $C = \sum C_{\om_1} \,\ee^{-\om_1 z}$,
$D C = \sum (-\om_1)C_{\om_1}C_{\om_2} \,\ee^{-(\om_1+\om_2) z}$,
$D^2 C = \ldots\,$, etc. The general formula is
\[
D^{s-1} C = (-1)^{s-1} \sum_{\om_1,\ldots,\om_s\in\Om} 
\Ga_{\om_1,\ldots,\om_s} C_{\om_1} \cdots C_{\om_s}\,
\ee^{-(\om_1+\cdots+\om_s)z},
\qquad s\ge1.
\]
We thus set, for every $\om\in\Om$,
\[
S_\om(t) \defeq \sum_{s\ge1} 
 \tfrac{(-1)^{s-1} t^s}{s!} \sum_{ \substack{ \om_1,\ldots,\om_s \in
    \Om \\ \om_1 + \cdots + \om_s = \om} }
\Ga_{\om_1,\ldots,\om_s}
C_{\om_1} \cdots C_{\om_s} 
\in \C[t]
\]
(observe that $S_\om(t)$ is a polynomial of degree $\le m$ if
$\om=\pm2\pi\I m$),
so that $S_\om(1) = S^-_\om$ and $S_\om(-1) = S^+_\om$,
and
\[
G_t(z) \defeq 
\sum_{s\ge1} \frac{t^s}{s!} D^{s-1}C =
\sum_{\om\in\Om} S_\om(t) \, \ee^{-\om z}
\in \C[t][[\ee^{\mp2\pi\I z}]].
\]
We leave it to the reader to check by induction the combinatorial identity
\[
D^s\ti\phi = \sum_{ \substack{ n\ge1,\  s_1,\ldots,s_n\ge1 \\ s_1+\cdots+s_n=s } }
\frac{s!}{s_1! \cdots s_n! n!} 
(D^{s_1-1}C) \cdots (D^{s_n-1}C) \pa^n\ti\phi,
\qquad s\ge1
\]
for any $\ti\phi\in\ti\gR_{2\pi\I\Z}\simp$, whence
$\exp(tD)\ti\phi = \ti\phi + \sum_{n\ge 1} \frac{1}{n!} (G_t)^n
\pa\ti\phi = \ti\phi\circ(\id+G_t)$.
\end{proof}


In view of Theorem~\ref{thmFromBEtoStokes}, we get
\begin{cor}
\[
\exp\big(C_{\I\R^-}(z) \pa\big) \ti\phi =
\ti\phi \circ h_*\upp,
\qquad
\exp\big(-C_{\I\R^+}(z) \pa\big) \ti\phi =
\ti\phi \circ h_*\low
\]
for every $\ti\phi\in\ti\gR_{2\pi\I\Z}\simp$.
\end{cor}


Since the Bridge Equation can be rephrased as
\[
\DD{}{\I\R^\pm} \ti u_* = C_{\I\R^\pm} \pa \ti u_*
\]
and the operators~$\DD{}{\I\R^\pm}$ and~$D_{\I\R^\pm}$ commute, we
obtain
\begin{cor}
\[
\exp(t\DD{}{\I\R^\pm}) \ti u_* = \exp(t C_{\I\R^\pm} \pa) \ti u_*,
\qquad t\in \C.
\]
In particular
\begin{align*}
\DD+{\I\R^-} \ti u_* &= \ti u_* \circ h_*\upp,
& \hspace{-4em}
\DD-{\I\R^-} \ti u_* &= \ti u_* \circ (h_*\upp)\ic,
\\[1ex]
\DD+{\I\R^+} \ti u_* &= \ti u_* \circ (h_*\low)\ic,
& \hspace{-4em}
\DD-{\I\R^+} \ti u_* &= \ti u_* \circ h_*\low.
\end{align*}
\end{cor}


Expanding the last equation, we get
\[
\om \in 2\pi\I\N^* \Imp
\De^+_\om \ti u_* = \sum_{ \substack{ n\ge1,\  \om_1,\ldots,\om_n\ge1 \\ \om_1+\cdots+\om_n=\om } }
\tfrac{1}{n!} S^+_{\om_1}\cdots S^+_{\om_n} \pa^n\ti u_*.
\]
We leave it to the reader to compute the formula for $\De^+_\om \ti
u_*$ when $\om \in 2\pi\I\N^*$, and the formulas for $\De^\pm_\om u_*$ when $\om \in -2\pi\I\N^*$.


\newpage


\noindent {\em Acknowledgements.}
We are grateful to Mich\`ele Loday-Richaud for her help with the pictures.
The research leading to these results has received funding from the
European Comunity's Seventh Framework Program (FP7/2007--2013) under Grant
Agreement n.~236346
and from the French National Research Agency under the reference ANR-12-BS01-0017.


\bigskip


\vspace{.6cm}

\noindent
David Sauzin\\[1ex]
CNRS UMI 3483 - Laboratorio Fibonacci \\
Collegio Puteano, Scuola Normale Superiore di Pisa \\
Piazza dei Cavalieri 3, 56126 Pisa, Italy\\
email:\,{\tt{david.sauzin@sns.it}}


\end{document}